\newif\ifejorlayout \ejorlayoutfalse
\newif\ifnumrefs \numrefsfalse
\ifejorlayout
  \ifnumrefs\documentclass[preprint,11pt,number,a4paper]{elsarticle}\else\documentclass[preprint,11pt,authoryear,a4paper]{elsarticle}\fi
\else
  \ifnumrefs\documentclass[preprint,12pt,number,a4paper]{elsarticle}\else\documentclass[preprint,12pt,authoryear,a4paper]{elsarticle}\fi
\fi

\IfFileExists{lmodern.sty}{\usepackage[T1]{fontenc}\usepackage{lmodern}}{\usepackage[OT1]{fontenc}}
\usepackage[utf8]{inputenc}
\usepackage{amsmath,amssymb,amsthm,mathtools,bm}
\usepackage{booktabs,multirow,threeparttable}
\usepackage{graphicx}
\usepackage{algorithm}
\usepackage{algpseudocode}
\usepackage{xcolor}
\usepackage{tikz}
\usetikzlibrary{positioning,arrows.meta,calc,fit,backgrounds}
\definecolor{alecblue}{RGB}{42,120,214}   % colour of ALEC in the figures
\usepackage{colortbl}
\usepackage{enumitem}
\usepackage[hidelinks]{hyperref}
\usepackage{lineno}
\usepackage{placeins}
\ifejorlayout\usepackage{setspace}\AtBeginDocument{\onehalfspacing}\fi
\makeatletter
\newcommand{\resdef}[2]{\expandafter\def\csname res@#1\endcsname{#2}}
\newcommand{\res}[1]{\ifcsname res@#1\endcsname\csname res@#1\endcsname\else\textbf{\textcolor{red}{??#1??}}\fi}
\makeatother
\resdef{ex-cr}{0.0166}
\resdef{ex-rank-evm}{Water$\succ$Coffee$\succ$Soda$\succ$Milk$\succ$Beer$\succ$Tea$\succ$Wine}
\resdef{ex-w0-Coffee}{0.1765}
\resdef{ex-w0-Wine}{0.0294}
\resdef{ex-w0-Tea}{0.0882}
\resdef{ex-w0-Beer}{0.0882}
\resdef{ex-w0-Soda}{0.1765}
\resdef{ex-w0-Milk}{0.0882}
\resdef{ex-w0-Water}{0.3529}
\resdef{ex-stop-lc}{4}
\resdef{ex-stop-b10}{2}
\resdef{ex-beta-0}{0.353}
\resdef{ex-beta-1}{0.143}
\resdef{ex-beta-2}{0.049}
\resdef{ex-beta-3}{0.037}
\resdef{ex-beta-4}{0.021}
\resdef{ex-beta-5}{0.006}
\resdef{ex-beta-6}{0.002}
\resdef{esm-RANDOM-hit-0}{50.2}
\resdef{esm-RANDOM-tau-0}{0.717}
\resdef{esm-RANDOM-l1-0}{0.420}
\resdef{esm-RANDOM-hit-1}{54.2}
\resdef{esm-RANDOM-tau-1}{0.745}
\resdef{esm-RANDOM-l1-1}{0.377}
\resdef{esm-RANDOM-hit-2}{57.0}
\resdef{esm-RANDOM-tau-2}{0.768}
\resdef{esm-RANDOM-l1-2}{0.339}
\resdef{esm-RANDOM-hit-3}{59.4}
\resdef{esm-RANDOM-tau-3}{0.779}
\resdef{esm-RANDOM-l1-3}{0.315}
\resdef{esm-RANDOM-hit-4}{62.0}
\resdef{esm-RANDOM-tau-4}{0.797}
\resdef{esm-RANDOM-l1-4}{0.292}
\resdef{esm-RANDOM-hit-5}{65.0}
\resdef{esm-RANDOM-tau-5}{0.815}
\resdef{esm-RANDOM-l1-5}{0.271}
\resdef{esm-RANDOM-hit-6}{67.4}
\resdef{esm-RANDOM-tau-6}{0.824}
\resdef{esm-RANDOM-l1-6}{0.254}
\resdef{esm-RANDOM-hit-7}{69.8}
\resdef{esm-RANDOM-tau-7}{0.835}
\resdef{esm-RANDOM-l1-7}{0.238}
\resdef{esm-RANDOM-hit-8}{73.2}
\resdef{esm-RANDOM-tau-8}{0.844}
\resdef{esm-RANDOM-l1-8}{0.226}
\resdef{esm-RANDOM-hit-9}{75.8}
\resdef{esm-RANDOM-tau-9}{0.852}
\resdef{esm-RANDOM-l1-9}{0.213}
\resdef{esm-RANDOM-hit-10}{78.4}
\resdef{esm-RANDOM-tau-10}{0.858}
\resdef{esm-RANDOM-l1-10}{0.199}
\resdef{esm-HARKER-PG-hit-0}{51.2}
\resdef{esm-HARKER-PG-tau-0}{0.715}
\resdef{esm-HARKER-PG-l1-0}{0.423}
\resdef{esm-HARKER-PG-hit-1}{64.8}
\resdef{esm-HARKER-PG-tau-1}{0.739}
\resdef{esm-HARKER-PG-l1-1}{0.346}
\resdef{esm-HARKER-PG-hit-2}{66.4}
\resdef{esm-HARKER-PG-tau-2}{0.750}
\resdef{esm-HARKER-PG-l1-2}{0.322}
\resdef{esm-HARKER-PG-hit-3}{70.4}
\resdef{esm-HARKER-PG-tau-3}{0.764}
\resdef{esm-HARKER-PG-l1-3}{0.309}
\resdef{esm-HARKER-PG-hit-4}{71.6}
\resdef{esm-HARKER-PG-tau-4}{0.772}
\resdef{esm-HARKER-PG-l1-4}{0.295}
\resdef{esm-HARKER-PG-hit-5}{72.4}
\resdef{esm-HARKER-PG-tau-5}{0.781}
\resdef{esm-HARKER-PG-l1-5}{0.286}
\resdef{esm-HARKER-PG-hit-6}{74.2}
\resdef{esm-HARKER-PG-tau-6}{0.787}
\resdef{esm-HARKER-PG-l1-6}{0.278}
\resdef{esm-HARKER-PG-hit-7}{75.6}
\resdef{esm-HARKER-PG-tau-7}{0.793}
\resdef{esm-HARKER-PG-l1-7}{0.271}
\resdef{esm-HARKER-PG-hit-8}{77.6}
\resdef{esm-HARKER-PG-tau-8}{0.798}
\resdef{esm-HARKER-PG-l1-8}{0.264}
\resdef{esm-HARKER-PG-hit-9}{77.2}
\resdef{esm-HARKER-PG-tau-9}{0.803}
\resdef{esm-HARKER-PG-l1-9}{0.259}
\resdef{esm-HARKER-PG-hit-10}{79.4}
\resdef{esm-HARKER-PG-tau-10}{0.808}
\resdef{esm-HARKER-PG-l1-10}{0.252}
\resdef{esm-HARKER-hit-0}{52.2}
\resdef{esm-HARKER-tau-0}{0.705}
\resdef{esm-HARKER-l1-0}{0.428}
\resdef{esm-HARKER-hit-1}{61.4}
\resdef{esm-HARKER-tau-1}{0.728}
\resdef{esm-HARKER-l1-1}{0.351}
\resdef{esm-HARKER-hit-2}{64.4}
\resdef{esm-HARKER-tau-2}{0.735}
\resdef{esm-HARKER-l1-2}{0.327}
\resdef{esm-HARKER-hit-3}{65.2}
\resdef{esm-HARKER-tau-3}{0.742}
\resdef{esm-HARKER-l1-3}{0.318}
\resdef{esm-HARKER-hit-4}{66.4}
\resdef{esm-HARKER-tau-4}{0.747}
\resdef{esm-HARKER-l1-4}{0.313}
\resdef{esm-HARKER-hit-5}{67.6}
\resdef{esm-HARKER-tau-5}{0.755}
\resdef{esm-HARKER-l1-5}{0.304}
\resdef{esm-HARKER-hit-6}{68.6}
\resdef{esm-HARKER-tau-6}{0.759}
\resdef{esm-HARKER-l1-6}{0.299}
\resdef{esm-HARKER-hit-7}{69.4}
\resdef{esm-HARKER-tau-7}{0.764}
\resdef{esm-HARKER-l1-7}{0.293}
\resdef{esm-HARKER-hit-8}{70.4}
\resdef{esm-HARKER-tau-8}{0.771}
\resdef{esm-HARKER-l1-8}{0.288}
\resdef{esm-HARKER-hit-9}{70.2}
\resdef{esm-HARKER-tau-9}{0.776}
\resdef{esm-HARKER-l1-9}{0.281}
\resdef{esm-HARKER-hit-10}{70.6}
\resdef{esm-HARKER-tau-10}{0.782}
\resdef{esm-HARKER-l1-10}{0.276}
\resdef{esm-HARKER-LOG-hit-0}{52.2}
\resdef{esm-HARKER-LOG-tau-0}{0.705}
\resdef{esm-HARKER-LOG-l1-0}{0.428}
\resdef{esm-HARKER-LOG-hit-1}{63.6}
\resdef{esm-HARKER-LOG-tau-1}{0.714}
\resdef{esm-HARKER-LOG-l1-1}{0.400}
\resdef{esm-HARKER-LOG-hit-2}{72.6}
\resdef{esm-HARKER-LOG-tau-2}{0.734}
\resdef{esm-HARKER-LOG-l1-2}{0.379}
\resdef{esm-HARKER-LOG-hit-3}{77.4}
\resdef{esm-HARKER-LOG-tau-3}{0.745}
\resdef{esm-HARKER-LOG-l1-3}{0.354}
\resdef{esm-HARKER-LOG-hit-4}{83.0}
\resdef{esm-HARKER-LOG-tau-4}{0.766}
\resdef{esm-HARKER-LOG-l1-4}{0.328}
\resdef{esm-HARKER-LOG-hit-5}{86.8}
\resdef{esm-HARKER-LOG-tau-5}{0.778}
\resdef{esm-HARKER-LOG-l1-5}{0.312}
\resdef{esm-HARKER-LOG-hit-6}{88.6}
\resdef{esm-HARKER-LOG-tau-6}{0.790}
\resdef{esm-HARKER-LOG-l1-6}{0.296}
\resdef{esm-HARKER-LOG-hit-7}{93.2}
\resdef{esm-HARKER-LOG-tau-7}{0.805}
\resdef{esm-HARKER-LOG-l1-7}{0.275}
\resdef{esm-HARKER-LOG-hit-8}{94.4}
\resdef{esm-HARKER-LOG-tau-8}{0.818}
\resdef{esm-HARKER-LOG-l1-8}{0.258}
\resdef{esm-HARKER-LOG-hit-9}{96.2}
\resdef{esm-HARKER-LOG-tau-9}{0.829}
\resdef{esm-HARKER-LOG-l1-9}{0.242}
\resdef{esm-HARKER-LOG-hit-10}{97.2}
\resdef{esm-HARKER-LOG-tau-10}{0.835}
\resdef{esm-HARKER-LOG-l1-10}{0.229}
\resdef{esm-AOPT-hit-0}{50.2}
\resdef{esm-AOPT-tau-0}{0.717}
\resdef{esm-AOPT-l1-0}{0.420}
\resdef{esm-AOPT-hit-1}{54.6}
\resdef{esm-AOPT-tau-1}{0.766}
\resdef{esm-AOPT-l1-1}{0.346}
\resdef{esm-AOPT-hit-2}{58.4}
\resdef{esm-AOPT-tau-2}{0.785}
\resdef{esm-AOPT-l1-2}{0.310}
\resdef{esm-AOPT-hit-3}{62.6}
\resdef{esm-AOPT-tau-3}{0.806}
\resdef{esm-AOPT-l1-3}{0.278}
\resdef{esm-AOPT-hit-4}{63.6}
\resdef{esm-AOPT-tau-4}{0.817}
\resdef{esm-AOPT-l1-4}{0.258}
\resdef{esm-AOPT-hit-5}{65.8}
\resdef{esm-AOPT-tau-5}{0.832}
\resdef{esm-AOPT-l1-5}{0.240}
\resdef{esm-AOPT-hit-6}{70.2}
\resdef{esm-AOPT-tau-6}{0.843}
\resdef{esm-AOPT-l1-6}{0.224}
\resdef{esm-AOPT-hit-7}{72.4}
\resdef{esm-AOPT-tau-7}{0.848}
\resdef{esm-AOPT-l1-7}{0.211}
\resdef{esm-AOPT-hit-8}{75.6}
\resdef{esm-AOPT-tau-8}{0.855}
\resdef{esm-AOPT-l1-8}{0.198}
\resdef{esm-AOPT-hit-9}{76.6}
\resdef{esm-AOPT-tau-9}{0.862}
\resdef{esm-AOPT-l1-9}{0.188}
\resdef{esm-AOPT-hit-10}{77.2}
\resdef{esm-AOPT-tau-10}{0.867}
\resdef{esm-AOPT-l1-10}{0.180}
\resdef{esm-FROZEN-hit-0}{50.2}
\resdef{esm-FROZEN-tau-0}{0.717}
\resdef{esm-FROZEN-l1-0}{0.420}
\resdef{esm-FROZEN-hit-1}{74.4}
\resdef{esm-FROZEN-tau-1}{0.758}
\resdef{esm-FROZEN-l1-1}{0.357}
\resdef{esm-FROZEN-hit-2}{81.0}
\resdef{esm-FROZEN-tau-2}{0.792}
\resdef{esm-FROZEN-l1-2}{0.310}
\resdef{esm-FROZEN-hit-3}{85.6}
\resdef{esm-FROZEN-tau-3}{0.824}
\resdef{esm-FROZEN-l1-3}{0.266}
\resdef{esm-FROZEN-hit-4}{86.4}
\resdef{esm-FROZEN-tau-4}{0.846}
\resdef{esm-FROZEN-l1-4}{0.236}
\resdef{esm-FROZEN-hit-5}{87.8}
\resdef{esm-FROZEN-tau-5}{0.858}
\resdef{esm-FROZEN-l1-5}{0.217}
\resdef{esm-FROZEN-hit-6}{89.6}
\resdef{esm-FROZEN-tau-6}{0.866}
\resdef{esm-FROZEN-l1-6}{0.201}
\resdef{esm-FROZEN-hit-7}{89.6}
\resdef{esm-FROZEN-tau-7}{0.870}
\resdef{esm-FROZEN-l1-7}{0.189}
\resdef{esm-FROZEN-hit-8}{90.4}
\resdef{esm-FROZEN-tau-8}{0.871}
\resdef{esm-FROZEN-l1-8}{0.177}
\resdef{esm-FROZEN-hit-9}{91.2}
\resdef{esm-FROZEN-tau-9}{0.876}
\resdef{esm-FROZEN-l1-9}{0.168}
\resdef{esm-FROZEN-hit-10}{91.2}
\resdef{esm-FROZEN-tau-10}{0.883}
\resdef{esm-FROZEN-l1-10}{0.161}
\resdef{esm-STAR-hit-0}{50.2}
\resdef{esm-STAR-tau-0}{0.717}
\resdef{esm-STAR-l1-0}{0.420}
\resdef{esm-STAR-hit-1}{74.4}
\resdef{esm-STAR-tau-1}{0.758}
\resdef{esm-STAR-l1-1}{0.357}
\resdef{esm-STAR-hit-2}{81.0}
\resdef{esm-STAR-tau-2}{0.792}
\resdef{esm-STAR-l1-2}{0.310}
\resdef{esm-STAR-hit-3}{85.6}
\resdef{esm-STAR-tau-3}{0.824}
\resdef{esm-STAR-l1-3}{0.266}
\resdef{esm-STAR-hit-4}{86.4}
\resdef{esm-STAR-tau-4}{0.846}
\resdef{esm-STAR-l1-4}{0.236}
\resdef{esm-STAR-hit-5}{87.8}
\resdef{esm-STAR-tau-5}{0.858}
\resdef{esm-STAR-l1-5}{0.217}
\resdef{esm-STAR-hit-6}{89.6}
\resdef{esm-STAR-tau-6}{0.866}
\resdef{esm-STAR-l1-6}{0.201}
\resdef{esm-STAR-hit-7}{90.4}
\resdef{esm-STAR-tau-7}{0.871}
\resdef{esm-STAR-l1-7}{0.188}
\resdef{esm-STAR-hit-8}{91.6}
\resdef{esm-STAR-tau-8}{0.871}
\resdef{esm-STAR-l1-8}{0.176}
\resdef{esm-STAR-hit-9}{92.8}
\resdef{esm-STAR-tau-9}{0.877}
\resdef{esm-STAR-l1-9}{0.168}
\resdef{esm-STAR-hit-10}{95.6}
\resdef{esm-STAR-tau-10}{0.885}
\resdef{esm-STAR-l1-10}{0.164}
\resdef{esm-ALEC-H-hit-0}{52.2}
\resdef{esm-ALEC-H-tau-0}{0.705}
\resdef{esm-ALEC-H-l1-0}{0.428}
\resdef{esm-ALEC-H-hit-1}{74.4}
\resdef{esm-ALEC-H-tau-1}{0.723}
\resdef{esm-ALEC-H-l1-1}{0.401}
\resdef{esm-ALEC-H-hit-2}{80.0}
\resdef{esm-ALEC-H-tau-2}{0.746}
\resdef{esm-ALEC-H-l1-2}{0.371}
\resdef{esm-ALEC-H-hit-3}{82.4}
\resdef{esm-ALEC-H-tau-3}{0.772}
\resdef{esm-ALEC-H-l1-3}{0.335}
\resdef{esm-ALEC-H-hit-4}{85.8}
\resdef{esm-ALEC-H-tau-4}{0.805}
\resdef{esm-ALEC-H-l1-4}{0.293}
\resdef{esm-ALEC-H-hit-5}{88.0}
\resdef{esm-ALEC-H-tau-5}{0.828}
\resdef{esm-ALEC-H-l1-5}{0.261}
\resdef{esm-ALEC-H-hit-6}{89.4}
\resdef{esm-ALEC-H-tau-6}{0.840}
\resdef{esm-ALEC-H-l1-6}{0.240}
\resdef{esm-ALEC-H-hit-7}{89.8}
\resdef{esm-ALEC-H-tau-7}{0.844}
\resdef{esm-ALEC-H-l1-7}{0.227}
\resdef{esm-ALEC-H-hit-8}{91.4}
\resdef{esm-ALEC-H-tau-8}{0.849}
\resdef{esm-ALEC-H-l1-8}{0.215}
\resdef{esm-ALEC-H-hit-9}{94.4}
\resdef{esm-ALEC-H-tau-9}{0.854}
\resdef{esm-ALEC-H-l1-9}{0.204}
\resdef{esm-ALEC-H-hit-10}{94.4}
\resdef{esm-ALEC-H-tau-10}{0.861}
\resdef{esm-ALEC-H-l1-10}{0.195}
\resdef{esm-ALEC-hit-0}{50.2}
\resdef{esm-ALEC-tau-0}{0.717}
\resdef{esm-ALEC-l1-0}{0.420}
\resdef{esm-ALEC-hit-1}{74.4}
\resdef{esm-ALEC-tau-1}{0.758}
\resdef{esm-ALEC-l1-1}{0.357}
\resdef{esm-ALEC-hit-2}{82.6}
\resdef{esm-ALEC-tau-2}{0.791}
\resdef{esm-ALEC-l1-2}{0.313}
\resdef{esm-ALEC-hit-3}{86.8}
\resdef{esm-ALEC-tau-3}{0.822}
\resdef{esm-ALEC-l1-3}{0.275}
\resdef{esm-ALEC-hit-4}{89.2}
\resdef{esm-ALEC-tau-4}{0.841}
\resdef{esm-ALEC-l1-4}{0.252}
\resdef{esm-ALEC-hit-5}{90.0}
\resdef{esm-ALEC-tau-5}{0.855}
\resdef{esm-ALEC-l1-5}{0.230}
\resdef{esm-ALEC-hit-6}{91.4}
\resdef{esm-ALEC-tau-6}{0.865}
\resdef{esm-ALEC-l1-6}{0.213}
\resdef{esm-ALEC-hit-7}{91.2}
\resdef{esm-ALEC-tau-7}{0.869}
\resdef{esm-ALEC-l1-7}{0.198}
\resdef{esm-ALEC-hit-8}{92.4}
\resdef{esm-ALEC-tau-8}{0.869}
\resdef{esm-ALEC-l1-8}{0.183}
\resdef{esm-ALEC-hit-9}{90.8}
\resdef{esm-ALEC-tau-9}{0.875}
\resdef{esm-ALEC-l1-9}{0.173}
\resdef{esm-ALEC-hit-10}{92.6}
\resdef{esm-ALEC-tau-10}{0.883}
\resdef{esm-ALEC-l1-10}{0.166}
\resdef{esm-ALEC-U-hit-0}{50.2}
\resdef{esm-ALEC-U-tau-0}{0.717}
\resdef{esm-ALEC-U-l1-0}{0.420}
\resdef{esm-ALEC-U-hit-1}{71.8}
\resdef{esm-ALEC-U-tau-1}{0.770}
\resdef{esm-ALEC-U-l1-1}{0.338}
\resdef{esm-ALEC-U-hit-2}{80.6}
\resdef{esm-ALEC-U-tau-2}{0.801}
\resdef{esm-ALEC-U-l1-2}{0.300}
\resdef{esm-ALEC-U-hit-3}{84.2}
\resdef{esm-ALEC-U-tau-3}{0.825}
\resdef{esm-ALEC-U-l1-3}{0.270}
\resdef{esm-ALEC-U-hit-4}{86.8}
\resdef{esm-ALEC-U-tau-4}{0.843}
\resdef{esm-ALEC-U-l1-4}{0.250}
\resdef{esm-ALEC-U-hit-5}{90.8}
\resdef{esm-ALEC-U-tau-5}{0.855}
\resdef{esm-ALEC-U-l1-5}{0.235}
\resdef{esm-ALEC-U-hit-6}{94.0}
\resdef{esm-ALEC-U-tau-6}{0.867}
\resdef{esm-ALEC-U-l1-6}{0.222}
\resdef{esm-ALEC-U-hit-7}{96.8}
\resdef{esm-ALEC-U-tau-7}{0.876}
\resdef{esm-ALEC-U-l1-7}{0.209}
\resdef{esm-ALEC-U-hit-8}{97.6}
\resdef{esm-ALEC-U-tau-8}{0.885}
\resdef{esm-ALEC-U-l1-8}{0.198}
\resdef{esm-ALEC-U-hit-9}{98.4}
\resdef{esm-ALEC-U-tau-9}{0.887}
\resdef{esm-ALEC-U-l1-9}{0.189}
\resdef{esm-ALEC-U-hit-10}{99.2}
\resdef{esm-ALEC-U-tau-10}{0.894}
\resdef{esm-ALEC-U-l1-10}{0.176}
\resdef{esm-p-RANDOM-1}{$p<0.001$}
\resdef{esm-p-HARKER-PG-1}{$p<0.001$}
\resdef{esm-p-HARKER-1}{$p<0.001$}
\resdef{esm-p-HARKER-LOG-1}{$p<0.001$}
\resdef{esm-p-AOPT-1}{$p<0.001$}
\resdef{esm-p-FROZEN-1}{$p=1.000$}
\resdef{esm-p-STAR-1}{$p=1.000$}
\resdef{esm-p-ALEC-U-1}{$p=0.182$}
\resdef{esm-p-ALEC-H-1}{$p=1.000$}
\resdef{esm-p-RANDOM-2}{$p<0.001$}
\resdef{esm-p-HARKER-PG-2}{$p<0.001$}
\resdef{esm-p-HARKER-2}{$p<0.001$}
\resdef{esm-p-HARKER-LOG-2}{$p<0.001$}
\resdef{esm-p-AOPT-2}{$p<0.001$}
\resdef{esm-p-FROZEN-2}{$p=0.185$}
\resdef{esm-p-STAR-2}{$p=0.185$}
\resdef{esm-p-ALEC-U-2}{$p=0.245$}
\resdef{esm-p-ALEC-H-2}{$p=0.024$}
\resdef{esm-p-RANDOM-3}{$p<0.001$}
\resdef{esm-p-HARKER-PG-3}{$p<0.001$}
\resdef{esm-p-HARKER-3}{$p<0.001$}
\resdef{esm-p-HARKER-LOG-3}{$p<0.001$}
\resdef{esm-p-AOPT-3}{$p<0.001$}
\resdef{esm-p-FROZEN-3}{$p=0.362$}
\resdef{esm-p-STAR-3}{$p=0.362$}
\resdef{esm-p-ALEC-U-3}{$p=0.117$}
\resdef{esm-p-ALEC-H-3}{$p<0.001$}
\resdef{esm-p-RANDOM-4}{$p<0.001$}
\resdef{esm-p-HARKER-PG-4}{$p<0.001$}
\resdef{esm-p-HARKER-4}{$p<0.001$}
\resdef{esm-p-HARKER-LOG-4}{$p<0.001$}
\resdef{esm-p-AOPT-4}{$p<0.001$}
\resdef{esm-p-FROZEN-4}{$p=0.020$}
\resdef{esm-p-STAR-4}{$p=0.020$}
\resdef{esm-p-ALEC-U-4}{$p=0.126$}
\resdef{esm-p-ALEC-H-4}{$p=0.002$}
\resdef{esm-p-RANDOM-5}{$p<0.001$}
\resdef{esm-p-HARKER-PG-5}{$p<0.001$}
\resdef{esm-p-HARKER-5}{$p<0.001$}
\resdef{esm-p-HARKER-LOG-5}{$p=0.056$}
\resdef{esm-p-AOPT-5}{$p<0.001$}
\resdef{esm-p-FROZEN-5}{$p=0.080$}
\resdef{esm-p-STAR-5}{$p=0.080$}
\resdef{esm-p-ALEC-U-5}{$p=0.665$}
\resdef{esm-p-ALEC-H-5}{$p=0.041$}
\resdef{esm-p-RANDOM-6}{$p<0.001$}
\resdef{esm-p-HARKER-PG-6}{$p<0.001$}
\resdef{esm-p-HARKER-6}{$p<0.001$}
\resdef{esm-p-HARKER-LOG-6}{$p=0.081$}
\resdef{esm-p-AOPT-6}{$p<0.001$}
\resdef{esm-p-FROZEN-6}{$p=0.150$}
\resdef{esm-p-STAR-6}{$p=0.150$}
\resdef{esm-p-ALEC-U-6}{$p=0.066$}
\resdef{esm-p-ALEC-H-6}{$p=0.041$}
\resdef{esm-p-RANDOM-7}{$p<0.001$}
\resdef{esm-p-HARKER-PG-7}{$p<0.001$}
\resdef{esm-p-HARKER-7}{$p<0.001$}
\resdef{esm-p-HARKER-LOG-7}{$p=0.220$}
\resdef{esm-p-AOPT-7}{$p<0.001$}
\resdef{esm-p-FROZEN-7}{$p=0.243$}
\resdef{esm-p-STAR-7}{$p=0.608$}
\resdef{esm-p-ALEC-U-7}{$p<0.001$}
\resdef{esm-p-ALEC-H-7}{$p=0.248$}
\resdef{esm-p-RANDOM-8}{$p<0.001$}
\resdef{esm-p-HARKER-PG-8}{$p<0.001$}
\resdef{esm-p-HARKER-8}{$p<0.001$}
\resdef{esm-p-HARKER-LOG-8}{$p=0.203$}
\resdef{esm-p-AOPT-8}{$p<0.001$}
\resdef{esm-p-FROZEN-8}{$p=0.121$}
\resdef{esm-p-STAR-8}{$p=0.585$}
\resdef{esm-p-ALEC-U-8}{$p<0.001$}
\resdef{esm-p-ALEC-H-8}{$p=0.473$}
\resdef{esm-p-RANDOM-9}{$p<0.001$}
\resdef{esm-p-HARKER-PG-9}{$p<0.001$}
\resdef{esm-p-HARKER-9}{$p<0.001$}
\resdef{esm-p-HARKER-LOG-9}{$p<0.001$}
\resdef{esm-p-AOPT-9}{$p<0.001$}
\resdef{esm-p-FROZEN-9}{$p=0.880$}
\resdef{esm-p-STAR-9}{$p=0.064$}
\resdef{esm-p-ALEC-U-9}{$p<0.001$}
\resdef{esm-p-ALEC-H-9}{$p=0.003$}
\resdef{esm-p-RANDOM-10}{$p<0.001$}
\resdef{esm-p-HARKER-PG-10}{$p<0.001$}
\resdef{esm-p-HARKER-10}{$p<0.001$}
\resdef{esm-p-HARKER-LOG-10}{$p=0.002$}
\resdef{esm-p-AOPT-10}{$p<0.001$}
\resdef{esm-p-FROZEN-10}{$p=0.401$}
\resdef{esm-p-STAR-10}{$p=0.004$}
\resdef{esm-p-ALEC-U-10}{$p<0.001$}
\resdef{esm-p-ALEC-H-10}{$p=0.188$}
\resdef{bev-RANDOM-hit-0}{87.6}
\resdef{bev-RANDOM-tau-0}{0.838}
\resdef{bev-RANDOM-l1-0}{0.141}
\resdef{bev-RANDOM-hit-1}{93.0}
\resdef{bev-RANDOM-tau-1}{0.847}
\resdef{bev-RANDOM-l1-1}{0.120}
\resdef{bev-RANDOM-hit-2}{95.2}
\resdef{bev-RANDOM-tau-2}{0.857}
\resdef{bev-RANDOM-l1-2}{0.101}
\resdef{bev-RANDOM-hit-3}{97.6}
\resdef{bev-RANDOM-tau-3}{0.874}
\resdef{bev-RANDOM-l1-3}{0.089}
\resdef{bev-RANDOM-hit-4}{98.6}
\resdef{bev-RANDOM-tau-4}{0.884}
\resdef{bev-RANDOM-l1-4}{0.079}
\resdef{bev-RANDOM-hit-5}{99.4}
\resdef{bev-RANDOM-tau-5}{0.893}
\resdef{bev-RANDOM-l1-5}{0.069}
\resdef{bev-RANDOM-hit-6}{100.0}
\resdef{bev-RANDOM-tau-6}{0.899}
\resdef{bev-RANDOM-l1-6}{0.061}
\resdef{bev-RANDOM-hit-7}{100.0}
\resdef{bev-RANDOM-tau-7}{0.912}
\resdef{bev-RANDOM-l1-7}{0.054}
\resdef{bev-RANDOM-hit-8}{100.0}
\resdef{bev-RANDOM-tau-8}{0.920}
\resdef{bev-RANDOM-l1-8}{0.049}
\resdef{bev-RANDOM-hit-9}{100.0}
\resdef{bev-RANDOM-tau-9}{0.928}
\resdef{bev-RANDOM-l1-9}{0.043}
\resdef{bev-RANDOM-hit-10}{100.0}
\resdef{bev-RANDOM-tau-10}{0.940}
\resdef{bev-RANDOM-l1-10}{0.037}
\resdef{bev-HARKER-PG-hit-0}{88.4}
\resdef{bev-HARKER-PG-tau-0}{0.829}
\resdef{bev-HARKER-PG-l1-0}{0.142}
\resdef{bev-HARKER-PG-hit-1}{92.2}
\resdef{bev-HARKER-PG-tau-1}{0.840}
\resdef{bev-HARKER-PG-l1-1}{0.114}
\resdef{bev-HARKER-PG-hit-2}{93.6}
\resdef{bev-HARKER-PG-tau-2}{0.861}
\resdef{bev-HARKER-PG-l1-2}{0.103}
\resdef{bev-HARKER-PG-hit-3}{93.2}
\resdef{bev-HARKER-PG-tau-3}{0.884}
\resdef{bev-HARKER-PG-l1-3}{0.096}
\resdef{bev-HARKER-PG-hit-4}{95.0}
\resdef{bev-HARKER-PG-tau-4}{0.903}
\resdef{bev-HARKER-PG-l1-4}{0.087}
\resdef{bev-HARKER-PG-hit-5}{95.8}
\resdef{bev-HARKER-PG-tau-5}{0.911}
\resdef{bev-HARKER-PG-l1-5}{0.082}
\resdef{bev-HARKER-PG-hit-6}{98.6}
\resdef{bev-HARKER-PG-tau-6}{0.910}
\resdef{bev-HARKER-PG-l1-6}{0.070}
\resdef{bev-HARKER-PG-hit-7}{99.8}
\resdef{bev-HARKER-PG-tau-7}{0.915}
\resdef{bev-HARKER-PG-l1-7}{0.058}
\resdef{bev-HARKER-PG-hit-8}{100.0}
\resdef{bev-HARKER-PG-tau-8}{0.925}
\resdef{bev-HARKER-PG-l1-8}{0.049}
\resdef{bev-HARKER-PG-hit-9}{100.0}
\resdef{bev-HARKER-PG-tau-9}{0.938}
\resdef{bev-HARKER-PG-l1-9}{0.043}
\resdef{bev-HARKER-PG-hit-10}{100.0}
\resdef{bev-HARKER-PG-tau-10}{0.949}
\resdef{bev-HARKER-PG-l1-10}{0.039}
\resdef{bev-HARKER-hit-0}{87.2}
\resdef{bev-HARKER-tau-0}{0.819}
\resdef{bev-HARKER-l1-0}{0.142}
\resdef{bev-HARKER-hit-1}{92.8}
\resdef{bev-HARKER-tau-1}{0.836}
\resdef{bev-HARKER-l1-1}{0.116}
\resdef{bev-HARKER-hit-2}{94.8}
\resdef{bev-HARKER-tau-2}{0.860}
\resdef{bev-HARKER-l1-2}{0.103}
\resdef{bev-HARKER-hit-3}{94.4}
\resdef{bev-HARKER-tau-3}{0.883}
\resdef{bev-HARKER-l1-3}{0.095}
\resdef{bev-HARKER-hit-4}{96.2}
\resdef{bev-HARKER-tau-4}{0.905}
\resdef{bev-HARKER-l1-4}{0.087}
\resdef{bev-HARKER-hit-5}{98.0}
\resdef{bev-HARKER-tau-5}{0.914}
\resdef{bev-HARKER-l1-5}{0.080}
\resdef{bev-HARKER-hit-6}{99.8}
\resdef{bev-HARKER-tau-6}{0.919}
\resdef{bev-HARKER-l1-6}{0.072}
\resdef{bev-HARKER-hit-7}{100.0}
\resdef{bev-HARKER-tau-7}{0.923}
\resdef{bev-HARKER-l1-7}{0.061}
\resdef{bev-HARKER-hit-8}{100.0}
\resdef{bev-HARKER-tau-8}{0.935}
\resdef{bev-HARKER-l1-8}{0.051}
\resdef{bev-HARKER-hit-9}{100.0}
\resdef{bev-HARKER-tau-9}{0.944}
\resdef{bev-HARKER-l1-9}{0.046}
\resdef{bev-HARKER-hit-10}{100.0}
\resdef{bev-HARKER-tau-10}{0.942}
\resdef{bev-HARKER-l1-10}{0.044}
\resdef{bev-HARKER-LOG-hit-0}{87.2}
\resdef{bev-HARKER-LOG-tau-0}{0.819}
\resdef{bev-HARKER-LOG-l1-0}{0.142}
\resdef{bev-HARKER-LOG-hit-1}{97.6}
\resdef{bev-HARKER-LOG-tau-1}{0.828}
\resdef{bev-HARKER-LOG-l1-1}{0.117}
\resdef{bev-HARKER-LOG-hit-2}{99.6}
\resdef{bev-HARKER-LOG-tau-2}{0.842}
\resdef{bev-HARKER-LOG-l1-2}{0.099}
\resdef{bev-HARKER-LOG-hit-3}{100.0}
\resdef{bev-HARKER-LOG-tau-3}{0.857}
\resdef{bev-HARKER-LOG-l1-3}{0.085}
\resdef{bev-HARKER-LOG-hit-4}{100.0}
\resdef{bev-HARKER-LOG-tau-4}{0.869}
\resdef{bev-HARKER-LOG-l1-4}{0.073}
\resdef{bev-HARKER-LOG-hit-5}{100.0}
\resdef{bev-HARKER-LOG-tau-5}{0.880}
\resdef{bev-HARKER-LOG-l1-5}{0.064}
\resdef{bev-HARKER-LOG-hit-6}{100.0}
\resdef{bev-HARKER-LOG-tau-6}{0.887}
\resdef{bev-HARKER-LOG-l1-6}{0.058}
\resdef{bev-HARKER-LOG-hit-7}{100.0}
\resdef{bev-HARKER-LOG-tau-7}{0.882}
\resdef{bev-HARKER-LOG-l1-7}{0.054}
\resdef{bev-HARKER-LOG-hit-8}{100.0}
\resdef{bev-HARKER-LOG-tau-8}{0.881}
\resdef{bev-HARKER-LOG-l1-8}{0.048}
\resdef{bev-HARKER-LOG-hit-9}{100.0}
\resdef{bev-HARKER-LOG-tau-9}{0.889}
\resdef{bev-HARKER-LOG-l1-9}{0.041}
\resdef{bev-HARKER-LOG-hit-10}{100.0}
\resdef{bev-HARKER-LOG-tau-10}{0.918}
\resdef{bev-HARKER-LOG-l1-10}{0.029}
\resdef{bev-AOPT-hit-0}{87.6}
\resdef{bev-AOPT-tau-0}{0.838}
\resdef{bev-AOPT-l1-0}{0.141}
\resdef{bev-AOPT-hit-1}{97.0}
\resdef{bev-AOPT-tau-1}{0.859}
\resdef{bev-AOPT-l1-1}{0.107}
\resdef{bev-AOPT-hit-2}{98.8}
\resdef{bev-AOPT-tau-2}{0.870}
\resdef{bev-AOPT-l1-2}{0.088}
\resdef{bev-AOPT-hit-3}{100.0}
\resdef{bev-AOPT-tau-3}{0.881}
\resdef{bev-AOPT-l1-3}{0.076}
\resdef{bev-AOPT-hit-4}{100.0}
\resdef{bev-AOPT-tau-4}{0.897}
\resdef{bev-AOPT-l1-4}{0.067}
\resdef{bev-AOPT-hit-5}{100.0}
\resdef{bev-AOPT-tau-5}{0.904}
\resdef{bev-AOPT-l1-5}{0.062}
\resdef{bev-AOPT-hit-6}{100.0}
\resdef{bev-AOPT-tau-6}{0.914}
\resdef{bev-AOPT-l1-6}{0.056}
\resdef{bev-AOPT-hit-7}{100.0}
\resdef{bev-AOPT-tau-7}{0.923}
\resdef{bev-AOPT-l1-7}{0.052}
\resdef{bev-AOPT-hit-8}{100.0}
\resdef{bev-AOPT-tau-8}{0.925}
\resdef{bev-AOPT-l1-8}{0.047}
\resdef{bev-AOPT-hit-9}{100.0}
\resdef{bev-AOPT-tau-9}{0.942}
\resdef{bev-AOPT-l1-9}{0.039}
\resdef{bev-AOPT-hit-10}{100.0}
\resdef{bev-AOPT-tau-10}{0.953}
\resdef{bev-AOPT-l1-10}{0.034}
\resdef{bev-FROZEN-hit-0}{87.6}
\resdef{bev-FROZEN-tau-0}{0.838}
\resdef{bev-FROZEN-l1-0}{0.141}
\resdef{bev-FROZEN-hit-1}{98.8}
\resdef{bev-FROZEN-tau-1}{0.842}
\resdef{bev-FROZEN-l1-1}{0.113}
\resdef{bev-FROZEN-hit-2}{99.8}
\resdef{bev-FROZEN-tau-2}{0.851}
\resdef{bev-FROZEN-l1-2}{0.095}
\resdef{bev-FROZEN-hit-3}{100.0}
\resdef{bev-FROZEN-tau-3}{0.872}
\resdef{bev-FROZEN-l1-3}{0.079}
\resdef{bev-FROZEN-hit-4}{100.0}
\resdef{bev-FROZEN-tau-4}{0.901}
\resdef{bev-FROZEN-l1-4}{0.067}
\resdef{bev-FROZEN-hit-5}{100.0}
\resdef{bev-FROZEN-tau-5}{0.912}
\resdef{bev-FROZEN-l1-5}{0.060}
\resdef{bev-FROZEN-hit-6}{100.0}
\resdef{bev-FROZEN-tau-6}{0.912}
\resdef{bev-FROZEN-l1-6}{0.057}
\resdef{bev-FROZEN-hit-7}{100.0}
\resdef{bev-FROZEN-tau-7}{0.912}
\resdef{bev-FROZEN-l1-7}{0.052}
\resdef{bev-FROZEN-hit-8}{100.0}
\resdef{bev-FROZEN-tau-8}{0.913}
\resdef{bev-FROZEN-l1-8}{0.046}
\resdef{bev-FROZEN-hit-9}{100.0}
\resdef{bev-FROZEN-tau-9}{0.925}
\resdef{bev-FROZEN-l1-9}{0.036}
\resdef{bev-FROZEN-hit-10}{100.0}
\resdef{bev-FROZEN-tau-10}{0.938}
\resdef{bev-FROZEN-l1-10}{0.027}
\resdef{bev-STAR-hit-0}{87.6}
\resdef{bev-STAR-tau-0}{0.838}
\resdef{bev-STAR-l1-0}{0.141}
\resdef{bev-STAR-hit-1}{98.8}
\resdef{bev-STAR-tau-1}{0.842}
\resdef{bev-STAR-l1-1}{0.113}
\resdef{bev-STAR-hit-2}{99.8}
\resdef{bev-STAR-tau-2}{0.851}
\resdef{bev-STAR-l1-2}{0.095}
\resdef{bev-STAR-hit-3}{100.0}
\resdef{bev-STAR-tau-3}{0.872}
\resdef{bev-STAR-l1-3}{0.079}
\resdef{bev-STAR-hit-4}{100.0}
\resdef{bev-STAR-tau-4}{0.901}
\resdef{bev-STAR-l1-4}{0.067}
\resdef{bev-STAR-hit-5}{100.0}
\resdef{bev-STAR-tau-5}{0.905}
\resdef{bev-STAR-l1-5}{0.062}
\resdef{bev-STAR-hit-6}{100.0}
\resdef{bev-STAR-tau-6}{0.901}
\resdef{bev-STAR-l1-6}{0.058}
\resdef{bev-STAR-hit-7}{100.0}
\resdef{bev-STAR-tau-7}{0.896}
\resdef{bev-STAR-l1-7}{0.054}
\resdef{bev-STAR-hit-8}{100.0}
\resdef{bev-STAR-tau-8}{0.895}
\resdef{bev-STAR-l1-8}{0.047}
\resdef{bev-STAR-hit-9}{100.0}
\resdef{bev-STAR-tau-9}{0.910}
\resdef{bev-STAR-l1-9}{0.037}
\resdef{bev-STAR-hit-10}{100.0}
\resdef{bev-STAR-tau-10}{0.935}
\resdef{bev-STAR-l1-10}{0.025}
\resdef{bev-ALEC-H-hit-0}{87.2}
\resdef{bev-ALEC-H-tau-0}{0.819}
\resdef{bev-ALEC-H-l1-0}{0.142}
\resdef{bev-ALEC-H-hit-1}{98.8}
\resdef{bev-ALEC-H-tau-1}{0.827}
\resdef{bev-ALEC-H-l1-1}{0.116}
\resdef{bev-ALEC-H-hit-2}{100.0}
\resdef{bev-ALEC-H-tau-2}{0.846}
\resdef{bev-ALEC-H-l1-2}{0.101}
\resdef{bev-ALEC-H-hit-3}{100.0}
\resdef{bev-ALEC-H-tau-3}{0.881}
\resdef{bev-ALEC-H-l1-3}{0.081}
\resdef{bev-ALEC-H-hit-4}{100.0}
\resdef{bev-ALEC-H-tau-4}{0.907}
\resdef{bev-ALEC-H-l1-4}{0.068}
\resdef{bev-ALEC-H-hit-5}{100.0}
\resdef{bev-ALEC-H-tau-5}{0.907}
\resdef{bev-ALEC-H-l1-5}{0.063}
\resdef{bev-ALEC-H-hit-6}{100.0}
\resdef{bev-ALEC-H-tau-6}{0.904}
\resdef{bev-ALEC-H-l1-6}{0.059}
\resdef{bev-ALEC-H-hit-7}{100.0}
\resdef{bev-ALEC-H-tau-7}{0.896}
\resdef{bev-ALEC-H-l1-7}{0.056}
\resdef{bev-ALEC-H-hit-8}{100.0}
\resdef{bev-ALEC-H-tau-8}{0.894}
\resdef{bev-ALEC-H-l1-8}{0.049}
\resdef{bev-ALEC-H-hit-9}{100.0}
\resdef{bev-ALEC-H-tau-9}{0.908}
\resdef{bev-ALEC-H-l1-9}{0.039}
\resdef{bev-ALEC-H-hit-10}{100.0}
\resdef{bev-ALEC-H-tau-10}{0.934}
\resdef{bev-ALEC-H-l1-10}{0.027}
\resdef{bev-ALEC-hit-0}{87.6}
\resdef{bev-ALEC-tau-0}{0.838}
\resdef{bev-ALEC-l1-0}{0.141}
\resdef{bev-ALEC-hit-1}{98.8}
\resdef{bev-ALEC-tau-1}{0.842}
\resdef{bev-ALEC-l1-1}{0.113}
\resdef{bev-ALEC-hit-2}{100.0}
\resdef{bev-ALEC-tau-2}{0.855}
\resdef{bev-ALEC-l1-2}{0.096}
\resdef{bev-ALEC-hit-3}{100.0}
\resdef{bev-ALEC-tau-3}{0.881}
\resdef{bev-ALEC-l1-3}{0.078}
\resdef{bev-ALEC-hit-4}{100.0}
\resdef{bev-ALEC-tau-4}{0.906}
\resdef{bev-ALEC-l1-4}{0.065}
\resdef{bev-ALEC-hit-5}{100.0}
\resdef{bev-ALEC-tau-5}{0.906}
\resdef{bev-ALEC-l1-5}{0.061}
\resdef{bev-ALEC-hit-6}{100.0}
\resdef{bev-ALEC-tau-6}{0.900}
\resdef{bev-ALEC-l1-6}{0.057}
\resdef{bev-ALEC-hit-7}{100.0}
\resdef{bev-ALEC-tau-7}{0.893}
\resdef{bev-ALEC-l1-7}{0.054}
\resdef{bev-ALEC-hit-8}{100.0}
\resdef{bev-ALEC-tau-8}{0.893}
\resdef{bev-ALEC-l1-8}{0.048}
\resdef{bev-ALEC-hit-9}{100.0}
\resdef{bev-ALEC-tau-9}{0.909}
\resdef{bev-ALEC-l1-9}{0.037}
\resdef{bev-ALEC-hit-10}{100.0}
\resdef{bev-ALEC-tau-10}{0.935}
\resdef{bev-ALEC-l1-10}{0.025}
\resdef{bev-ALEC-U-hit-0}{87.6}
\resdef{bev-ALEC-U-tau-0}{0.838}
\resdef{bev-ALEC-U-l1-0}{0.141}
\resdef{bev-ALEC-U-hit-1}{97.4}
\resdef{bev-ALEC-U-tau-1}{0.849}
\resdef{bev-ALEC-U-l1-1}{0.112}
\resdef{bev-ALEC-U-hit-2}{99.6}
\resdef{bev-ALEC-U-tau-2}{0.854}
\resdef{bev-ALEC-U-l1-2}{0.093}
\resdef{bev-ALEC-U-hit-3}{100.0}
\resdef{bev-ALEC-U-tau-3}{0.862}
\resdef{bev-ALEC-U-l1-3}{0.079}
\resdef{bev-ALEC-U-hit-4}{100.0}
\resdef{bev-ALEC-U-tau-4}{0.875}
\resdef{bev-ALEC-U-l1-4}{0.067}
\resdef{bev-ALEC-U-hit-5}{100.0}
\resdef{bev-ALEC-U-tau-5}{0.880}
\resdef{bev-ALEC-U-l1-5}{0.063}
\resdef{bev-ALEC-U-hit-6}{100.0}
\resdef{bev-ALEC-U-tau-6}{0.886}
\resdef{bev-ALEC-U-l1-6}{0.058}
\resdef{bev-ALEC-U-hit-7}{100.0}
\resdef{bev-ALEC-U-tau-7}{0.890}
\resdef{bev-ALEC-U-l1-7}{0.053}
\resdef{bev-ALEC-U-hit-8}{100.0}
\resdef{bev-ALEC-U-tau-8}{0.901}
\resdef{bev-ALEC-U-l1-8}{0.046}
\resdef{bev-ALEC-U-hit-9}{100.0}
\resdef{bev-ALEC-U-tau-9}{0.914}
\resdef{bev-ALEC-U-l1-9}{0.037}
\resdef{bev-ALEC-U-hit-10}{100.0}
\resdef{bev-ALEC-U-tau-10}{0.936}
\resdef{bev-ALEC-U-l1-10}{0.025}
\resdef{bev-p-RANDOM-1}{$p<0.001$}
\resdef{bev-p-HARKER-PG-1}{$p<0.001$}
\resdef{bev-p-HARKER-1}{$p<0.001$}
\resdef{bev-p-HARKER-LOG-1}{$p=0.109$}
\resdef{bev-p-AOPT-1}{$p=0.064$}
\resdef{bev-p-FROZEN-1}{$p=1.000$}
\resdef{bev-p-STAR-1}{$p=1.000$}
\resdef{bev-p-ALEC-U-1}{$p=0.118$}
\resdef{bev-p-ALEC-H-1}{$p=1.000$}
\resdef{bev-p-RANDOM-2}{$p<0.001$}
\resdef{bev-p-HARKER-PG-2}{$p<0.001$}
\resdef{bev-p-HARKER-2}{$p<0.001$}
\resdef{bev-p-HARKER-LOG-2}{$p=0.500$}
\resdef{bev-p-AOPT-2}{$p=0.031$}
\resdef{bev-p-FROZEN-2}{$p=1.000$}
\resdef{bev-p-STAR-2}{$p=1.000$}
\resdef{bev-p-ALEC-U-2}{$p=0.500$}
\resdef{bev-p-ALEC-H-2}{$p=1.000$}
\resdef{bev-p-RANDOM-3}{$p<0.001$}
\resdef{bev-p-HARKER-PG-3}{$p<0.001$}
\resdef{bev-p-HARKER-3}{$p<0.001$}
\resdef{bev-p-HARKER-LOG-3}{$p=1.000$}
\resdef{bev-p-AOPT-3}{$p=1.000$}
\resdef{bev-p-FROZEN-3}{$p=1.000$}
\resdef{bev-p-STAR-3}{$p=1.000$}
\resdef{bev-p-ALEC-U-3}{$p=1.000$}
\resdef{bev-p-ALEC-H-3}{$p=1.000$}
\resdef{bev-p-RANDOM-4}{$p=0.016$}
\resdef{bev-p-HARKER-PG-4}{$p<0.001$}
\resdef{bev-p-HARKER-4}{$p<0.001$}
\resdef{bev-p-HARKER-LOG-4}{$p=1.000$}
\resdef{bev-p-AOPT-4}{$p=1.000$}
\resdef{bev-p-FROZEN-4}{$p=1.000$}
\resdef{bev-p-STAR-4}{$p=1.000$}
\resdef{bev-p-ALEC-U-4}{$p=1.000$}
\resdef{bev-p-ALEC-H-4}{$p=1.000$}
\resdef{bev-p-RANDOM-5}{$p=0.250$}
\resdef{bev-p-HARKER-PG-5}{$p<0.001$}
\resdef{bev-p-HARKER-5}{$p=0.002$}
\resdef{bev-p-HARKER-LOG-5}{$p=1.000$}
\resdef{bev-p-AOPT-5}{$p=1.000$}
\resdef{bev-p-FROZEN-5}{$p=1.000$}
\resdef{bev-p-STAR-5}{$p=1.000$}
\resdef{bev-p-ALEC-U-5}{$p=1.000$}
\resdef{bev-p-ALEC-H-5}{$p=1.000$}
\resdef{bev-p-RANDOM-6}{$p=1.000$}
\resdef{bev-p-HARKER-PG-6}{$p=0.016$}
\resdef{bev-p-HARKER-6}{$p=1.000$}
\resdef{bev-p-HARKER-LOG-6}{$p=1.000$}
\resdef{bev-p-AOPT-6}{$p=1.000$}
\resdef{bev-p-FROZEN-6}{$p=1.000$}
\resdef{bev-p-STAR-6}{$p=1.000$}
\resdef{bev-p-ALEC-U-6}{$p=1.000$}
\resdef{bev-p-ALEC-H-6}{$p=1.000$}
\resdef{bev-p-RANDOM-7}{$p=1.000$}
\resdef{bev-p-HARKER-PG-7}{$p=1.000$}
\resdef{bev-p-HARKER-7}{$p=1.000$}
\resdef{bev-p-HARKER-LOG-7}{$p=1.000$}
\resdef{bev-p-AOPT-7}{$p=1.000$}
\resdef{bev-p-FROZEN-7}{$p=1.000$}
\resdef{bev-p-STAR-7}{$p=1.000$}
\resdef{bev-p-ALEC-U-7}{$p=1.000$}
\resdef{bev-p-ALEC-H-7}{$p=1.000$}
\resdef{bev-p-RANDOM-8}{$p=1.000$}
\resdef{bev-p-HARKER-PG-8}{$p=1.000$}
\resdef{bev-p-HARKER-8}{$p=1.000$}
\resdef{bev-p-HARKER-LOG-8}{$p=1.000$}
\resdef{bev-p-AOPT-8}{$p=1.000$}
\resdef{bev-p-FROZEN-8}{$p=1.000$}
\resdef{bev-p-STAR-8}{$p=1.000$}
\resdef{bev-p-ALEC-U-8}{$p=1.000$}
\resdef{bev-p-ALEC-H-8}{$p=1.000$}
\resdef{bev-p-RANDOM-9}{$p=1.000$}
\resdef{bev-p-HARKER-PG-9}{$p=1.000$}
\resdef{bev-p-HARKER-9}{$p=1.000$}
\resdef{bev-p-HARKER-LOG-9}{$p=1.000$}
\resdef{bev-p-AOPT-9}{$p=1.000$}
\resdef{bev-p-FROZEN-9}{$p=1.000$}
\resdef{bev-p-STAR-9}{$p=1.000$}
\resdef{bev-p-ALEC-U-9}{$p=1.000$}
\resdef{bev-p-ALEC-H-9}{$p=1.000$}
\resdef{bev-p-RANDOM-10}{$p=1.000$}
\resdef{bev-p-HARKER-PG-10}{$p=1.000$}
\resdef{bev-p-HARKER-10}{$p=1.000$}
\resdef{bev-p-HARKER-LOG-10}{$p=1.000$}
\resdef{bev-p-AOPT-10}{$p=1.000$}
\resdef{bev-p-FROZEN-10}{$p=1.000$}
\resdef{bev-p-STAR-10}{$p=1.000$}
\resdef{bev-p-ALEC-U-10}{$p=1.000$}
\resdef{bev-p-ALEC-H-10}{$p=1.000$}
\resdef{R-ref}{500}
\resdef{esmdiag-ALEC-fails}{37}
\resdef{esmdiag-ALEC-fails-F}{36}
\resdef{esmdiag-ALEC-degC-all}{9.5}
\resdef{esmdiag-ALEC-degF-all}{5.2}
\resdef{esmdiag-ALEC-degC-fail}{9.5}
\resdef{esmdiag-ALEC-degF-fail}{4.6}
\resdef{esmdiag-ALEC-U-fails}{4}
\resdef{esmdiag-ALEC-U-fails-F}{4}
\resdef{esmdiag-ALEC-U-degC-all}{7.2}
\resdef{esmdiag-ALEC-U-degF-all}{6.5}
\resdef{esmdiag-ALEC-U-degC-fail}{3.8}
\resdef{esmdiag-ALEC-U-degF-fail}{6.8}
\resdef{esmdiag-R}{500}
\resdef{alec2n-worst-ci}{[$-3.8$, $-0.8$]}
\resdef{alec2n-worst-n}{15}
\resdef{alec2n-worst-sigma}{0.5}
\resdef{alec2n-regret-max}{3.51}
\resdef{full-regret-max}{3.34}
\resdef{alec2n-ratio-min}{96.5}
\resdef{alec2n-diff-min}{-2.3}
\resdef{alec2n-diff-max}{0.4}
\resdef{alec2n-diff-min-abs}{2.3}
\resdef{alec-regret-n-max}{5.7}
\resdef{gapz-n7-b0-ci-h}{[$-0.5$, $+3.3$]}
\resdef{gapz-n7-b0-ci-n}{[$-2.0$, $+1.1$]}
\resdef{gapz-n7-b0-N}{1205}
\resdef{gapz-n7-b0-full}{60.5}
\resdef{gapz-n7-b0-diff-h}{1.4}
\resdef{gapz-n7-b0-diff-n}{-0.4}
\resdef{gapz-n7-b0-RANDOM}{55.4}
\resdef{gapz-n7-b0-HARKER-LOG}{58.0}
\resdef{gapz-n7-b0-FROZEN}{59.4}
\resdef{gapz-n7-b0-ALEC}{59.0}
\resdef{gapz-n7-b0-ALEC-U}{58.3}
\resdef{gapz-n7-b1-ci-h}{[$-2.5$, $+1.9$]}
\resdef{gapz-n7-b1-ci-n}{[$-0.5$, $+2.7$]}
\resdef{gapz-n7-b1-N}{657}
\resdef{gapz-n7-b1-full}{89.0}
\resdef{gapz-n7-b1-diff-h}{-0.3}
\resdef{gapz-n7-b1-diff-n}{1.1}
\resdef{gapz-n7-b1-RANDOM}{78.8}
\resdef{gapz-n7-b1-HARKER-LOG}{85.8}
\resdef{gapz-n7-b1-FROZEN}{87.1}
\resdef{gapz-n7-b1-ALEC}{88.1}
\resdef{gapz-n7-b1-ALEC-U}{86.9}
\resdef{gapz-n7-b2-ci-h}{[$-1.9$, $+0.6$]}
\resdef{gapz-n7-b2-ci-n}{[$-0.4$, $+1.7$]}
\resdef{gapz-n7-b2-N}{668}
\resdef{gapz-n7-b2-full}{99.3}
\resdef{gapz-n7-b2-diff-h}{-0.6}
\resdef{gapz-n7-b2-diff-n}{0.6}
\resdef{gapz-n7-b2-RANDOM}{95.2}
\resdef{gapz-n7-b2-HARKER-LOG}{97.6}
\resdef{gapz-n7-b2-FROZEN}{98.4}
\resdef{gapz-n7-b2-ALEC}{99.0}
\resdef{gapz-n7-b2-ALEC-U}{97.6}
\resdef{gapz-n7-b3-ci-h}{[$-1.0$, $+1.0$]}
\resdef{gapz-n7-b3-ci-n}{[$-0.8$, $+0.8$]}
\resdef{gapz-n7-b3-N}{470}
\resdef{gapz-n7-b3-full}{100.0}
\resdef{gapz-n7-b3-diff-h}{0.0}
\resdef{gapz-n7-b3-diff-n}{0.0}
\resdef{gapz-n7-b3-RANDOM}{99.4}
\resdef{gapz-n7-b3-HARKER-LOG}{100.0}
\resdef{gapz-n7-b3-FROZEN}{100.0}
\resdef{gapz-n7-b3-ALEC}{100.0}
\resdef{gapz-n7-b3-ALEC-U}{100.0}
\resdef{gapz-n11-b0-ci-h}{[$+2.2$, $+6.1$]}
\resdef{gapz-n11-b0-ci-n}{[$-0.7$, $+3.3$]}
\resdef{gapz-n11-b0-N}{1613}
\resdef{gapz-n11-b0-full}{61.7}
\resdef{gapz-n11-b0-diff-h}{4.2}
\resdef{gapz-n11-b0-diff-n}{1.3}
\resdef{gapz-n11-b0-RANDOM}{49.2}
\resdef{gapz-n11-b0-HARKER-LOG}{57.2}
\resdef{gapz-n11-b0-FROZEN}{55.4}
\resdef{gapz-n11-b0-ALEC}{56.7}
\resdef{gapz-n11-b0-ALEC-U}{57.0}
\resdef{gapz-n11-b1-ci-h}{[$+2.4$, $+7.5$]}
\resdef{gapz-n11-b1-ci-n}{[$-0.7$, $+4.0$]}
\resdef{gapz-n11-b1-N}{671}
\resdef{gapz-n11-b1-full}{92.7}
\resdef{gapz-n11-b1-diff-h}{4.9}
\resdef{gapz-n11-b1-diff-n}{1.6}
\resdef{gapz-n11-b1-RANDOM}{75.0}
\resdef{gapz-n11-b1-HARKER-LOG}{88.7}
\resdef{gapz-n11-b1-FROZEN}{85.1}
\resdef{gapz-n11-b1-ALEC}{86.7}
\resdef{gapz-n11-b1-ALEC-U}{88.5}
\resdef{gapz-n11-b2-ci-h}{[$+0.3$, $+4.6$]}
\resdef{gapz-n11-b2-ci-n}{[$-1.3$, $+1.3$]}
\resdef{gapz-n11-b2-N}{474}
\resdef{gapz-n11-b2-full}{99.8}
\resdef{gapz-n11-b2-diff-h}{2.3}
\resdef{gapz-n11-b2-diff-n}{0.0}
\resdef{gapz-n11-b2-RANDOM}{95.1}
\resdef{gapz-n11-b2-HARKER-LOG}{98.5}
\resdef{gapz-n11-b2-FROZEN}{98.7}
\resdef{gapz-n11-b2-ALEC}{98.7}
\resdef{gapz-n11-b2-ALEC-U}{98.5}
\resdef{gapz-n11-b3-ci-h}{[$-1.7$, $+1.7$]}
\resdef{gapz-n11-b3-ci-n}{[$-1.2$, $+2.3$]}
\resdef{gapz-n11-b3-N}{242}
\resdef{gapz-n11-b3-full}{100.0}
\resdef{gapz-n11-b3-diff-h}{0.0}
\resdef{gapz-n11-b3-diff-n}{0.4}
\resdef{gapz-n11-b3-RANDOM}{99.6}
\resdef{gapz-n11-b3-HARKER-LOG}{100.0}
\resdef{gapz-n11-b3-FROZEN}{99.6}
\resdef{gapz-n11-b3-ALEC}{100.0}
\resdef{gapz-n11-b3-ALEC-U}{100.0}
\resdef{gapz-n15-b0-ci-h}{[$+1.6$, $+5.6$]}
\resdef{gapz-n15-b0-ci-n}{[$+0.9$, $+5.3$]}
\resdef{gapz-n15-b0-N}{1919}
\resdef{gapz-n15-b0-full}{64.5}
\resdef{gapz-n15-b0-diff-h}{3.6}
\resdef{gapz-n15-b0-diff-n}{3.1}
\resdef{gapz-n15-b0-RANDOM}{47.5}
\resdef{gapz-n15-b0-HARKER-LOG}{55.8}
\resdef{gapz-n15-b0-FROZEN}{52.7}
\resdef{gapz-n15-b0-ALEC}{55.8}
\resdef{gapz-n15-b0-ALEC-U}{57.2}
\resdef{gapz-n15-b1-ci-h}{[$+4.1$, $+10.3$]}
\resdef{gapz-n15-b1-ci-n}{[$+2.5$, $+8.3$]}
\resdef{gapz-n15-b1-N}{612}
\resdef{gapz-n15-b1-full}{95.4}
\resdef{gapz-n15-b1-diff-h}{7.2}
\resdef{gapz-n15-b1-diff-n}{5.4}
\resdef{gapz-n15-b1-RANDOM}{73.0}
\resdef{gapz-n15-b1-HARKER-LOG}{86.6}
\resdef{gapz-n15-b1-FROZEN}{80.7}
\resdef{gapz-n15-b1-ALEC}{86.1}
\resdef{gapz-n15-b1-ALEC-U}{87.7}
\resdef{gapz-n15-b2-ci-h}{[$+1.5$, $+6.9$]}
\resdef{gapz-n15-b2-ci-n}{[$+0.1$, $+4.8$]}
\resdef{gapz-n15-b2-N}{347}
\resdef{gapz-n15-b2-full}{100.0}
\resdef{gapz-n15-b2-diff-h}{4.0}
\resdef{gapz-n15-b2-diff-n}{2.3}
\resdef{gapz-n15-b2-RANDOM}{90.5}
\resdef{gapz-n15-b2-HARKER-LOG}{98.6}
\resdef{gapz-n15-b2-FROZEN}{96.8}
\resdef{gapz-n15-b2-ALEC}{99.1}
\resdef{gapz-n15-b2-ALEC-U}{98.8}
\resdef{gapz-n15-b3-ci-h}{[$-3.1$, $+3.1$]}
\resdef{gapz-n15-b3-ci-n}{[$-3.1$, $+3.1$]}
\resdef{gapz-n15-b3-N}{122}
\resdef{gapz-n15-b3-full}{100.0}
\resdef{gapz-n15-b3-diff-h}{0.0}
\resdef{gapz-n15-b3-diff-n}{0.0}
\resdef{gapz-n15-b3-RANDOM}{100.0}
\resdef{gapz-n15-b3-HARKER-LOG}{100.0}
\resdef{gapz-n15-b3-FROZEN}{100.0}
\resdef{gapz-n15-b3-ALEC}{100.0}
\resdef{gapz-n15-b3-ALEC-U}{100.0}
\resdef{gap-expected-n15}{0.137}
\resdef{gap-expected-n7}{0.275}
\resdef{n7-s015-tol-save-add}{20.7}
\resdef{n7-s015-tol-save-tot}{9.9}
\resdef{n7-s015-tol-regdiff}{0.00}
\resdef{n7-s015-stop-LC-capped}{0.0}
\resdef{n7-s015-stop-B10-ALEC-capped}{0.0}
\resdef{n7-s015-stop-B05-ALECU-capped}{0.0}
\resdef{n7-s015-stop-B10-ALECU-capped}{0.0}
\resdef{n7-s015-stop-B20-ALECU-capped}{0.0}
\resdef{n7-s015-stop-B10-ALECU-s06-capped}{0.0}
\resdef{n7-s015-stop-B10-ALECU-cap-capped}{0.0}
\resdef{n7-s015-stop-B10-ALECU-oracle-capped}{0.0}
\resdef{n7-s015-stop-E02-ALEC-capped}{0.0}
\resdef{n7-s015-stop-E05-ALEC-capped}{0.0}
\resdef{n7-s015-stop-B10-ALEC-cap-capped}{0.0}
\resdef{n7-s015-stop-E02-ALEC-cap-capped}{0.0}
\resdef{n7-s015-stop-E05-ALEC-cap-capped}{0.0}
\resdef{n7-s015-stop-E05-ALEC-cap-s06-capped}{0.0}
\resdef{n7-s030-tol-save-add}{14.8}
\resdef{n7-s030-tol-save-tot}{7.9}
\resdef{n7-s030-tol-regdiff}{0.01}
\resdef{n7-s030-stop-LC-capped}{0.0}
\resdef{n7-s030-stop-B10-ALEC-capped}{0.0}
\resdef{n7-s030-stop-B05-ALECU-capped}{0.0}
\resdef{n7-s030-stop-B10-ALECU-capped}{0.0}
\resdef{n7-s030-stop-B20-ALECU-capped}{0.0}
\resdef{n7-s030-stop-B10-ALECU-s06-capped}{0.0}
\resdef{n7-s030-stop-B10-ALECU-cap-capped}{0.0}
\resdef{n7-s030-stop-B10-ALECU-oracle-capped}{0.0}
\resdef{n7-s030-stop-E02-ALEC-capped}{0.0}
\resdef{n7-s030-stop-E05-ALEC-capped}{0.0}
\resdef{n7-s030-stop-B10-ALEC-cap-capped}{0.0}
\resdef{n7-s030-stop-E02-ALEC-cap-capped}{0.0}
\resdef{n7-s030-stop-E05-ALEC-cap-capped}{0.0}
\resdef{n7-s030-stop-E05-ALEC-cap-s06-capped}{0.0}
\resdef{n7-s050-tol-save-add}{11.2}
\resdef{n7-s050-tol-save-tot}{6.3}
\resdef{n7-s050-tol-regdiff}{0.22}
\resdef{n7-s050-stop-LC-capped}{0.0}
\resdef{n7-s050-stop-B10-ALEC-capped}{0.0}
\resdef{n7-s050-stop-B05-ALECU-capped}{0.0}
\resdef{n7-s050-stop-B10-ALECU-capped}{0.0}
\resdef{n7-s050-stop-B20-ALECU-capped}{0.0}
\resdef{n7-s050-stop-B10-ALECU-s06-capped}{0.0}
\resdef{n7-s050-stop-B10-ALECU-cap-capped}{0.0}
\resdef{n7-s050-stop-B10-ALECU-oracle-capped}{0.0}
\resdef{n7-s050-stop-E02-ALEC-capped}{0.0}
\resdef{n7-s050-stop-E05-ALEC-capped}{0.0}
\resdef{n7-s050-stop-B10-ALEC-cap-capped}{0.0}
\resdef{n7-s050-stop-E02-ALEC-cap-capped}{0.0}
\resdef{n7-s050-stop-E05-ALEC-cap-capped}{0.0}
\resdef{n7-s050-stop-E05-ALEC-cap-s06-capped}{0.0}
\resdef{n11-s015-tol-save-add}{40.0}
\resdef{n11-s015-tol-save-tot}{25.3}
\resdef{n11-s015-tol-regdiff}{0.02}
\resdef{n11-s015-stop-LC-capped}{0.0}
\resdef{n11-s015-stop-B10-ALEC-capped}{0.0}
\resdef{n11-s015-stop-B05-ALECU-capped}{0.0}
\resdef{n11-s015-stop-B10-ALECU-capped}{0.0}
\resdef{n11-s015-stop-B20-ALECU-capped}{0.0}
\resdef{n11-s015-stop-B10-ALECU-s06-capped}{0.0}
\resdef{n11-s015-stop-B10-ALECU-cap-capped}{31.8}
\resdef{n11-s015-stop-B10-ALECU-oracle-capped}{0.0}
\resdef{n11-s015-stop-E02-ALEC-capped}{0.0}
\resdef{n11-s015-stop-E05-ALEC-capped}{0.0}
\resdef{n11-s015-stop-B10-ALEC-cap-capped}{33.1}
\resdef{n11-s015-stop-E02-ALEC-cap-capped}{25.2}
\resdef{n11-s015-stop-E05-ALEC-cap-capped}{12.5}
\resdef{n11-s015-stop-E05-ALEC-cap-s06-capped}{31.6}
\resdef{n11-s030-tol-save-add}{21.2}
\resdef{n11-s030-tol-save-tot}{14.8}
\resdef{n11-s030-tol-regdiff}{-0.01}
\resdef{n11-s030-stop-LC-capped}{0.0}
\resdef{n11-s030-stop-B10-ALEC-capped}{0.0}
\resdef{n11-s030-stop-B05-ALECU-capped}{0.0}
\resdef{n11-s030-stop-B10-ALECU-capped}{0.0}
\resdef{n11-s030-stop-B20-ALECU-capped}{0.0}
\resdef{n11-s030-stop-B10-ALECU-s06-capped}{0.0}
\resdef{n11-s030-stop-B10-ALECU-cap-capped}{46.6}
\resdef{n11-s030-stop-B10-ALECU-oracle-capped}{0.0}
\resdef{n11-s030-stop-E02-ALEC-capped}{0.0}
\resdef{n11-s030-stop-E05-ALEC-capped}{0.0}
\resdef{n11-s030-stop-B10-ALEC-cap-capped}{50.4}
\resdef{n11-s030-stop-E02-ALEC-cap-capped}{44.0}
\resdef{n11-s030-stop-E05-ALEC-cap-capped}{36.0}
\resdef{n11-s030-stop-E05-ALEC-cap-s06-capped}{49.2}
\resdef{n11-s050-tol-save-add}{13.8}
\resdef{n11-s050-tol-save-tot}{10.0}
\resdef{n11-s050-tol-regdiff}{0.09}
\resdef{n11-s050-stop-LC-capped}{0.0}
\resdef{n11-s050-stop-B10-ALEC-capped}{0.0}
\resdef{n11-s050-stop-B05-ALECU-capped}{0.0}
\resdef{n11-s050-stop-B10-ALECU-capped}{0.0}
\resdef{n11-s050-stop-B20-ALECU-capped}{0.0}
\resdef{n11-s050-stop-B10-ALECU-s06-capped}{0.0}
\resdef{n11-s050-stop-B10-ALECU-cap-capped}{53.5}
\resdef{n11-s050-stop-B10-ALECU-oracle-capped}{0.0}
\resdef{n11-s050-stop-E02-ALEC-capped}{0.0}
\resdef{n11-s050-stop-E05-ALEC-capped}{0.0}
\resdef{n11-s050-stop-B10-ALEC-cap-capped}{60.4}
\resdef{n11-s050-stop-E02-ALEC-cap-capped}{56.1}
\resdef{n11-s050-stop-E05-ALEC-cap-capped}{49.6}
\resdef{n11-s050-stop-E05-ALEC-cap-s06-capped}{64.4}
\resdef{n15-s015-tol-save-add}{51.3}
\resdef{n15-s015-tol-save-tot}{37.6}
\resdef{n15-s015-tol-regdiff}{0.03}
\resdef{n15-s015-stop-LC-capped}{0.0}
\resdef{n15-s015-stop-B10-ALEC-capped}{0.0}
\resdef{n15-s015-stop-B05-ALECU-capped}{0.0}
\resdef{n15-s015-stop-B10-ALECU-capped}{0.0}
\resdef{n15-s015-stop-B20-ALECU-capped}{0.0}
\resdef{n15-s015-stop-B10-ALECU-s06-capped}{0.0}
\resdef{n15-s015-stop-B10-ALECU-cap-capped}{36.9}
\resdef{n15-s015-stop-B10-ALECU-oracle-capped}{0.0}
\resdef{n15-s015-stop-E02-ALEC-capped}{0.0}
\resdef{n15-s015-stop-E05-ALEC-capped}{0.0}
\resdef{n15-s015-stop-B10-ALEC-cap-capped}{39.1}
\resdef{n15-s015-stop-E02-ALEC-cap-capped}{27.8}
\resdef{n15-s015-stop-E05-ALEC-cap-capped}{13.0}
\resdef{n15-s015-stop-E05-ALEC-cap-s06-capped}{29.6}
\resdef{n15-s030-tol-save-add}{25.4}
\resdef{n15-s030-tol-save-tot}{20.0}
\resdef{n15-s030-tol-regdiff}{0.07}
\resdef{n15-s030-stop-LC-capped}{0.0}
\resdef{n15-s030-stop-B10-ALEC-capped}{0.0}
\resdef{n15-s030-stop-B05-ALECU-capped}{0.0}
\resdef{n15-s030-stop-B10-ALECU-capped}{0.0}
\resdef{n15-s030-stop-B20-ALECU-capped}{0.0}
\resdef{n15-s030-stop-B10-ALECU-s06-capped}{0.0}
\resdef{n15-s030-stop-B10-ALECU-cap-capped}{50.9}
\resdef{n15-s030-stop-B10-ALECU-oracle-capped}{0.0}
\resdef{n15-s030-stop-E02-ALEC-capped}{0.0}
\resdef{n15-s030-stop-E05-ALEC-capped}{0.0}
\resdef{n15-s030-stop-B10-ALEC-cap-capped}{57.2}
\resdef{n15-s030-stop-E02-ALEC-cap-capped}{51.0}
\resdef{n15-s030-stop-E05-ALEC-cap-capped}{40.6}
\resdef{n15-s030-stop-E05-ALEC-cap-s06-capped}{51.2}
\resdef{n15-s050-tol-save-add}{13.5}
\resdef{n15-s050-tol-save-tot}{10.9}
\resdef{n15-s050-tol-regdiff}{0.15}
\resdef{n15-s050-stop-LC-capped}{0.0}
\resdef{n15-s050-stop-B10-ALEC-capped}{0.0}
\resdef{n15-s050-stop-B05-ALECU-capped}{0.0}
\resdef{n15-s050-stop-B10-ALECU-capped}{0.0}
\resdef{n15-s050-stop-B20-ALECU-capped}{0.0}
\resdef{n15-s050-stop-B10-ALECU-s06-capped}{0.0}
\resdef{n15-s050-stop-B10-ALECU-cap-capped}{62.6}
\resdef{n15-s050-stop-B10-ALECU-oracle-capped}{0.0}
\resdef{n15-s050-stop-E02-ALEC-capped}{0.0}
\resdef{n15-s050-stop-E05-ALEC-capped}{0.0}
\resdef{n15-s050-stop-B10-ALEC-cap-capped}{70.4}
\resdef{n15-s050-stop-E02-ALEC-cap-capped}{66.5}
\resdef{n15-s050-stop-E05-ALEC-cap-capped}{60.4}
\resdef{n15-s050-stop-E05-ALEC-cap-s06-capped}{71.7}
\resdef{errearly-E05-ALEC-cap-s015-min}{1.7}
\resdef{errearly-E05-ALEC-cap-s015-max}{2.1}
\resdef{errearly-E05-ALEC-cap-s06-s015-min}{0.1}
\resdef{errearly-E05-ALEC-cap-s06-s015-max}{0.4}
\resdef{errearly-B10-ALECU-s06-s015-min}{1.0}
\resdef{errearly-B10-ALECU-s06-s015-max}{1.6}
\resdef{errearly-B10-ALECU-s015-min}{2.7}
\resdef{errearly-B10-ALECU-s015-max}{4.1}
\resdef{errearly-E05-ALEC-cap-s030-min}{4.2}
\resdef{errearly-E05-ALEC-cap-s030-max}{7.2}
\resdef{errearly-E05-ALEC-cap-s06-s030-min}{1.8}
\resdef{errearly-E05-ALEC-cap-s06-s030-max}{3.5}
\resdef{errearly-B10-ALECU-s06-s030-min}{3.4}
\resdef{errearly-B10-ALECU-s06-s030-max}{5.9}
\resdef{errearly-B10-ALECU-s030-min}{7.3}
\resdef{errearly-B10-ALECU-s030-max}{10.6}
\resdef{errearly-E05-ALEC-cap-s050-min}{17.1}
\resdef{errearly-E05-ALEC-cap-s050-max}{17.7}
\resdef{errearly-E05-ALEC-cap-s06-s050-min}{7.6}
\resdef{errearly-E05-ALEC-cap-s06-s050-max}{9.5}
\resdef{errearly-B10-ALECU-s06-s050-min}{11.1}
\resdef{errearly-B10-ALECU-s06-s050-max}{12.3}
\resdef{errearly-B10-ALECU-s050-min}{21.3}
\resdef{errearly-B10-ALECU-s050-max}{21.4}
\resdef{s06-extra-q-min}{3.5}
\resdef{s06-extra-q-max}{6.2}
\resdef{s06-extra-qrel-min}{16}
\resdef{s06-extra-qrel-max}{58}
\resdef{s06-errearly-max}{9.5}
\resdef{s06-errci-max}{13.5}
\resdef{prot-errearly-max-s03}{24.3}
\resdef{prot-save2n-min}{24}
\resdef{prot-save2n-max}{60}
\resdef{eps-s050-errearly-min}{16.4}
\resdef{eps-s050-errearly-max}{20.9}
\resdef{prot-s03-saveadd-min}{24}
\resdef{prot-s03-saveadd-max}{60}
\resdef{prot-s03-savetot-min}{16}
\resdef{prot-s03-savetot-max}{40}
\resdef{prot-s03-stop-min}{40}
\resdef{prot-s03-stop-max}{88}
\resdef{prot-s06-saveadd-min}{11}
\resdef{prot-s06-saveadd-max}{36}
\resdef{prot-s06-savetot-min}{8}
\resdef{prot-s06-savetot-max}{24}
\resdef{prot-s06-stop-min}{28}
\resdef{prot-s06-stop-max}{70}
\resdef{n11-s015-fix-s06-m}{15}
\resdef{n11-s015-fix-s06-dhit}{$+1.5$}
\resdef{n11-s015-fix-s06-dhitci}{[$+0.5$, $+2.6$]}
\resdef{n11-s015-fix-s06-hitp}{88.8}
\resdef{n11-s015-fix-s06-hitf}{87.3}
\resdef{n11-s015-fix-s06-regp}{0.43}
\resdef{n11-s015-fix-s06-regf}{0.52}
\resdef{n11-s015-fix-s06-dreg}{$-0.09$}
\resdef{n11-s015-fix-s06-dregci}{[$-0.15$, $-0.03$]}
\resdef{n11-s030-fix-s06-m}{17}
\resdef{n11-s030-fix-s06-dhit}{$+0.2$}
\resdef{n11-s030-fix-s06-dhitci}{[$-0.8$, $+1.2$]}
\resdef{n11-s030-fix-s06-hitp}{76.6}
\resdef{n11-s030-fix-s06-hitf}{76.4}
\resdef{n11-s030-fix-s06-regp}{1.59}
\resdef{n11-s030-fix-s06-regf}{1.60}
\resdef{n11-s030-fix-s06-dreg}{$-0.02$}
\resdef{n11-s030-fix-s06-dregci}{[$-0.13$, $+0.10$]}
\resdef{n11-s050-fix-s06-m}{19}
\resdef{n11-s050-fix-s06-dhit}{$+0.2$}
\resdef{n11-s050-fix-s06-dhitci}{[$-1.0$, $+1.4$]}
\resdef{n11-s050-fix-s06-hitp}{67.7}
\resdef{n11-s050-fix-s06-hitf}{67.5}
\resdef{n11-s050-fix-s06-regp}{3.47}
\resdef{n11-s050-fix-s06-regf}{3.59}
\resdef{n11-s050-fix-s06-dreg}{$-0.12$}
\resdef{n11-s050-fix-s06-dregci}{[$-0.35$, $+0.10$]}
\resdef{n15-s015-fix-s06-m}{21}
\resdef{n15-s015-fix-s06-dhit}{$+2.0$}
\resdef{n15-s015-fix-s06-dhitci}{[$+0.7$, $+3.4$]}
\resdef{n15-s015-fix-s06-hitp}{86.0}
\resdef{n15-s015-fix-s06-hitf}{84.0}
\resdef{n15-s015-fix-s06-regp}{0.42}
\resdef{n15-s015-fix-s06-regf}{0.55}
\resdef{n15-s015-fix-s06-dreg}{$-0.12$}
\resdef{n15-s015-fix-s06-dregci}{[$-0.21$, $-0.04$]}
\resdef{n15-s030-fix-s06-m}{24}
\resdef{n15-s030-fix-s06-dhit}{$-0.2$}
\resdef{n15-s030-fix-s06-dhitci}{[$-1.7$, $+1.3$]}
\resdef{n15-s030-fix-s06-hitp}{74.9}
\resdef{n15-s030-fix-s06-hitf}{75.1}
\resdef{n15-s030-fix-s06-regp}{1.71}
\resdef{n15-s030-fix-s06-regf}{1.65}
\resdef{n15-s030-fix-s06-dreg}{$+0.06$}
\resdef{n15-s030-fix-s06-dregci}{[$-0.10$, $+0.22$]}
\resdef{n15-s050-fix-s06-m}{27}
\resdef{n15-s050-fix-s06-dhit}{$+0.5$}
\resdef{n15-s050-fix-s06-dhitci}{[$-0.9$, $+1.9$]}
\resdef{n15-s050-fix-s06-hitp}{64.9}
\resdef{n15-s050-fix-s06-hitf}{64.4}
\resdef{n15-s050-fix-s06-regp}{3.57}
\resdef{n15-s050-fix-s06-regf}{3.72}
\resdef{n15-s050-fix-s06-dreg}{$-0.15$}
\resdef{n15-s050-fix-s06-dregci}{[$-0.41$, $+0.12$]}
\resdef{n11-s015-fix-s03-m}{9}
\resdef{n11-s015-fix-s03-dhit}{$+3.3$}
\resdef{n11-s015-fix-s03-dhitci}{[$+1.5$, $+5.2$]}
\resdef{n11-s015-fix-s03-hitp}{88.3}
\resdef{n11-s015-fix-s03-hitf}{85.0}
\resdef{n11-s015-fix-s03-regp}{0.45}
\resdef{n11-s015-fix-s03-regf}{0.76}
\resdef{n11-s015-fix-s03-dreg}{$-0.31$}
\resdef{n11-s015-fix-s03-dregci}{[$-0.45$, $-0.17$]}
\resdef{n11-s030-fix-s03-m}{13}
\resdef{n11-s030-fix-s03-dhit}{$+2.3$}
\resdef{n11-s030-fix-s03-dhitci}{[$+0.4$, $+4.2$]}
\resdef{n11-s030-fix-s03-hitp}{77.6}
\resdef{n11-s030-fix-s03-hitf}{75.3}
\resdef{n11-s030-fix-s03-regp}{1.54}
\resdef{n11-s030-fix-s03-regf}{1.88}
\resdef{n11-s030-fix-s03-dreg}{$-0.33$}
\resdef{n11-s030-fix-s03-dregci}{[$-0.56$, $-0.11$]}
\resdef{n11-s050-fix-s03-m}{15}
\resdef{n11-s050-fix-s03-dhit}{$+1.9$}
\resdef{n11-s050-fix-s03-dhitci}{[$-0.2$, $+4.0$]}
\resdef{n11-s050-fix-s03-hitp}{66.6}
\resdef{n11-s050-fix-s03-hitf}{64.7}
\resdef{n11-s050-fix-s03-regp}{3.94}
\resdef{n11-s050-fix-s03-regf}{4.23}
\resdef{n11-s050-fix-s03-dreg}{$-0.29$}
\resdef{n11-s050-fix-s03-dregci}{[$-0.76$, $+0.18$]}
\resdef{n15-s015-fix-s03-m}{14}
\resdef{n15-s015-fix-s03-dhit}{$+4.7$}
\resdef{n15-s015-fix-s03-dhitci}{[$+2.6$, $+6.9$]}
\resdef{n15-s015-fix-s03-hitp}{85.3}
\resdef{n15-s015-fix-s03-hitf}{80.6}
\resdef{n15-s015-fix-s03-regp}{0.46}
\resdef{n15-s015-fix-s03-regf}{0.88}
\resdef{n15-s015-fix-s03-dreg}{$-0.42$}
\resdef{n15-s015-fix-s03-dregci}{[$-0.57$, $-0.28$]}
\resdef{n15-s030-fix-s03-m}{20}
\resdef{n15-s030-fix-s03-dhit}{$+0.8$}
\resdef{n15-s030-fix-s03-dhitci}{[$-1.3$, $+2.9$]}
\resdef{n15-s030-fix-s03-hitp}{74.2}
\resdef{n15-s030-fix-s03-hitf}{73.4}
\resdef{n15-s030-fix-s03-regp}{1.84}
\resdef{n15-s030-fix-s03-regf}{1.93}
\resdef{n15-s030-fix-s03-dreg}{$-0.09$}
\resdef{n15-s030-fix-s03-dregci}{[$-0.34$, $+0.16$]}
\resdef{n15-s050-fix-s03-m}{23}
\resdef{n15-s050-fix-s03-dhit}{$+1.2$}
\resdef{n15-s050-fix-s03-dhitci}{[$-1.1$, $+3.5$]}
\resdef{n15-s050-fix-s03-hitp}{63.0}
\resdef{n15-s050-fix-s03-hitf}{61.8}
\resdef{n15-s050-fix-s03-regp}{3.99}
\resdef{n15-s050-fix-s03-regf}{4.16}
\resdef{n15-s050-fix-s03-dreg}{$-0.17$}
\resdef{n15-s050-fix-s03-dregci}{[$-0.58$, $+0.24$]}
\resdef{fix-s03-hitpos}{6}
\resdef{fix-s06-hitpos}{5}
\resdef{fix-hit-holm}{3}
\resdef{fix-s03-dhit-min}{$+0.8$}
\resdef{fix-s03-dhit-max}{$+4.7$}
\resdef{fix-s03-hit-holm}{2}
\resdef{fix-s03-reg-holm}{3}
\resdef{fix-s06-dhit-min}{$-0.2$}
\resdef{fix-s06-dhit-max}{$+2.0$}
\resdef{fix-s06-hit-holm}{1}
\resdef{fix-s06-reg-holm}{0}
\resdef{fix-reg-holm}{3}
\resdef{fix-rows}{12}
\resdef{fix-hitci-pos}{5}
\resdef{fix-regci-neg}{5}
\resdef{fix-dhit-min}{$-0.2$}
\resdef{fix-hitpoint-pos}{11}
\resdef{fix-regpoint-neg}{11}
\resdef{fix-dhit-max}{$+4.7$}
\resdef{ss-holm-q}{4}
\resdef{ss-holm-hit}{2}
\resdef{ss-holm-reg}{1}
\resdef{ss-holm-exc}{1}
\resdef{ss-holm-stop}{1}
\resdef{n11-s015-ss-alec-q}{14.03}
\resdef{n11-s015-ss-alec-hit}{88.8}
\resdef{n11-s015-ss-alec-reg}{0.43}
\resdef{n11-s015-ss-alec-exc}{3.2}
\resdef{n11-s015-ss-alec-stop}{68.4}
\resdef{n11-s015-ss-alec-err}{0.4}
\resdef{n11-s015-ss-alec-errci}{[0.1, 1.3]}
\resdef{n11-s015-ss-alec-ne}{684}
\resdef{n11-s015-ss-hlog-q}{13.93}
\resdef{n11-s015-ss-hlog-hit}{86.9}
\resdef{n11-s015-ss-hlog-reg}{0.51}
\resdef{n11-s015-ss-hlog-exc}{3.6}
\resdef{n11-s015-ss-hlog-stop}{68.1}
\resdef{n11-s015-ss-hlog-err}{0.3}
\resdef{n11-s015-ss-hlog-errci}{[0.1, 1.1]}
\resdef{n11-s015-ss-hlog-ne}{681}
\resdef{n11-s030-ss-alec-q}{16.38}
\resdef{n11-s030-ss-alec-hit}{76.6}
\resdef{n11-s030-ss-alec-reg}{1.59}
\resdef{n11-s030-ss-alec-exc}{11.6}
\resdef{n11-s030-ss-alec-stop}{50.8}
\resdef{n11-s030-ss-alec-err}{1.8}
\resdef{n11-s030-ss-alec-errci}{[0.9, 3.3]}
\resdef{n11-s030-ss-alec-ne}{508}
\resdef{n11-s030-ss-hlog-q}{15.85}
\resdef{n11-s030-ss-hlog-hit}{76.5}
\resdef{n11-s030-ss-hlog-reg}{1.62}
\resdef{n11-s030-ss-hlog-exc}{11.7}
\resdef{n11-s030-ss-hlog-stop}{51.7}
\resdef{n11-s030-ss-hlog-err}{2.5}
\resdef{n11-s030-ss-hlog-errci}{[1.5, 4.3]}
\resdef{n11-s030-ss-hlog-ne}{517}
\resdef{n11-s050-ss-alec-q}{18.19}
\resdef{n11-s050-ss-alec-hit}{67.7}
\resdef{n11-s050-ss-alec-reg}{3.47}
\resdef{n11-s050-ss-alec-exc}{22.9}
\resdef{n11-s050-ss-alec-stop}{35.6}
\resdef{n11-s050-ss-alec-err}{7.6}
\resdef{n11-s050-ss-alec-errci}{[5.3, 10.8]}
\resdef{n11-s050-ss-alec-ne}{356}
\resdef{n11-s050-ss-hlog-q}{17.48}
\resdef{n11-s050-ss-hlog-hit}{66.3}
\resdef{n11-s050-ss-hlog-reg}{3.78}
\resdef{n11-s050-ss-hlog-exc}{24.6}
\resdef{n11-s050-ss-hlog-stop}{38.3}
\resdef{n11-s050-ss-hlog-err}{10.7}
\resdef{n11-s050-ss-hlog-errci}{[8.0, 14.2]}
\resdef{n11-s050-ss-hlog-ne}{383}
\resdef{n15-s015-ss-alec-q}{20.22}
\resdef{n15-s015-ss-alec-hit}{86.0}
\resdef{n15-s015-ss-alec-reg}{0.42}
\resdef{n15-s015-ss-alec-exc}{3.0}
\resdef{n15-s015-ss-alec-stop}{70.4}
\resdef{n15-s015-ss-alec-err}{0.1}
\resdef{n15-s015-ss-alec-errci}{[0.0, 0.8]}
\resdef{n15-s015-ss-alec-ne}{704}
\resdef{n15-s015-ss-hlog-q}{20.24}
\resdef{n15-s015-ss-hlog-hit}{85.3}
\resdef{n15-s015-ss-hlog-reg}{0.46}
\resdef{n15-s015-ss-hlog-exc}{3.3}
\resdef{n15-s015-ss-hlog-stop}{64.9}
\resdef{n15-s015-ss-hlog-err}{0.2}
\resdef{n15-s015-ss-hlog-errci}{[0.0, 0.9]}
\resdef{n15-s015-ss-hlog-ne}{649}
\resdef{n15-s030-ss-alec-q}{23.73}
\resdef{n15-s030-ss-alec-hit}{74.9}
\resdef{n15-s030-ss-alec-reg}{1.71}
\resdef{n15-s030-ss-alec-exc}{14.2}
\resdef{n15-s030-ss-alec-stop}{48.8}
\resdef{n15-s030-ss-alec-err}{3.5}
\resdef{n15-s030-ss-alec-errci}{[2.2, 5.5]}
\resdef{n15-s030-ss-alec-ne}{488}
\resdef{n15-s030-ss-hlog-q}{22.77}
\resdef{n15-s030-ss-hlog-hit}{74.1}
\resdef{n15-s030-ss-hlog-reg}{1.88}
\resdef{n15-s030-ss-hlog-exc}{15.5}
\resdef{n15-s030-ss-hlog-stop}{48.9}
\resdef{n15-s030-ss-hlog-err}{5.5}
\resdef{n15-s030-ss-hlog-errci}{[3.8, 7.9]}
\resdef{n15-s030-ss-hlog-ne}{489}
\resdef{n15-s050-ss-alec-q}{26.57}
\resdef{n15-s050-ss-alec-hit}{64.9}
\resdef{n15-s050-ss-alec-reg}{3.57}
\resdef{n15-s050-ss-alec-exc}{23.3}
\resdef{n15-s050-ss-alec-stop}{28.3}
\resdef{n15-s050-ss-alec-err}{9.5}
\resdef{n15-s050-ss-alec-errci}{[6.6, 13.5]}
\resdef{n15-s050-ss-alec-ne}{283}
\resdef{n15-s050-ss-hlog-q}{25.44}
\resdef{n15-s050-ss-hlog-hit}{59.9}
\resdef{n15-s050-ss-hlog-reg}{4.28}
\resdef{n15-s050-ss-hlog-exc}{27.4}
\resdef{n15-s050-ss-hlog-stop}{31.0}
\resdef{n15-s050-ss-hlog-err}{14.5}
\resdef{n15-s050-ss-hlog-errci}{[11.0, 18.9]}
\resdef{n15-s050-ss-hlog-ne}{310}
\resdef{n11-s015-ss-dq}{$+0.10$}
\resdef{n11-s015-ss-dqci}{[$-0.04$, $+0.25$]}
\resdef{n11-s015-ss-dhit}{$+1.9$}
\resdef{n11-s015-ss-dhitci}{[$+0.8$, $+3.1$]}
\resdef{n11-s015-ss-dreg}{$-0.08$}
\resdef{n11-s015-ss-dregci}{[$-0.15$, $-0.01$]}
\resdef{n11-s015-ss-dexc}{$-0.4$}
\resdef{n11-s015-ss-dexcci}{[$-1.2$, $+0.4$]}
\resdef{n11-s015-ss-dstop}{$+0.3$}
\resdef{n11-s015-ss-dstopci}{[$-1.2$, $+1.8$]}
\resdef{n11-s030-ss-dq}{$+0.53$}
\resdef{n11-s030-ss-dqci}{[$+0.30$, $+0.76$]}
\resdef{n11-s030-ss-dhit}{$+0.1$}
\resdef{n11-s030-ss-dhitci}{[$-1.4$, $+1.6$]}
\resdef{n11-s030-ss-dreg}{$-0.03$}
\resdef{n11-s030-ss-dregci}{[$-0.21$, $+0.14$]}
\resdef{n11-s030-ss-dexc}{$-0.1$}
\resdef{n11-s030-ss-dexcci}{[$-1.5$, $+1.3$]}
\resdef{n11-s030-ss-dstop}{$-0.9$}
\resdef{n11-s030-ss-dstopci}{[$-2.9$, $+1.2$]}
\resdef{n11-s050-ss-dq}{$+0.71$}
\resdef{n11-s050-ss-dqci}{[$+0.40$, $+1.02$]}
\resdef{n11-s050-ss-dhit}{$+1.4$}
\resdef{n11-s050-ss-dhitci}{[$-0.5$, $+3.3$]}
\resdef{n11-s050-ss-dreg}{$-0.31$}
\resdef{n11-s050-ss-dregci}{[$-0.64$, $+0.02$]}
\resdef{n11-s050-ss-dexc}{$-1.7$}
\resdef{n11-s050-ss-dexcci}{[$-3.5$, $+0.1$]}
\resdef{n11-s050-ss-dstop}{$-2.7$}
\resdef{n11-s050-ss-dstopci}{[$-5.0$, $-0.4$]}
\resdef{n15-s015-ss-dq}{$-0.02$}
\resdef{n15-s015-ss-dqci}{[$-0.28$, $+0.24$]}
\resdef{n15-s015-ss-dhit}{$+0.7$}
\resdef{n15-s015-ss-dhitci}{[$-1.0$, $+2.4$]}
\resdef{n15-s015-ss-dreg}{$-0.04$}
\resdef{n15-s015-ss-dregci}{[$-0.14$, $+0.06$]}
\resdef{n15-s015-ss-dexc}{$-0.3$}
\resdef{n15-s015-ss-dexcci}{[$-1.4$, $+0.8$]}
\resdef{n15-s015-ss-dstop}{$+5.5$}
\resdef{n15-s015-ss-dstopci}{[$+3.4$, $+7.6$]}
\resdef{n15-s030-ss-dq}{$+0.95$}
\resdef{n15-s030-ss-dqci}{[$+0.59$, $+1.32$]}
\resdef{n15-s030-ss-dhit}{$+0.8$}
\resdef{n15-s030-ss-dhitci}{[$-1.3$, $+2.9$]}
\resdef{n15-s030-ss-dreg}{$-0.17$}
\resdef{n15-s030-ss-dregci}{[$-0.39$, $+0.05$]}
\resdef{n15-s030-ss-dexc}{$-1.3$}
\resdef{n15-s030-ss-dexcci}{[$-3.2$, $+0.5$]}
\resdef{n15-s030-ss-dstop}{$-0.1$}
\resdef{n15-s030-ss-dstopci}{[$-2.6$, $+2.4$]}
\resdef{n15-s050-ss-dq}{$+1.13$}
\resdef{n15-s050-ss-dqci}{[$+0.73$, $+1.53$]}
\resdef{n15-s050-ss-dhit}{$+5.0$}
\resdef{n15-s050-ss-dhitci}{[$+2.5$, $+7.5$]}
\resdef{n15-s050-ss-dreg}{$-0.71$}
\resdef{n15-s050-ss-dregci}{[$-1.14$, $-0.27$]}
\resdef{n15-s050-ss-dexc}{$-4.1$}
\resdef{n15-s050-ss-dexcci}{[$-6.3$, $-1.9$]}
\resdef{n15-s050-ss-dstop}{$-2.7$}
\resdef{n15-s050-ss-dstopci}{[$-5.2$, $-0.2$]}
\resdef{ss-dq-min}{$-0.02$}
\resdef{ss-dq-max}{$+1.13$}
\resdef{ss-dhit-min}{$+0.1$}
\resdef{ss-dhit-max}{$+5.0$}
\resdef{ss-dexc-min}{$-4.1$}
\resdef{ss-dexc-max}{$-0.1$}
\resdef{ss-hitpos}{6}
\resdef{ss-excneg}{6}
\resdef{n11-s030-rd-cont-q}{16.38}
\resdef{n11-s030-rd-cont-hit}{76.6}
\resdef{n11-s030-rd-cont-reg}{1.59}
\resdef{n11-s030-rd-cont-exc}{11.6}
\resdef{n11-s030-rd-cont-stop}{50.8}
\resdef{n11-s030-rd-cont-err}{1.8}
\resdef{n11-s030-rd-cont-errci}{[0.9, 3.3]}
\resdef{n11-s030-rd-round-q}{16.77}
\resdef{n11-s030-rd-round-hit}{75.2}
\resdef{n11-s030-rd-round-reg}{2.06}
\resdef{n11-s030-rd-round-exc}{13.6}
\resdef{n11-s030-rd-round-stop}{45.1}
\resdef{n11-s030-rd-round-err}{2.4}
\resdef{n11-s030-rd-round-errci}{[1.4, 4.3]}
\resdef{n11-s050-rd-cont-q}{18.19}
\resdef{n11-s050-rd-cont-hit}{67.7}
\resdef{n11-s050-rd-cont-reg}{3.47}
\resdef{n11-s050-rd-cont-exc}{22.9}
\resdef{n11-s050-rd-cont-stop}{35.6}
\resdef{n11-s050-rd-cont-err}{7.6}
\resdef{n11-s050-rd-cont-errci}{[5.3, 10.8]}
\resdef{n11-s050-rd-round-q}{18.12}
\resdef{n11-s050-rd-round-hit}{66.6}
\resdef{n11-s050-rd-round-reg}{3.74}
\resdef{n11-s050-rd-round-exc}{24.0}
\resdef{n11-s050-rd-round-stop}{35.0}
\resdef{n11-s050-rd-round-err}{8.9}
\resdef{n11-s050-rd-round-errci}{[6.3, 12.3]}
\resdef{n11-s030-rd-dq}{$+0.40$}
\resdef{n11-s030-rd-dqci}{[$+0.14$, $+0.65$]}
\resdef{n11-s030-rd-dhit}{$-1.4$}
\resdef{n11-s030-rd-dhitci}{[$-3.5$, $+0.6$]}
\resdef{n11-s030-rd-dreg}{$+0.47$}
\resdef{n11-s030-rd-dregci}{[$+0.21$, $+0.73$]}
\resdef{n11-s030-rd-dexc}{$+2.0$}
\resdef{n11-s030-rd-dexcci}{[$+0.3$, $+3.8$]}
\resdef{n11-s030-rd-dstop}{$-5.7$}
\resdef{n11-s030-rd-dstopci}{[$-8.2$, $-3.2$]}
\resdef{n11-s050-rd-dq}{$-0.07$}
\resdef{n11-s050-rd-dqci}{[$-0.31$, $+0.17$]}
\resdef{n11-s050-rd-dhit}{$-1.1$}
\resdef{n11-s050-rd-dhitci}{[$-2.8$, $+0.6$]}
\resdef{n11-s050-rd-dreg}{$+0.28$}
\resdef{n11-s050-rd-dregci}{[$-0.02$, $+0.57$]}
\resdef{n11-s050-rd-dexc}{$+1.1$}
\resdef{n11-s050-rd-dexcci}{[$-0.4$, $+2.7$]}
\resdef{n11-s050-rd-dstop}{$-0.6$}
\resdef{n11-s050-rd-dstopci}{[$-2.7$, $+1.5$]}
\resdef{rd-R}{1000}
\resdef{dw-R}{1000}
\resdef{dw-strict-min}{31.4}
\resdef{dw-strict-max}{87.0}
\resdef{dw-agree-min}{82.7}
\resdef{dw-agree-max}{96.1}
\resdef{dw-other-agree-min}{71.7}
\resdef{dw-other-agree-max}{87.7}
\resdef{dw-cert-max}{21.6}
\resdef{dw-cert-max-n}{7}
\resdef{dw-cert-max-s}{0.15}
\resdef{dw-cert-max-n11}{0.2}
\resdef{dw-cert-rough-max}{16.0}
\resdef{dw-wrong}{0}
\resdef{n7-s015-dw-direct}{87.0}
\resdef{n7-s015-dw-weak}{0.0}
\resdef{n7-s015-dw-agreedirect}{96.1}
\resdef{n7-s015-dw-agreeweak}{--}
\resdef{n7-s015-dw-agreeother}{87.7}
\resdef{n7-s015-dw-certdirect}{21.6}
\resdef{n7-s030-dw-direct}{76.7}
\resdef{n7-s030-dw-weak}{0.0}
\resdef{n7-s030-dw-agreedirect}{91.9}
\resdef{n7-s030-dw-agreeweak}{--}
\resdef{n7-s030-dw-agreeother}{82.0}
\resdef{n7-s030-dw-certdirect}{2.7}
\resdef{n7-s050-dw-direct}{60.5}
\resdef{n7-s050-dw-weak}{0.0}
\resdef{n7-s050-dw-agreedirect}{90.9}
\resdef{n7-s050-dw-agreeweak}{--}
\resdef{n7-s050-dw-agreeother}{77.5}
\resdef{n7-s050-dw-certdirect}{0.2}
\resdef{n11-s015-dw-direct}{82.7}
\resdef{n11-s015-dw-weak}{0.0}
\resdef{n11-s015-dw-agreedirect}{93.0}
\resdef{n11-s015-dw-agreeweak}{--}
\resdef{n11-s015-dw-agreeother}{80.3}
\resdef{n11-s015-dw-certdirect}{0.2}
\resdef{n11-s030-dw-direct}{60.1}
\resdef{n11-s030-dw-weak}{0.0}
\resdef{n11-s030-dw-agreedirect}{89.7}
\resdef{n11-s030-dw-agreeweak}{--}
\resdef{n11-s030-dw-agreeother}{71.7}
\resdef{n11-s030-dw-certdirect}{0.0}
\resdef{n11-s050-dw-direct}{46.8}
\resdef{n11-s050-dw-weak}{0.0}
\resdef{n11-s050-dw-agreedirect}{86.8}
\resdef{n11-s050-dw-agreeweak}{--}
\resdef{n11-s050-dw-agreeother}{71.8}
\resdef{n11-s050-dw-certdirect}{0.0}
\resdef{n15-s015-dw-direct}{74.4}
\resdef{n15-s015-dw-weak}{0.0}
\resdef{n15-s015-dw-agreedirect}{88.6}
\resdef{n15-s015-dw-agreeweak}{--}
\resdef{n15-s015-dw-agreeother}{78.1}
\resdef{n15-s015-dw-certdirect}{0.0}
\resdef{n15-s030-dw-direct}{53.7}
\resdef{n15-s030-dw-weak}{0.0}
\resdef{n15-s030-dw-agreedirect}{82.7}
\resdef{n15-s030-dw-agreeweak}{--}
\resdef{n15-s030-dw-agreeother}{73.0}
\resdef{n15-s030-dw-certdirect}{0.0}
\resdef{n15-s050-dw-direct}{31.4}
\resdef{n15-s050-dw-weak}{0.0}
\resdef{n15-s050-dw-agreedirect}{87.6}
\resdef{n15-s050-dw-agreeweak}{--}
\resdef{n15-s050-dw-agreeother}{74.1}
\resdef{n15-s050-dw-certdirect}{0.0}
\resdef{n11-s030-saaty-dw-direct}{16.0}
\resdef{n11-s030-saaty-dw-weak}{82.2}
\resdef{n11-s030-saaty-dw-agreedirect}{100.0}
\resdef{n11-s030-saaty-dw-agreeweak}{78.8}
\resdef{n11-s030-saaty-dw-agreeother}{61.1}
\resdef{n11-s030-saaty-dw-certdirect}{0.0}
\resdef{dw-gci-n7-a3}{0.113}
\resdef{dw-gci-n7-a9}{0.451}
\resdef{dw-gci-n11-a3}{0.033}
\resdef{dw-gci-n11-a9}{0.131}
\resdef{dw-gci-n15-a3}{0.015}
\resdef{dw-gci-n15-a9}{0.061}
\resdef{pool-n7-h-diff}{0.4}
\resdef{pool-n7-h-ci}{[$-0.6$, $+1.3$]}
\resdef{strat-perm-n7-h}{0.36}
\resdef{strat-n7-h-b1b0}{-1.7}
\resdef{strat-n7-h-b1b0-ci}{[$-4.6$, $+1.1$]}
\resdef{strat-n7-h-b1b2}{0.3}
\resdef{strat-n7-h-b1b2-ci}{[$-2.1$, $+2.7$]}
\resdef{pool-n7-n-diff}{0.2}
\resdef{pool-n7-n-ci}{[$-0.5$, $+0.9$]}
\resdef{strat-perm-n7-n}{0.44}
\resdef{strat-n7-n-b1b0}{1.5}
\resdef{strat-n7-n-b1b0-ci}{[$-0.7$, $+3.6$]}
\resdef{strat-n7-n-b1b2}{0.5}
\resdef{strat-n7-n-b1b2-ci}{[$-1.2$, $+2.2$]}
\resdef{pool-n11-h-diff}{3.7}
\resdef{pool-n11-h-ci}{[$+2.4$, $+5.0$]}
\resdef{strat-perm-n11-h}{0.31}
\resdef{strat-n11-h-b1b0}{0.8}
\resdef{strat-n11-h-b1b0-ci}{[$-2.5$, $+4.0$]}
\resdef{strat-n11-h-b1b2}{2.6}
\resdef{strat-n11-h-b1b2-ci}{[$-0.6$, $+5.8$]}
\resdef{pool-n11-n-diff}{1.1}
\resdef{pool-n11-n-ci}{[$-0.1$, $+2.3$]}
\resdef{strat-perm-n11-n}{0.82}
\resdef{strat-n11-n-b1b0}{0.3}
\resdef{strat-n11-n-b1b0-ci}{[$-2.7$, $+3.4$]}
\resdef{strat-n11-n-b1b2}{1.6}
\resdef{strat-n11-n-b1b2-ci}{[$-0.9$, $+4.2$]}
\resdef{pool-n15-h-diff}{4.2}
\resdef{pool-n15-h-ci}{[$+2.8$, $+5.7$]}
\resdef{strat-perm-n15-h}{0.11}
\resdef{strat-n15-h-b1b0}{3.6}
\resdef{strat-n15-h-b1b0-ci}{[$-0.1$, $+7.3$]}
\resdef{strat-n15-h-b1b2}{3.2}
\resdef{strat-n15-h-b1b2-ci}{[$-0.8$, $+7.2$]}
\resdef{pool-n15-n-diff}{3.3}
\resdef{pool-n15-n-ci}{[$+1.8$, $+4.9$]}
\resdef{strat-perm-n15-n}{0.46}
\resdef{strat-n15-n-b1b0}{2.3}
\resdef{strat-n15-n-b1b0-ci}{[$-1.3$, $+5.9$]}
\resdef{strat-n15-n-b1b2}{3.1}
\resdef{strat-n15-n-b1b2-ci}{[$-0.5$, $+6.7$]}
\resdef{bayes-matched-n7-s030-a05-err}{1.9}
\resdef{bayes-matched-n7-s030-a05-hi}{2.5}
\resdef{bayes-matched-n7-s030-a05-stop}{59}
\resdef{bayes-flat-n7-s030-a05-err}{2.3}
\resdef{bayes-flat-n7-s030-a05-hi}{3.0}
\resdef{bayes-flat-n7-s030-a05-stop}{59}
\resdef{bayes-matched-n7-s030-a10-err}{3.7}
\resdef{bayes-matched-n7-s030-a10-hi}{4.5}
\resdef{bayes-matched-n7-s030-a10-stop}{69}
\resdef{bayes-flat-n7-s030-a10-err}{4.1}
\resdef{bayes-flat-n7-s030-a10-hi}{4.9}
\resdef{bayes-flat-n7-s030-a10-stop}{68}
\resdef{bayes-matched-n7-s030-a20-err}{8.3}
\resdef{bayes-matched-n7-s030-a20-hi}{9.3}
\resdef{bayes-matched-n7-s030-a20-stop}{81}
\resdef{bayes-flat-n7-s030-a20-err}{9.0}
\resdef{bayes-flat-n7-s030-a20-hi}{10.0}
\resdef{bayes-flat-n7-s030-a20-stop}{81}
\resdef{bayes-matched-n7-s050-a05-err}{3.2}
\resdef{bayes-matched-n7-s050-a05-hi}{4.2}
\resdef{bayes-matched-n7-s050-a05-stop}{36}
\resdef{bayes-flat-n7-s050-a05-err}{4.0}
\resdef{bayes-flat-n7-s050-a05-hi}{5.2}
\resdef{bayes-flat-n7-s050-a05-stop}{37}
\resdef{bayes-matched-n7-s050-a10-err}{5.5}
\resdef{bayes-matched-n7-s050-a10-hi}{6.6}
\resdef{bayes-matched-n7-s050-a10-stop}{47}
\resdef{bayes-flat-n7-s050-a10-err}{6.8}
\resdef{bayes-flat-n7-s050-a10-hi}{8.0}
\resdef{bayes-flat-n7-s050-a10-stop}{48}
\resdef{bayes-matched-n7-s050-a20-err}{11.4}
\resdef{bayes-matched-n7-s050-a20-hi}{12.7}
\resdef{bayes-matched-n7-s050-a20-stop}{64}
\resdef{bayes-flat-n7-s050-a20-err}{13.0}
\resdef{bayes-flat-n7-s050-a20-hi}{14.3}
\resdef{bayes-flat-n7-s050-a20-stop}{66}
\resdef{bayes-matched-n11-s030-a05-err}{2.4}
\resdef{bayes-matched-n11-s030-a05-hi}{3.1}
\resdef{bayes-matched-n11-s030-a05-stop}{63}
\resdef{bayes-flat-n11-s030-a05-err}{2.5}
\resdef{bayes-flat-n11-s030-a05-hi}{3.2}
\resdef{bayes-flat-n11-s030-a05-stop}{62}
\resdef{bayes-matched-n11-s030-a10-err}{4.6}
\resdef{bayes-matched-n11-s030-a10-hi}{5.4}
\resdef{bayes-matched-n11-s030-a10-stop}{73}
\resdef{bayes-flat-n11-s030-a10-err}{5.0}
\resdef{bayes-flat-n11-s030-a10-hi}{5.8}
\resdef{bayes-flat-n11-s030-a10-stop}{72}
\resdef{bayes-matched-n11-s030-a20-err}{9.6}
\resdef{bayes-matched-n11-s030-a20-hi}{10.6}
\resdef{bayes-matched-n11-s030-a20-stop}{85}
\resdef{bayes-flat-n11-s030-a20-err}{10.3}
\resdef{bayes-flat-n11-s030-a20-hi}{11.4}
\resdef{bayes-flat-n11-s030-a20-stop}{85}
\resdef{bayes-matched-n11-s050-a05-err}{3.1}
\resdef{bayes-matched-n11-s050-a05-hi}{4.1}
\resdef{bayes-matched-n11-s050-a05-stop}{41}
\resdef{bayes-flat-n11-s050-a05-err}{3.2}
\resdef{bayes-flat-n11-s050-a05-hi}{4.2}
\resdef{bayes-flat-n11-s050-a05-stop}{39}
\resdef{bayes-matched-n11-s050-a10-err}{6.4}
\resdef{bayes-matched-n11-s050-a10-hi}{7.5}
\resdef{bayes-matched-n11-s050-a10-stop}{52}
\resdef{bayes-flat-n11-s050-a10-err}{6.3}
\resdef{bayes-flat-n11-s050-a10-hi}{7.5}
\resdef{bayes-flat-n11-s050-a10-stop}{50}
\resdef{bayes-matched-n11-s050-a20-err}{11.9}
\resdef{bayes-matched-n11-s050-a20-hi}{13.1}
\resdef{bayes-matched-n11-s050-a20-stop}{68}
\resdef{bayes-flat-n11-s050-a20-err}{12.2}
\resdef{bayes-flat-n11-s050-a20-hi}{13.5}
\resdef{bayes-flat-n11-s050-a20-stop}{66}
\resdef{bayes-R}{4000}
\resdef{bayes-tau}{0.634}
\resdef{bayes-matched-maxexcess}{-1.8}
\resdef{time-ALEC-n7}{0.04}
\resdef{time-ALEC-n7-int}{0}
\resdef{time-FROZEN-n7}{0.04}
\resdef{time-FROZEN-n7-int}{0}
\resdef{time-AOPT-n7}{0.06}
\resdef{time-AOPT-n7-int}{0}
\resdef{time-ALEC-U-n7}{0.19}
\resdef{time-ALEC-U-n7-int}{0}
\resdef{time-HARKER-n7}{0.14}
\resdef{time-HARKER-n7-int}{0}
\resdef{time-HARKER-LOG-n7}{0.14}
\resdef{time-HARKER-LOG-n7-int}{0}
\resdef{time-HARKER-PG-n7}{0.62}
\resdef{time-HARKER-PG-n7-int}{1}
\resdef{time-ALEC-n11}{0.06}
\resdef{time-ALEC-n11-int}{0}
\resdef{time-FROZEN-n11}{0.07}
\resdef{time-FROZEN-n11-int}{0}
\resdef{time-AOPT-n11}{0.12}
\resdef{time-AOPT-n11-int}{0}
\resdef{time-ALEC-U-n11}{0.41}
\resdef{time-ALEC-U-n11-int}{0}
\resdef{time-HARKER-n11}{0.48}
\resdef{time-HARKER-n11-int}{0}
\resdef{time-HARKER-LOG-n11}{0.30}
\resdef{time-HARKER-LOG-n11-int}{0}
\resdef{time-HARKER-PG-n11}{5.47}
\resdef{time-HARKER-PG-n11-int}{5}
\resdef{time-ALEC-n15}{0.06}
\resdef{time-ALEC-n15-int}{0}
\resdef{time-FROZEN-n15}{0.08}
\resdef{time-FROZEN-n15-int}{0}
\resdef{time-AOPT-n15}{0.11}
\resdef{time-AOPT-n15-int}{0}
\resdef{time-ALEC-U-n15}{0.50}
\resdef{time-ALEC-U-n15-int}{0}
\resdef{time-HARKER-n15}{0.67}
\resdef{time-HARKER-n15-int}{1}
\resdef{time-HARKER-LOG-n15}{0.69}
\resdef{time-HARKER-LOG-n15-int}{1}
\resdef{time-HARKER-PG-n15}{40.46}
\resdef{time-HARKER-PG-n15-int}{40}
\resdef{time-ALEC-min}{0.04}
\resdef{time-ALEC-max}{0.06}
\resdef{time-FROZEN-min}{0.04}
\resdef{time-FROZEN-max}{0.08}
\resdef{time-AOPT-min}{0.06}
\resdef{time-AOPT-max}{0.12}
\resdef{time-ALEC-U-min}{0.19}
\resdef{time-ALEC-U-max}{0.50}
\resdef{time-HARKER-min}{0.14}
\resdef{time-HARKER-max}{0.67}
\resdef{time-HARKER-LOG-min}{0.14}
\resdef{time-HARKER-LOG-max}{0.69}
\resdef{time-HARKER-PG-min}{0.62}
\resdef{time-HARKER-PG-max}{40.46}
\resdef{cr-s015-mean}{0.006}
\resdef{crok-s015-min}{100}
\resdef{cr-s030-mean}{0.025}
\resdef{crok-s030-min}{100}
\resdef{cr-s050-mean}{0.073}
\resdef{crok-s050-min}{83}
\resdef{n7-s015-full-hit}{91.4}
\resdef{n7-s015-tot-h}{11}
\resdef{n7-s015-frac-h}{52}
\resdef{n7-s015-RANDOM-hit-h}{85.4}
\resdef{n7-s015-RANDOM-hitllsm-h}{85.4}
\resdef{n7-s015-RANDOM-tau-h}{0.912}
\resdef{n7-s015-RANDOM-reg-h}{1.05}
\resdef{n7-s015-RANDOM-agree-h}{90.4}
\resdef{n7-s015-p-RANDOM-h}{$p=0.001$}
\resdef{n7-s015-p-RANDOM-h-llsm}{$p=0.001$}
\resdef{n7-s015-HARKER-PG-hit-h}{88.0}
\resdef{n7-s015-HARKER-PG-hitllsm-h}{87.4}
\resdef{n7-s015-HARKER-PG-tau-h}{0.912}
\resdef{n7-s015-HARKER-PG-reg-h}{0.88}
\resdef{n7-s015-HARKER-PG-agree-h}{90.9}
\resdef{n7-s015-p-HARKER-PG-h}{$p=0.530$}
\resdef{n7-s015-p-HARKER-PG-h-llsm}{$p=0.203$}
\resdef{n7-s015-HARKER-hit-h}{87.2}
\resdef{n7-s015-HARKER-hitllsm-h}{87.3}
\resdef{n7-s015-HARKER-tau-h}{0.913}
\resdef{n7-s015-HARKER-reg-h}{0.88}
\resdef{n7-s015-HARKER-agree-h}{90.9}
\resdef{n7-s015-p-HARKER-h}{$p=0.133$}
\resdef{n7-s015-p-HARKER-h-llsm}{$p=0.161$}
\resdef{n7-s015-HARKER-LOG-hit-h}{89.8}
\resdef{n7-s015-HARKER-LOG-hitllsm-h}{89.6}
\resdef{n7-s015-HARKER-LOG-tau-h}{0.912}
\resdef{n7-s015-HARKER-LOG-reg-h}{0.60}
\resdef{n7-s015-HARKER-LOG-agree-h}{94.3}
\resdef{n7-s015-p-HARKER-LOG-h}{$p=0.222$}
\resdef{n7-s015-p-HARKER-LOG-h-llsm}{$p=0.328$}
\resdef{n7-s015-AOPT-hit-h}{87.6}
\resdef{n7-s015-AOPT-hitllsm-h}{87.6}
\resdef{n7-s015-AOPT-tau-h}{0.922}
\resdef{n7-s015-AOPT-reg-h}{0.82}
\resdef{n7-s015-AOPT-agree-h}{90.9}
\resdef{n7-s015-p-AOPT-h}{$p=0.289$}
\resdef{n7-s015-p-AOPT-h-llsm}{$p=0.289$}
\resdef{n7-s015-FROZEN-hit-h}{89.1}
\resdef{n7-s015-FROZEN-hitllsm-h}{89.1}
\resdef{n7-s015-FROZEN-tau-h}{0.914}
\resdef{n7-s015-FROZEN-reg-h}{0.70}
\resdef{n7-s015-FROZEN-agree-h}{93.1}
\resdef{n7-s015-p-FROZEN-h}{$p=0.644$}
\resdef{n7-s015-p-FROZEN-h-llsm}{$p=0.644$}
\resdef{n7-s015-STAR-hit-h}{88.9}
\resdef{n7-s015-STAR-hitllsm-h}{88.9}
\resdef{n7-s015-STAR-tau-h}{0.913}
\resdef{n7-s015-STAR-reg-h}{0.71}
\resdef{n7-s015-STAR-agree-h}{93.3}
\resdef{n7-s015-p-STAR-h}{$p=0.871$}
\resdef{n7-s015-p-STAR-h-llsm}{$p=0.871$}
\resdef{n7-s015-ALEC-hit-h}{88.7}
\resdef{n7-s015-ALEC-hitllsm-h}{88.7}
\resdef{n7-s015-ALEC-tau-h}{0.914}
\resdef{n7-s015-ALEC-reg-h}{0.72}
\resdef{n7-s015-ALEC-agree-h}{93.6}
\resdef{n7-s015-ALEC-U-hit-h}{89.0}
\resdef{n7-s015-ALEC-U-hitllsm-h}{89.0}
\resdef{n7-s015-ALEC-U-tau-h}{0.918}
\resdef{n7-s015-ALEC-U-reg-h}{0.62}
\resdef{n7-s015-ALEC-U-agree-h}{94.3}
\resdef{n7-s015-p-ALEC-U-h}{$p=0.771$}
\resdef{n7-s015-p-ALEC-U-h-llsm}{$p=0.771$}
\resdef{n7-s015-tot-n}{14}
\resdef{n7-s015-frac-n}{67}
\resdef{n7-s015-RANDOM-hit-n}{86.8}
\resdef{n7-s015-RANDOM-hitllsm-n}{86.8}
\resdef{n7-s015-RANDOM-tau-n}{0.929}
\resdef{n7-s015-RANDOM-reg-n}{0.88}
\resdef{n7-s015-RANDOM-agree-n}{91.2}
\resdef{n7-s015-p-RANDOM-n}{$p<0.001$}
\resdef{n7-s015-p-RANDOM-n-llsm}{$p<0.001$}
\resdef{n7-s015-HARKER-PG-hit-n}{90.2}
\resdef{n7-s015-HARKER-PG-hitllsm-n}{90.3}
\resdef{n7-s015-HARKER-PG-tau-n}{0.929}
\resdef{n7-s015-HARKER-PG-reg-n}{0.55}
\resdef{n7-s015-HARKER-PG-agree-n}{94.5}
\resdef{n7-s015-p-HARKER-PG-n}{$p=0.126$}
\resdef{n7-s015-p-HARKER-PG-n-llsm}{$p=0.169$}
\resdef{n7-s015-HARKER-hit-n}{89.3}
\resdef{n7-s015-HARKER-hitllsm-n}{89.5}
\resdef{n7-s015-HARKER-tau-n}{0.927}
\resdef{n7-s015-HARKER-reg-n}{0.61}
\resdef{n7-s015-HARKER-agree-n}{93.9}
\resdef{n7-s015-p-HARKER-n}{$p=0.008$}
\resdef{n7-s015-p-HARKER-n-llsm}{$p=0.016$}
\resdef{n7-s015-HARKER-LOG-hit-n}{90.4}
\resdef{n7-s015-HARKER-LOG-hitllsm-n}{90.2}
\resdef{n7-s015-HARKER-LOG-tau-n}{0.927}
\resdef{n7-s015-HARKER-LOG-reg-n}{0.50}
\resdef{n7-s015-HARKER-LOG-agree-n}{96.0}
\resdef{n7-s015-p-HARKER-LOG-n}{$p=0.143$}
\resdef{n7-s015-p-HARKER-LOG-n-llsm}{$p=0.073$}
\resdef{n7-s015-AOPT-hit-n}{89.4}
\resdef{n7-s015-AOPT-hitllsm-n}{89.4}
\resdef{n7-s015-AOPT-tau-n}{0.934}
\resdef{n7-s015-AOPT-reg-n}{0.57}
\resdef{n7-s015-AOPT-agree-n}{93.4}
\resdef{n7-s015-p-AOPT-n}{$p=0.019$}
\resdef{n7-s015-p-AOPT-n-llsm}{$p=0.019$}
\resdef{n7-s015-FROZEN-hit-n}{91.7}
\resdef{n7-s015-FROZEN-hitllsm-n}{91.7}
\resdef{n7-s015-FROZEN-tau-n}{0.930}
\resdef{n7-s015-FROZEN-reg-n}{0.41}
\resdef{n7-s015-FROZEN-agree-n}{97.7}
\resdef{n7-s015-p-FROZEN-n}{$p=0.607$}
\resdef{n7-s015-p-FROZEN-n-llsm}{$p=0.607$}
\resdef{n7-s015-STAR-hit-n}{91.6}
\resdef{n7-s015-STAR-hitllsm-n}{91.6}
\resdef{n7-s015-STAR-tau-n}{0.930}
\resdef{n7-s015-STAR-reg-n}{0.39}
\resdef{n7-s015-STAR-agree-n}{98.3}
\resdef{n7-s015-p-STAR-n}{$p=0.754$}
\resdef{n7-s015-p-STAR-n-llsm}{$p=0.754$}
\resdef{n7-s015-ALEC-hit-n}{91.4}
\resdef{n7-s015-ALEC-hitllsm-n}{91.4}
\resdef{n7-s015-ALEC-tau-n}{0.930}
\resdef{n7-s015-ALEC-reg-n}{0.41}
\resdef{n7-s015-ALEC-agree-n}{98.4}
\resdef{n7-s015-ALEC-U-hit-n}{90.2}
\resdef{n7-s015-ALEC-U-hitllsm-n}{90.2}
\resdef{n7-s015-ALEC-U-tau-n}{0.928}
\resdef{n7-s015-ALEC-U-reg-n}{0.52}
\resdef{n7-s015-ALEC-U-agree-n}{95.8}
\resdef{n7-s015-p-ALEC-U-n}{$p=0.065$}
\resdef{n7-s015-p-ALEC-U-n-llsm}{$p=0.065$}
\resdef{n7-s015-tot-2n}{21}
\resdef{n7-s015-frac-2n}{100}
\resdef{n7-s015-RANDOM-hit-2n}{91.4}
\resdef{n7-s015-RANDOM-hitllsm-2n}{91.4}
\resdef{n7-s015-RANDOM-tau-2n}{0.946}
\resdef{n7-s015-RANDOM-reg-2n}{0.40}
\resdef{n7-s015-RANDOM-agree-2n}{100.0}
\resdef{n7-s015-p-RANDOM-2n}{$p=1.000$}
\resdef{n7-s015-p-RANDOM-2n-llsm}{$p=1.000$}
\resdef{n7-s015-HARKER-PG-hit-2n}{91.8}
\resdef{n7-s015-HARKER-PG-hitllsm-2n}{91.4}
\resdef{n7-s015-HARKER-PG-tau-2n}{0.947}
\resdef{n7-s015-HARKER-PG-reg-2n}{0.36}
\resdef{n7-s015-HARKER-PG-agree-2n}{99.2}
\resdef{n7-s015-p-HARKER-PG-2n}{$p=0.289$}
\resdef{n7-s015-p-HARKER-PG-2n-llsm}{$p=1.000$}
\resdef{n7-s015-HARKER-hit-2n}{91.8}
\resdef{n7-s015-HARKER-hitllsm-2n}{91.4}
\resdef{n7-s015-HARKER-tau-2n}{0.947}
\resdef{n7-s015-HARKER-reg-2n}{0.36}
\resdef{n7-s015-HARKER-agree-2n}{99.2}
\resdef{n7-s015-p-HARKER-2n}{$p=0.289$}
\resdef{n7-s015-p-HARKER-2n-llsm}{$p=1.000$}
\resdef{n7-s015-HARKER-LOG-hit-2n}{91.8}
\resdef{n7-s015-HARKER-LOG-hitllsm-2n}{91.4}
\resdef{n7-s015-HARKER-LOG-tau-2n}{0.947}
\resdef{n7-s015-HARKER-LOG-reg-2n}{0.36}
\resdef{n7-s015-HARKER-LOG-agree-2n}{99.2}
\resdef{n7-s015-p-HARKER-LOG-2n}{$p=0.289$}
\resdef{n7-s015-p-HARKER-LOG-2n-llsm}{$p=1.000$}
\resdef{n7-s015-AOPT-hit-2n}{91.4}
\resdef{n7-s015-AOPT-hitllsm-2n}{91.4}
\resdef{n7-s015-AOPT-tau-2n}{0.946}
\resdef{n7-s015-AOPT-reg-2n}{0.40}
\resdef{n7-s015-AOPT-agree-2n}{100.0}
\resdef{n7-s015-p-AOPT-2n}{$p=1.000$}
\resdef{n7-s015-p-AOPT-2n-llsm}{$p=1.000$}
\resdef{n7-s015-FROZEN-hit-2n}{91.4}
\resdef{n7-s015-FROZEN-hitllsm-2n}{91.4}
\resdef{n7-s015-FROZEN-tau-2n}{0.946}
\resdef{n7-s015-FROZEN-reg-2n}{0.40}
\resdef{n7-s015-FROZEN-agree-2n}{100.0}
\resdef{n7-s015-p-FROZEN-2n}{$p=1.000$}
\resdef{n7-s015-p-FROZEN-2n-llsm}{$p=1.000$}
\resdef{n7-s015-STAR-hit-2n}{91.4}
\resdef{n7-s015-STAR-hitllsm-2n}{91.4}
\resdef{n7-s015-STAR-tau-2n}{0.946}
\resdef{n7-s015-STAR-reg-2n}{0.40}
\resdef{n7-s015-STAR-agree-2n}{100.0}
\resdef{n7-s015-p-STAR-2n}{$p=1.000$}
\resdef{n7-s015-p-STAR-2n-llsm}{$p=1.000$}
\resdef{n7-s015-ALEC-hit-2n}{91.4}
\resdef{n7-s015-ALEC-hitllsm-2n}{91.4}
\resdef{n7-s015-ALEC-tau-2n}{0.946}
\resdef{n7-s015-ALEC-reg-2n}{0.40}
\resdef{n7-s015-ALEC-agree-2n}{100.0}
\resdef{n7-s015-ALEC-U-hit-2n}{91.4}
\resdef{n7-s015-ALEC-U-hitllsm-2n}{91.4}
\resdef{n7-s015-ALEC-U-tau-2n}{0.946}
\resdef{n7-s015-ALEC-U-reg-2n}{0.40}
\resdef{n7-s015-ALEC-U-agree-2n}{100.0}
\resdef{n7-s015-p-ALEC-U-2n}{$p=1.000$}
\resdef{n7-s015-p-ALEC-U-2n-llsm}{$p=1.000$}
\resdef{n7-s015-stop-LC-q}{4.1}
\resdef{n7-s015-stop-LC-hit}{88.7}
\resdef{n7-s015-stop-LC-agree}{95.0}
\resdef{n7-s015-stop-LC-early}{100.0}
\resdef{n7-s015-stop-LC-exceed}{11.3}
\resdef{n7-s015-stop-LC-qmax}{11}
\resdef{n7-s015-stop-LC-qtot}{11.1}
\resdef{n7-s015-stop-LC-errearly}{11.3}
\resdef{n7-s015-stop-LC-errci}{[9.5, 13.4]}
\resdef{n7-s015-stop-LC-regret}{0.64}
\resdef{n7-s015-stop-LC-qround}{4}
\resdef{n7-s015-stop-LC-fixedhit}{88.7}
\resdef{n7-s015-stop-LC-fixedreg}{0.72}
\resdef{n7-s015-stop-LC-save2n}{71}
\resdef{n7-s015-stop-B10-ALEC-q}{6.4}
\resdef{n7-s015-stop-B10-ALEC-hit}{91.4}
\resdef{n7-s015-stop-B10-ALEC-agree}{100.0}
\resdef{n7-s015-stop-B10-ALEC-early}{69.7}
\resdef{n7-s015-stop-B10-ALEC-exceed}{8.6}
\resdef{n7-s015-stop-B10-ALEC-qmax}{14}
\resdef{n7-s015-stop-B10-ALEC-qtot}{13.4}
\resdef{n7-s015-stop-B10-ALEC-errearly}{1.0}
\resdef{n7-s015-stop-B10-ALEC-errci}{[0.5, 2.1]}
\resdef{n7-s015-stop-B10-ALEC-regret}{0.40}
\resdef{n7-s015-stop-B10-ALEC-qround}{6}
\resdef{n7-s015-stop-B10-ALEC-fixedhit}{90.7}
\resdef{n7-s015-stop-B10-ALEC-fixedreg}{0.46}
\resdef{n7-s015-stop-B10-ALEC-save2n}{54}
\resdef{n7-s015-stop-B05-ALECU-q}{7.6}
\resdef{n7-s015-stop-B05-ALECU-hit}{91.4}
\resdef{n7-s015-stop-B05-ALECU-agree}{100.0}
\resdef{n7-s015-stop-B05-ALECU-early}{60.9}
\resdef{n7-s015-stop-B05-ALECU-exceed}{8.6}
\resdef{n7-s015-stop-B05-ALECU-qmax}{14}
\resdef{n7-s015-stop-B05-ALECU-qtot}{14.6}
\resdef{n7-s015-stop-B05-ALECU-errearly}{0.3}
\resdef{n7-s015-stop-B05-ALECU-errci}{[0.1, 1.2]}
\resdef{n7-s015-stop-B05-ALECU-regret}{0.40}
\resdef{n7-s015-stop-B05-ALECU-qround}{8}
\resdef{n7-s015-stop-B05-ALECU-fixedhit}{91.3}
\resdef{n7-s015-stop-B05-ALECU-fixedreg}{0.40}
\resdef{n7-s015-stop-B05-ALECU-save2n}{45}
\resdef{n7-s015-stop-B10-ALECU-q}{6.3}
\resdef{n7-s015-stop-B10-ALECU-hit}{91.4}
\resdef{n7-s015-stop-B10-ALECU-agree}{100.0}
\resdef{n7-s015-stop-B10-ALECU-early}{70.2}
\resdef{n7-s015-stop-B10-ALECU-exceed}{8.6}
\resdef{n7-s015-stop-B10-ALECU-qmax}{14}
\resdef{n7-s015-stop-B10-ALECU-qtot}{13.3}
\resdef{n7-s015-stop-B10-ALECU-errearly}{0.9}
\resdef{n7-s015-stop-B10-ALECU-errci}{[0.4, 1.9]}
\resdef{n7-s015-stop-B10-ALECU-regret}{0.40}
\resdef{n7-s015-stop-B10-ALECU-qround}{6}
\resdef{n7-s015-stop-B10-ALECU-fixedhit}{90.7}
\resdef{n7-s015-stop-B10-ALECU-fixedreg}{0.46}
\resdef{n7-s015-stop-B10-ALECU-save2n}{55}
\resdef{n7-s015-stop-B20-ALECU-q}{4.6}
\resdef{n7-s015-stop-B20-ALECU-hit}{91.3}
\resdef{n7-s015-stop-B20-ALECU-agree}{99.9}
\resdef{n7-s015-stop-B20-ALECU-early}{80.3}
\resdef{n7-s015-stop-B20-ALECU-exceed}{8.7}
\resdef{n7-s015-stop-B20-ALECU-qmax}{14}
\resdef{n7-s015-stop-B20-ALECU-qtot}{11.6}
\resdef{n7-s015-stop-B20-ALECU-errearly}{2.9}
\resdef{n7-s015-stop-B20-ALECU-errci}{[1.9, 4.3]}
\resdef{n7-s015-stop-B20-ALECU-regret}{0.41}
\resdef{n7-s015-stop-B20-ALECU-qround}{5}
\resdef{n7-s015-stop-B20-ALECU-fixedhit}{89.2}
\resdef{n7-s015-stop-B20-ALECU-fixedreg}{0.63}
\resdef{n7-s015-stop-B20-ALECU-save2n}{67}
\resdef{n7-s015-stop-B10-ALECU-s06-q}{9.2}
\resdef{n7-s015-stop-B10-ALECU-s06-hit}{91.4}
\resdef{n7-s015-stop-B10-ALECU-s06-agree}{100.0}
\resdef{n7-s015-stop-B10-ALECU-s06-early}{54.8}
\resdef{n7-s015-stop-B10-ALECU-s06-exceed}{8.6}
\resdef{n7-s015-stop-B10-ALECU-s06-qmax}{14}
\resdef{n7-s015-stop-B10-ALECU-s06-qtot}{16.2}
\resdef{n7-s015-stop-B10-ALECU-s06-errearly}{0.2}
\resdef{n7-s015-stop-B10-ALECU-s06-errci}{[0.0, 1.0]}
\resdef{n7-s015-stop-B10-ALECU-s06-regret}{0.40}
\resdef{n7-s015-stop-B10-ALECU-s06-qround}{9}
\resdef{n7-s015-stop-B10-ALECU-s06-fixedhit}{91.7}
\resdef{n7-s015-stop-B10-ALECU-s06-fixedreg}{0.39}
\resdef{n7-s015-stop-B10-ALECU-s06-save2n}{34}
\resdef{n7-s015-stop-B10-ALECU-cap-q}{6.3}
\resdef{n7-s015-stop-B10-ALECU-cap-hit}{91.4}
\resdef{n7-s015-stop-B10-ALECU-cap-agree}{100.0}
\resdef{n7-s015-stop-B10-ALECU-cap-early}{70.2}
\resdef{n7-s015-stop-B10-ALECU-cap-exceed}{8.6}
\resdef{n7-s015-stop-B10-ALECU-cap-qmax}{14}
\resdef{n7-s015-stop-B10-ALECU-cap-qtot}{13.3}
\resdef{n7-s015-stop-B10-ALECU-cap-errearly}{0.9}
\resdef{n7-s015-stop-B10-ALECU-cap-errci}{[0.4, 1.9]}
\resdef{n7-s015-stop-B10-ALECU-cap-regret}{0.40}
\resdef{n7-s015-stop-B10-ALECU-cap-qround}{6}
\resdef{n7-s015-stop-B10-ALECU-cap-fixedhit}{90.7}
\resdef{n7-s015-stop-B10-ALECU-cap-fixedreg}{0.46}
\resdef{n7-s015-stop-B10-ALECU-cap-save2n}{55}
\resdef{n7-s015-stop-B10-ALECU-oracle-q}{3.6}
\resdef{n7-s015-stop-B10-ALECU-oracle-hit}{90.3}
\resdef{n7-s015-stop-B10-ALECU-oracle-agree}{98.9}
\resdef{n7-s015-stop-B10-ALECU-oracle-early}{79.9}
\resdef{n7-s015-stop-B10-ALECU-oracle-exceed}{9.7}
\resdef{n7-s015-stop-B10-ALECU-oracle-qmax}{14}
\resdef{n7-s015-stop-B10-ALECU-oracle-qtot}{10.6}
\resdef{n7-s015-stop-B10-ALECU-oracle-errearly}{4.3}
\resdef{n7-s015-stop-B10-ALECU-oracle-errci}{[3.1, 5.9]}
\resdef{n7-s015-stop-B10-ALECU-oracle-regret}{0.53}
\resdef{n7-s015-stop-B10-ALECU-oracle-qround}{4}
\resdef{n7-s015-stop-B10-ALECU-oracle-fixedhit}{88.7}
\resdef{n7-s015-stop-B10-ALECU-oracle-fixedreg}{0.72}
\resdef{n7-s015-stop-B10-ALECU-oracle-save2n}{74}
\resdef{n7-s015-stop-E02-ALEC-q}{5.9}
\resdef{n7-s015-stop-E02-ALEC-hit}{91.4}
\resdef{n7-s015-stop-E02-ALEC-agree}{100.0}
\resdef{n7-s015-stop-E02-ALEC-early}{73.8}
\resdef{n7-s015-stop-E02-ALEC-exceed}{6.6}
\resdef{n7-s015-stop-E02-ALEC-qmax}{14}
\resdef{n7-s015-stop-E02-ALEC-qtot}{12.9}
\resdef{n7-s015-stop-E02-ALEC-errearly}{0.8}
\resdef{n7-s015-stop-E02-ALEC-errci}{[0.4, 1.8]}
\resdef{n7-s015-stop-E02-ALEC-regret}{0.40}
\resdef{n7-s015-stop-E02-ALEC-qround}{6}
\resdef{n7-s015-stop-E02-ALEC-fixedhit}{90.7}
\resdef{n7-s015-stop-E02-ALEC-fixedreg}{0.46}
\resdef{n7-s015-stop-E02-ALEC-save2n}{58}
\resdef{n7-s015-stop-E05-ALEC-q}{5.1}
\resdef{n7-s015-stop-E05-ALEC-hit}{91.4}
\resdef{n7-s015-stop-E05-ALEC-agree}{100.0}
\resdef{n7-s015-stop-E05-ALEC-early}{81.4}
\resdef{n7-s015-stop-E05-ALEC-exceed}{3.2}
\resdef{n7-s015-stop-E05-ALEC-qmax}{14}
\resdef{n7-s015-stop-E05-ALEC-qtot}{12.1}
\resdef{n7-s015-stop-E05-ALEC-errearly}{0.6}
\resdef{n7-s015-stop-E05-ALEC-errci}{[0.3, 1.4]}
\resdef{n7-s015-stop-E05-ALEC-regret}{0.40}
\resdef{n7-s015-stop-E05-ALEC-qround}{5}
\resdef{n7-s015-stop-E05-ALEC-fixedhit}{89.2}
\resdef{n7-s015-stop-E05-ALEC-fixedreg}{0.63}
\resdef{n7-s015-stop-E05-ALEC-save2n}{64}
\resdef{n7-s015-stop-B10-ALEC-cap-q}{6.4}
\resdef{n7-s015-stop-B10-ALEC-cap-hit}{91.4}
\resdef{n7-s015-stop-B10-ALEC-cap-agree}{100.0}
\resdef{n7-s015-stop-B10-ALEC-cap-early}{69.7}
\resdef{n7-s015-stop-B10-ALEC-cap-exceed}{8.6}
\resdef{n7-s015-stop-B10-ALEC-cap-qmax}{14}
\resdef{n7-s015-stop-B10-ALEC-cap-qtot}{13.4}
\resdef{n7-s015-stop-B10-ALEC-cap-errearly}{1.0}
\resdef{n7-s015-stop-B10-ALEC-cap-errci}{[0.5, 2.1]}
\resdef{n7-s015-stop-B10-ALEC-cap-regret}{0.40}
\resdef{n7-s015-stop-B10-ALEC-cap-qround}{6}
\resdef{n7-s015-stop-B10-ALEC-cap-fixedhit}{90.7}
\resdef{n7-s015-stop-B10-ALEC-cap-fixedreg}{0.46}
\resdef{n7-s015-stop-B10-ALEC-cap-save2n}{54}
\resdef{n7-s015-stop-E02-ALEC-cap-q}{5.9}
\resdef{n7-s015-stop-E02-ALEC-cap-hit}{91.4}
\resdef{n7-s015-stop-E02-ALEC-cap-agree}{100.0}
\resdef{n7-s015-stop-E02-ALEC-cap-early}{73.8}
\resdef{n7-s015-stop-E02-ALEC-cap-exceed}{6.6}
\resdef{n7-s015-stop-E02-ALEC-cap-qmax}{14}
\resdef{n7-s015-stop-E02-ALEC-cap-qtot}{12.9}
\resdef{n7-s015-stop-E02-ALEC-cap-errearly}{0.8}
\resdef{n7-s015-stop-E02-ALEC-cap-errci}{[0.4, 1.8]}
\resdef{n7-s015-stop-E02-ALEC-cap-regret}{0.40}
\resdef{n7-s015-stop-E02-ALEC-cap-qround}{6}
\resdef{n7-s015-stop-E02-ALEC-cap-fixedhit}{90.7}
\resdef{n7-s015-stop-E02-ALEC-cap-fixedreg}{0.46}
\resdef{n7-s015-stop-E02-ALEC-cap-save2n}{58}
\resdef{n7-s015-stop-E05-ALEC-cap-q}{5.1}
\resdef{n7-s015-stop-E05-ALEC-cap-hit}{91.4}
\resdef{n7-s015-stop-E05-ALEC-cap-agree}{100.0}
\resdef{n7-s015-stop-E05-ALEC-cap-early}{81.4}
\resdef{n7-s015-stop-E05-ALEC-cap-exceed}{3.2}
\resdef{n7-s015-stop-E05-ALEC-cap-qmax}{14}
\resdef{n7-s015-stop-E05-ALEC-cap-qtot}{12.1}
\resdef{n7-s015-stop-E05-ALEC-cap-errearly}{0.6}
\resdef{n7-s015-stop-E05-ALEC-cap-errci}{[0.3, 1.4]}
\resdef{n7-s015-stop-E05-ALEC-cap-regret}{0.40}
\resdef{n7-s015-stop-E05-ALEC-cap-qround}{5}
\resdef{n7-s015-stop-E05-ALEC-cap-fixedhit}{89.2}
\resdef{n7-s015-stop-E05-ALEC-cap-fixedreg}{0.63}
\resdef{n7-s015-stop-E05-ALEC-cap-save2n}{64}
\resdef{n7-s015-stop-E05-ALEC-cap-s06-q}{8.3}
\resdef{n7-s015-stop-E05-ALEC-cap-s06-hit}{91.4}
\resdef{n7-s015-stop-E05-ALEC-cap-s06-agree}{100.0}
\resdef{n7-s015-stop-E05-ALEC-cap-s06-early}{65.1}
\resdef{n7-s015-stop-E05-ALEC-cap-s06-exceed}{3.2}
\resdef{n7-s015-stop-E05-ALEC-cap-s06-qmax}{14}
\resdef{n7-s015-stop-E05-ALEC-cap-s06-qtot}{15.3}
\resdef{n7-s015-stop-E05-ALEC-cap-s06-errearly}{0.2}
\resdef{n7-s015-stop-E05-ALEC-cap-s06-errci}{[0.0, 0.9]}
\resdef{n7-s015-stop-E05-ALEC-cap-s06-regret}{0.40}
\resdef{n7-s015-stop-E05-ALEC-cap-s06-qround}{8}
\resdef{n7-s015-stop-E05-ALEC-cap-s06-fixedhit}{91.3}
\resdef{n7-s015-stop-E05-ALEC-cap-s06-fixedreg}{0.40}
\resdef{n7-s015-stop-E05-ALEC-cap-s06-save2n}{41}
\resdef{n7-s030-full-hit}{80.8}
\resdef{n7-s030-tot-h}{11}
\resdef{n7-s030-frac-h}{52}
\resdef{n7-s030-RANDOM-hit-h}{72.8}
\resdef{n7-s030-RANDOM-hitllsm-h}{72.8}
\resdef{n7-s030-RANDOM-tau-h}{0.824}
\resdef{n7-s030-RANDOM-reg-h}{3.52}
\resdef{n7-s030-RANDOM-agree-h}{80.7}
\resdef{n7-s030-p-RANDOM-h}{$p=0.005$}
\resdef{n7-s030-p-RANDOM-h-llsm}{$p=0.005$}
\resdef{n7-s030-HARKER-PG-hit-h}{75.0}
\resdef{n7-s030-HARKER-PG-hitllsm-h}{75.4}
\resdef{n7-s030-HARKER-PG-tau-h}{0.819}
\resdef{n7-s030-HARKER-PG-reg-h}{2.80}
\resdef{n7-s030-HARKER-PG-agree-h}{84.4}
\resdef{n7-s030-p-HARKER-PG-h}{$p=0.176$}
\resdef{n7-s030-p-HARKER-PG-h-llsm}{$p=0.310$}
\resdef{n7-s030-HARKER-hit-h}{77.0}
\resdef{n7-s030-HARKER-hitllsm-h}{77.3}
\resdef{n7-s030-HARKER-tau-h}{0.823}
\resdef{n7-s030-HARKER-reg-h}{2.43}
\resdef{n7-s030-HARKER-agree-h}{85.6}
\resdef{n7-s030-p-HARKER-h}{$p=0.936$}
\resdef{n7-s030-p-HARKER-h-llsm}{$p=0.747$}
\resdef{n7-s030-HARKER-LOG-hit-h}{77.1}
\resdef{n7-s030-HARKER-LOG-hitllsm-h}{78.3}
\resdef{n7-s030-HARKER-LOG-tau-h}{0.817}
\resdef{n7-s030-HARKER-LOG-reg-h}{2.49}
\resdef{n7-s030-HARKER-LOG-agree-h}{87.3}
\resdef{n7-s030-p-HARKER-LOG-h}{$p=0.832$}
\resdef{n7-s030-p-HARKER-LOG-h-llsm}{$p=0.119$}
\resdef{n7-s030-AOPT-hit-h}{75.1}
\resdef{n7-s030-AOPT-hitllsm-h}{75.1}
\resdef{n7-s030-AOPT-tau-h}{0.835}
\resdef{n7-s030-AOPT-reg-h}{2.80}
\resdef{n7-s030-AOPT-agree-h}{82.0}
\resdef{n7-s030-p-AOPT-h}{$p=0.224$}
\resdef{n7-s030-p-AOPT-h-llsm}{$p=0.224$}
\resdef{n7-s030-FROZEN-hit-h}{77.0}
\resdef{n7-s030-FROZEN-hitllsm-h}{77.0}
\resdef{n7-s030-FROZEN-tau-h}{0.820}
\resdef{n7-s030-FROZEN-reg-h}{2.59}
\resdef{n7-s030-FROZEN-agree-h}{87.0}
\resdef{n7-s030-p-FROZEN-h}{$p=0.905$}
\resdef{n7-s030-p-FROZEN-h-llsm}{$p=0.905$}
\resdef{n7-s030-STAR-hit-h}{77.0}
\resdef{n7-s030-STAR-hitllsm-h}{77.0}
\resdef{n7-s030-STAR-tau-h}{0.821}
\resdef{n7-s030-STAR-reg-h}{2.57}
\resdef{n7-s030-STAR-agree-h}{87.1}
\resdef{n7-s030-p-STAR-h}{$p=0.902$}
\resdef{n7-s030-p-STAR-h-llsm}{$p=0.902$}
\resdef{n7-s030-ALEC-hit-h}{76.8}
\resdef{n7-s030-ALEC-hitllsm-h}{76.8}
\resdef{n7-s030-ALEC-tau-h}{0.821}
\resdef{n7-s030-ALEC-reg-h}{2.47}
\resdef{n7-s030-ALEC-agree-h}{88.3}
\resdef{n7-s030-ALEC-U-hit-h}{76.7}
\resdef{n7-s030-ALEC-U-hitllsm-h}{76.7}
\resdef{n7-s030-ALEC-U-tau-h}{0.823}
\resdef{n7-s030-ALEC-U-reg-h}{2.48}
\resdef{n7-s030-ALEC-U-agree-h}{88.2}
\resdef{n7-s030-p-ALEC-U-h}{$p=1.000$}
\resdef{n7-s030-p-ALEC-U-h-llsm}{$p=1.000$}
\resdef{n7-s030-tot-n}{14}
\resdef{n7-s030-frac-n}{67}
\resdef{n7-s030-RANDOM-hit-n}{74.9}
\resdef{n7-s030-RANDOM-hitllsm-n}{74.9}
\resdef{n7-s030-RANDOM-tau-n}{0.853}
\resdef{n7-s030-RANDOM-reg-n}{2.77}
\resdef{n7-s030-RANDOM-agree-n}{86.7}
\resdef{n7-s030-p-RANDOM-n}{$p<0.001$}
\resdef{n7-s030-p-RANDOM-n-llsm}{$p<0.001$}
\resdef{n7-s030-HARKER-PG-hit-n}{77.6}
\resdef{n7-s030-HARKER-PG-hitllsm-n}{77.9}
\resdef{n7-s030-HARKER-PG-tau-n}{0.851}
\resdef{n7-s030-HARKER-PG-reg-n}{2.22}
\resdef{n7-s030-HARKER-PG-agree-n}{89.3}
\resdef{n7-s030-p-HARKER-PG-n}{$p=0.061$}
\resdef{n7-s030-p-HARKER-PG-n-llsm}{$p=0.117$}
\resdef{n7-s030-HARKER-hit-n}{77.5}
\resdef{n7-s030-HARKER-hitllsm-n}{77.2}
\resdef{n7-s030-HARKER-tau-n}{0.853}
\resdef{n7-s030-HARKER-reg-n}{2.14}
\resdef{n7-s030-HARKER-agree-n}{90.2}
\resdef{n7-s030-p-HARKER-n}{$p=0.042$}
\resdef{n7-s030-p-HARKER-n-llsm}{$p=0.021$}
\resdef{n7-s030-HARKER-LOG-hit-n}{78.5}
\resdef{n7-s030-HARKER-LOG-hitllsm-n}{78.8}
\resdef{n7-s030-HARKER-LOG-tau-n}{0.844}
\resdef{n7-s030-HARKER-LOG-reg-n}{1.92}
\resdef{n7-s030-HARKER-LOG-agree-n}{93.6}
\resdef{n7-s030-p-HARKER-LOG-n}{$p=0.275$}
\resdef{n7-s030-p-HARKER-LOG-n-llsm}{$p=0.435$}
\resdef{n7-s030-AOPT-hit-n}{77.1}
\resdef{n7-s030-AOPT-hitllsm-n}{77.1}
\resdef{n7-s030-AOPT-tau-n}{0.858}
\resdef{n7-s030-AOPT-reg-n}{2.44}
\resdef{n7-s030-AOPT-agree-n}{87.8}
\resdef{n7-s030-p-AOPT-n}{$p=0.025$}
\resdef{n7-s030-p-AOPT-n-llsm}{$p=0.025$}
\resdef{n7-s030-FROZEN-hit-n}{79.8}
\resdef{n7-s030-FROZEN-hitllsm-n}{79.8}
\resdef{n7-s030-FROZEN-tau-n}{0.849}
\resdef{n7-s030-FROZEN-reg-n}{1.90}
\resdef{n7-s030-FROZEN-agree-n}{95.0}
\resdef{n7-s030-p-FROZEN-n}{$p=0.728$}
\resdef{n7-s030-p-FROZEN-n-llsm}{$p=0.728$}
\resdef{n7-s030-STAR-hit-n}{80.2}
\resdef{n7-s030-STAR-hitllsm-n}{80.2}
\resdef{n7-s030-STAR-tau-n}{0.848}
\resdef{n7-s030-STAR-reg-n}{1.75}
\resdef{n7-s030-STAR-agree-n}{96.2}
\resdef{n7-s030-p-STAR-n}{$p=0.143$}
\resdef{n7-s030-p-STAR-n-llsm}{$p=0.143$}
\resdef{n7-s030-ALEC-hit-n}{79.5}
\resdef{n7-s030-ALEC-hitllsm-n}{79.5}
\resdef{n7-s030-ALEC-tau-n}{0.848}
\resdef{n7-s030-ALEC-reg-n}{1.81}
\resdef{n7-s030-ALEC-agree-n}{96.2}
\resdef{n7-s030-ALEC-U-hit-n}{79.2}
\resdef{n7-s030-ALEC-U-hitllsm-n}{79.2}
\resdef{n7-s030-ALEC-U-tau-n}{0.851}
\resdef{n7-s030-ALEC-U-reg-n}{1.89}
\resdef{n7-s030-ALEC-U-agree-n}{92.8}
\resdef{n7-s030-p-ALEC-U-n}{$p=0.784$}
\resdef{n7-s030-p-ALEC-U-n-llsm}{$p=0.784$}
\resdef{n7-s030-tot-2n}{21}
\resdef{n7-s030-frac-2n}{100}
\resdef{n7-s030-RANDOM-hit-2n}{80.8}
\resdef{n7-s030-RANDOM-hitllsm-2n}{80.8}
\resdef{n7-s030-RANDOM-tau-2n}{0.887}
\resdef{n7-s030-RANDOM-reg-2n}{1.61}
\resdef{n7-s030-RANDOM-agree-2n}{100.0}
\resdef{n7-s030-p-RANDOM-2n}{$p=1.000$}
\resdef{n7-s030-p-RANDOM-2n-llsm}{$p=1.000$}
\resdef{n7-s030-HARKER-PG-hit-2n}{81.7}
\resdef{n7-s030-HARKER-PG-hitllsm-2n}{80.8}
\resdef{n7-s030-HARKER-PG-tau-2n}{0.887}
\resdef{n7-s030-HARKER-PG-reg-2n}{1.51}
\resdef{n7-s030-HARKER-PG-agree-2n}{97.4}
\resdef{n7-s030-p-HARKER-PG-2n}{$p=0.093$}
\resdef{n7-s030-p-HARKER-PG-2n-llsm}{$p=1.000$}
\resdef{n7-s030-HARKER-hit-2n}{81.7}
\resdef{n7-s030-HARKER-hitllsm-2n}{80.8}
\resdef{n7-s030-HARKER-tau-2n}{0.887}
\resdef{n7-s030-HARKER-reg-2n}{1.51}
\resdef{n7-s030-HARKER-agree-2n}{97.4}
\resdef{n7-s030-p-HARKER-2n}{$p=0.093$}
\resdef{n7-s030-p-HARKER-2n-llsm}{$p=1.000$}
\resdef{n7-s030-HARKER-LOG-hit-2n}{81.7}
\resdef{n7-s030-HARKER-LOG-hitllsm-2n}{80.8}
\resdef{n7-s030-HARKER-LOG-tau-2n}{0.887}
\resdef{n7-s030-HARKER-LOG-reg-2n}{1.51}
\resdef{n7-s030-HARKER-LOG-agree-2n}{97.4}
\resdef{n7-s030-p-HARKER-LOG-2n}{$p=0.093$}
\resdef{n7-s030-p-HARKER-LOG-2n-llsm}{$p=1.000$}
\resdef{n7-s030-AOPT-hit-2n}{80.8}
\resdef{n7-s030-AOPT-hitllsm-2n}{80.8}
\resdef{n7-s030-AOPT-tau-2n}{0.887}
\resdef{n7-s030-AOPT-reg-2n}{1.61}
\resdef{n7-s030-AOPT-agree-2n}{100.0}
\resdef{n7-s030-p-AOPT-2n}{$p=1.000$}
\resdef{n7-s030-p-AOPT-2n-llsm}{$p=1.000$}
\resdef{n7-s030-FROZEN-hit-2n}{80.8}
\resdef{n7-s030-FROZEN-hitllsm-2n}{80.8}
\resdef{n7-s030-FROZEN-tau-2n}{0.887}
\resdef{n7-s030-FROZEN-reg-2n}{1.61}
\resdef{n7-s030-FROZEN-agree-2n}{100.0}
\resdef{n7-s030-p-FROZEN-2n}{$p=1.000$}
\resdef{n7-s030-p-FROZEN-2n-llsm}{$p=1.000$}
\resdef{n7-s030-STAR-hit-2n}{80.8}
\resdef{n7-s030-STAR-hitllsm-2n}{80.8}
\resdef{n7-s030-STAR-tau-2n}{0.887}
\resdef{n7-s030-STAR-reg-2n}{1.61}
\resdef{n7-s030-STAR-agree-2n}{100.0}
\resdef{n7-s030-p-STAR-2n}{$p=1.000$}
\resdef{n7-s030-p-STAR-2n-llsm}{$p=1.000$}
\resdef{n7-s030-ALEC-hit-2n}{80.8}
\resdef{n7-s030-ALEC-hitllsm-2n}{80.8}
\resdef{n7-s030-ALEC-tau-2n}{0.887}
\resdef{n7-s030-ALEC-reg-2n}{1.61}
\resdef{n7-s030-ALEC-agree-2n}{100.0}
\resdef{n7-s030-ALEC-U-hit-2n}{80.8}
\resdef{n7-s030-ALEC-U-hitllsm-2n}{80.8}
\resdef{n7-s030-ALEC-U-tau-2n}{0.887}
\resdef{n7-s030-ALEC-U-reg-2n}{1.61}
\resdef{n7-s030-ALEC-U-agree-2n}{100.0}
\resdef{n7-s030-p-ALEC-U-2n}{$p=1.000$}
\resdef{n7-s030-p-ALEC-U-2n-llsm}{$p=1.000$}
\resdef{n7-s030-stop-LC-q}{4.3}
\resdef{n7-s030-stop-LC-hit}{77.5}
\resdef{n7-s030-stop-LC-agree}{89.6}
\resdef{n7-s030-stop-LC-early}{100.0}
\resdef{n7-s030-stop-LC-exceed}{22.5}
\resdef{n7-s030-stop-LC-qmax}{9}
\resdef{n7-s030-stop-LC-qtot}{11.3}
\resdef{n7-s030-stop-LC-errearly}{22.5}
\resdef{n7-s030-stop-LC-errci}{[20.0, 25.2]}
\resdef{n7-s030-stop-LC-regret}{2.28}
\resdef{n7-s030-stop-LC-qround}{4}
\resdef{n7-s030-stop-LC-fixedhit}{76.8}
\resdef{n7-s030-stop-LC-fixedreg}{2.47}
\resdef{n7-s030-stop-LC-save2n}{69}
\resdef{n7-s030-stop-B10-ALEC-q}{8.0}
\resdef{n7-s030-stop-B10-ALEC-hit}{80.6}
\resdef{n7-s030-stop-B10-ALEC-agree}{99.8}
\resdef{n7-s030-stop-B10-ALEC-early}{56.7}
\resdef{n7-s030-stop-B10-ALEC-exceed}{19.4}
\resdef{n7-s030-stop-B10-ALEC-qmax}{14}
\resdef{n7-s030-stop-B10-ALEC-qtot}{15.0}
\resdef{n7-s030-stop-B10-ALEC-errearly}{7.9}
\resdef{n7-s030-stop-B10-ALEC-errci}{[6.0, 10.5]}
\resdef{n7-s030-stop-B10-ALEC-regret}{1.65}
\resdef{n7-s030-stop-B10-ALEC-qround}{8}
\resdef{n7-s030-stop-B10-ALEC-fixedhit}{80.5}
\resdef{n7-s030-stop-B10-ALEC-fixedreg}{1.72}
\resdef{n7-s030-stop-B10-ALEC-save2n}{43}
\resdef{n7-s030-stop-B05-ALECU-q}{9.4}
\resdef{n7-s030-stop-B05-ALECU-hit}{80.7}
\resdef{n7-s030-stop-B05-ALECU-agree}{99.9}
\resdef{n7-s030-stop-B05-ALECU-early}{46.3}
\resdef{n7-s030-stop-B05-ALECU-exceed}{19.3}
\resdef{n7-s030-stop-B05-ALECU-qmax}{14}
\resdef{n7-s030-stop-B05-ALECU-qtot}{16.4}
\resdef{n7-s030-stop-B05-ALECU-errearly}{4.8}
\resdef{n7-s030-stop-B05-ALECU-errci}{[3.2, 7.1]}
\resdef{n7-s030-stop-B05-ALECU-regret}{1.62}
\resdef{n7-s030-stop-B05-ALECU-qround}{9}
\resdef{n7-s030-stop-B05-ALECU-fixedhit}{81.1}
\resdef{n7-s030-stop-B05-ALECU-fixedreg}{1.59}
\resdef{n7-s030-stop-B05-ALECU-save2n}{33}
\resdef{n7-s030-stop-B10-ALECU-q}{7.8}
\resdef{n7-s030-stop-B10-ALECU-hit}{80.9}
\resdef{n7-s030-stop-B10-ALECU-agree}{99.7}
\resdef{n7-s030-stop-B10-ALECU-early}{57.4}
\resdef{n7-s030-stop-B10-ALECU-exceed}{19.1}
\resdef{n7-s030-stop-B10-ALECU-qmax}{14}
\resdef{n7-s030-stop-B10-ALECU-qtot}{14.8}
\resdef{n7-s030-stop-B10-ALECU-errearly}{8.2}
\resdef{n7-s030-stop-B10-ALECU-errci}{[6.2, 10.7]}
\resdef{n7-s030-stop-B10-ALECU-regret}{1.61}
\resdef{n7-s030-stop-B10-ALECU-qround}{8}
\resdef{n7-s030-stop-B10-ALECU-fixedhit}{80.5}
\resdef{n7-s030-stop-B10-ALECU-fixedreg}{1.72}
\resdef{n7-s030-stop-B10-ALECU-save2n}{45}
\resdef{n7-s030-stop-B20-ALECU-q}{5.3}
\resdef{n7-s030-stop-B20-ALECU-hit}{80.5}
\resdef{n7-s030-stop-B20-ALECU-agree}{97.6}
\resdef{n7-s030-stop-B20-ALECU-early}{75.7}
\resdef{n7-s030-stop-B20-ALECU-exceed}{19.5}
\resdef{n7-s030-stop-B20-ALECU-qmax}{14}
\resdef{n7-s030-stop-B20-ALECU-qtot}{12.3}
\resdef{n7-s030-stop-B20-ALECU-errearly}{13.1}
\resdef{n7-s030-stop-B20-ALECU-errci}{[10.9, 15.7]}
\resdef{n7-s030-stop-B20-ALECU-regret}{1.78}
\resdef{n7-s030-stop-B20-ALECU-qround}{5}
\resdef{n7-s030-stop-B20-ALECU-fixedhit}{77.4}
\resdef{n7-s030-stop-B20-ALECU-fixedreg}{2.23}
\resdef{n7-s030-stop-B20-ALECU-save2n}{62}
\resdef{n7-s030-stop-B10-ALECU-s06-q}{10.1}
\resdef{n7-s030-stop-B10-ALECU-s06-hit}{80.8}
\resdef{n7-s030-stop-B10-ALECU-s06-agree}{100.0}
\resdef{n7-s030-stop-B10-ALECU-s06-early}{43.3}
\resdef{n7-s030-stop-B10-ALECU-s06-exceed}{19.2}
\resdef{n7-s030-stop-B10-ALECU-s06-qmax}{14}
\resdef{n7-s030-stop-B10-ALECU-s06-qtot}{17.1}
\resdef{n7-s030-stop-B10-ALECU-s06-errearly}{2.5}
\resdef{n7-s030-stop-B10-ALECU-s06-errci}{[1.4, 4.5]}
\resdef{n7-s030-stop-B10-ALECU-s06-regret}{1.61}
\resdef{n7-s030-stop-B10-ALECU-s06-qround}{10}
\resdef{n7-s030-stop-B10-ALECU-s06-fixedhit}{80.7}
\resdef{n7-s030-stop-B10-ALECU-s06-fixedreg}{1.63}
\resdef{n7-s030-stop-B10-ALECU-s06-save2n}{28}
\resdef{n7-s030-stop-B10-ALECU-cap-q}{7.8}
\resdef{n7-s030-stop-B10-ALECU-cap-hit}{80.9}
\resdef{n7-s030-stop-B10-ALECU-cap-agree}{99.7}
\resdef{n7-s030-stop-B10-ALECU-cap-early}{57.4}
\resdef{n7-s030-stop-B10-ALECU-cap-exceed}{19.1}
\resdef{n7-s030-stop-B10-ALECU-cap-qmax}{14}
\resdef{n7-s030-stop-B10-ALECU-cap-qtot}{14.8}
\resdef{n7-s030-stop-B10-ALECU-cap-errearly}{8.2}
\resdef{n7-s030-stop-B10-ALECU-cap-errci}{[6.2, 10.7]}
\resdef{n7-s030-stop-B10-ALECU-cap-regret}{1.61}
\resdef{n7-s030-stop-B10-ALECU-cap-qround}{8}
\resdef{n7-s030-stop-B10-ALECU-cap-fixedhit}{80.5}
\resdef{n7-s030-stop-B10-ALECU-cap-fixedreg}{1.72}
\resdef{n7-s030-stop-B10-ALECU-cap-save2n}{45}
\resdef{n7-s030-stop-B10-ALECU-oracle-q}{6.8}
\resdef{n7-s030-stop-B10-ALECU-oracle-hit}{80.4}
\resdef{n7-s030-stop-B10-ALECU-oracle-agree}{98.8}
\resdef{n7-s030-stop-B10-ALECU-oracle-early}{61.0}
\resdef{n7-s030-stop-B10-ALECU-oracle-exceed}{19.6}
\resdef{n7-s030-stop-B10-ALECU-oracle-qmax}{14}
\resdef{n7-s030-stop-B10-ALECU-oracle-qtot}{13.8}
\resdef{n7-s030-stop-B10-ALECU-oracle-errearly}{9.5}
\resdef{n7-s030-stop-B10-ALECU-oracle-errci}{[7.4, 12.1]}
\resdef{n7-s030-stop-B10-ALECU-oracle-regret}{1.73}
\resdef{n7-s030-stop-B10-ALECU-oracle-qround}{7}
\resdef{n7-s030-stop-B10-ALECU-oracle-fixedhit}{79.5}
\resdef{n7-s030-stop-B10-ALECU-oracle-fixedreg}{1.81}
\resdef{n7-s030-stop-B10-ALECU-oracle-save2n}{52}
\resdef{n7-s030-stop-E02-ALEC-q}{7.5}
\resdef{n7-s030-stop-E02-ALEC-hit}{80.6}
\resdef{n7-s030-stop-E02-ALEC-agree}{99.7}
\resdef{n7-s030-stop-E02-ALEC-early}{61.7}
\resdef{n7-s030-stop-E02-ALEC-exceed}{16.2}
\resdef{n7-s030-stop-E02-ALEC-qmax}{14}
\resdef{n7-s030-stop-E02-ALEC-qtot}{14.5}
\resdef{n7-s030-stop-E02-ALEC-errearly}{7.3}
\resdef{n7-s030-stop-E02-ALEC-errci}{[5.5, 9.6]}
\resdef{n7-s030-stop-E02-ALEC-regret}{1.68}
\resdef{n7-s030-stop-E02-ALEC-qround}{8}
\resdef{n7-s030-stop-E02-ALEC-fixedhit}{80.5}
\resdef{n7-s030-stop-E02-ALEC-fixedreg}{1.72}
\resdef{n7-s030-stop-E02-ALEC-save2n}{46}
\resdef{n7-s030-stop-E05-ALEC-q}{6.8}
\resdef{n7-s030-stop-E05-ALEC-hit}{80.8}
\resdef{n7-s030-stop-E05-ALEC-agree}{99.5}
\resdef{n7-s030-stop-E05-ALEC-early}{67.5}
\resdef{n7-s030-stop-E05-ALEC-exceed}{11.4}
\resdef{n7-s030-stop-E05-ALEC-qmax}{14}
\resdef{n7-s030-stop-E05-ALEC-qtot}{13.8}
\resdef{n7-s030-stop-E05-ALEC-errearly}{5.8}
\resdef{n7-s030-stop-E05-ALEC-errci}{[4.3, 7.8]}
\resdef{n7-s030-stop-E05-ALEC-regret}{1.66}
\resdef{n7-s030-stop-E05-ALEC-qround}{7}
\resdef{n7-s030-stop-E05-ALEC-fixedhit}{79.5}
\resdef{n7-s030-stop-E05-ALEC-fixedreg}{1.81}
\resdef{n7-s030-stop-E05-ALEC-save2n}{51}
\resdef{n7-s030-stop-B10-ALEC-cap-q}{8.0}
\resdef{n7-s030-stop-B10-ALEC-cap-hit}{80.6}
\resdef{n7-s030-stop-B10-ALEC-cap-agree}{99.8}
\resdef{n7-s030-stop-B10-ALEC-cap-early}{56.7}
\resdef{n7-s030-stop-B10-ALEC-cap-exceed}{19.4}
\resdef{n7-s030-stop-B10-ALEC-cap-qmax}{14}
\resdef{n7-s030-stop-B10-ALEC-cap-qtot}{15.0}
\resdef{n7-s030-stop-B10-ALEC-cap-errearly}{7.9}
\resdef{n7-s030-stop-B10-ALEC-cap-errci}{[6.0, 10.5]}
\resdef{n7-s030-stop-B10-ALEC-cap-regret}{1.65}
\resdef{n7-s030-stop-B10-ALEC-cap-qround}{8}
\resdef{n7-s030-stop-B10-ALEC-cap-fixedhit}{80.5}
\resdef{n7-s030-stop-B10-ALEC-cap-fixedreg}{1.72}
\resdef{n7-s030-stop-B10-ALEC-cap-save2n}{43}
\resdef{n7-s030-stop-E02-ALEC-cap-q}{7.5}
\resdef{n7-s030-stop-E02-ALEC-cap-hit}{80.6}
\resdef{n7-s030-stop-E02-ALEC-cap-agree}{99.7}
\resdef{n7-s030-stop-E02-ALEC-cap-early}{61.7}
\resdef{n7-s030-stop-E02-ALEC-cap-exceed}{16.2}
\resdef{n7-s030-stop-E02-ALEC-cap-qmax}{14}
\resdef{n7-s030-stop-E02-ALEC-cap-qtot}{14.5}
\resdef{n7-s030-stop-E02-ALEC-cap-errearly}{7.3}
\resdef{n7-s030-stop-E02-ALEC-cap-errci}{[5.5, 9.6]}
\resdef{n7-s030-stop-E02-ALEC-cap-regret}{1.68}
\resdef{n7-s030-stop-E02-ALEC-cap-qround}{8}
\resdef{n7-s030-stop-E02-ALEC-cap-fixedhit}{80.5}
\resdef{n7-s030-stop-E02-ALEC-cap-fixedreg}{1.72}
\resdef{n7-s030-stop-E02-ALEC-cap-save2n}{46}
\resdef{n7-s030-stop-E05-ALEC-cap-q}{6.8}
\resdef{n7-s030-stop-E05-ALEC-cap-hit}{80.8}
\resdef{n7-s030-stop-E05-ALEC-cap-agree}{99.5}
\resdef{n7-s030-stop-E05-ALEC-cap-early}{67.5}
\resdef{n7-s030-stop-E05-ALEC-cap-exceed}{11.4}
\resdef{n7-s030-stop-E05-ALEC-cap-qmax}{14}
\resdef{n7-s030-stop-E05-ALEC-cap-qtot}{13.8}
\resdef{n7-s030-stop-E05-ALEC-cap-errearly}{5.8}
\resdef{n7-s030-stop-E05-ALEC-cap-errci}{[4.3, 7.8]}
\resdef{n7-s030-stop-E05-ALEC-cap-regret}{1.66}
\resdef{n7-s030-stop-E05-ALEC-cap-qround}{7}
\resdef{n7-s030-stop-E05-ALEC-cap-fixedhit}{79.5}
\resdef{n7-s030-stop-E05-ALEC-cap-fixedreg}{1.81}
\resdef{n7-s030-stop-E05-ALEC-cap-save2n}{51}
\resdef{n7-s030-stop-E05-ALEC-cap-s06-q}{9.5}
\resdef{n7-s030-stop-E05-ALEC-cap-s06-hit}{80.8}
\resdef{n7-s030-stop-E05-ALEC-cap-s06-agree}{100.0}
\resdef{n7-s030-stop-E05-ALEC-cap-s06-early}{52.0}
\resdef{n7-s030-stop-E05-ALEC-cap-s06-exceed}{11.3}
\resdef{n7-s030-stop-E05-ALEC-cap-s06-qmax}{14}
\resdef{n7-s030-stop-E05-ALEC-cap-s06-qtot}{16.5}
\resdef{n7-s030-stop-E05-ALEC-cap-s06-errearly}{1.5}
\resdef{n7-s030-stop-E05-ALEC-cap-s06-errci}{[0.8, 3.0]}
\resdef{n7-s030-stop-E05-ALEC-cap-s06-regret}{1.61}
\resdef{n7-s030-stop-E05-ALEC-cap-s06-qround}{9}
\resdef{n7-s030-stop-E05-ALEC-cap-s06-fixedhit}{81.1}
\resdef{n7-s030-stop-E05-ALEC-cap-s06-fixedreg}{1.59}
\resdef{n7-s030-stop-E05-ALEC-cap-s06-save2n}{32}
\resdef{n7-s050-full-hit}{72.5}
\resdef{n7-s050-tot-h}{11}
\resdef{n7-s050-frac-h}{52}
\resdef{n7-s050-RANDOM-hit-h}{63.8}
\resdef{n7-s050-RANDOM-hitllsm-h}{63.8}
\resdef{n7-s050-RANDOM-tau-h}{0.731}
\resdef{n7-s050-RANDOM-reg-h}{6.95}
\resdef{n7-s050-RANDOM-agree-h}{73.2}
\resdef{n7-s050-p-RANDOM-h}{$p=0.026$}
\resdef{n7-s050-p-RANDOM-h-llsm}{$p=0.026$}
\resdef{n7-s050-HARKER-PG-hit-h}{65.4}
\resdef{n7-s050-HARKER-PG-hitllsm-h}{65.2}
\resdef{n7-s050-HARKER-PG-tau-h}{0.728}
\resdef{n7-s050-HARKER-PG-reg-h}{6.56}
\resdef{n7-s050-HARKER-PG-agree-h}{72.9}
\resdef{n7-s050-p-HARKER-PG-h}{$p=0.264$}
\resdef{n7-s050-p-HARKER-PG-h-llsm}{$p=0.200$}
\resdef{n7-s050-HARKER-hit-h}{63.7}
\resdef{n7-s050-HARKER-hitllsm-h}{65.5}
\resdef{n7-s050-HARKER-tau-h}{0.716}
\resdef{n7-s050-HARKER-reg-h}{6.63}
\resdef{n7-s050-HARKER-agree-h}{72.8}
\resdef{n7-s050-p-HARKER-h}{$p=0.027$}
\resdef{n7-s050-p-HARKER-h-llsm}{$p=0.282$}
\resdef{n7-s050-HARKER-LOG-hit-h}{66.3}
\resdef{n7-s050-HARKER-LOG-hitllsm-h}{67.2}
\resdef{n7-s050-HARKER-LOG-tau-h}{0.697}
\resdef{n7-s050-HARKER-LOG-reg-h}{5.61}
\resdef{n7-s050-HARKER-LOG-agree-h}{80.6}
\resdef{n7-s050-p-HARKER-LOG-h}{$p=0.512$}
\resdef{n7-s050-p-HARKER-LOG-h-llsm}{$p=1.000$}
\resdef{n7-s050-AOPT-hit-h}{65.4}
\resdef{n7-s050-AOPT-hitllsm-h}{65.4}
\resdef{n7-s050-AOPT-tau-h}{0.747}
\resdef{n7-s050-AOPT-reg-h}{5.83}
\resdef{n7-s050-AOPT-agree-h}{74.1}
\resdef{n7-s050-p-AOPT-h}{$p=0.238$}
\resdef{n7-s050-p-AOPT-h-llsm}{$p=0.238$}
\resdef{n7-s050-FROZEN-hit-h}{65.5}
\resdef{n7-s050-FROZEN-hitllsm-h}{65.5}
\resdef{n7-s050-FROZEN-tau-h}{0.716}
\resdef{n7-s050-FROZEN-reg-h}{5.76}
\resdef{n7-s050-FROZEN-agree-h}{77.6}
\resdef{n7-s050-p-FROZEN-h}{$p=0.086$}
\resdef{n7-s050-p-FROZEN-h-llsm}{$p=0.086$}
\resdef{n7-s050-STAR-hit-h}{65.4}
\resdef{n7-s050-STAR-hitllsm-h}{65.4}
\resdef{n7-s050-STAR-tau-h}{0.716}
\resdef{n7-s050-STAR-reg-h}{5.77}
\resdef{n7-s050-STAR-agree-h}{77.2}
\resdef{n7-s050-p-STAR-h}{$p=0.060$}
\resdef{n7-s050-p-STAR-h-llsm}{$p=0.060$}
\resdef{n7-s050-ALEC-hit-h}{67.2}
\resdef{n7-s050-ALEC-hitllsm-h}{67.2}
\resdef{n7-s050-ALEC-tau-h}{0.714}
\resdef{n7-s050-ALEC-reg-h}{5.26}
\resdef{n7-s050-ALEC-agree-h}{81.0}
\resdef{n7-s050-ALEC-U-hit-h}{67.9}
\resdef{n7-s050-ALEC-U-hitllsm-h}{67.9}
\resdef{n7-s050-ALEC-U-tau-h}{0.725}
\resdef{n7-s050-ALEC-U-reg-h}{4.94}
\resdef{n7-s050-ALEC-U-agree-h}{81.8}
\resdef{n7-s050-p-ALEC-U-h}{$p=0.576$}
\resdef{n7-s050-p-ALEC-U-h-llsm}{$p=0.576$}
\resdef{n7-s050-tot-n}{14}
\resdef{n7-s050-frac-n}{67}
\resdef{n7-s050-RANDOM-hit-n}{67.2}
\resdef{n7-s050-RANDOM-hitllsm-n}{67.2}
\resdef{n7-s050-RANDOM-tau-n}{0.768}
\resdef{n7-s050-RANDOM-reg-n}{5.49}
\resdef{n7-s050-RANDOM-agree-n}{80.7}
\resdef{n7-s050-p-RANDOM-n}{$p=0.004$}
\resdef{n7-s050-p-RANDOM-n-llsm}{$p=0.004$}
\resdef{n7-s050-HARKER-PG-hit-n}{68.5}
\resdef{n7-s050-HARKER-PG-hitllsm-n}{69.0}
\resdef{n7-s050-HARKER-PG-tau-n}{0.767}
\resdef{n7-s050-HARKER-PG-reg-n}{4.86}
\resdef{n7-s050-HARKER-PG-agree-n}{80.6}
\resdef{n7-s050-p-HARKER-PG-n}{$p=0.044$}
\resdef{n7-s050-p-HARKER-PG-n-llsm}{$p=0.097$}
\resdef{n7-s050-HARKER-hit-n}{68.7}
\resdef{n7-s050-HARKER-hitllsm-n}{69.7}
\resdef{n7-s050-HARKER-tau-n}{0.762}
\resdef{n7-s050-HARKER-reg-n}{4.83}
\resdef{n7-s050-HARKER-agree-n}{80.7}
\resdef{n7-s050-p-HARKER-n}{$p=0.061$}
\resdef{n7-s050-p-HARKER-n-llsm}{$p=0.251$}
\resdef{n7-s050-HARKER-LOG-hit-n}{69.6}
\resdef{n7-s050-HARKER-LOG-hitllsm-n}{72.1}
\resdef{n7-s050-HARKER-LOG-tau-n}{0.743}
\resdef{n7-s050-HARKER-LOG-reg-n}{4.29}
\resdef{n7-s050-HARKER-LOG-agree-n}{87.9}
\resdef{n7-s050-p-HARKER-LOG-n}{$p=0.145$}
\resdef{n7-s050-p-HARKER-LOG-n-llsm}{$p=0.368$}
\resdef{n7-s050-AOPT-hit-n}{67.7}
\resdef{n7-s050-AOPT-hitllsm-n}{67.7}
\resdef{n7-s050-AOPT-tau-n}{0.776}
\resdef{n7-s050-AOPT-reg-n}{4.93}
\resdef{n7-s050-AOPT-agree-n}{80.7}
\resdef{n7-s050-p-AOPT-n}{$p=0.010$}
\resdef{n7-s050-p-AOPT-n-llsm}{$p=0.010$}
\resdef{n7-s050-FROZEN-hit-n}{70.0}
\resdef{n7-s050-FROZEN-hitllsm-n}{70.0}
\resdef{n7-s050-FROZEN-tau-n}{0.761}
\resdef{n7-s050-FROZEN-reg-n}{4.25}
\resdef{n7-s050-FROZEN-agree-n}{89.0}
\resdef{n7-s050-p-FROZEN-n}{$p=0.201$}
\resdef{n7-s050-p-FROZEN-n-llsm}{$p=0.201$}
\resdef{n7-s050-STAR-hit-n}{70.4}
\resdef{n7-s050-STAR-hitllsm-n}{70.4}
\resdef{n7-s050-STAR-tau-n}{0.761}
\resdef{n7-s050-STAR-reg-n}{4.07}
\resdef{n7-s050-STAR-agree-n}{90.9}
\resdef{n7-s050-p-STAR-n}{$p=0.229$}
\resdef{n7-s050-p-STAR-n-llsm}{$p=0.229$}
\resdef{n7-s050-ALEC-hit-n}{71.2}
\resdef{n7-s050-ALEC-hitllsm-n}{71.2}
\resdef{n7-s050-ALEC-tau-n}{0.762}
\resdef{n7-s050-ALEC-reg-n}{3.83}
\resdef{n7-s050-ALEC-agree-n}{92.3}
\resdef{n7-s050-ALEC-U-hit-n}{70.2}
\resdef{n7-s050-ALEC-U-hitllsm-n}{70.2}
\resdef{n7-s050-ALEC-U-tau-n}{0.768}
\resdef{n7-s050-ALEC-U-reg-n}{4.19}
\resdef{n7-s050-ALEC-U-agree-n}{88.9}
\resdef{n7-s050-p-ALEC-U-n}{$p=0.348$}
\resdef{n7-s050-p-ALEC-U-n-llsm}{$p=0.348$}
\resdef{n7-s050-tot-2n}{21}
\resdef{n7-s050-frac-2n}{100}
\resdef{n7-s050-RANDOM-hit-2n}{72.5}
\resdef{n7-s050-RANDOM-hitllsm-2n}{72.5}
\resdef{n7-s050-RANDOM-tau-2n}{0.821}
\resdef{n7-s050-RANDOM-reg-2n}{3.48}
\resdef{n7-s050-RANDOM-agree-2n}{100.0}
\resdef{n7-s050-p-RANDOM-2n}{$p=1.000$}
\resdef{n7-s050-p-RANDOM-2n-llsm}{$p=1.000$}
\resdef{n7-s050-HARKER-PG-hit-2n}{72.5}
\resdef{n7-s050-HARKER-PG-hitllsm-2n}{72.5}
\resdef{n7-s050-HARKER-PG-tau-2n}{0.815}
\resdef{n7-s050-HARKER-PG-reg-2n}{3.50}
\resdef{n7-s050-HARKER-PG-agree-2n}{94.9}
\resdef{n7-s050-p-HARKER-PG-2n}{$p=1.000$}
\resdef{n7-s050-p-HARKER-PG-2n-llsm}{$p=1.000$}
\resdef{n7-s050-HARKER-hit-2n}{72.5}
\resdef{n7-s050-HARKER-hitllsm-2n}{72.5}
\resdef{n7-s050-HARKER-tau-2n}{0.815}
\resdef{n7-s050-HARKER-reg-2n}{3.50}
\resdef{n7-s050-HARKER-agree-2n}{94.9}
\resdef{n7-s050-p-HARKER-2n}{$p=1.000$}
\resdef{n7-s050-p-HARKER-2n-llsm}{$p=1.000$}
\resdef{n7-s050-HARKER-LOG-hit-2n}{72.5}
\resdef{n7-s050-HARKER-LOG-hitllsm-2n}{72.5}
\resdef{n7-s050-HARKER-LOG-tau-2n}{0.815}
\resdef{n7-s050-HARKER-LOG-reg-2n}{3.50}
\resdef{n7-s050-HARKER-LOG-agree-2n}{94.9}
\resdef{n7-s050-p-HARKER-LOG-2n}{$p=1.000$}
\resdef{n7-s050-p-HARKER-LOG-2n-llsm}{$p=1.000$}
\resdef{n7-s050-AOPT-hit-2n}{72.5}
\resdef{n7-s050-AOPT-hitllsm-2n}{72.5}
\resdef{n7-s050-AOPT-tau-2n}{0.821}
\resdef{n7-s050-AOPT-reg-2n}{3.48}
\resdef{n7-s050-AOPT-agree-2n}{100.0}
\resdef{n7-s050-p-AOPT-2n}{$p=1.000$}
\resdef{n7-s050-p-AOPT-2n-llsm}{$p=1.000$}
\resdef{n7-s050-FROZEN-hit-2n}{72.5}
\resdef{n7-s050-FROZEN-hitllsm-2n}{72.5}
\resdef{n7-s050-FROZEN-tau-2n}{0.821}
\resdef{n7-s050-FROZEN-reg-2n}{3.48}
\resdef{n7-s050-FROZEN-agree-2n}{100.0}
\resdef{n7-s050-p-FROZEN-2n}{$p=1.000$}
\resdef{n7-s050-p-FROZEN-2n-llsm}{$p=1.000$}
\resdef{n7-s050-STAR-hit-2n}{72.5}
\resdef{n7-s050-STAR-hitllsm-2n}{72.5}
\resdef{n7-s050-STAR-tau-2n}{0.821}
\resdef{n7-s050-STAR-reg-2n}{3.48}
\resdef{n7-s050-STAR-agree-2n}{100.0}
\resdef{n7-s050-p-STAR-2n}{$p=1.000$}
\resdef{n7-s050-p-STAR-2n-llsm}{$p=1.000$}
\resdef{n7-s050-ALEC-hit-2n}{72.5}
\resdef{n7-s050-ALEC-hitllsm-2n}{72.5}
\resdef{n7-s050-ALEC-tau-2n}{0.821}
\resdef{n7-s050-ALEC-reg-2n}{3.48}
\resdef{n7-s050-ALEC-agree-2n}{100.0}
\resdef{n7-s050-ALEC-U-hit-2n}{72.5}
\resdef{n7-s050-ALEC-U-hitllsm-2n}{72.5}
\resdef{n7-s050-ALEC-U-tau-2n}{0.821}
\resdef{n7-s050-ALEC-U-reg-2n}{3.48}
\resdef{n7-s050-ALEC-U-agree-2n}{100.0}
\resdef{n7-s050-p-ALEC-U-2n}{$p=1.000$}
\resdef{n7-s050-p-ALEC-U-2n-llsm}{$p=1.000$}
\resdef{n7-s050-stop-LC-q}{4.6}
\resdef{n7-s050-stop-LC-hit}{70.6}
\resdef{n7-s050-stop-LC-agree}{85.6}
\resdef{n7-s050-stop-LC-early}{100.0}
\resdef{n7-s050-stop-LC-exceed}{29.4}
\resdef{n7-s050-stop-LC-qmax}{10}
\resdef{n7-s050-stop-LC-qtot}{11.6}
\resdef{n7-s050-stop-LC-errearly}{29.4}
\resdef{n7-s050-stop-LC-errci}{[26.7, 32.3]}
\resdef{n7-s050-stop-LC-regret}{4.32}
\resdef{n7-s050-stop-LC-qround}{5}
\resdef{n7-s050-stop-LC-fixedhit}{69.0}
\resdef{n7-s050-stop-LC-fixedreg}{4.88}
\resdef{n7-s050-stop-LC-save2n}{67}
\resdef{n7-s050-stop-B10-ALEC-q}{9.0}
\resdef{n7-s050-stop-B10-ALEC-hit}{72.1}
\resdef{n7-s050-stop-B10-ALEC-agree}{97.9}
\resdef{n7-s050-stop-B10-ALEC-early}{47.4}
\resdef{n7-s050-stop-B10-ALEC-exceed}{27.9}
\resdef{n7-s050-stop-B10-ALEC-qmax}{14}
\resdef{n7-s050-stop-B10-ALEC-qtot}{16.0}
\resdef{n7-s050-stop-B10-ALEC-errearly}{14.8}
\resdef{n7-s050-stop-B10-ALEC-errci}{[11.9, 18.2]}
\resdef{n7-s050-stop-B10-ALEC-regret}{3.70}
\resdef{n7-s050-stop-B10-ALEC-qround}{9}
\resdef{n7-s050-stop-B10-ALEC-fixedhit}{71.5}
\resdef{n7-s050-stop-B10-ALEC-fixedreg}{3.62}
\resdef{n7-s050-stop-B10-ALEC-save2n}{36}
\resdef{n7-s050-stop-B05-ALECU-q}{10.5}
\resdef{n7-s050-stop-B05-ALECU-hit}{72.4}
\resdef{n7-s050-stop-B05-ALECU-agree}{99.1}
\resdef{n7-s050-stop-B05-ALECU-early}{34.6}
\resdef{n7-s050-stop-B05-ALECU-exceed}{27.6}
\resdef{n7-s050-stop-B05-ALECU-qmax}{14}
\resdef{n7-s050-stop-B05-ALECU-qtot}{17.5}
\resdef{n7-s050-stop-B05-ALECU-errearly}{9.5}
\resdef{n7-s050-stop-B05-ALECU-errci}{[6.9, 13.1]}
\resdef{n7-s050-stop-B05-ALECU-regret}{3.55}
\resdef{n7-s050-stop-B05-ALECU-qround}{10}
\resdef{n7-s050-stop-B05-ALECU-fixedhit}{72.0}
\resdef{n7-s050-stop-B05-ALECU-fixedreg}{3.53}
\resdef{n7-s050-stop-B05-ALECU-save2n}{25}
\resdef{n7-s050-stop-B10-ALECU-q}{8.6}
\resdef{n7-s050-stop-B10-ALECU-hit}{72.2}
\resdef{n7-s050-stop-B10-ALECU-agree}{97.3}
\resdef{n7-s050-stop-B10-ALECU-early}{49.8}
\resdef{n7-s050-stop-B10-ALECU-exceed}{27.8}
\resdef{n7-s050-stop-B10-ALECU-qmax}{14}
\resdef{n7-s050-stop-B10-ALECU-qtot}{15.6}
\resdef{n7-s050-stop-B10-ALECU-errearly}{15.1}
\resdef{n7-s050-stop-B10-ALECU-errci}{[12.2, 18.5]}
\resdef{n7-s050-stop-B10-ALECU-regret}{3.75}
\resdef{n7-s050-stop-B10-ALECU-qround}{9}
\resdef{n7-s050-stop-B10-ALECU-fixedhit}{71.5}
\resdef{n7-s050-stop-B10-ALECU-fixedreg}{3.62}
\resdef{n7-s050-stop-B10-ALECU-save2n}{39}
\resdef{n7-s050-stop-B20-ALECU-q}{5.6}
\resdef{n7-s050-stop-B20-ALECU-hit}{69.9}
\resdef{n7-s050-stop-B20-ALECU-agree}{91.4}
\resdef{n7-s050-stop-B20-ALECU-early}{71.3}
\resdef{n7-s050-stop-B20-ALECU-exceed}{30.1}
\resdef{n7-s050-stop-B20-ALECU-qmax}{14}
\resdef{n7-s050-stop-B20-ALECU-qtot}{12.6}
\resdef{n7-s050-stop-B20-ALECU-errearly}{23.0}
\resdef{n7-s050-stop-B20-ALECU-errci}{[20.1, 26.2]}
\resdef{n7-s050-stop-B20-ALECU-regret}{4.56}
\resdef{n7-s050-stop-B20-ALECU-qround}{6}
\resdef{n7-s050-stop-B20-ALECU-fixedhit}{69.7}
\resdef{n7-s050-stop-B20-ALECU-fixedreg}{4.42}
\resdef{n7-s050-stop-B20-ALECU-save2n}{60}
\resdef{n7-s050-stop-B10-ALECU-s06-q}{10.9}
\resdef{n7-s050-stop-B10-ALECU-s06-hit}{72.3}
\resdef{n7-s050-stop-B10-ALECU-s06-agree}{99.8}
\resdef{n7-s050-stop-B10-ALECU-s06-early}{34.2}
\resdef{n7-s050-stop-B10-ALECU-s06-exceed}{27.7}
\resdef{n7-s050-stop-B10-ALECU-s06-qmax}{14}
\resdef{n7-s050-stop-B10-ALECU-s06-qtot}{17.9}
\resdef{n7-s050-stop-B10-ALECU-s06-errearly}{9.1}
\resdef{n7-s050-stop-B10-ALECU-s06-errci}{[6.5, 12.6]}
\resdef{n7-s050-stop-B10-ALECU-s06-regret}{3.50}
\resdef{n7-s050-stop-B10-ALECU-s06-qround}{11}
\resdef{n7-s050-stop-B10-ALECU-s06-fixedhit}{72.4}
\resdef{n7-s050-stop-B10-ALECU-s06-fixedreg}{3.51}
\resdef{n7-s050-stop-B10-ALECU-s06-save2n}{22}
\resdef{n7-s050-stop-B10-ALECU-cap-q}{8.6}
\resdef{n7-s050-stop-B10-ALECU-cap-hit}{72.2}
\resdef{n7-s050-stop-B10-ALECU-cap-agree}{97.3}
\resdef{n7-s050-stop-B10-ALECU-cap-early}{49.8}
\resdef{n7-s050-stop-B10-ALECU-cap-exceed}{27.8}
\resdef{n7-s050-stop-B10-ALECU-cap-qmax}{14}
\resdef{n7-s050-stop-B10-ALECU-cap-qtot}{15.6}
\resdef{n7-s050-stop-B10-ALECU-cap-errearly}{15.1}
\resdef{n7-s050-stop-B10-ALECU-cap-errci}{[12.2, 18.5]}
\resdef{n7-s050-stop-B10-ALECU-cap-regret}{3.75}
\resdef{n7-s050-stop-B10-ALECU-cap-qround}{9}
\resdef{n7-s050-stop-B10-ALECU-cap-fixedhit}{71.5}
\resdef{n7-s050-stop-B10-ALECU-cap-fixedreg}{3.62}
\resdef{n7-s050-stop-B10-ALECU-cap-save2n}{39}
\resdef{n7-s050-stop-B10-ALECU-oracle-q}{9.1}
\resdef{n7-s050-stop-B10-ALECU-oracle-hit}{72.2}
\resdef{n7-s050-stop-B10-ALECU-oracle-agree}{98.4}
\resdef{n7-s050-stop-B10-ALECU-oracle-early}{44.1}
\resdef{n7-s050-stop-B10-ALECU-oracle-exceed}{27.8}
\resdef{n7-s050-stop-B10-ALECU-oracle-qmax}{14}
\resdef{n7-s050-stop-B10-ALECU-oracle-qtot}{16.1}
\resdef{n7-s050-stop-B10-ALECU-oracle-errearly}{11.6}
\resdef{n7-s050-stop-B10-ALECU-oracle-errci}{[8.9, 14.9]}
\resdef{n7-s050-stop-B10-ALECU-oracle-regret}{3.67}
\resdef{n7-s050-stop-B10-ALECU-oracle-qround}{9}
\resdef{n7-s050-stop-B10-ALECU-oracle-fixedhit}{71.5}
\resdef{n7-s050-stop-B10-ALECU-oracle-fixedreg}{3.62}
\resdef{n7-s050-stop-B10-ALECU-oracle-save2n}{35}
\resdef{n7-s050-stop-E02-ALEC-q}{8.6}
\resdef{n7-s050-stop-E02-ALEC-hit}{71.7}
\resdef{n7-s050-stop-E02-ALEC-agree}{97.5}
\resdef{n7-s050-stop-E02-ALEC-early}{50.8}
\resdef{n7-s050-stop-E02-ALEC-exceed}{25.0}
\resdef{n7-s050-stop-E02-ALEC-qmax}{14}
\resdef{n7-s050-stop-E02-ALEC-qtot}{15.6}
\resdef{n7-s050-stop-E02-ALEC-errearly}{13.8}
\resdef{n7-s050-stop-E02-ALEC-errci}{[11.1, 17.0]}
\resdef{n7-s050-stop-E02-ALEC-regret}{3.80}
\resdef{n7-s050-stop-E02-ALEC-qround}{9}
\resdef{n7-s050-stop-E02-ALEC-fixedhit}{71.5}
\resdef{n7-s050-stop-E02-ALEC-fixedreg}{3.62}
\resdef{n7-s050-stop-E02-ALEC-save2n}{39}
\resdef{n7-s050-stop-E05-ALEC-q}{8.0}
\resdef{n7-s050-stop-E05-ALEC-hit}{71.4}
\resdef{n7-s050-stop-E05-ALEC-agree}{97.0}
\resdef{n7-s050-stop-E05-ALEC-early}{56.3}
\resdef{n7-s050-stop-E05-ALEC-exceed}{21.5}
\resdef{n7-s050-stop-E05-ALEC-qmax}{14}
\resdef{n7-s050-stop-E05-ALEC-qtot}{15.0}
\resdef{n7-s050-stop-E05-ALEC-errearly}{13.1}
\resdef{n7-s050-stop-E05-ALEC-errci}{[10.6, 16.2]}
\resdef{n7-s050-stop-E05-ALEC-regret}{3.92}
\resdef{n7-s050-stop-E05-ALEC-qround}{8}
\resdef{n7-s050-stop-E05-ALEC-fixedhit}{71.8}
\resdef{n7-s050-stop-E05-ALEC-fixedreg}{3.67}
\resdef{n7-s050-stop-E05-ALEC-save2n}{43}
\resdef{n7-s050-stop-B10-ALEC-cap-q}{9.0}
\resdef{n7-s050-stop-B10-ALEC-cap-hit}{72.1}
\resdef{n7-s050-stop-B10-ALEC-cap-agree}{97.9}
\resdef{n7-s050-stop-B10-ALEC-cap-early}{47.4}
\resdef{n7-s050-stop-B10-ALEC-cap-exceed}{27.9}
\resdef{n7-s050-stop-B10-ALEC-cap-qmax}{14}
\resdef{n7-s050-stop-B10-ALEC-cap-qtot}{16.0}
\resdef{n7-s050-stop-B10-ALEC-cap-errearly}{14.8}
\resdef{n7-s050-stop-B10-ALEC-cap-errci}{[11.9, 18.2]}
\resdef{n7-s050-stop-B10-ALEC-cap-regret}{3.70}
\resdef{n7-s050-stop-B10-ALEC-cap-qround}{9}
\resdef{n7-s050-stop-B10-ALEC-cap-fixedhit}{71.5}
\resdef{n7-s050-stop-B10-ALEC-cap-fixedreg}{3.62}
\resdef{n7-s050-stop-B10-ALEC-cap-save2n}{36}
\resdef{n7-s050-stop-E02-ALEC-cap-q}{8.6}
\resdef{n7-s050-stop-E02-ALEC-cap-hit}{71.7}
\resdef{n7-s050-stop-E02-ALEC-cap-agree}{97.5}
\resdef{n7-s050-stop-E02-ALEC-cap-early}{50.8}
\resdef{n7-s050-stop-E02-ALEC-cap-exceed}{25.0}
\resdef{n7-s050-stop-E02-ALEC-cap-qmax}{14}
\resdef{n7-s050-stop-E02-ALEC-cap-qtot}{15.6}
\resdef{n7-s050-stop-E02-ALEC-cap-errearly}{13.8}
\resdef{n7-s050-stop-E02-ALEC-cap-errci}{[11.1, 17.0]}
\resdef{n7-s050-stop-E02-ALEC-cap-regret}{3.80}
\resdef{n7-s050-stop-E02-ALEC-cap-qround}{9}
\resdef{n7-s050-stop-E02-ALEC-cap-fixedhit}{71.5}
\resdef{n7-s050-stop-E02-ALEC-cap-fixedreg}{3.62}
\resdef{n7-s050-stop-E02-ALEC-cap-save2n}{39}
\resdef{n7-s050-stop-E05-ALEC-cap-q}{8.0}
\resdef{n7-s050-stop-E05-ALEC-cap-hit}{71.4}
\resdef{n7-s050-stop-E05-ALEC-cap-agree}{97.0}
\resdef{n7-s050-stop-E05-ALEC-cap-early}{56.3}
\resdef{n7-s050-stop-E05-ALEC-cap-exceed}{21.5}
\resdef{n7-s050-stop-E05-ALEC-cap-qmax}{14}
\resdef{n7-s050-stop-E05-ALEC-cap-qtot}{15.0}
\resdef{n7-s050-stop-E05-ALEC-cap-errearly}{13.1}
\resdef{n7-s050-stop-E05-ALEC-cap-errci}{[10.6, 16.2]}
\resdef{n7-s050-stop-E05-ALEC-cap-regret}{3.92}
\resdef{n7-s050-stop-E05-ALEC-cap-qround}{8}
\resdef{n7-s050-stop-E05-ALEC-cap-fixedhit}{71.8}
\resdef{n7-s050-stop-E05-ALEC-cap-fixedreg}{3.67}
\resdef{n7-s050-stop-E05-ALEC-cap-save2n}{43}
\resdef{n7-s050-stop-E05-ALEC-cap-s06-q}{10.5}
\resdef{n7-s050-stop-E05-ALEC-cap-s06-hit}{72.4}
\resdef{n7-s050-stop-E05-ALEC-cap-s06-agree}{99.9}
\resdef{n7-s050-stop-E05-ALEC-cap-s06-early}{39.6}
\resdef{n7-s050-stop-E05-ALEC-cap-s06-exceed}{20.4}
\resdef{n7-s050-stop-E05-ALEC-cap-s06-qmax}{14}
\resdef{n7-s050-stop-E05-ALEC-cap-s06-qtot}{17.5}
\resdef{n7-s050-stop-E05-ALEC-cap-s06-errearly}{6.1}
\resdef{n7-s050-stop-E05-ALEC-cap-s06-errci}{[4.1, 8.9]}
\resdef{n7-s050-stop-E05-ALEC-cap-s06-regret}{3.49}
\resdef{n7-s050-stop-E05-ALEC-cap-s06-qround}{11}
\resdef{n7-s050-stop-E05-ALEC-cap-s06-fixedhit}{72.4}
\resdef{n7-s050-stop-E05-ALEC-cap-s06-fixedreg}{3.51}
\resdef{n7-s050-stop-E05-ALEC-cap-s06-save2n}{25}
\resdef{n11-s015-full-hit}{88.9}
\resdef{n11-s015-tot-h}{17}
\resdef{n11-s015-frac-h}{31}
\resdef{n11-s015-RANDOM-hit-h}{78.1}
\resdef{n11-s015-RANDOM-hitllsm-h}{78.1}
\resdef{n11-s015-RANDOM-tau-h}{0.899}
\resdef{n11-s015-RANDOM-reg-h}{1.68}
\resdef{n11-s015-RANDOM-agree-h}{81.7}
\resdef{n11-s015-p-RANDOM-h}{$p<0.001$}
\resdef{n11-s015-p-RANDOM-h-llsm}{$p<0.001$}
\resdef{n11-s015-HARKER-PG-hit-h}{78.4}
\resdef{n11-s015-HARKER-PG-hitllsm-h}{78.4}
\resdef{n11-s015-HARKER-PG-tau-h}{0.890}
\resdef{n11-s015-HARKER-PG-reg-h}{1.66}
\resdef{n11-s015-HARKER-PG-agree-h}{82.4}
\resdef{n11-s015-p-HARKER-PG-h}{$p<0.001$}
\resdef{n11-s015-p-HARKER-PG-h-llsm}{$p<0.001$}
\resdef{n11-s015-HARKER-hit-h}{78.3}
\resdef{n11-s015-HARKER-hitllsm-h}{78.2}
\resdef{n11-s015-HARKER-tau-h}{0.894}
\resdef{n11-s015-HARKER-reg-h}{1.59}
\resdef{n11-s015-HARKER-agree-h}{82.3}
\resdef{n11-s015-p-HARKER-h}{$p<0.001$}
\resdef{n11-s015-p-HARKER-h-llsm}{$p<0.001$}
\resdef{n11-s015-HARKER-LOG-hit-h}{82.5}
\resdef{n11-s015-HARKER-LOG-hitllsm-h}{82.7}
\resdef{n11-s015-HARKER-LOG-tau-h}{0.894}
\resdef{n11-s015-HARKER-LOG-reg-h}{1.01}
\resdef{n11-s015-HARKER-LOG-agree-h}{88.7}
\resdef{n11-s015-p-HARKER-LOG-h}{$p=0.223$}
\resdef{n11-s015-p-HARKER-LOG-h-llsm}{$p=0.294$}
\resdef{n11-s015-AOPT-hit-h}{80.4}
\resdef{n11-s015-AOPT-hitllsm-h}{80.4}
\resdef{n11-s015-AOPT-tau-h}{0.908}
\resdef{n11-s015-AOPT-reg-h}{1.38}
\resdef{n11-s015-AOPT-agree-h}{83.0}
\resdef{n11-s015-p-AOPT-h}{$p=0.012$}
\resdef{n11-s015-p-AOPT-h-llsm}{$p=0.012$}
\resdef{n11-s015-FROZEN-hit-h}{83.1}
\resdef{n11-s015-FROZEN-hitllsm-h}{83.1}
\resdef{n11-s015-FROZEN-tau-h}{0.899}
\resdef{n11-s015-FROZEN-reg-h}{0.98}
\resdef{n11-s015-FROZEN-agree-h}{88.6}
\resdef{n11-s015-p-FROZEN-h}{$p=0.477$}
\resdef{n11-s015-p-FROZEN-h-llsm}{$p=0.477$}
\resdef{n11-s015-STAR-hit-h}{83.0}
\resdef{n11-s015-STAR-hitllsm-h}{83.0}
\resdef{n11-s015-STAR-tau-h}{0.899}
\resdef{n11-s015-STAR-reg-h}{0.99}
\resdef{n11-s015-STAR-agree-h}{88.5}
\resdef{n11-s015-p-STAR-h}{$p=0.410$}
\resdef{n11-s015-p-STAR-h-llsm}{$p=0.410$}
\resdef{n11-s015-ALEC-hit-h}{83.8}
\resdef{n11-s015-ALEC-hitllsm-h}{83.8}
\resdef{n11-s015-ALEC-tau-h}{0.897}
\resdef{n11-s015-ALEC-reg-h}{0.93}
\resdef{n11-s015-ALEC-agree-h}{88.6}
\resdef{n11-s015-ALEC-U-hit-h}{83.8}
\resdef{n11-s015-ALEC-U-hitllsm-h}{83.8}
\resdef{n11-s015-ALEC-U-tau-h}{0.894}
\resdef{n11-s015-ALEC-U-reg-h}{0.92}
\resdef{n11-s015-ALEC-U-agree-h}{89.3}
\resdef{n11-s015-p-ALEC-U-h}{$p=1.000$}
\resdef{n11-s015-p-ALEC-U-h-llsm}{$p=1.000$}
\resdef{n11-s015-tot-n}{22}
\resdef{n11-s015-frac-n}{40}
\resdef{n11-s015-RANDOM-hit-n}{80.2}
\resdef{n11-s015-RANDOM-hitllsm-n}{80.2}
\resdef{n11-s015-RANDOM-tau-n}{0.917}
\resdef{n11-s015-RANDOM-reg-n}{1.29}
\resdef{n11-s015-RANDOM-agree-n}{84.4}
\resdef{n11-s015-p-RANDOM-n}{$p<0.001$}
\resdef{n11-s015-p-RANDOM-n-llsm}{$p<0.001$}
\resdef{n11-s015-HARKER-PG-hit-n}{83.0}
\resdef{n11-s015-HARKER-PG-hitllsm-n}{82.4}
\resdef{n11-s015-HARKER-PG-tau-n}{0.908}
\resdef{n11-s015-HARKER-PG-reg-n}{1.11}
\resdef{n11-s015-HARKER-PG-agree-n}{87.4}
\resdef{n11-s015-p-HARKER-PG-n}{$p=0.056$}
\resdef{n11-s015-p-HARKER-PG-n-llsm}{$p=0.016$}
\resdef{n11-s015-HARKER-hit-n}{82.3}
\resdef{n11-s015-HARKER-hitllsm-n}{82.2}
\resdef{n11-s015-HARKER-tau-n}{0.911}
\resdef{n11-s015-HARKER-reg-n}{1.09}
\resdef{n11-s015-HARKER-agree-n}{86.8}
\resdef{n11-s015-p-HARKER-n}{$p=0.010$}
\resdef{n11-s015-p-HARKER-n-llsm}{$p=0.009$}
\resdef{n11-s015-HARKER-LOG-hit-n}{85.1}
\resdef{n11-s015-HARKER-LOG-hitllsm-n}{85.1}
\resdef{n11-s015-HARKER-LOG-tau-n}{0.911}
\resdef{n11-s015-HARKER-LOG-reg-n}{0.68}
\resdef{n11-s015-HARKER-LOG-agree-n}{92.0}
\resdef{n11-s015-p-HARKER-LOG-n}{$p=0.899$}
\resdef{n11-s015-p-HARKER-LOG-n-llsm}{$p=0.899$}
\resdef{n11-s015-AOPT-hit-n}{82.7}
\resdef{n11-s015-AOPT-hitllsm-n}{82.7}
\resdef{n11-s015-AOPT-tau-n}{0.923}
\resdef{n11-s015-AOPT-reg-n}{1.02}
\resdef{n11-s015-AOPT-agree-n}{85.8}
\resdef{n11-s015-p-AOPT-n}{$p=0.027$}
\resdef{n11-s015-p-AOPT-n-llsm}{$p=0.027$}
\resdef{n11-s015-FROZEN-hit-n}{85.3}
\resdef{n11-s015-FROZEN-hitllsm-n}{85.3}
\resdef{n11-s015-FROZEN-tau-n}{0.914}
\resdef{n11-s015-FROZEN-reg-n}{0.75}
\resdef{n11-s015-FROZEN-agree-n}{91.1}
\resdef{n11-s015-p-FROZEN-n}{$p=1.000$}
\resdef{n11-s015-p-FROZEN-n-llsm}{$p=1.000$}
\resdef{n11-s015-STAR-hit-n}{85.9}
\resdef{n11-s015-STAR-hitllsm-n}{85.9}
\resdef{n11-s015-STAR-tau-n}{0.913}
\resdef{n11-s015-STAR-reg-n}{0.69}
\resdef{n11-s015-STAR-agree-n}{92.2}
\resdef{n11-s015-p-STAR-n}{$p=0.405$}
\resdef{n11-s015-p-STAR-n-llsm}{$p=0.405$}
\resdef{n11-s015-ALEC-hit-n}{85.3}
\resdef{n11-s015-ALEC-hitllsm-n}{85.3}
\resdef{n11-s015-ALEC-tau-n}{0.913}
\resdef{n11-s015-ALEC-reg-n}{0.70}
\resdef{n11-s015-ALEC-agree-n}{92.9}
\resdef{n11-s015-ALEC-U-hit-n}{85.8}
\resdef{n11-s015-ALEC-U-hitllsm-n}{85.8}
\resdef{n11-s015-ALEC-U-tau-n}{0.909}
\resdef{n11-s015-ALEC-U-reg-n}{0.58}
\resdef{n11-s015-ALEC-U-agree-n}{93.8}
\resdef{n11-s015-p-ALEC-U-n}{$p=0.597$}
\resdef{n11-s015-p-ALEC-U-n-llsm}{$p=0.597$}
\resdef{n11-s015-tot-2n}{33}
\resdef{n11-s015-frac-2n}{60}
\resdef{n11-s015-RANDOM-hit-2n}{86.3}
\resdef{n11-s015-RANDOM-hitllsm-2n}{86.3}
\resdef{n11-s015-RANDOM-tau-2n}{0.936}
\resdef{n11-s015-RANDOM-reg-2n}{0.67}
\resdef{n11-s015-RANDOM-agree-2n}{91.1}
\resdef{n11-s015-p-RANDOM-2n}{$p=0.009$}
\resdef{n11-s015-p-RANDOM-2n-llsm}{$p=0.009$}
\resdef{n11-s015-HARKER-PG-hit-2n}{--}
\resdef{n11-s015-HARKER-PG-hitllsm-2n}{--}
\resdef{n11-s015-HARKER-PG-tau-2n}{--}
\resdef{n11-s015-HARKER-PG-reg-2n}{--}
\resdef{n11-s015-HARKER-PG-agree-2n}{--}
\resdef{n11-s015-HARKER-hit-2n}{85.9}
\resdef{n11-s015-HARKER-hitllsm-2n}{86.5}
\resdef{n11-s015-HARKER-tau-2n}{0.935}
\resdef{n11-s015-HARKER-reg-2n}{0.62}
\resdef{n11-s015-HARKER-agree-2n}{92.5}
\resdef{n11-s015-p-HARKER-2n}{$p<0.001$}
\resdef{n11-s015-p-HARKER-2n-llsm}{$p=0.007$}
\resdef{n11-s015-HARKER-LOG-hit-2n}{86.7}
\resdef{n11-s015-HARKER-LOG-hitllsm-2n}{86.9}
\resdef{n11-s015-HARKER-LOG-tau-2n}{0.931}
\resdef{n11-s015-HARKER-LOG-reg-2n}{0.52}
\resdef{n11-s015-HARKER-LOG-agree-2n}{96.4}
\resdef{n11-s015-p-HARKER-LOG-2n}{$p<0.001$}
\resdef{n11-s015-p-HARKER-LOG-2n-llsm}{$p<0.001$}
\resdef{n11-s015-AOPT-hit-2n}{86.1}
\resdef{n11-s015-AOPT-hitllsm-2n}{86.1}
\resdef{n11-s015-AOPT-tau-2n}{0.940}
\resdef{n11-s015-AOPT-reg-2n}{0.65}
\resdef{n11-s015-AOPT-agree-2n}{90.7}
\resdef{n11-s015-p-AOPT-2n}{$p=0.005$}
\resdef{n11-s015-p-AOPT-2n-llsm}{$p=0.005$}
\resdef{n11-s015-FROZEN-hit-2n}{88.4}
\resdef{n11-s015-FROZEN-hitllsm-2n}{88.4}
\resdef{n11-s015-FROZEN-tau-2n}{0.935}
\resdef{n11-s015-FROZEN-reg-2n}{0.47}
\resdef{n11-s015-FROZEN-agree-2n}{99.2}
\resdef{n11-s015-p-FROZEN-2n}{$p=0.289$}
\resdef{n11-s015-p-FROZEN-2n-llsm}{$p=0.289$}
\resdef{n11-s015-STAR-hit-2n}{88.8}
\resdef{n11-s015-STAR-hitllsm-2n}{88.8}
\resdef{n11-s015-STAR-tau-2n}{0.933}
\resdef{n11-s015-STAR-reg-2n}{0.42}
\resdef{n11-s015-STAR-agree-2n}{99.8}
\resdef{n11-s015-p-STAR-2n}{$p=1.000$}
\resdef{n11-s015-p-STAR-2n-llsm}{$p=1.000$}
\resdef{n11-s015-ALEC-hit-2n}{88.8}
\resdef{n11-s015-ALEC-hitllsm-2n}{88.8}
\resdef{n11-s015-ALEC-tau-2n}{0.933}
\resdef{n11-s015-ALEC-reg-2n}{0.43}
\resdef{n11-s015-ALEC-agree-2n}{99.8}
\resdef{n11-s015-ALEC-U-hit-2n}{88.8}
\resdef{n11-s015-ALEC-U-hitllsm-2n}{88.8}
\resdef{n11-s015-ALEC-U-tau-2n}{0.934}
\resdef{n11-s015-ALEC-U-reg-2n}{0.44}
\resdef{n11-s015-ALEC-U-agree-2n}{98.4}
\resdef{n11-s015-p-ALEC-U-2n}{$p=1.000$}
\resdef{n11-s015-p-ALEC-U-2n-llsm}{$p=1.000$}
\resdef{n11-s015-stop-LC-q}{8.5}
\resdef{n11-s015-stop-LC-hit}{84.7}
\resdef{n11-s015-stop-LC-agree}{90.8}
\resdef{n11-s015-stop-LC-early}{100.0}
\resdef{n11-s015-stop-LC-exceed}{15.3}
\resdef{n11-s015-stop-LC-qmax}{18}
\resdef{n11-s015-stop-LC-qtot}{19.5}
\resdef{n11-s015-stop-LC-errearly}{15.3}
\resdef{n11-s015-stop-LC-errci}{[13.2, 17.7]}
\resdef{n11-s015-stop-LC-regret}{0.76}
\resdef{n11-s015-stop-LC-qround}{8}
\resdef{n11-s015-stop-LC-fixedhit}{84.9}
\resdef{n11-s015-stop-LC-fixedreg}{0.79}
\resdef{n11-s015-stop-LC-save2n}{62}
\resdef{n11-s015-stop-B10-ALEC-q}{19.0}
\resdef{n11-s015-stop-B10-ALEC-hit}{88.9}
\resdef{n11-s015-stop-B10-ALEC-agree}{100.0}
\resdef{n11-s015-stop-B10-ALEC-early}{69.3}
\resdef{n11-s015-stop-B10-ALEC-exceed}{11.1}
\resdef{n11-s015-stop-B10-ALEC-qmax}{44}
\resdef{n11-s015-stop-B10-ALEC-qtot}{30.0}
\resdef{n11-s015-stop-B10-ALEC-errearly}{2.5}
\resdef{n11-s015-stop-B10-ALEC-errci}{[1.5, 3.9]}
\resdef{n11-s015-stop-B10-ALEC-regret}{0.43}
\resdef{n11-s015-stop-B10-ALEC-qround}{19}
\resdef{n11-s015-stop-B10-ALEC-fixedhit}{88.8}
\resdef{n11-s015-stop-B10-ALEC-fixedreg}{0.43}
\resdef{n11-s015-stop-B10-ALEC-save2n}{14}
\resdef{n11-s015-stop-B05-ALECU-q}{22.2}
\resdef{n11-s015-stop-B05-ALECU-hit}{88.9}
\resdef{n11-s015-stop-B05-ALECU-agree}{100.0}
\resdef{n11-s015-stop-B05-ALECU-early}{60.9}
\resdef{n11-s015-stop-B05-ALECU-exceed}{11.1}
\resdef{n11-s015-stop-B05-ALECU-qmax}{44}
\resdef{n11-s015-stop-B05-ALECU-qtot}{33.2}
\resdef{n11-s015-stop-B05-ALECU-errearly}{1.0}
\resdef{n11-s015-stop-B05-ALECU-errci}{[0.5, 2.1]}
\resdef{n11-s015-stop-B05-ALECU-regret}{0.43}
\resdef{n11-s015-stop-B05-ALECU-qround}{22}
\resdef{n11-s015-stop-B05-ALECU-fixedhit}{88.8}
\resdef{n11-s015-stop-B05-ALECU-fixedreg}{0.43}
\resdef{n11-s015-stop-B05-ALECU-save2n}{-1}
\resdef{n11-s015-stop-B10-ALECU-q}{18.1}
\resdef{n11-s015-stop-B10-ALECU-hit}{88.9}
\resdef{n11-s015-stop-B10-ALECU-agree}{100.0}
\resdef{n11-s015-stop-B10-ALECU-early}{70.0}
\resdef{n11-s015-stop-B10-ALECU-exceed}{11.1}
\resdef{n11-s015-stop-B10-ALECU-qmax}{44}
\resdef{n11-s015-stop-B10-ALECU-qtot}{29.1}
\resdef{n11-s015-stop-B10-ALECU-errearly}{2.7}
\resdef{n11-s015-stop-B10-ALECU-errci}{[1.7, 4.2]}
\resdef{n11-s015-stop-B10-ALECU-regret}{0.43}
\resdef{n11-s015-stop-B10-ALECU-qround}{18}
\resdef{n11-s015-stop-B10-ALECU-fixedhit}{89.0}
\resdef{n11-s015-stop-B10-ALECU-fixedreg}{0.41}
\resdef{n11-s015-stop-B10-ALECU-save2n}{18}
\resdef{n11-s015-stop-B20-ALECU-q}{12.5}
\resdef{n11-s015-stop-B20-ALECU-hit}{88.8}
\resdef{n11-s015-stop-B20-ALECU-agree}{99.7}
\resdef{n11-s015-stop-B20-ALECU-early}{82.0}
\resdef{n11-s015-stop-B20-ALECU-exceed}{11.2}
\resdef{n11-s015-stop-B20-ALECU-qmax}{44}
\resdef{n11-s015-stop-B20-ALECU-qtot}{23.5}
\resdef{n11-s015-stop-B20-ALECU-errearly}{6.0}
\resdef{n11-s015-stop-B20-ALECU-errci}{[4.5, 7.8]}
\resdef{n11-s015-stop-B20-ALECU-regret}{0.43}
\resdef{n11-s015-stop-B20-ALECU-qround}{12}
\resdef{n11-s015-stop-B20-ALECU-fixedhit}{85.9}
\resdef{n11-s015-stop-B20-ALECU-fixedreg}{0.61}
\resdef{n11-s015-stop-B20-ALECU-save2n}{43}
\resdef{n11-s015-stop-B10-ALECU-s06-q}{24.8}
\resdef{n11-s015-stop-B10-ALECU-s06-hit}{88.9}
\resdef{n11-s015-stop-B10-ALECU-s06-agree}{100.0}
\resdef{n11-s015-stop-B10-ALECU-s06-early}{59.1}
\resdef{n11-s015-stop-B10-ALECU-s06-exceed}{11.1}
\resdef{n11-s015-stop-B10-ALECU-s06-qmax}{44}
\resdef{n11-s015-stop-B10-ALECU-s06-qtot}{35.8}
\resdef{n11-s015-stop-B10-ALECU-s06-errearly}{1.0}
\resdef{n11-s015-stop-B10-ALECU-s06-errci}{[0.5, 2.2]}
\resdef{n11-s015-stop-B10-ALECU-s06-regret}{0.43}
\resdef{n11-s015-stop-B10-ALECU-s06-save2n}{-13}
\resdef{n11-s015-stop-B10-ALECU-cap-q}{11.3}
\resdef{n11-s015-stop-B10-ALECU-cap-hit}{88.8}
\resdef{n11-s015-stop-B10-ALECU-cap-agree}{98.4}
\resdef{n11-s015-stop-B10-ALECU-cap-early}{68.2}
\resdef{n11-s015-stop-B10-ALECU-cap-exceed}{11.2}
\resdef{n11-s015-stop-B10-ALECU-cap-qmax}{22}
\resdef{n11-s015-stop-B10-ALECU-cap-qtot}{22.3}
\resdef{n11-s015-stop-B10-ALECU-cap-errearly}{2.6}
\resdef{n11-s015-stop-B10-ALECU-cap-errci}{[1.7, 4.1]}
\resdef{n11-s015-stop-B10-ALECU-cap-regret}{0.44}
\resdef{n11-s015-stop-B10-ALECU-cap-qround}{11}
\resdef{n11-s015-stop-B10-ALECU-cap-fixedhit}{85.3}
\resdef{n11-s015-stop-B10-ALECU-cap-fixedreg}{0.70}
\resdef{n11-s015-stop-B10-ALECU-cap-save2n}{48}
\resdef{n11-s015-stop-B10-ALECU-oracle-q}{12.0}
\resdef{n11-s015-stop-B10-ALECU-oracle-hit}{87.5}
\resdef{n11-s015-stop-B10-ALECU-oracle-agree}{98.1}
\resdef{n11-s015-stop-B10-ALECU-oracle-early}{77.9}
\resdef{n11-s015-stop-B10-ALECU-oracle-exceed}{12.5}
\resdef{n11-s015-stop-B10-ALECU-oracle-qmax}{44}
\resdef{n11-s015-stop-B10-ALECU-oracle-qtot}{23.0}
\resdef{n11-s015-stop-B10-ALECU-oracle-errearly}{6.9}
\resdef{n11-s015-stop-B10-ALECU-oracle-errci}{[5.4, 8.9]}
\resdef{n11-s015-stop-B10-ALECU-oracle-regret}{0.54}
\resdef{n11-s015-stop-B10-ALECU-oracle-qround}{12}
\resdef{n11-s015-stop-B10-ALECU-oracle-fixedhit}{85.9}
\resdef{n11-s015-stop-B10-ALECU-oracle-fixedreg}{0.61}
\resdef{n11-s015-stop-B10-ALECU-oracle-save2n}{45}
\resdef{n11-s015-stop-E02-ALEC-q}{16.2}
\resdef{n11-s015-stop-E02-ALEC-hit}{88.9}
\resdef{n11-s015-stop-E02-ALEC-agree}{100.0}
\resdef{n11-s015-stop-E02-ALEC-early}{77.1}
\resdef{n11-s015-stop-E02-ALEC-exceed}{6.1}
\resdef{n11-s015-stop-E02-ALEC-qmax}{44}
\resdef{n11-s015-stop-E02-ALEC-qtot}{27.2}
\resdef{n11-s015-stop-E02-ALEC-errearly}{1.4}
\resdef{n11-s015-stop-E02-ALEC-errci}{[0.8, 2.5]}
\resdef{n11-s015-stop-E02-ALEC-regret}{0.43}
\resdef{n11-s015-stop-E02-ALEC-qround}{16}
\resdef{n11-s015-stop-E02-ALEC-fixedhit}{88.7}
\resdef{n11-s015-stop-E02-ALEC-fixedreg}{0.43}
\resdef{n11-s015-stop-E02-ALEC-save2n}{27}
\resdef{n11-s015-stop-E05-ALEC-q}{11.4}
\resdef{n11-s015-stop-E05-ALEC-hit}{88.4}
\resdef{n11-s015-stop-E05-ALEC-agree}{99.3}
\resdef{n11-s015-stop-E05-ALEC-early}{89.3}
\resdef{n11-s015-stop-E05-ALEC-exceed}{3.4}
\resdef{n11-s015-stop-E05-ALEC-qmax}{44}
\resdef{n11-s015-stop-E05-ALEC-qtot}{22.4}
\resdef{n11-s015-stop-E05-ALEC-errearly}{1.9}
\resdef{n11-s015-stop-E05-ALEC-errci}{[1.2, 3.0]}
\resdef{n11-s015-stop-E05-ALEC-regret}{0.44}
\resdef{n11-s015-stop-E05-ALEC-qround}{11}
\resdef{n11-s015-stop-E05-ALEC-fixedhit}{85.3}
\resdef{n11-s015-stop-E05-ALEC-fixedreg}{0.70}
\resdef{n11-s015-stop-E05-ALEC-save2n}{48}
\resdef{n11-s015-stop-B10-ALEC-cap-q}{12.0}
\resdef{n11-s015-stop-B10-ALEC-cap-hit}{88.8}
\resdef{n11-s015-stop-B10-ALEC-cap-agree}{99.8}
\resdef{n11-s015-stop-B10-ALEC-cap-early}{66.9}
\resdef{n11-s015-stop-B10-ALEC-cap-exceed}{11.2}
\resdef{n11-s015-stop-B10-ALEC-cap-qmax}{22}
\resdef{n11-s015-stop-B10-ALEC-cap-qtot}{23.0}
\resdef{n11-s015-stop-B10-ALEC-cap-errearly}{2.5}
\resdef{n11-s015-stop-B10-ALEC-cap-errci}{[1.6, 4.0]}
\resdef{n11-s015-stop-B10-ALEC-cap-regret}{0.43}
\resdef{n11-s015-stop-B10-ALEC-cap-qround}{12}
\resdef{n11-s015-stop-B10-ALEC-cap-fixedhit}{85.9}
\resdef{n11-s015-stop-B10-ALEC-cap-fixedreg}{0.61}
\resdef{n11-s015-stop-B10-ALEC-cap-save2n}{45}
\resdef{n11-s015-stop-E02-ALEC-cap-q}{10.9}
\resdef{n11-s015-stop-E02-ALEC-cap-hit}{88.8}
\resdef{n11-s015-stop-E02-ALEC-cap-agree}{99.8}
\resdef{n11-s015-stop-E02-ALEC-cap-early}{74.8}
\resdef{n11-s015-stop-E02-ALEC-cap-exceed}{6.2}
\resdef{n11-s015-stop-E02-ALEC-cap-qmax}{22}
\resdef{n11-s015-stop-E02-ALEC-cap-qtot}{21.9}
\resdef{n11-s015-stop-E02-ALEC-cap-errearly}{1.5}
\resdef{n11-s015-stop-E02-ALEC-cap-errci}{[0.8, 2.6]}
\resdef{n11-s015-stop-E02-ALEC-cap-regret}{0.43}
\resdef{n11-s015-stop-E02-ALEC-cap-qround}{11}
\resdef{n11-s015-stop-E02-ALEC-cap-fixedhit}{85.3}
\resdef{n11-s015-stop-E02-ALEC-cap-fixedreg}{0.70}
\resdef{n11-s015-stop-E02-ALEC-cap-save2n}{50}
\resdef{n11-s015-stop-E05-ALEC-cap-q}{8.9}
\resdef{n11-s015-stop-E05-ALEC-cap-hit}{88.3}
\resdef{n11-s015-stop-E05-ALEC-cap-agree}{99.1}
\resdef{n11-s015-stop-E05-ALEC-cap-early}{87.5}
\resdef{n11-s015-stop-E05-ALEC-cap-exceed}{3.4}
\resdef{n11-s015-stop-E05-ALEC-cap-qmax}{22}
\resdef{n11-s015-stop-E05-ALEC-cap-qtot}{19.9}
\resdef{n11-s015-stop-E05-ALEC-cap-errearly}{1.7}
\resdef{n11-s015-stop-E05-ALEC-cap-errci}{[1.0, 2.8]}
\resdef{n11-s015-stop-E05-ALEC-cap-regret}{0.45}
\resdef{n11-s015-stop-E05-ALEC-cap-qround}{9}
\resdef{n11-s015-stop-E05-ALEC-cap-fixedhit}{85.0}
\resdef{n11-s015-stop-E05-ALEC-cap-fixedreg}{0.76}
\resdef{n11-s015-stop-E05-ALEC-cap-save2n}{60}
\resdef{n11-s015-stop-E05-ALEC-cap-s06-q}{14.0}
\resdef{n11-s015-stop-E05-ALEC-cap-s06-hit}{88.8}
\resdef{n11-s015-stop-E05-ALEC-cap-s06-agree}{99.8}
\resdef{n11-s015-stop-E05-ALEC-cap-s06-early}{68.4}
\resdef{n11-s015-stop-E05-ALEC-cap-s06-exceed}{3.2}
\resdef{n11-s015-stop-E05-ALEC-cap-s06-qmax}{22}
\resdef{n11-s015-stop-E05-ALEC-cap-s06-qtot}{25.0}
\resdef{n11-s015-stop-E05-ALEC-cap-s06-errearly}{0.4}
\resdef{n11-s015-stop-E05-ALEC-cap-s06-errci}{[0.1, 1.3]}
\resdef{n11-s015-stop-E05-ALEC-cap-s06-regret}{0.43}
\resdef{n11-s015-stop-E05-ALEC-cap-s06-qround}{14}
\resdef{n11-s015-stop-E05-ALEC-cap-s06-fixedhit}{86.3}
\resdef{n11-s015-stop-E05-ALEC-cap-s06-fixedreg}{0.56}
\resdef{n11-s015-stop-E05-ALEC-cap-s06-save2n}{36}
\resdef{n11-s030-full-hit}{76.6}
\resdef{n11-s030-tot-h}{17}
\resdef{n11-s030-frac-h}{31}
\resdef{n11-s030-RANDOM-hit-h}{63.7}
\resdef{n11-s030-RANDOM-hitllsm-h}{63.7}
\resdef{n11-s030-RANDOM-tau-h}{0.811}
\resdef{n11-s030-RANDOM-reg-h}{4.86}
\resdef{n11-s030-RANDOM-agree-h}{67.3}
\resdef{n11-s030-p-RANDOM-h}{$p<0.001$}
\resdef{n11-s030-p-RANDOM-h-llsm}{$p<0.001$}
\resdef{n11-s030-HARKER-PG-hit-h}{64.8}
\resdef{n11-s030-HARKER-PG-hitllsm-h}{65.0}
\resdef{n11-s030-HARKER-PG-tau-h}{0.795}
\resdef{n11-s030-HARKER-PG-reg-h}{4.79}
\resdef{n11-s030-HARKER-PG-agree-h}{68.5}
\resdef{n11-s030-p-HARKER-PG-h}{$p<0.001$}
\resdef{n11-s030-p-HARKER-PG-h-llsm}{$p<0.001$}
\resdef{n11-s030-HARKER-hit-h}{65.7}
\resdef{n11-s030-HARKER-hitllsm-h}{64.6}
\resdef{n11-s030-HARKER-tau-h}{0.794}
\resdef{n11-s030-HARKER-reg-h}{4.54}
\resdef{n11-s030-HARKER-agree-h}{69.0}
\resdef{n11-s030-p-HARKER-h}{$p<0.001$}
\resdef{n11-s030-p-HARKER-h-llsm}{$p<0.001$}
\resdef{n11-s030-HARKER-LOG-hit-h}{70.8}
\resdef{n11-s030-HARKER-LOG-hitllsm-h}{71.3}
\resdef{n11-s030-HARKER-LOG-tau-h}{0.783}
\resdef{n11-s030-HARKER-LOG-reg-h}{3.03}
\resdef{n11-s030-HARKER-LOG-agree-h}{76.1}
\resdef{n11-s030-p-HARKER-LOG-h}{$p=0.039$}
\resdef{n11-s030-p-HARKER-LOG-h-llsm}{$p=0.083$}
\resdef{n11-s030-AOPT-hit-h}{64.5}
\resdef{n11-s030-AOPT-hitllsm-h}{64.5}
\resdef{n11-s030-AOPT-tau-h}{0.825}
\resdef{n11-s030-AOPT-reg-h}{4.51}
\resdef{n11-s030-AOPT-agree-h}{68.5}
\resdef{n11-s030-p-AOPT-h}{$p<0.001$}
\resdef{n11-s030-p-AOPT-h-llsm}{$p<0.001$}
\resdef{n11-s030-FROZEN-hit-h}{67.0}
\resdef{n11-s030-FROZEN-hitllsm-h}{67.0}
\resdef{n11-s030-FROZEN-tau-h}{0.801}
\resdef{n11-s030-FROZEN-reg-h}{4.02}
\resdef{n11-s030-FROZEN-agree-h}{72.1}
\resdef{n11-s030-p-FROZEN-h}{$p<0.001$}
\resdef{n11-s030-p-FROZEN-h-llsm}{$p<0.001$}
\resdef{n11-s030-STAR-hit-h}{67.1}
\resdef{n11-s030-STAR-hitllsm-h}{67.1}
\resdef{n11-s030-STAR-tau-h}{0.802}
\resdef{n11-s030-STAR-reg-h}{4.02}
\resdef{n11-s030-STAR-agree-h}{72.0}
\resdef{n11-s030-p-STAR-h}{$p<0.001$}
\resdef{n11-s030-p-STAR-h-llsm}{$p<0.001$}
\resdef{n11-s030-ALEC-hit-h}{73.6}
\resdef{n11-s030-ALEC-hitllsm-h}{73.6}
\resdef{n11-s030-ALEC-tau-h}{0.802}
\resdef{n11-s030-ALEC-reg-h}{2.51}
\resdef{n11-s030-ALEC-agree-h}{78.5}
\resdef{n11-s030-ALEC-U-hit-h}{71.2}
\resdef{n11-s030-ALEC-U-hitllsm-h}{71.2}
\resdef{n11-s030-ALEC-U-tau-h}{0.804}
\resdef{n11-s030-ALEC-U-reg-h}{2.84}
\resdef{n11-s030-ALEC-U-agree-h}{78.4}
\resdef{n11-s030-p-ALEC-U-h}{$p=0.057$}
\resdef{n11-s030-p-ALEC-U-h-llsm}{$p=0.057$}
\resdef{n11-s030-tot-n}{22}
\resdef{n11-s030-frac-n}{40}
\resdef{n11-s030-RANDOM-hit-n}{66.3}
\resdef{n11-s030-RANDOM-hitllsm-n}{66.3}
\resdef{n11-s030-RANDOM-tau-n}{0.840}
\resdef{n11-s030-RANDOM-reg-n}{4.01}
\resdef{n11-s030-RANDOM-agree-n}{72.7}
\resdef{n11-s030-p-RANDOM-n}{$p<0.001$}
\resdef{n11-s030-p-RANDOM-n-llsm}{$p<0.001$}
\resdef{n11-s030-HARKER-PG-hit-n}{68.7}
\resdef{n11-s030-HARKER-PG-hitllsm-n}{70.9}
\resdef{n11-s030-HARKER-PG-tau-n}{0.830}
\resdef{n11-s030-HARKER-PG-reg-n}{3.48}
\resdef{n11-s030-HARKER-PG-agree-n}{73.7}
\resdef{n11-s030-p-HARKER-PG-n}{$p<0.001$}
\resdef{n11-s030-p-HARKER-PG-n-llsm}{$p=0.029$}
\resdef{n11-s030-HARKER-hit-n}{71.0}
\resdef{n11-s030-HARKER-hitllsm-n}{70.5}
\resdef{n11-s030-HARKER-tau-n}{0.830}
\resdef{n11-s030-HARKER-reg-n}{3.14}
\resdef{n11-s030-HARKER-agree-n}{76.1}
\resdef{n11-s030-p-HARKER-n}{$p=0.033$}
\resdef{n11-s030-p-HARKER-n-llsm}{$p=0.015$}
\resdef{n11-s030-HARKER-LOG-hit-n}{74.5}
\resdef{n11-s030-HARKER-LOG-hitllsm-n}{75.7}
\resdef{n11-s030-HARKER-LOG-tau-n}{0.822}
\resdef{n11-s030-HARKER-LOG-reg-n}{2.09}
\resdef{n11-s030-HARKER-LOG-agree-n}{83.5}
\resdef{n11-s030-p-HARKER-LOG-n}{$p=0.928$}
\resdef{n11-s030-p-HARKER-LOG-n-llsm}{$p=0.219$}
\resdef{n11-s030-AOPT-hit-n}{69.0}
\resdef{n11-s030-AOPT-hitllsm-n}{69.0}
\resdef{n11-s030-AOPT-tau-n}{0.853}
\resdef{n11-s030-AOPT-reg-n}{3.33}
\resdef{n11-s030-AOPT-agree-n}{75.5}
\resdef{n11-s030-p-AOPT-n}{$p<0.001$}
\resdef{n11-s030-p-AOPT-n-llsm}{$p<0.001$}
\resdef{n11-s030-FROZEN-hit-n}{73.0}
\resdef{n11-s030-FROZEN-hitllsm-n}{73.0}
\resdef{n11-s030-FROZEN-tau-n}{0.838}
\resdef{n11-s030-FROZEN-reg-n}{2.43}
\resdef{n11-s030-FROZEN-agree-n}{80.3}
\resdef{n11-s030-p-FROZEN-n}{$p=0.279$}
\resdef{n11-s030-p-FROZEN-n-llsm}{$p=0.279$}
\resdef{n11-s030-STAR-hit-n}{73.7}
\resdef{n11-s030-STAR-hitllsm-n}{73.7}
\resdef{n11-s030-STAR-tau-n}{0.836}
\resdef{n11-s030-STAR-reg-n}{2.33}
\resdef{n11-s030-STAR-agree-n}{82.4}
\resdef{n11-s030-p-STAR-n}{$p=0.550$}
\resdef{n11-s030-p-STAR-n-llsm}{$p=0.550$}
\resdef{n11-s030-ALEC-hit-n}{74.3}
\resdef{n11-s030-ALEC-hitllsm-n}{74.3}
\resdef{n11-s030-ALEC-tau-n}{0.834}
\resdef{n11-s030-ALEC-reg-n}{2.16}
\resdef{n11-s030-ALEC-agree-n}{83.0}
\resdef{n11-s030-ALEC-U-hit-n}{74.2}
\resdef{n11-s030-ALEC-U-hitllsm-n}{74.2}
\resdef{n11-s030-ALEC-U-tau-n}{0.831}
\resdef{n11-s030-ALEC-U-reg-n}{2.13}
\resdef{n11-s030-ALEC-U-agree-n}{82.8}
\resdef{n11-s030-p-ALEC-U-n}{$p=1.000$}
\resdef{n11-s030-p-ALEC-U-n-llsm}{$p=1.000$}
\resdef{n11-s030-tot-2n}{33}
\resdef{n11-s030-frac-2n}{60}
\resdef{n11-s030-RANDOM-hit-2n}{70.9}
\resdef{n11-s030-RANDOM-hitllsm-2n}{70.9}
\resdef{n11-s030-RANDOM-tau-2n}{0.877}
\resdef{n11-s030-RANDOM-reg-2n}{2.47}
\resdef{n11-s030-RANDOM-agree-2n}{81.8}
\resdef{n11-s030-p-RANDOM-2n}{$p<0.001$}
\resdef{n11-s030-p-RANDOM-2n-llsm}{$p<0.001$}
\resdef{n11-s030-HARKER-PG-hit-2n}{--}
\resdef{n11-s030-HARKER-PG-hitllsm-2n}{--}
\resdef{n11-s030-HARKER-PG-tau-2n}{--}
\resdef{n11-s030-HARKER-PG-reg-2n}{--}
\resdef{n11-s030-HARKER-PG-agree-2n}{--}
\resdef{n11-s030-HARKER-hit-2n}{74.4}
\resdef{n11-s030-HARKER-hitllsm-2n}{74.8}
\resdef{n11-s030-HARKER-tau-2n}{0.875}
\resdef{n11-s030-HARKER-reg-2n}{2.23}
\resdef{n11-s030-HARKER-agree-2n}{83.8}
\resdef{n11-s030-p-HARKER-2n}{$p=0.076$}
\resdef{n11-s030-p-HARKER-2n-llsm}{$p=0.145$}
\resdef{n11-s030-HARKER-LOG-hit-2n}{76.7}
\resdef{n11-s030-HARKER-LOG-hitllsm-2n}{76.4}
\resdef{n11-s030-HARKER-LOG-tau-2n}{0.868}
\resdef{n11-s030-HARKER-LOG-reg-2n}{1.63}
\resdef{n11-s030-HARKER-LOG-agree-2n}{90.8}
\resdef{n11-s030-p-HARKER-LOG-2n}{$p=1.000$}
\resdef{n11-s030-p-HARKER-LOG-2n-llsm}{$p=0.897$}
\resdef{n11-s030-AOPT-hit-2n}{73.9}
\resdef{n11-s030-AOPT-hitllsm-2n}{73.9}
\resdef{n11-s030-AOPT-tau-2n}{0.882}
\resdef{n11-s030-AOPT-reg-2n}{2.21}
\resdef{n11-s030-AOPT-agree-2n}{83.9}
\resdef{n11-s030-p-AOPT-2n}{$p=0.028$}
\resdef{n11-s030-p-AOPT-2n-llsm}{$p=0.028$}
\resdef{n11-s030-FROZEN-hit-2n}{76.7}
\resdef{n11-s030-FROZEN-hitllsm-2n}{76.7}
\resdef{n11-s030-FROZEN-tau-2n}{0.872}
\resdef{n11-s030-FROZEN-reg-2n}{1.72}
\resdef{n11-s030-FROZEN-agree-2n}{94.6}
\resdef{n11-s030-p-FROZEN-2n}{$p=1.000$}
\resdef{n11-s030-p-FROZEN-2n-llsm}{$p=1.000$}
\resdef{n11-s030-STAR-hit-2n}{76.6}
\resdef{n11-s030-STAR-hitllsm-2n}{76.6}
\resdef{n11-s030-STAR-tau-2n}{0.872}
\resdef{n11-s030-STAR-reg-2n}{1.59}
\resdef{n11-s030-STAR-agree-2n}{98.8}
\resdef{n11-s030-p-STAR-2n}{$p=1.000$}
\resdef{n11-s030-p-STAR-2n-llsm}{$p=1.000$}
\resdef{n11-s030-ALEC-hit-2n}{76.6}
\resdef{n11-s030-ALEC-hitllsm-2n}{76.6}
\resdef{n11-s030-ALEC-tau-2n}{0.872}
\resdef{n11-s030-ALEC-reg-2n}{1.59}
\resdef{n11-s030-ALEC-agree-2n}{99.0}
\resdef{n11-s030-ALEC-U-hit-2n}{76.6}
\resdef{n11-s030-ALEC-U-hitllsm-2n}{76.6}
\resdef{n11-s030-ALEC-U-tau-2n}{0.872}
\resdef{n11-s030-ALEC-U-reg-2n}{1.57}
\resdef{n11-s030-ALEC-U-agree-2n}{94.3}
\resdef{n11-s030-p-ALEC-U-2n}{$p=1.000$}
\resdef{n11-s030-p-ALEC-U-2n-llsm}{$p=1.000$}
\resdef{n11-s030-stop-LC-q}{9.1}
\resdef{n11-s030-stop-LC-hit}{75.2}
\resdef{n11-s030-stop-LC-agree}{82.5}
\resdef{n11-s030-stop-LC-early}{100.0}
\resdef{n11-s030-stop-LC-exceed}{24.8}
\resdef{n11-s030-stop-LC-qmax}{19}
\resdef{n11-s030-stop-LC-qtot}{20.1}
\resdef{n11-s030-stop-LC-errearly}{24.8}
\resdef{n11-s030-stop-LC-errci}{[22.2, 27.6]}
\resdef{n11-s030-stop-LC-regret}{2.06}
\resdef{n11-s030-stop-LC-qround}{9}
\resdef{n11-s030-stop-LC-fixedhit}{72.3}
\resdef{n11-s030-stop-LC-fixedreg}{2.49}
\resdef{n11-s030-stop-LC-save2n}{59}
\resdef{n11-s030-stop-B10-ALEC-q}{25.7}
\resdef{n11-s030-stop-B10-ALEC-hit}{77.2}
\resdef{n11-s030-stop-B10-ALEC-agree}{99.2}
\resdef{n11-s030-stop-B10-ALEC-early}{51.9}
\resdef{n11-s030-stop-B10-ALEC-exceed}{22.8}
\resdef{n11-s030-stop-B10-ALEC-qmax}{44}
\resdef{n11-s030-stop-B10-ALEC-qtot}{36.7}
\resdef{n11-s030-stop-B10-ALEC-errearly}{5.4}
\resdef{n11-s030-stop-B10-ALEC-errci}{[3.8, 7.7]}
\resdef{n11-s030-stop-B10-ALEC-regret}{1.50}
\resdef{n11-s030-stop-B10-ALEC-save2n}{-17}
\resdef{n11-s030-stop-B05-ALECU-q}{28.6}
\resdef{n11-s030-stop-B05-ALECU-hit}{76.4}
\resdef{n11-s030-stop-B05-ALECU-agree}{99.6}
\resdef{n11-s030-stop-B05-ALECU-early}{43.8}
\resdef{n11-s030-stop-B05-ALECU-exceed}{23.6}
\resdef{n11-s030-stop-B05-ALECU-qmax}{44}
\resdef{n11-s030-stop-B05-ALECU-qtot}{39.6}
\resdef{n11-s030-stop-B05-ALECU-errearly}{3.9}
\resdef{n11-s030-stop-B05-ALECU-errci}{[2.4, 6.1]}
\resdef{n11-s030-stop-B05-ALECU-regret}{1.59}
\resdef{n11-s030-stop-B05-ALECU-save2n}{-30}
\resdef{n11-s030-stop-B10-ALECU-q}{23.9}
\resdef{n11-s030-stop-B10-ALECU-hit}{76.7}
\resdef{n11-s030-stop-B10-ALECU-agree}{98.9}
\resdef{n11-s030-stop-B10-ALECU-early}{54.7}
\resdef{n11-s030-stop-B10-ALECU-exceed}{23.3}
\resdef{n11-s030-stop-B10-ALECU-qmax}{44}
\resdef{n11-s030-stop-B10-ALECU-qtot}{34.9}
\resdef{n11-s030-stop-B10-ALECU-errearly}{7.3}
\resdef{n11-s030-stop-B10-ALECU-errci}{[5.4, 9.8]}
\resdef{n11-s030-stop-B10-ALECU-regret}{1.58}
\resdef{n11-s030-stop-B10-ALECU-save2n}{-9}
\resdef{n11-s030-stop-B20-ALECU-q}{16.6}
\resdef{n11-s030-stop-B20-ALECU-hit}{76.2}
\resdef{n11-s030-stop-B20-ALECU-agree}{96.0}
\resdef{n11-s030-stop-B20-ALECU-early}{72.0}
\resdef{n11-s030-stop-B20-ALECU-exceed}{23.8}
\resdef{n11-s030-stop-B20-ALECU-qmax}{44}
\resdef{n11-s030-stop-B20-ALECU-qtot}{27.6}
\resdef{n11-s030-stop-B20-ALECU-errearly}{14.2}
\resdef{n11-s030-stop-B20-ALECU-errci}{[11.8, 16.9]}
\resdef{n11-s030-stop-B20-ALECU-regret}{1.65}
\resdef{n11-s030-stop-B20-ALECU-qround}{17}
\resdef{n11-s030-stop-B20-ALECU-fixedhit}{76.4}
\resdef{n11-s030-stop-B20-ALECU-fixedreg}{1.60}
\resdef{n11-s030-stop-B20-ALECU-save2n}{24}
\resdef{n11-s030-stop-B10-ALECU-s06-q}{29.3}
\resdef{n11-s030-stop-B10-ALECU-s06-hit}{76.7}
\resdef{n11-s030-stop-B10-ALECU-s06-agree}{99.9}
\resdef{n11-s030-stop-B10-ALECU-s06-early}{44.6}
\resdef{n11-s030-stop-B10-ALECU-s06-exceed}{23.3}
\resdef{n11-s030-stop-B10-ALECU-s06-qmax}{44}
\resdef{n11-s030-stop-B10-ALECU-s06-qtot}{40.3}
\resdef{n11-s030-stop-B10-ALECU-s06-errearly}{3.4}
\resdef{n11-s030-stop-B10-ALECU-s06-errci}{[2.0, 5.5]}
\resdef{n11-s030-stop-B10-ALECU-s06-regret}{1.54}
\resdef{n11-s030-stop-B10-ALECU-s06-save2n}{-33}
\resdef{n11-s030-stop-B10-ALECU-cap-q}{13.8}
\resdef{n11-s030-stop-B10-ALECU-cap-hit}{76.7}
\resdef{n11-s030-stop-B10-ALECU-cap-agree}{93.4}
\resdef{n11-s030-stop-B10-ALECU-cap-early}{53.4}
\resdef{n11-s030-stop-B10-ALECU-cap-exceed}{23.3}
\resdef{n11-s030-stop-B10-ALECU-cap-qmax}{22}
\resdef{n11-s030-stop-B10-ALECU-cap-qtot}{24.8}
\resdef{n11-s030-stop-B10-ALECU-cap-errearly}{7.1}
\resdef{n11-s030-stop-B10-ALECU-cap-errci}{[5.2, 9.6]}
\resdef{n11-s030-stop-B10-ALECU-cap-regret}{1.58}
\resdef{n11-s030-stop-B10-ALECU-cap-qround}{14}
\resdef{n11-s030-stop-B10-ALECU-cap-fixedhit}{75.3}
\resdef{n11-s030-stop-B10-ALECU-cap-fixedreg}{1.79}
\resdef{n11-s030-stop-B10-ALECU-cap-save2n}{37}
\resdef{n11-s030-stop-B10-ALECU-oracle-q}{22.1}
\resdef{n11-s030-stop-B10-ALECU-oracle-hit}{75.8}
\resdef{n11-s030-stop-B10-ALECU-oracle-agree}{97.5}
\resdef{n11-s030-stop-B10-ALECU-oracle-early}{56.3}
\resdef{n11-s030-stop-B10-ALECU-oracle-exceed}{24.2}
\resdef{n11-s030-stop-B10-ALECU-oracle-qmax}{44}
\resdef{n11-s030-stop-B10-ALECU-oracle-qtot}{33.1}
\resdef{n11-s030-stop-B10-ALECU-oracle-errearly}{10.3}
\resdef{n11-s030-stop-B10-ALECU-oracle-errci}{[8.1, 13.1]}
\resdef{n11-s030-stop-B10-ALECU-oracle-regret}{1.86}
\resdef{n11-s030-stop-B10-ALECU-oracle-qround}{22}
\resdef{n11-s030-stop-B10-ALECU-oracle-fixedhit}{76.6}
\resdef{n11-s030-stop-B10-ALECU-oracle-fixedreg}{1.59}
\resdef{n11-s030-stop-B10-ALECU-oracle-save2n}{-0}
\resdef{n11-s030-stop-E02-ALEC-q}{23.5}
\resdef{n11-s030-stop-E02-ALEC-hit}{77.3}
\resdef{n11-s030-stop-E02-ALEC-agree}{98.9}
\resdef{n11-s030-stop-E02-ALEC-early}{57.8}
\resdef{n11-s030-stop-E02-ALEC-exceed}{17.9}
\resdef{n11-s030-stop-E02-ALEC-qmax}{44}
\resdef{n11-s030-stop-E02-ALEC-qtot}{34.5}
\resdef{n11-s030-stop-E02-ALEC-errearly}{4.8}
\resdef{n11-s030-stop-E02-ALEC-errci}{[3.4, 6.9]}
\resdef{n11-s030-stop-E02-ALEC-regret}{1.51}
\resdef{n11-s030-stop-E02-ALEC-save2n}{-7}
\resdef{n11-s030-stop-E05-ALEC-q}{20.3}
\resdef{n11-s030-stop-E05-ALEC-hit}{77.7}
\resdef{n11-s030-stop-E05-ALEC-agree}{98.1}
\resdef{n11-s030-stop-E05-ALEC-early}{66.8}
\resdef{n11-s030-stop-E05-ALEC-exceed}{11.1}
\resdef{n11-s030-stop-E05-ALEC-qmax}{44}
\resdef{n11-s030-stop-E05-ALEC-qtot}{31.3}
\resdef{n11-s030-stop-E05-ALEC-errearly}{4.0}
\resdef{n11-s030-stop-E05-ALEC-errci}{[2.8, 5.8]}
\resdef{n11-s030-stop-E05-ALEC-regret}{1.50}
\resdef{n11-s030-stop-E05-ALEC-qround}{20}
\resdef{n11-s030-stop-E05-ALEC-fixedhit}{76.5}
\resdef{n11-s030-stop-E05-ALEC-fixedreg}{1.62}
\resdef{n11-s030-stop-E05-ALEC-save2n}{8}
\resdef{n11-s030-stop-B10-ALEC-cap-q}{15.0}
\resdef{n11-s030-stop-B10-ALEC-cap-hit}{77.2}
\resdef{n11-s030-stop-B10-ALEC-cap-agree}{98.2}
\resdef{n11-s030-stop-B10-ALEC-cap-early}{49.6}
\resdef{n11-s030-stop-B10-ALEC-cap-exceed}{22.8}
\resdef{n11-s030-stop-B10-ALEC-cap-qmax}{22}
\resdef{n11-s030-stop-B10-ALEC-cap-qtot}{26.0}
\resdef{n11-s030-stop-B10-ALEC-cap-errearly}{5.4}
\resdef{n11-s030-stop-B10-ALEC-cap-errci}{[3.8, 7.8]}
\resdef{n11-s030-stop-B10-ALEC-cap-regret}{1.54}
\resdef{n11-s030-stop-B10-ALEC-cap-qround}{15}
\resdef{n11-s030-stop-B10-ALEC-cap-fixedhit}{76.1}
\resdef{n11-s030-stop-B10-ALEC-cap-fixedreg}{1.70}
\resdef{n11-s030-stop-B10-ALEC-cap-save2n}{32}
\resdef{n11-s030-stop-E02-ALEC-cap-q}{14.0}
\resdef{n11-s030-stop-E02-ALEC-cap-hit}{77.3}
\resdef{n11-s030-stop-E02-ALEC-cap-agree}{97.9}
\resdef{n11-s030-stop-E02-ALEC-cap-early}{56.0}
\resdef{n11-s030-stop-E02-ALEC-cap-exceed}{18.0}
\resdef{n11-s030-stop-E02-ALEC-cap-qmax}{22}
\resdef{n11-s030-stop-E02-ALEC-cap-qtot}{25.0}
\resdef{n11-s030-stop-E02-ALEC-cap-errearly}{5.0}
\resdef{n11-s030-stop-E02-ALEC-cap-errci}{[3.5, 7.1]}
\resdef{n11-s030-stop-E02-ALEC-cap-regret}{1.55}
\resdef{n11-s030-stop-E02-ALEC-cap-qround}{14}
\resdef{n11-s030-stop-E02-ALEC-cap-fixedhit}{75.3}
\resdef{n11-s030-stop-E02-ALEC-cap-fixedreg}{1.79}
\resdef{n11-s030-stop-E02-ALEC-cap-save2n}{36}
\resdef{n11-s030-stop-E05-ALEC-cap-q}{12.8}
\resdef{n11-s030-stop-E05-ALEC-cap-hit}{77.6}
\resdef{n11-s030-stop-E05-ALEC-cap-agree}{97.2}
\resdef{n11-s030-stop-E05-ALEC-cap-early}{64.0}
\resdef{n11-s030-stop-E05-ALEC-cap-exceed}{11.3}
\resdef{n11-s030-stop-E05-ALEC-cap-qmax}{22}
\resdef{n11-s030-stop-E05-ALEC-cap-qtot}{23.8}
\resdef{n11-s030-stop-E05-ALEC-cap-errearly}{4.2}
\resdef{n11-s030-stop-E05-ALEC-cap-errci}{[2.9, 6.1]}
\resdef{n11-s030-stop-E05-ALEC-cap-regret}{1.54}
\resdef{n11-s030-stop-E05-ALEC-cap-qround}{13}
\resdef{n11-s030-stop-E05-ALEC-cap-fixedhit}{75.3}
\resdef{n11-s030-stop-E05-ALEC-cap-fixedreg}{1.88}
\resdef{n11-s030-stop-E05-ALEC-cap-save2n}{42}
\resdef{n11-s030-stop-E05-ALEC-cap-s06-q}{16.4}
\resdef{n11-s030-stop-E05-ALEC-cap-s06-hit}{76.6}
\resdef{n11-s030-stop-E05-ALEC-cap-s06-agree}{99.0}
\resdef{n11-s030-stop-E05-ALEC-cap-s06-early}{50.8}
\resdef{n11-s030-stop-E05-ALEC-cap-s06-exceed}{11.6}
\resdef{n11-s030-stop-E05-ALEC-cap-s06-qmax}{22}
\resdef{n11-s030-stop-E05-ALEC-cap-s06-qtot}{27.4}
\resdef{n11-s030-stop-E05-ALEC-cap-s06-errearly}{1.8}
\resdef{n11-s030-stop-E05-ALEC-cap-s06-errci}{[0.9, 3.3]}
\resdef{n11-s030-stop-E05-ALEC-cap-s06-regret}{1.59}
\resdef{n11-s030-stop-E05-ALEC-cap-s06-qround}{16}
\resdef{n11-s030-stop-E05-ALEC-cap-s06-fixedhit}{76.2}
\resdef{n11-s030-stop-E05-ALEC-cap-s06-fixedreg}{1.60}
\resdef{n11-s030-stop-E05-ALEC-cap-s06-save2n}{26}
\resdef{n11-s050-full-hit}{67.7}
\resdef{n11-s050-tot-h}{17}
\resdef{n11-s050-frac-h}{31}
\resdef{n11-s050-RANDOM-hit-h}{51.3}
\resdef{n11-s050-RANDOM-hitllsm-h}{51.3}
\resdef{n11-s050-RANDOM-tau-h}{0.701}
\resdef{n11-s050-RANDOM-reg-h}{9.19}
\resdef{n11-s050-RANDOM-agree-h}{56.3}
\resdef{n11-s050-p-RANDOM-h}{$p<0.001$}
\resdef{n11-s050-p-RANDOM-h-llsm}{$p<0.001$}
\resdef{n11-s050-HARKER-PG-hit-h}{49.2}
\resdef{n11-s050-HARKER-PG-hitllsm-h}{49.7}
\resdef{n11-s050-HARKER-PG-tau-h}{0.674}
\resdef{n11-s050-HARKER-PG-reg-h}{9.72}
\resdef{n11-s050-HARKER-PG-agree-h}{56.4}
\resdef{n11-s050-p-HARKER-PG-h}{$p<0.001$}
\resdef{n11-s050-p-HARKER-PG-h-llsm}{$p<0.001$}
\resdef{n11-s050-HARKER-hit-h}{52.7}
\resdef{n11-s050-HARKER-hitllsm-h}{52.6}
\resdef{n11-s050-HARKER-tau-h}{0.671}
\resdef{n11-s050-HARKER-reg-h}{8.45}
\resdef{n11-s050-HARKER-agree-h}{59.8}
\resdef{n11-s050-p-HARKER-h}{$p=0.002$}
\resdef{n11-s050-p-HARKER-h-llsm}{$p=0.001$}
\resdef{n11-s050-HARKER-LOG-hit-h}{58.2}
\resdef{n11-s050-HARKER-LOG-hitllsm-h}{58.4}
\resdef{n11-s050-HARKER-LOG-tau-h}{0.642}
\resdef{n11-s050-HARKER-LOG-reg-h}{6.19}
\resdef{n11-s050-HARKER-LOG-agree-h}{68.7}
\resdef{n11-s050-p-HARKER-LOG-h}{$p=1.000$}
\resdef{n11-s050-p-HARKER-LOG-h-llsm}{$p=0.939$}
\resdef{n11-s050-AOPT-hit-h}{51.6}
\resdef{n11-s050-AOPT-hitllsm-h}{51.6}
\resdef{n11-s050-AOPT-tau-h}{0.720}
\resdef{n11-s050-AOPT-reg-h}{8.77}
\resdef{n11-s050-AOPT-agree-h}{58.1}
\resdef{n11-s050-p-AOPT-h}{$p<0.001$}
\resdef{n11-s050-p-AOPT-h-llsm}{$p<0.001$}
\resdef{n11-s050-FROZEN-hit-h}{54.4}
\resdef{n11-s050-FROZEN-hitllsm-h}{54.4}
\resdef{n11-s050-FROZEN-tau-h}{0.692}
\resdef{n11-s050-FROZEN-reg-h}{7.53}
\resdef{n11-s050-FROZEN-agree-h}{65.0}
\resdef{n11-s050-p-FROZEN-h}{$p=0.003$}
\resdef{n11-s050-p-FROZEN-h-llsm}{$p=0.003$}
\resdef{n11-s050-STAR-hit-h}{54.6}
\resdef{n11-s050-STAR-hitllsm-h}{54.6}
\resdef{n11-s050-STAR-tau-h}{0.692}
\resdef{n11-s050-STAR-reg-h}{7.54}
\resdef{n11-s050-STAR-agree-h}{65.1}
\resdef{n11-s050-p-STAR-h}{$p=0.004$}
\resdef{n11-s050-p-STAR-h-llsm}{$p=0.004$}
\resdef{n11-s050-ALEC-hit-h}{58.2}
\resdef{n11-s050-ALEC-hitllsm-h}{58.2}
\resdef{n11-s050-ALEC-tau-h}{0.689}
\resdef{n11-s050-ALEC-reg-h}{5.97}
\resdef{n11-s050-ALEC-agree-h}{69.6}
\resdef{n11-s050-ALEC-U-hit-h}{60.2}
\resdef{n11-s050-ALEC-U-hitllsm-h}{60.2}
\resdef{n11-s050-ALEC-U-tau-h}{0.692}
\resdef{n11-s050-ALEC-U-reg-h}{5.68}
\resdef{n11-s050-ALEC-U-agree-h}{69.2}
\resdef{n11-s050-p-ALEC-U-h}{$p=0.140$}
\resdef{n11-s050-p-ALEC-U-h-llsm}{$p=0.140$}
\resdef{n11-s050-tot-n}{22}
\resdef{n11-s050-frac-n}{40}
\resdef{n11-s050-RANDOM-hit-n}{52.3}
\resdef{n11-s050-RANDOM-hitllsm-n}{52.3}
\resdef{n11-s050-RANDOM-tau-n}{0.748}
\resdef{n11-s050-RANDOM-reg-n}{8.15}
\resdef{n11-s050-RANDOM-agree-n}{60.9}
\resdef{n11-s050-p-RANDOM-n}{$p<0.001$}
\resdef{n11-s050-p-RANDOM-n-llsm}{$p<0.001$}
\resdef{n11-s050-HARKER-PG-hit-n}{54.7}
\resdef{n11-s050-HARKER-PG-hitllsm-n}{55.5}
\resdef{n11-s050-HARKER-PG-tau-n}{0.727}
\resdef{n11-s050-HARKER-PG-reg-n}{7.55}
\resdef{n11-s050-HARKER-PG-agree-n}{65.4}
\resdef{n11-s050-p-HARKER-PG-n}{$p<0.001$}
\resdef{n11-s050-p-HARKER-PG-n-llsm}{$p<0.001$}
\resdef{n11-s050-HARKER-hit-n}{56.2}
\resdef{n11-s050-HARKER-hitllsm-n}{56.3}
\resdef{n11-s050-HARKER-tau-n}{0.721}
\resdef{n11-s050-HARKER-reg-n}{7.20}
\resdef{n11-s050-HARKER-agree-n}{66.9}
\resdef{n11-s050-p-HARKER-n}{$p=0.003$}
\resdef{n11-s050-p-HARKER-n-llsm}{$p=0.004$}
\resdef{n11-s050-HARKER-LOG-hit-n}{63.0}
\resdef{n11-s050-HARKER-LOG-hitllsm-n}{64.2}
\resdef{n11-s050-HARKER-LOG-tau-n}{0.689}
\resdef{n11-s050-HARKER-LOG-reg-n}{4.70}
\resdef{n11-s050-HARKER-LOG-agree-n}{76.8}
\resdef{n11-s050-p-HARKER-LOG-n}{$p=0.197$}
\resdef{n11-s050-p-HARKER-LOG-n-llsm}{$p=0.028$}
\resdef{n11-s050-AOPT-hit-n}{56.3}
\resdef{n11-s050-AOPT-hitllsm-n}{56.3}
\resdef{n11-s050-AOPT-tau-n}{0.762}
\resdef{n11-s050-AOPT-reg-n}{6.70}
\resdef{n11-s050-AOPT-agree-n}{64.9}
\resdef{n11-s050-p-AOPT-n}{$p=0.005$}
\resdef{n11-s050-p-AOPT-n-llsm}{$p=0.005$}
\resdef{n11-s050-FROZEN-hit-n}{59.1}
\resdef{n11-s050-FROZEN-hitllsm-n}{59.1}
\resdef{n11-s050-FROZEN-tau-n}{0.740}
\resdef{n11-s050-FROZEN-reg-n}{6.15}
\resdef{n11-s050-FROZEN-agree-n}{72.0}
\resdef{n11-s050-p-FROZEN-n}{$p=0.150$}
\resdef{n11-s050-p-FROZEN-n-llsm}{$p=0.150$}
\resdef{n11-s050-STAR-hit-n}{60.7}
\resdef{n11-s050-STAR-hitllsm-n}{60.7}
\resdef{n11-s050-STAR-tau-n}{0.739}
\resdef{n11-s050-STAR-reg-n}{5.52}
\resdef{n11-s050-STAR-agree-n}{75.6}
\resdef{n11-s050-p-STAR-n}{$p=0.764$}
\resdef{n11-s050-p-STAR-n-llsm}{$p=0.764$}
\resdef{n11-s050-ALEC-hit-n}{61.1}
\resdef{n11-s050-ALEC-hitllsm-n}{61.1}
\resdef{n11-s050-ALEC-tau-n}{0.739}
\resdef{n11-s050-ALEC-reg-n}{4.98}
\resdef{n11-s050-ALEC-agree-n}{76.7}
\resdef{n11-s050-ALEC-U-hit-n}{62.2}
\resdef{n11-s050-ALEC-U-hitllsm-n}{62.2}
\resdef{n11-s050-ALEC-U-tau-n}{0.740}
\resdef{n11-s050-ALEC-U-reg-n}{4.71}
\resdef{n11-s050-ALEC-U-agree-n}{76.8}
\resdef{n11-s050-p-ALEC-U-n}{$p=0.460$}
\resdef{n11-s050-p-ALEC-U-n-llsm}{$p=0.460$}
\resdef{n11-s050-tot-2n}{33}
\resdef{n11-s050-frac-2n}{60}
\resdef{n11-s050-RANDOM-hit-2n}{61.7}
\resdef{n11-s050-RANDOM-hitllsm-2n}{61.7}
\resdef{n11-s050-RANDOM-tau-2n}{0.805}
\resdef{n11-s050-RANDOM-reg-2n}{5.17}
\resdef{n11-s050-RANDOM-agree-2n}{75.1}
\resdef{n11-s050-p-RANDOM-2n}{$p<0.001$}
\resdef{n11-s050-p-RANDOM-2n-llsm}{$p<0.001$}
\resdef{n11-s050-HARKER-PG-hit-2n}{--}
\resdef{n11-s050-HARKER-PG-hitllsm-2n}{--}
\resdef{n11-s050-HARKER-PG-tau-2n}{--}
\resdef{n11-s050-HARKER-PG-reg-2n}{--}
\resdef{n11-s050-HARKER-PG-agree-2n}{--}
\resdef{n11-s050-HARKER-hit-2n}{62.4}
\resdef{n11-s050-HARKER-hitllsm-2n}{62.3}
\resdef{n11-s050-HARKER-tau-2n}{0.788}
\resdef{n11-s050-HARKER-reg-2n}{4.80}
\resdef{n11-s050-HARKER-agree-2n}{78.4}
\resdef{n11-s050-p-HARKER-2n}{$p<0.001$}
\resdef{n11-s050-p-HARKER-2n-llsm}{$p<0.001$}
\resdef{n11-s050-HARKER-LOG-hit-2n}{66.3}
\resdef{n11-s050-HARKER-LOG-hitllsm-2n}{66.4}
\resdef{n11-s050-HARKER-LOG-tau-2n}{0.779}
\resdef{n11-s050-HARKER-LOG-reg-2n}{3.82}
\resdef{n11-s050-HARKER-LOG-agree-2n}{86.8}
\resdef{n11-s050-p-HARKER-LOG-2n}{$p=0.351$}
\resdef{n11-s050-p-HARKER-LOG-2n-llsm}{$p=0.337$}
\resdef{n11-s050-AOPT-hit-2n}{62.0}
\resdef{n11-s050-AOPT-hitllsm-2n}{62.0}
\resdef{n11-s050-AOPT-tau-2n}{0.810}
\resdef{n11-s050-AOPT-reg-2n}{4.89}
\resdef{n11-s050-AOPT-agree-2n}{76.7}
\resdef{n11-s050-p-AOPT-2n}{$p<0.001$}
\resdef{n11-s050-p-AOPT-2n-llsm}{$p<0.001$}
\resdef{n11-s050-FROZEN-hit-2n}{63.2}
\resdef{n11-s050-FROZEN-hitllsm-2n}{63.2}
\resdef{n11-s050-FROZEN-tau-2n}{0.799}
\resdef{n11-s050-FROZEN-reg-2n}{4.53}
\resdef{n11-s050-FROZEN-agree-2n}{87.2}
\resdef{n11-s050-p-FROZEN-2n}{$p<0.001$}
\resdef{n11-s050-p-FROZEN-2n-llsm}{$p<0.001$}
\resdef{n11-s050-STAR-hit-2n}{67.1}
\resdef{n11-s050-STAR-hitllsm-2n}{67.1}
\resdef{n11-s050-STAR-tau-2n}{0.800}
\resdef{n11-s050-STAR-reg-2n}{3.53}
\resdef{n11-s050-STAR-agree-2n}{96.9}
\resdef{n11-s050-p-STAR-2n}{$p=0.549$}
\resdef{n11-s050-p-STAR-2n-llsm}{$p=0.549$}
\resdef{n11-s050-ALEC-hit-2n}{67.4}
\resdef{n11-s050-ALEC-hitllsm-2n}{67.4}
\resdef{n11-s050-ALEC-tau-2n}{0.800}
\resdef{n11-s050-ALEC-reg-2n}{3.45}
\resdef{n11-s050-ALEC-agree-2n}{97.6}
\resdef{n11-s050-ALEC-U-hit-2n}{66.2}
\resdef{n11-s050-ALEC-U-hitllsm-2n}{66.2}
\resdef{n11-s050-ALEC-U-tau-2n}{0.799}
\resdef{n11-s050-ALEC-U-reg-2n}{3.73}
\resdef{n11-s050-ALEC-U-agree-2n}{89.7}
\resdef{n11-s050-p-ALEC-U-2n}{$p=0.213$}
\resdef{n11-s050-p-ALEC-U-2n-llsm}{$p=0.213$}
\resdef{n11-s050-stop-LC-q}{9.8}
\resdef{n11-s050-stop-LC-hit}{62.1}
\resdef{n11-s050-stop-LC-agree}{78.8}
\resdef{n11-s050-stop-LC-early}{100.0}
\resdef{n11-s050-stop-LC-exceed}{37.9}
\resdef{n11-s050-stop-LC-qmax}{22}
\resdef{n11-s050-stop-LC-qtot}{20.8}
\resdef{n11-s050-stop-LC-errearly}{37.9}
\resdef{n11-s050-stop-LC-errci}{[34.9, 40.9]}
\resdef{n11-s050-stop-LC-regret}{4.63}
\resdef{n11-s050-stop-LC-qround}{10}
\resdef{n11-s050-stop-LC-fixedhit}{60.7}
\resdef{n11-s050-stop-LC-fixedreg}{5.32}
\resdef{n11-s050-stop-LC-save2n}{56}
\resdef{n11-s050-stop-B10-ALEC-q}{29.0}
\resdef{n11-s050-stop-B10-ALEC-hit}{66.5}
\resdef{n11-s050-stop-B10-ALEC-agree}{96.6}
\resdef{n11-s050-stop-B10-ALEC-early}{43.0}
\resdef{n11-s050-stop-B10-ALEC-exceed}{33.5}
\resdef{n11-s050-stop-B10-ALEC-qmax}{44}
\resdef{n11-s050-stop-B10-ALEC-qtot}{40.0}
\resdef{n11-s050-stop-B10-ALEC-errearly}{19.5}
\resdef{n11-s050-stop-B10-ALEC-errci}{[16.1, 23.5]}
\resdef{n11-s050-stop-B10-ALEC-regret}{3.74}
\resdef{n11-s050-stop-B10-ALEC-save2n}{-32}
\resdef{n11-s050-stop-B05-ALECU-q}{32.2}
\resdef{n11-s050-stop-B05-ALECU-hit}{67.9}
\resdef{n11-s050-stop-B05-ALECU-agree}{98.8}
\resdef{n11-s050-stop-B05-ALECU-early}{33.9}
\resdef{n11-s050-stop-B05-ALECU-exceed}{32.1}
\resdef{n11-s050-stop-B05-ALECU-qmax}{44}
\resdef{n11-s050-stop-B05-ALECU-qtot}{43.2}
\resdef{n11-s050-stop-B05-ALECU-errearly}{13.3}
\resdef{n11-s050-stop-B05-ALECU-errci}{[10.1, 17.3]}
\resdef{n11-s050-stop-B05-ALECU-regret}{3.37}
\resdef{n11-s050-stop-B05-ALECU-save2n}{-47}
\resdef{n11-s050-stop-B10-ALECU-q}{26.2}
\resdef{n11-s050-stop-B10-ALECU-hit}{66.8}
\resdef{n11-s050-stop-B10-ALECU-agree}{95.4}
\resdef{n11-s050-stop-B10-ALECU-early}{48.8}
\resdef{n11-s050-stop-B10-ALECU-exceed}{33.2}
\resdef{n11-s050-stop-B10-ALECU-qmax}{44}
\resdef{n11-s050-stop-B10-ALECU-qtot}{37.2}
\resdef{n11-s050-stop-B10-ALECU-errearly}{21.3}
\resdef{n11-s050-stop-B10-ALECU-errci}{[17.9, 25.2]}
\resdef{n11-s050-stop-B10-ALECU-regret}{3.76}
\resdef{n11-s050-stop-B10-ALECU-save2n}{-19}
\resdef{n11-s050-stop-B20-ALECU-q}{17.6}
\resdef{n11-s050-stop-B20-ALECU-hit}{63.8}
\resdef{n11-s050-stop-B20-ALECU-agree}{88.1}
\resdef{n11-s050-stop-B20-ALECU-early}{69.1}
\resdef{n11-s050-stop-B20-ALECU-exceed}{36.2}
\resdef{n11-s050-stop-B20-ALECU-qmax}{44}
\resdef{n11-s050-stop-B20-ALECU-qtot}{28.6}
\resdef{n11-s050-stop-B20-ALECU-errearly}{29.5}
\resdef{n11-s050-stop-B20-ALECU-errci}{[26.2, 33.0]}
\resdef{n11-s050-stop-B20-ALECU-regret}{4.80}
\resdef{n11-s050-stop-B20-ALECU-qround}{18}
\resdef{n11-s050-stop-B20-ALECU-fixedhit}{67.6}
\resdef{n11-s050-stop-B20-ALECU-fixedreg}{3.55}
\resdef{n11-s050-stop-B20-ALECU-save2n}{20}
\resdef{n11-s050-stop-B10-ALECU-s06-q}{32.2}
\resdef{n11-s050-stop-B10-ALECU-s06-hit}{67.7}
\resdef{n11-s050-stop-B10-ALECU-s06-agree}{99.4}
\resdef{n11-s050-stop-B10-ALECU-s06-early}{35.9}
\resdef{n11-s050-stop-B10-ALECU-s06-exceed}{32.3}
\resdef{n11-s050-stop-B10-ALECU-s06-qmax}{44}
\resdef{n11-s050-stop-B10-ALECU-s06-qtot}{43.2}
\resdef{n11-s050-stop-B10-ALECU-s06-errearly}{11.1}
\resdef{n11-s050-stop-B10-ALECU-s06-errci}{[8.3, 14.8]}
\resdef{n11-s050-stop-B10-ALECU-s06-regret}{3.35}
\resdef{n11-s050-stop-B10-ALECU-s06-save2n}{-47}
\resdef{n11-s050-stop-B10-ALECU-cap-q}{14.8}
\resdef{n11-s050-stop-B10-ALECU-cap-hit}{66.1}
\resdef{n11-s050-stop-B10-ALECU-cap-agree}{86.3}
\resdef{n11-s050-stop-B10-ALECU-cap-early}{46.5}
\resdef{n11-s050-stop-B10-ALECU-cap-exceed}{33.9}
\resdef{n11-s050-stop-B10-ALECU-cap-qmax}{22}
\resdef{n11-s050-stop-B10-ALECU-cap-qtot}{25.8}
\resdef{n11-s050-stop-B10-ALECU-cap-errearly}{22.2}
\resdef{n11-s050-stop-B10-ALECU-cap-errci}{[18.6, 26.1]}
\resdef{n11-s050-stop-B10-ALECU-cap-regret}{3.99}
\resdef{n11-s050-stop-B10-ALECU-cap-qround}{15}
\resdef{n11-s050-stop-B10-ALECU-cap-fixedhit}{64.7}
\resdef{n11-s050-stop-B10-ALECU-cap-fixedreg}{4.23}
\resdef{n11-s050-stop-B10-ALECU-cap-save2n}{33}
\resdef{n11-s050-stop-B10-ALECU-oracle-q}{30.1}
\resdef{n11-s050-stop-B10-ALECU-oracle-hit}{67.3}
\resdef{n11-s050-stop-B10-ALECU-oracle-agree}{98.1}
\resdef{n11-s050-stop-B10-ALECU-oracle-early}{38.6}
\resdef{n11-s050-stop-B10-ALECU-oracle-exceed}{32.7}
\resdef{n11-s050-stop-B10-ALECU-oracle-qmax}{44}
\resdef{n11-s050-stop-B10-ALECU-oracle-qtot}{41.1}
\resdef{n11-s050-stop-B10-ALECU-oracle-errearly}{14.2}
\resdef{n11-s050-stop-B10-ALECU-oracle-errci}{[11.1, 18.1]}
\resdef{n11-s050-stop-B10-ALECU-oracle-regret}{3.46}
\resdef{n11-s050-stop-B10-ALECU-oracle-save2n}{-37}
\resdef{n11-s050-stop-E02-ALEC-q}{27.5}
\resdef{n11-s050-stop-E02-ALEC-hit}{67.1}
\resdef{n11-s050-stop-E02-ALEC-agree}{95.8}
\resdef{n11-s050-stop-E02-ALEC-early}{46.8}
\resdef{n11-s050-stop-E02-ALEC-exceed}{29.0}
\resdef{n11-s050-stop-E02-ALEC-qmax}{44}
\resdef{n11-s050-stop-E02-ALEC-qtot}{38.5}
\resdef{n11-s050-stop-E02-ALEC-errearly}{17.1}
\resdef{n11-s050-stop-E02-ALEC-errci}{[14.0, 20.8]}
\resdef{n11-s050-stop-E02-ALEC-regret}{3.72}
\resdef{n11-s050-stop-E02-ALEC-save2n}{-25}
\resdef{n11-s050-stop-E05-ALEC-q}{25.0}
\resdef{n11-s050-stop-E05-ALEC-hit}{66.9}
\resdef{n11-s050-stop-E05-ALEC-agree}{94.6}
\resdef{n11-s050-stop-E05-ALEC-early}{53.6}
\resdef{n11-s050-stop-E05-ALEC-exceed}{24.1}
\resdef{n11-s050-stop-E05-ALEC-qmax}{44}
\resdef{n11-s050-stop-E05-ALEC-qtot}{36.0}
\resdef{n11-s050-stop-E05-ALEC-errearly}{16.4}
\resdef{n11-s050-stop-E05-ALEC-errci}{[13.5, 19.8]}
\resdef{n11-s050-stop-E05-ALEC-regret}{3.84}
\resdef{n11-s050-stop-E05-ALEC-save2n}{-14}
\resdef{n11-s050-stop-B10-ALEC-cap-q}{16.2}
\resdef{n11-s050-stop-B10-ALEC-cap-hit}{66.2}
\resdef{n11-s050-stop-B10-ALEC-cap-agree}{94.2}
\resdef{n11-s050-stop-B10-ALEC-cap-early}{39.6}
\resdef{n11-s050-stop-B10-ALEC-cap-exceed}{33.8}
\resdef{n11-s050-stop-B10-ALEC-cap-qmax}{22}
\resdef{n11-s050-stop-B10-ALEC-cap-qtot}{27.2}
\resdef{n11-s050-stop-B10-ALEC-cap-errearly}{20.5}
\resdef{n11-s050-stop-B10-ALEC-cap-errci}{[16.8, 24.7]}
\resdef{n11-s050-stop-B10-ALEC-cap-regret}{3.85}
\resdef{n11-s050-stop-B10-ALEC-cap-qround}{16}
\resdef{n11-s050-stop-B10-ALEC-cap-fixedhit}{66.4}
\resdef{n11-s050-stop-B10-ALEC-cap-fixedreg}{3.78}
\resdef{n11-s050-stop-B10-ALEC-cap-save2n}{26}
\resdef{n11-s050-stop-E02-ALEC-cap-q}{15.6}
\resdef{n11-s050-stop-E02-ALEC-cap-hit}{66.8}
\resdef{n11-s050-stop-E02-ALEC-cap-agree}{93.4}
\resdef{n11-s050-stop-E02-ALEC-cap-early}{43.9}
\resdef{n11-s050-stop-E02-ALEC-cap-exceed}{29.3}
\resdef{n11-s050-stop-E02-ALEC-cap-qmax}{22}
\resdef{n11-s050-stop-E02-ALEC-cap-qtot}{26.6}
\resdef{n11-s050-stop-E02-ALEC-cap-errearly}{18.0}
\resdef{n11-s050-stop-E02-ALEC-cap-errci}{[14.7, 21.9]}
\resdef{n11-s050-stop-E02-ALEC-cap-regret}{3.82}
\resdef{n11-s050-stop-E02-ALEC-cap-qround}{16}
\resdef{n11-s050-stop-E02-ALEC-cap-fixedhit}{66.4}
\resdef{n11-s050-stop-E02-ALEC-cap-fixedreg}{3.78}
\resdef{n11-s050-stop-E02-ALEC-cap-save2n}{29}
\resdef{n11-s050-stop-E05-ALEC-cap-q}{14.6}
\resdef{n11-s050-stop-E05-ALEC-cap-hit}{66.6}
\resdef{n11-s050-stop-E05-ALEC-cap-agree}{92.2}
\resdef{n11-s050-stop-E05-ALEC-cap-early}{50.4}
\resdef{n11-s050-stop-E05-ALEC-cap-exceed}{24.6}
\resdef{n11-s050-stop-E05-ALEC-cap-qmax}{22}
\resdef{n11-s050-stop-E05-ALEC-cap-qtot}{25.6}
\resdef{n11-s050-stop-E05-ALEC-cap-errearly}{17.1}
\resdef{n11-s050-stop-E05-ALEC-cap-errci}{[14.0, 20.6]}
\resdef{n11-s050-stop-E05-ALEC-cap-regret}{3.94}
\resdef{n11-s050-stop-E05-ALEC-cap-qround}{15}
\resdef{n11-s050-stop-E05-ALEC-cap-fixedhit}{64.7}
\resdef{n11-s050-stop-E05-ALEC-cap-fixedreg}{4.23}
\resdef{n11-s050-stop-E05-ALEC-cap-save2n}{33}
\resdef{n11-s050-stop-E05-ALEC-cap-s06-q}{18.2}
\resdef{n11-s050-stop-E05-ALEC-cap-s06-hit}{67.7}
\resdef{n11-s050-stop-E05-ALEC-cap-s06-agree}{97.0}
\resdef{n11-s050-stop-E05-ALEC-cap-s06-early}{35.6}
\resdef{n11-s050-stop-E05-ALEC-cap-s06-exceed}{22.9}
\resdef{n11-s050-stop-E05-ALEC-cap-s06-qmax}{22}
\resdef{n11-s050-stop-E05-ALEC-cap-s06-qtot}{29.2}
\resdef{n11-s050-stop-E05-ALEC-cap-s06-errearly}{7.6}
\resdef{n11-s050-stop-E05-ALEC-cap-s06-errci}{[5.3, 10.8]}
\resdef{n11-s050-stop-E05-ALEC-cap-s06-regret}{3.47}
\resdef{n11-s050-stop-E05-ALEC-cap-s06-qround}{18}
\resdef{n11-s050-stop-E05-ALEC-cap-s06-fixedhit}{67.6}
\resdef{n11-s050-stop-E05-ALEC-cap-s06-fixedreg}{3.55}
\resdef{n11-s050-stop-E05-ALEC-cap-s06-save2n}{17}
\resdef{n15-s015-full-hit}{85.5}
\resdef{n15-s015-tot-h}{23}
\resdef{n15-s015-frac-h}{22}
\resdef{n15-s015-RANDOM-hit-h}{69.2}
\resdef{n15-s015-RANDOM-hitllsm-h}{69.2}
\resdef{n15-s015-RANDOM-tau-h}{0.895}
\resdef{n15-s015-RANDOM-reg-h}{2.44}
\resdef{n15-s015-RANDOM-agree-h}{72.5}
\resdef{n15-s015-p-RANDOM-h}{$p<0.001$}
\resdef{n15-s015-p-RANDOM-h-llsm}{$p<0.001$}
\resdef{n15-s015-HARKER-PG-hit-h}{72.7}
\resdef{n15-s015-HARKER-PG-hitllsm-h}{73.7}
\resdef{n15-s015-HARKER-PG-tau-h}{0.883}
\resdef{n15-s015-HARKER-PG-reg-h}{1.98}
\resdef{n15-s015-HARKER-PG-agree-h}{75.2}
\resdef{n15-s015-p-HARKER-PG-h}{$p<0.001$}
\resdef{n15-s015-p-HARKER-PG-h-llsm}{$p=0.003$}
\resdef{n15-s015-HARKER-hit-h}{73.9}
\resdef{n15-s015-HARKER-hitllsm-h}{74.1}
\resdef{n15-s015-HARKER-tau-h}{0.887}
\resdef{n15-s015-HARKER-reg-h}{1.80}
\resdef{n15-s015-HARKER-agree-h}{76.9}
\resdef{n15-s015-p-HARKER-h}{$p=0.004$}
\resdef{n15-s015-p-HARKER-h-llsm}{$p=0.006$}
\resdef{n15-s015-HARKER-LOG-hit-h}{77.7}
\resdef{n15-s015-HARKER-LOG-hitllsm-h}{78.3}
\resdef{n15-s015-HARKER-LOG-tau-h}{0.884}
\resdef{n15-s015-HARKER-LOG-reg-h}{1.25}
\resdef{n15-s015-HARKER-LOG-agree-h}{81.6}
\resdef{n15-s015-p-HARKER-LOG-h}{$p=0.587$}
\resdef{n15-s015-p-HARKER-LOG-h-llsm}{$p=0.939$}
\resdef{n15-s015-AOPT-hit-h}{73.4}
\resdef{n15-s015-AOPT-hitllsm-h}{73.4}
\resdef{n15-s015-AOPT-tau-h}{0.904}
\resdef{n15-s015-AOPT-reg-h}{1.87}
\resdef{n15-s015-AOPT-agree-h}{76.1}
\resdef{n15-s015-p-AOPT-h}{$p<0.001$}
\resdef{n15-s015-p-AOPT-h-llsm}{$p<0.001$}
\resdef{n15-s015-FROZEN-hit-h}{75.2}
\resdef{n15-s015-FROZEN-hitllsm-h}{75.2}
\resdef{n15-s015-FROZEN-tau-h}{0.890}
\resdef{n15-s015-FROZEN-reg-h}{1.55}
\resdef{n15-s015-FROZEN-agree-h}{79.0}
\resdef{n15-s015-p-FROZEN-h}{$p=0.002$}
\resdef{n15-s015-p-FROZEN-h-llsm}{$p=0.002$}
\resdef{n15-s015-STAR-hit-h}{75.2}
\resdef{n15-s015-STAR-hitllsm-h}{75.2}
\resdef{n15-s015-STAR-tau-h}{0.890}
\resdef{n15-s015-STAR-reg-h}{1.55}
\resdef{n15-s015-STAR-agree-h}{79.0}
\resdef{n15-s015-p-STAR-h}{$p=0.002$}
\resdef{n15-s015-p-STAR-h-llsm}{$p=0.002$}
\resdef{n15-s015-ALEC-hit-h}{78.5}
\resdef{n15-s015-ALEC-hitllsm-h}{78.5}
\resdef{n15-s015-ALEC-tau-h}{0.889}
\resdef{n15-s015-ALEC-reg-h}{1.24}
\resdef{n15-s015-ALEC-agree-h}{82.7}
\resdef{n15-s015-ALEC-U-hit-h}{80.3}
\resdef{n15-s015-ALEC-U-hitllsm-h}{80.3}
\resdef{n15-s015-ALEC-U-tau-h}{0.886}
\resdef{n15-s015-ALEC-U-reg-h}{0.93}
\resdef{n15-s015-ALEC-U-agree-h}{83.5}
\resdef{n15-s015-p-ALEC-U-h}{$p=0.145$}
\resdef{n15-s015-p-ALEC-U-h-llsm}{$p=0.145$}
\resdef{n15-s015-tot-n}{30}
\resdef{n15-s015-frac-n}{29}
\resdef{n15-s015-RANDOM-hit-n}{74.2}
\resdef{n15-s015-RANDOM-hitllsm-n}{74.2}
\resdef{n15-s015-RANDOM-tau-n}{0.914}
\resdef{n15-s015-RANDOM-reg-n}{1.59}
\resdef{n15-s015-RANDOM-agree-n}{78.0}
\resdef{n15-s015-p-RANDOM-n}{$p<0.001$}
\resdef{n15-s015-p-RANDOM-n-llsm}{$p<0.001$}
\resdef{n15-s015-HARKER-PG-hit-n}{77.7}
\resdef{n15-s015-HARKER-PG-hitllsm-n}{79.1}
\resdef{n15-s015-HARKER-PG-tau-n}{0.904}
\resdef{n15-s015-HARKER-PG-reg-n}{1.30}
\resdef{n15-s015-HARKER-PG-agree-n}{81.1}
\resdef{n15-s015-p-HARKER-PG-n}{$p=0.003$}
\resdef{n15-s015-p-HARKER-PG-n-llsm}{$p=0.050$}
\resdef{n15-s015-HARKER-hit-n}{78.7}
\resdef{n15-s015-HARKER-hitllsm-n}{78.6}
\resdef{n15-s015-HARKER-tau-n}{0.908}
\resdef{n15-s015-HARKER-reg-n}{1.14}
\resdef{n15-s015-HARKER-agree-n}{82.8}
\resdef{n15-s015-p-HARKER-n}{$p=0.023$}
\resdef{n15-s015-p-HARKER-n-llsm}{$p=0.019$}
\resdef{n15-s015-HARKER-LOG-hit-n}{82.4}
\resdef{n15-s015-HARKER-LOG-hitllsm-n}{82.8}
\resdef{n15-s015-HARKER-LOG-tau-n}{0.905}
\resdef{n15-s015-HARKER-LOG-reg-n}{0.78}
\resdef{n15-s015-HARKER-LOG-agree-n}{86.7}
\resdef{n15-s015-p-HARKER-LOG-n}{$p=0.857$}
\resdef{n15-s015-p-HARKER-LOG-n-llsm}{$p=0.576$}
\resdef{n15-s015-AOPT-hit-n}{77.5}
\resdef{n15-s015-AOPT-hitllsm-n}{77.5}
\resdef{n15-s015-AOPT-tau-n}{0.921}
\resdef{n15-s015-AOPT-reg-n}{1.34}
\resdef{n15-s015-AOPT-agree-n}{80.5}
\resdef{n15-s015-p-AOPT-n}{$p=0.001$}
\resdef{n15-s015-p-AOPT-n-llsm}{$p=0.001$}
\resdef{n15-s015-FROZEN-hit-n}{79.9}
\resdef{n15-s015-FROZEN-hitllsm-n}{79.9}
\resdef{n15-s015-FROZEN-tau-n}{0.912}
\resdef{n15-s015-FROZEN-reg-n}{1.01}
\resdef{n15-s015-FROZEN-agree-n}{84.3}
\resdef{n15-s015-p-FROZEN-n}{$p=0.053$}
\resdef{n15-s015-p-FROZEN-n-llsm}{$p=0.053$}
\resdef{n15-s015-STAR-hit-n}{82.1}
\resdef{n15-s015-STAR-hitllsm-n}{82.1}
\resdef{n15-s015-STAR-tau-n}{0.913}
\resdef{n15-s015-STAR-reg-n}{0.80}
\resdef{n15-s015-STAR-agree-n}{86.4}
\resdef{n15-s015-p-STAR-n}{$p=1.000$}
\resdef{n15-s015-p-STAR-n-llsm}{$p=1.000$}
\resdef{n15-s015-ALEC-hit-n}{82.1}
\resdef{n15-s015-ALEC-hitllsm-n}{82.1}
\resdef{n15-s015-ALEC-tau-n}{0.911}
\resdef{n15-s015-ALEC-reg-n}{0.72}
\resdef{n15-s015-ALEC-agree-n}{87.2}
\resdef{n15-s015-ALEC-U-hit-n}{83.3}
\resdef{n15-s015-ALEC-U-hitllsm-n}{83.3}
\resdef{n15-s015-ALEC-U-tau-n}{0.902}
\resdef{n15-s015-ALEC-U-reg-n}{0.66}
\resdef{n15-s015-ALEC-U-agree-n}{88.7}
\resdef{n15-s015-p-ALEC-U-n}{$p=0.311$}
\resdef{n15-s015-p-ALEC-U-n-llsm}{$p=0.311$}
\resdef{n15-s015-tot-2n}{45}
\resdef{n15-s015-frac-2n}{43}
\resdef{n15-s015-RANDOM-hit-2n}{78.8}
\resdef{n15-s015-RANDOM-hitllsm-2n}{78.8}
\resdef{n15-s015-RANDOM-tau-2n}{0.935}
\resdef{n15-s015-RANDOM-reg-2n}{1.05}
\resdef{n15-s015-RANDOM-agree-2n}{83.5}
\resdef{n15-s015-p-RANDOM-2n}{$p<0.001$}
\resdef{n15-s015-p-RANDOM-2n-llsm}{$p<0.001$}
\resdef{n15-s015-HARKER-PG-hit-2n}{--}
\resdef{n15-s015-HARKER-PG-hitllsm-2n}{--}
\resdef{n15-s015-HARKER-PG-tau-2n}{--}
\resdef{n15-s015-HARKER-PG-reg-2n}{--}
\resdef{n15-s015-HARKER-PG-agree-2n}{--}
\resdef{n15-s015-HARKER-hit-2n}{83.3}
\resdef{n15-s015-HARKER-hitllsm-2n}{83.0}
\resdef{n15-s015-HARKER-tau-2n}{0.931}
\resdef{n15-s015-HARKER-reg-2n}{0.63}
\resdef{n15-s015-HARKER-agree-2n}{89.8}
\resdef{n15-s015-p-HARKER-2n}{$p=0.010$}
\resdef{n15-s015-p-HARKER-2n-llsm}{$p=0.004$}
\resdef{n15-s015-HARKER-LOG-hit-2n}{85.6}
\resdef{n15-s015-HARKER-LOG-hitllsm-2n}{85.3}
\resdef{n15-s015-HARKER-LOG-tau-2n}{0.928}
\resdef{n15-s015-HARKER-LOG-reg-2n}{0.45}
\resdef{n15-s015-HARKER-LOG-agree-2n}{92.6}
\resdef{n15-s015-p-HARKER-LOG-2n}{$p=0.810$}
\resdef{n15-s015-p-HARKER-LOG-2n-llsm}{$p=0.550$}
\resdef{n15-s015-AOPT-hit-2n}{81.3}
\resdef{n15-s015-AOPT-hitllsm-2n}{81.3}
\resdef{n15-s015-AOPT-tau-2n}{0.940}
\resdef{n15-s015-AOPT-reg-2n}{0.82}
\resdef{n15-s015-AOPT-agree-2n}{87.4}
\resdef{n15-s015-p-AOPT-2n}{$p<0.001$}
\resdef{n15-s015-p-AOPT-2n-llsm}{$p<0.001$}
\resdef{n15-s015-FROZEN-hit-2n}{84.8}
\resdef{n15-s015-FROZEN-hitllsm-2n}{84.8}
\resdef{n15-s015-FROZEN-tau-2n}{0.931}
\resdef{n15-s015-FROZEN-reg-2n}{0.55}
\resdef{n15-s015-FROZEN-agree-2n}{94.7}
\resdef{n15-s015-p-FROZEN-2n}{$p=0.144$}
\resdef{n15-s015-p-FROZEN-2n-llsm}{$p=0.144$}
\resdef{n15-s015-STAR-hit-2n}{85.7}
\resdef{n15-s015-STAR-hitllsm-2n}{85.7}
\resdef{n15-s015-STAR-tau-2n}{0.930}
\resdef{n15-s015-STAR-reg-2n}{0.43}
\resdef{n15-s015-STAR-agree-2n}{97.7}
\resdef{n15-s015-p-STAR-2n}{$p=0.791$}
\resdef{n15-s015-p-STAR-2n-llsm}{$p=0.791$}
\resdef{n15-s015-ALEC-hit-2n}{85.9}
\resdef{n15-s015-ALEC-hitllsm-2n}{85.9}
\resdef{n15-s015-ALEC-tau-2n}{0.930}
\resdef{n15-s015-ALEC-reg-2n}{0.42}
\resdef{n15-s015-ALEC-agree-2n}{98.5}
\resdef{n15-s015-ALEC-U-hit-2n}{84.8}
\resdef{n15-s015-ALEC-U-hitllsm-2n}{84.8}
\resdef{n15-s015-ALEC-U-tau-2n}{0.928}
\resdef{n15-s015-ALEC-U-reg-2n}{0.46}
\resdef{n15-s015-ALEC-U-agree-2n}{94.9}
\resdef{n15-s015-p-ALEC-U-2n}{$p=0.144$}
\resdef{n15-s015-p-ALEC-U-2n-llsm}{$p=0.144$}
\resdef{n15-s015-stop-LC-q}{13.1}
\resdef{n15-s015-stop-LC-hit}{80.3}
\resdef{n15-s015-stop-LC-agree}{85.9}
\resdef{n15-s015-stop-LC-early}{100.0}
\resdef{n15-s015-stop-LC-exceed}{19.7}
\resdef{n15-s015-stop-LC-qmax}{25}
\resdef{n15-s015-stop-LC-qtot}{28.1}
\resdef{n15-s015-stop-LC-errearly}{19.7}
\resdef{n15-s015-stop-LC-errci}{[17.4, 22.3]}
\resdef{n15-s015-stop-LC-regret}{0.93}
\resdef{n15-s015-stop-LC-qround}{13}
\resdef{n15-s015-stop-LC-fixedhit}{79.2}
\resdef{n15-s015-stop-LC-fixedreg}{1.04}
\resdef{n15-s015-stop-LC-save2n}{56}
\resdef{n15-s015-stop-B10-ALEC-q}{41.1}
\resdef{n15-s015-stop-B10-ALEC-hit}{85.5}
\resdef{n15-s015-stop-B10-ALEC-agree}{99.8}
\resdef{n15-s015-stop-B10-ALEC-early}{64.2}
\resdef{n15-s015-stop-B10-ALEC-exceed}{14.5}
\resdef{n15-s015-stop-B10-ALEC-qmax}{90}
\resdef{n15-s015-stop-B10-ALEC-qtot}{56.1}
\resdef{n15-s015-stop-B10-ALEC-errearly}{3.6}
\resdef{n15-s015-stop-B10-ALEC-errci}{[2.4, 5.3]}
\resdef{n15-s015-stop-B10-ALEC-regret}{0.43}
\resdef{n15-s015-stop-B10-ALEC-save2n}{-37}
\resdef{n15-s015-stop-B05-ALECU-q}{48.1}
\resdef{n15-s015-stop-B05-ALECU-hit}{85.5}
\resdef{n15-s015-stop-B05-ALECU-agree}{100.0}
\resdef{n15-s015-stop-B05-ALECU-early}{54.3}
\resdef{n15-s015-stop-B05-ALECU-exceed}{14.5}
\resdef{n15-s015-stop-B05-ALECU-qmax}{90}
\resdef{n15-s015-stop-B05-ALECU-qtot}{63.1}
\resdef{n15-s015-stop-B05-ALECU-errearly}{1.1}
\resdef{n15-s015-stop-B05-ALECU-errci}{[0.5, 2.4]}
\resdef{n15-s015-stop-B05-ALECU-regret}{0.42}
\resdef{n15-s015-stop-B05-ALECU-save2n}{-60}
\resdef{n15-s015-stop-B10-ALECU-q}{38.2}
\resdef{n15-s015-stop-B10-ALECU-hit}{85.5}
\resdef{n15-s015-stop-B10-ALECU-agree}{99.6}
\resdef{n15-s015-stop-B10-ALECU-early}{66.1}
\resdef{n15-s015-stop-B10-ALECU-exceed}{14.5}
\resdef{n15-s015-stop-B10-ALECU-qmax}{90}
\resdef{n15-s015-stop-B10-ALECU-qtot}{53.2}
\resdef{n15-s015-stop-B10-ALECU-errearly}{4.1}
\resdef{n15-s015-stop-B10-ALECU-errci}{[2.8, 5.9]}
\resdef{n15-s015-stop-B10-ALECU-regret}{0.42}
\resdef{n15-s015-stop-B10-ALECU-save2n}{-27}
\resdef{n15-s015-stop-B20-ALECU-q}{24.9}
\resdef{n15-s015-stop-B20-ALECU-hit}{84.8}
\resdef{n15-s015-stop-B20-ALECU-agree}{98.0}
\resdef{n15-s015-stop-B20-ALECU-early}{80.6}
\resdef{n15-s015-stop-B20-ALECU-exceed}{15.2}
\resdef{n15-s015-stop-B20-ALECU-qmax}{90}
\resdef{n15-s015-stop-B20-ALECU-qtot}{39.9}
\resdef{n15-s015-stop-B20-ALECU-errearly}{7.7}
\resdef{n15-s015-stop-B20-ALECU-errci}{[6.0, 9.7]}
\resdef{n15-s015-stop-B20-ALECU-regret}{0.49}
\resdef{n15-s015-stop-B20-ALECU-qround}{25}
\resdef{n15-s015-stop-B20-ALECU-fixedhit}{85.4}
\resdef{n15-s015-stop-B20-ALECU-fixedreg}{0.48}
\resdef{n15-s015-stop-B20-ALECU-save2n}{17}
\resdef{n15-s015-stop-B10-ALECU-s06-q}{50.2}
\resdef{n15-s015-stop-B10-ALECU-s06-hit}{85.5}
\resdef{n15-s015-stop-B10-ALECU-s06-agree}{100.0}
\resdef{n15-s015-stop-B10-ALECU-s06-early}{56.5}
\resdef{n15-s015-stop-B10-ALECU-s06-exceed}{14.5}
\resdef{n15-s015-stop-B10-ALECU-s06-qmax}{90}
\resdef{n15-s015-stop-B10-ALECU-s06-qtot}{65.2}
\resdef{n15-s015-stop-B10-ALECU-s06-errearly}{1.6}
\resdef{n15-s015-stop-B10-ALECU-s06-errci}{[0.8, 3.0]}
\resdef{n15-s015-stop-B10-ALECU-s06-regret}{0.42}
\resdef{n15-s015-stop-B10-ALECU-s06-save2n}{-67}
\resdef{n15-s015-stop-B10-ALECU-cap-q}{17.4}
\resdef{n15-s015-stop-B10-ALECU-cap-hit}{84.9}
\resdef{n15-s015-stop-B10-ALECU-cap-agree}{94.6}
\resdef{n15-s015-stop-B10-ALECU-cap-early}{63.1}
\resdef{n15-s015-stop-B10-ALECU-cap-exceed}{15.1}
\resdef{n15-s015-stop-B10-ALECU-cap-qmax}{30}
\resdef{n15-s015-stop-B10-ALECU-cap-qtot}{32.4}
\resdef{n15-s015-stop-B10-ALECU-cap-errearly}{4.1}
\resdef{n15-s015-stop-B10-ALECU-cap-errci}{[2.8, 6.0]}
\resdef{n15-s015-stop-B10-ALECU-cap-regret}{0.46}
\resdef{n15-s015-stop-B10-ALECU-cap-qround}{17}
\resdef{n15-s015-stop-B10-ALECU-cap-fixedhit}{81.7}
\resdef{n15-s015-stop-B10-ALECU-cap-fixedreg}{0.73}
\resdef{n15-s015-stop-B10-ALECU-cap-save2n}{42}
\resdef{n15-s015-stop-B10-ALECU-oracle-q}{26.3}
\resdef{n15-s015-stop-B10-ALECU-oracle-hit}{84.8}
\resdef{n15-s015-stop-B10-ALECU-oracle-agree}{95.7}
\resdef{n15-s015-stop-B10-ALECU-oracle-early}{75.5}
\resdef{n15-s015-stop-B10-ALECU-oracle-exceed}{15.2}
\resdef{n15-s015-stop-B10-ALECU-oracle-qmax}{90}
\resdef{n15-s015-stop-B10-ALECU-oracle-qtot}{41.3}
\resdef{n15-s015-stop-B10-ALECU-oracle-errearly}{8.6}
\resdef{n15-s015-stop-B10-ALECU-oracle-errci}{[6.8, 10.8]}
\resdef{n15-s015-stop-B10-ALECU-oracle-regret}{0.59}
\resdef{n15-s015-stop-B10-ALECU-oracle-qround}{26}
\resdef{n15-s015-stop-B10-ALECU-oracle-fixedhit}{85.2}
\resdef{n15-s015-stop-B10-ALECU-oracle-fixedreg}{0.47}
\resdef{n15-s015-stop-B10-ALECU-oracle-save2n}{12}
\resdef{n15-s015-stop-E02-ALEC-q}{32.2}
\resdef{n15-s015-stop-E02-ALEC-hit}{85.5}
\resdef{n15-s015-stop-E02-ALEC-agree}{99.6}
\resdef{n15-s015-stop-E02-ALEC-early}{75.5}
\resdef{n15-s015-stop-E02-ALEC-exceed}{7.8}
\resdef{n15-s015-stop-E02-ALEC-qmax}{90}
\resdef{n15-s015-stop-E02-ALEC-qtot}{47.2}
\resdef{n15-s015-stop-E02-ALEC-errearly}{1.5}
\resdef{n15-s015-stop-E02-ALEC-errci}{[0.8, 2.6]}
\resdef{n15-s015-stop-E02-ALEC-regret}{0.43}
\resdef{n15-s015-stop-E02-ALEC-save2n}{-7}
\resdef{n15-s015-stop-E05-ALEC-q}{20.0}
\resdef{n15-s015-stop-E05-ALEC-hit}{84.8}
\resdef{n15-s015-stop-E05-ALEC-agree}{98.1}
\resdef{n15-s015-stop-E05-ALEC-early}{90.7}
\resdef{n15-s015-stop-E05-ALEC-exceed}{3.5}
\resdef{n15-s015-stop-E05-ALEC-qmax}{90}
\resdef{n15-s015-stop-E05-ALEC-qtot}{35.0}
\resdef{n15-s015-stop-E05-ALEC-errearly}{2.3}
\resdef{n15-s015-stop-E05-ALEC-errci}{[1.5, 3.5]}
\resdef{n15-s015-stop-E05-ALEC-regret}{0.46}
\resdef{n15-s015-stop-E05-ALEC-qround}{20}
\resdef{n15-s015-stop-E05-ALEC-fixedhit}{84.3}
\resdef{n15-s015-stop-E05-ALEC-fixedreg}{0.56}
\resdef{n15-s015-stop-E05-ALEC-save2n}{33}
\resdef{n15-s015-stop-B10-ALEC-cap-q}{19.0}
\resdef{n15-s015-stop-B10-ALEC-cap-hit}{86.0}
\resdef{n15-s015-stop-B10-ALEC-cap-agree}{98.4}
\resdef{n15-s015-stop-B10-ALEC-cap-early}{60.9}
\resdef{n15-s015-stop-B10-ALEC-cap-exceed}{14.0}
\resdef{n15-s015-stop-B10-ALEC-cap-qmax}{30}
\resdef{n15-s015-stop-B10-ALEC-cap-qtot}{34.0}
\resdef{n15-s015-stop-B10-ALEC-cap-errearly}{3.3}
\resdef{n15-s015-stop-B10-ALEC-cap-errci}{[2.1, 5.0]}
\resdef{n15-s015-stop-B10-ALEC-cap-regret}{0.42}
\resdef{n15-s015-stop-B10-ALEC-cap-qround}{19}
\resdef{n15-s015-stop-B10-ALEC-cap-fixedhit}{84.1}
\resdef{n15-s015-stop-B10-ALEC-cap-fixedreg}{0.56}
\resdef{n15-s015-stop-B10-ALEC-cap-save2n}{37}
\resdef{n15-s015-stop-E02-ALEC-cap-q}{17.0}
\resdef{n15-s015-stop-E02-ALEC-cap-hit}{86.0}
\resdef{n15-s015-stop-E02-ALEC-cap-agree}{98.2}
\resdef{n15-s015-stop-E02-ALEC-cap-early}{72.2}
\resdef{n15-s015-stop-E02-ALEC-cap-exceed}{7.6}
\resdef{n15-s015-stop-E02-ALEC-cap-qmax}{30}
\resdef{n15-s015-stop-E02-ALEC-cap-qtot}{32.0}
\resdef{n15-s015-stop-E02-ALEC-cap-errearly}{1.4}
\resdef{n15-s015-stop-E02-ALEC-cap-errci}{[0.8, 2.5]}
\resdef{n15-s015-stop-E02-ALEC-cap-regret}{0.42}
\resdef{n15-s015-stop-E02-ALEC-cap-qround}{17}
\resdef{n15-s015-stop-E02-ALEC-cap-fixedhit}{81.7}
\resdef{n15-s015-stop-E02-ALEC-cap-fixedreg}{0.73}
\resdef{n15-s015-stop-E02-ALEC-cap-save2n}{43}
\resdef{n15-s015-stop-E05-ALEC-cap-q}{14.0}
\resdef{n15-s015-stop-E05-ALEC-cap-hit}{85.3}
\resdef{n15-s015-stop-E05-ALEC-cap-agree}{96.9}
\resdef{n15-s015-stop-E05-ALEC-cap-early}{87.0}
\resdef{n15-s015-stop-E05-ALEC-cap-exceed}{3.6}
\resdef{n15-s015-stop-E05-ALEC-cap-qmax}{30}
\resdef{n15-s015-stop-E05-ALEC-cap-qtot}{29.0}
\resdef{n15-s015-stop-E05-ALEC-cap-errearly}{2.1}
\resdef{n15-s015-stop-E05-ALEC-cap-errci}{[1.3, 3.2]}
\resdef{n15-s015-stop-E05-ALEC-cap-regret}{0.46}
\resdef{n15-s015-stop-E05-ALEC-cap-qround}{14}
\resdef{n15-s015-stop-E05-ALEC-cap-fixedhit}{80.6}
\resdef{n15-s015-stop-E05-ALEC-cap-fixedreg}{0.88}
\resdef{n15-s015-stop-E05-ALEC-cap-save2n}{53}
\resdef{n15-s015-stop-E05-ALEC-cap-s06-q}{20.2}
\resdef{n15-s015-stop-E05-ALEC-cap-s06-hit}{86.0}
\resdef{n15-s015-stop-E05-ALEC-cap-s06-agree}{98.4}
\resdef{n15-s015-stop-E05-ALEC-cap-s06-early}{70.4}
\resdef{n15-s015-stop-E05-ALEC-cap-s06-exceed}{3.0}
\resdef{n15-s015-stop-E05-ALEC-cap-s06-qmax}{30}
\resdef{n15-s015-stop-E05-ALEC-cap-s06-qtot}{35.2}
\resdef{n15-s015-stop-E05-ALEC-cap-s06-errearly}{0.1}
\resdef{n15-s015-stop-E05-ALEC-cap-s06-errci}{[0.0, 0.8]}
\resdef{n15-s015-stop-E05-ALEC-cap-s06-regret}{0.42}
\resdef{n15-s015-stop-E05-ALEC-cap-s06-qround}{20}
\resdef{n15-s015-stop-E05-ALEC-cap-s06-fixedhit}{84.3}
\resdef{n15-s015-stop-E05-ALEC-cap-s06-fixedreg}{0.56}
\resdef{n15-s015-stop-E05-ALEC-cap-s06-save2n}{33}
\resdef{n15-s030-full-hit}{76.1}
\resdef{n15-s030-tot-h}{23}
\resdef{n15-s030-frac-h}{22}
\resdef{n15-s030-RANDOM-hit-h}{53.3}
\resdef{n15-s030-RANDOM-hitllsm-h}{53.3}
\resdef{n15-s030-RANDOM-tau-h}{0.792}
\resdef{n15-s030-RANDOM-reg-h}{6.45}
\resdef{n15-s030-RANDOM-agree-h}{56.8}
\resdef{n15-s030-p-RANDOM-h}{$p<0.001$}
\resdef{n15-s030-p-RANDOM-h-llsm}{$p<0.001$}
\resdef{n15-s030-HARKER-PG-hit-h}{53.9}
\resdef{n15-s030-HARKER-PG-hitllsm-h}{55.9}
\resdef{n15-s030-HARKER-PG-tau-h}{0.770}
\resdef{n15-s030-HARKER-PG-reg-h}{6.50}
\resdef{n15-s030-HARKER-PG-agree-h}{61.0}
\resdef{n15-s030-p-HARKER-PG-h}{$p<0.001$}
\resdef{n15-s030-p-HARKER-PG-h-llsm}{$p<0.001$}
\resdef{n15-s030-HARKER-hit-h}{57.9}
\resdef{n15-s030-HARKER-hitllsm-h}{57.4}
\resdef{n15-s030-HARKER-tau-h}{0.775}
\resdef{n15-s030-HARKER-reg-h}{4.97}
\resdef{n15-s030-HARKER-agree-h}{62.8}
\resdef{n15-s030-p-HARKER-h}{$p<0.001$}
\resdef{n15-s030-p-HARKER-h-llsm}{$p<0.001$}
\resdef{n15-s030-HARKER-LOG-hit-h}{64.3}
\resdef{n15-s030-HARKER-LOG-hitllsm-h}{65.1}
\resdef{n15-s030-HARKER-LOG-tau-h}{0.749}
\resdef{n15-s030-HARKER-LOG-reg-h}{3.58}
\resdef{n15-s030-HARKER-LOG-agree-h}{69.9}
\resdef{n15-s030-p-HARKER-LOG-h}{$p=1.000$}
\resdef{n15-s030-p-HARKER-LOG-h-llsm}{$p=0.666$}
\resdef{n15-s030-AOPT-hit-h}{56.2}
\resdef{n15-s030-AOPT-hitllsm-h}{56.2}
\resdef{n15-s030-AOPT-tau-h}{0.814}
\resdef{n15-s030-AOPT-reg-h}{5.33}
\resdef{n15-s030-AOPT-agree-h}{61.6}
\resdef{n15-s030-p-AOPT-h}{$p<0.001$}
\resdef{n15-s030-p-AOPT-h-llsm}{$p<0.001$}
\resdef{n15-s030-FROZEN-hit-h}{59.3}
\resdef{n15-s030-FROZEN-hitllsm-h}{59.3}
\resdef{n15-s030-FROZEN-tau-h}{0.785}
\resdef{n15-s030-FROZEN-reg-h}{4.72}
\resdef{n15-s030-FROZEN-agree-h}{64.6}
\resdef{n15-s030-p-FROZEN-h}{$p<0.001$}
\resdef{n15-s030-p-FROZEN-h-llsm}{$p<0.001$}
\resdef{n15-s030-STAR-hit-h}{59.3}
\resdef{n15-s030-STAR-hitllsm-h}{59.3}
\resdef{n15-s030-STAR-tau-h}{0.785}
\resdef{n15-s030-STAR-reg-h}{4.72}
\resdef{n15-s030-STAR-agree-h}{64.6}
\resdef{n15-s030-p-STAR-h}{$p<0.001$}
\resdef{n15-s030-p-STAR-h-llsm}{$p<0.001$}
\resdef{n15-s030-ALEC-hit-h}{64.4}
\resdef{n15-s030-ALEC-hitllsm-h}{64.4}
\resdef{n15-s030-ALEC-tau-h}{0.783}
\resdef{n15-s030-ALEC-reg-h}{3.40}
\resdef{n15-s030-ALEC-agree-h}{70.7}
\resdef{n15-s030-ALEC-U-hit-h}{63.7}
\resdef{n15-s030-ALEC-U-hitllsm-h}{63.7}
\resdef{n15-s030-ALEC-U-tau-h}{0.780}
\resdef{n15-s030-ALEC-U-reg-h}{3.51}
\resdef{n15-s030-ALEC-U-agree-h}{70.6}
\resdef{n15-s030-p-ALEC-U-h}{$p=0.648$}
\resdef{n15-s030-p-ALEC-U-h-llsm}{$p=0.648$}
\resdef{n15-s030-tot-n}{30}
\resdef{n15-s030-frac-n}{29}
\resdef{n15-s030-RANDOM-hit-n}{61.4}
\resdef{n15-s030-RANDOM-hitllsm-n}{61.4}
\resdef{n15-s030-RANDOM-tau-n}{0.831}
\resdef{n15-s030-RANDOM-reg-n}{4.81}
\resdef{n15-s030-RANDOM-agree-n}{64.3}
\resdef{n15-s030-p-RANDOM-n}{$p<0.001$}
\resdef{n15-s030-p-RANDOM-n-llsm}{$p<0.001$}
\resdef{n15-s030-HARKER-PG-hit-n}{61.9}
\resdef{n15-s030-HARKER-PG-hitllsm-n}{63.4}
\resdef{n15-s030-HARKER-PG-tau-n}{0.808}
\resdef{n15-s030-HARKER-PG-reg-n}{4.12}
\resdef{n15-s030-HARKER-PG-agree-n}{70.1}
\resdef{n15-s030-p-HARKER-PG-n}{$p<0.001$}
\resdef{n15-s030-p-HARKER-PG-n-llsm}{$p=0.005$}
\resdef{n15-s030-HARKER-hit-n}{62.1}
\resdef{n15-s030-HARKER-hitllsm-n}{63.4}
\resdef{n15-s030-HARKER-tau-n}{0.817}
\resdef{n15-s030-HARKER-reg-n}{3.75}
\resdef{n15-s030-HARKER-agree-n}{70.1}
\resdef{n15-s030-p-HARKER-n}{$p<0.001$}
\resdef{n15-s030-p-HARKER-n-llsm}{$p=0.005$}
\resdef{n15-s030-HARKER-LOG-hit-n}{70.7}
\resdef{n15-s030-HARKER-LOG-hitllsm-n}{71.7}
\resdef{n15-s030-HARKER-LOG-tau-n}{0.792}
\resdef{n15-s030-HARKER-LOG-reg-n}{2.36}
\resdef{n15-s030-HARKER-LOG-agree-n}{77.7}
\resdef{n15-s030-p-HARKER-LOG-n}{$p=0.061$}
\resdef{n15-s030-p-HARKER-LOG-n-llsm}{$p=0.009$}
\resdef{n15-s030-AOPT-hit-n}{60.6}
\resdef{n15-s030-AOPT-hitllsm-n}{60.6}
\resdef{n15-s030-AOPT-tau-n}{0.845}
\resdef{n15-s030-AOPT-reg-n}{4.02}
\resdef{n15-s030-AOPT-agree-n}{65.7}
\resdef{n15-s030-p-AOPT-n}{$p<0.001$}
\resdef{n15-s030-p-AOPT-n-llsm}{$p<0.001$}
\resdef{n15-s030-FROZEN-hit-n}{66.5}
\resdef{n15-s030-FROZEN-hitllsm-n}{66.5}
\resdef{n15-s030-FROZEN-tau-n}{0.824}
\resdef{n15-s030-FROZEN-reg-n}{3.25}
\resdef{n15-s030-FROZEN-agree-n}{73.1}
\resdef{n15-s030-p-FROZEN-n}{$p=0.229$}
\resdef{n15-s030-p-FROZEN-n-llsm}{$p=0.229$}
\resdef{n15-s030-STAR-hit-n}{69.6}
\resdef{n15-s030-STAR-hitllsm-n}{69.6}
\resdef{n15-s030-STAR-tau-n}{0.823}
\resdef{n15-s030-STAR-reg-n}{2.69}
\resdef{n15-s030-STAR-agree-n}{74.2}
\resdef{n15-s030-p-STAR-n}{$p=0.250$}
\resdef{n15-s030-p-STAR-n-llsm}{$p=0.250$}
\resdef{n15-s030-ALEC-hit-n}{68.2}
\resdef{n15-s030-ALEC-hitllsm-n}{68.2}
\resdef{n15-s030-ALEC-tau-n}{0.822}
\resdef{n15-s030-ALEC-reg-n}{2.65}
\resdef{n15-s030-ALEC-agree-n}{78.0}
\resdef{n15-s030-ALEC-U-hit-n}{69.7}
\resdef{n15-s030-ALEC-U-hitllsm-n}{69.7}
\resdef{n15-s030-ALEC-U-tau-n}{0.813}
\resdef{n15-s030-ALEC-U-reg-n}{2.48}
\resdef{n15-s030-ALEC-U-agree-n}{78.6}
\resdef{n15-s030-p-ALEC-U-n}{$p=0.255$}
\resdef{n15-s030-p-ALEC-U-n-llsm}{$p=0.255$}
\resdef{n15-s030-tot-2n}{45}
\resdef{n15-s030-frac-2n}{43}
\resdef{n15-s030-RANDOM-hit-2n}{66.9}
\resdef{n15-s030-RANDOM-hitllsm-2n}{66.9}
\resdef{n15-s030-RANDOM-tau-2n}{0.871}
\resdef{n15-s030-RANDOM-reg-2n}{3.08}
\resdef{n15-s030-RANDOM-agree-2n}{73.0}
\resdef{n15-s030-p-RANDOM-2n}{$p<0.001$}
\resdef{n15-s030-p-RANDOM-2n-llsm}{$p<0.001$}
\resdef{n15-s030-HARKER-PG-hit-2n}{--}
\resdef{n15-s030-HARKER-PG-hitllsm-2n}{--}
\resdef{n15-s030-HARKER-PG-tau-2n}{--}
\resdef{n15-s030-HARKER-PG-reg-2n}{--}
\resdef{n15-s030-HARKER-PG-agree-2n}{--}
\resdef{n15-s030-HARKER-hit-2n}{69.4}
\resdef{n15-s030-HARKER-hitllsm-2n}{70.4}
\resdef{n15-s030-HARKER-tau-2n}{0.861}
\resdef{n15-s030-HARKER-reg-2n}{2.48}
\resdef{n15-s030-HARKER-agree-2n}{78.1}
\resdef{n15-s030-p-HARKER-2n}{$p<0.001$}
\resdef{n15-s030-p-HARKER-2n-llsm}{$p<0.001$}
\resdef{n15-s030-HARKER-LOG-hit-2n}{73.7}
\resdef{n15-s030-HARKER-LOG-hitllsm-2n}{74.3}
\resdef{n15-s030-HARKER-LOG-tau-2n}{0.846}
\resdef{n15-s030-HARKER-LOG-reg-2n}{1.88}
\resdef{n15-s030-HARKER-LOG-agree-2n}{85.8}
\resdef{n15-s030-p-HARKER-LOG-2n}{$p=0.207$}
\resdef{n15-s030-p-HARKER-LOG-2n-llsm}{$p=0.452$}
\resdef{n15-s030-AOPT-hit-2n}{66.1}
\resdef{n15-s030-AOPT-hitllsm-2n}{66.1}
\resdef{n15-s030-AOPT-tau-2n}{0.878}
\resdef{n15-s030-AOPT-reg-2n}{3.02}
\resdef{n15-s030-AOPT-agree-2n}{73.5}
\resdef{n15-s030-p-AOPT-2n}{$p<0.001$}
\resdef{n15-s030-p-AOPT-2n-llsm}{$p<0.001$}
\resdef{n15-s030-FROZEN-hit-2n}{72.0}
\resdef{n15-s030-FROZEN-hitllsm-2n}{72.0}
\resdef{n15-s030-FROZEN-tau-2n}{0.863}
\resdef{n15-s030-FROZEN-reg-2n}{2.14}
\resdef{n15-s030-FROZEN-agree-2n}{86.3}
\resdef{n15-s030-p-FROZEN-2n}{$p=0.002$}
\resdef{n15-s030-p-FROZEN-2n-llsm}{$p=0.002$}
\resdef{n15-s030-STAR-hit-2n}{75.0}
\resdef{n15-s030-STAR-hitllsm-2n}{75.0}
\resdef{n15-s030-STAR-tau-2n}{0.862}
\resdef{n15-s030-STAR-reg-2n}{1.71}
\resdef{n15-s030-STAR-agree-2n}{94.2}
\resdef{n15-s030-p-STAR-2n}{$p=0.868$}
\resdef{n15-s030-p-STAR-2n-llsm}{$p=0.868$}
\resdef{n15-s030-ALEC-hit-2n}{75.2}
\resdef{n15-s030-ALEC-hitllsm-2n}{75.2}
\resdef{n15-s030-ALEC-tau-2n}{0.862}
\resdef{n15-s030-ALEC-reg-2n}{1.69}
\resdef{n15-s030-ALEC-agree-2n}{96.4}
\resdef{n15-s030-ALEC-U-hit-2n}{74.7}
\resdef{n15-s030-ALEC-U-hitllsm-2n}{74.7}
\resdef{n15-s030-ALEC-U-tau-2n}{0.857}
\resdef{n15-s030-ALEC-U-reg-2n}{1.74}
\resdef{n15-s030-ALEC-U-agree-2n}{88.5}
\resdef{n15-s030-p-ALEC-U-2n}{$p=0.679$}
\resdef{n15-s030-p-ALEC-U-2n-llsm}{$p=0.679$}
\resdef{n15-s030-stop-LC-q}{14.0}
\resdef{n15-s030-stop-LC-hit}{69.9}
\resdef{n15-s030-stop-LC-agree}{78.2}
\resdef{n15-s030-stop-LC-early}{100.0}
\resdef{n15-s030-stop-LC-exceed}{30.1}
\resdef{n15-s030-stop-LC-qmax}{30}
\resdef{n15-s030-stop-LC-qtot}{29.0}
\resdef{n15-s030-stop-LC-errearly}{30.1}
\resdef{n15-s030-stop-LC-errci}{[27.3, 33.0]}
\resdef{n15-s030-stop-LC-regret}{2.42}
\resdef{n15-s030-stop-LC-qround}{14}
\resdef{n15-s030-stop-LC-fixedhit}{68.3}
\resdef{n15-s030-stop-LC-fixedreg}{2.64}
\resdef{n15-s030-stop-LC-save2n}{53}
\resdef{n15-s030-stop-B10-ALEC-q}{55.6}
\resdef{n15-s030-stop-B10-ALEC-hit}{75.5}
\resdef{n15-s030-stop-B10-ALEC-agree}{98.9}
\resdef{n15-s030-stop-B10-ALEC-early}{46.5}
\resdef{n15-s030-stop-B10-ALEC-exceed}{24.5}
\resdef{n15-s030-stop-B10-ALEC-qmax}{90}
\resdef{n15-s030-stop-B10-ALEC-qtot}{70.6}
\resdef{n15-s030-stop-B10-ALEC-errearly}{9.9}
\resdef{n15-s030-stop-B10-ALEC-errci}{[7.5, 12.9]}
\resdef{n15-s030-stop-B10-ALEC-regret}{1.64}
\resdef{n15-s030-stop-B10-ALEC-save2n}{-85}
\resdef{n15-s030-stop-B05-ALECU-q}{60.4}
\resdef{n15-s030-stop-B05-ALECU-hit}{76.2}
\resdef{n15-s030-stop-B05-ALECU-agree}{99.7}
\resdef{n15-s030-stop-B05-ALECU-early}{39.6}
\resdef{n15-s030-stop-B05-ALECU-exceed}{23.8}
\resdef{n15-s030-stop-B05-ALECU-qmax}{90}
\resdef{n15-s030-stop-B05-ALECU-qtot}{75.4}
\resdef{n15-s030-stop-B05-ALECU-errearly}{6.3}
\resdef{n15-s030-stop-B05-ALECU-errci}{[4.3, 9.2]}
\resdef{n15-s030-stop-B05-ALECU-regret}{1.56}
\resdef{n15-s030-stop-B05-ALECU-save2n}{-101}
\resdef{n15-s030-stop-B10-ALECU-q}{49.2}
\resdef{n15-s030-stop-B10-ALECU-hit}{75.7}
\resdef{n15-s030-stop-B10-ALECU-agree}{98.6}
\resdef{n15-s030-stop-B10-ALECU-early}{52.6}
\resdef{n15-s030-stop-B10-ALECU-exceed}{24.3}
\resdef{n15-s030-stop-B10-ALECU-qmax}{90}
\resdef{n15-s030-stop-B10-ALECU-qtot}{64.2}
\resdef{n15-s030-stop-B10-ALECU-errearly}{10.6}
\resdef{n15-s030-stop-B10-ALECU-errci}{[8.3, 13.6]}
\resdef{n15-s030-stop-B10-ALECU-regret}{1.58}
\resdef{n15-s030-stop-B10-ALECU-save2n}{-64}
\resdef{n15-s030-stop-B20-ALECU-q}{32.5}
\resdef{n15-s030-stop-B20-ALECU-hit}{72.8}
\resdef{n15-s030-stop-B20-ALECU-agree}{93.0}
\resdef{n15-s030-stop-B20-ALECU-early}{71.8}
\resdef{n15-s030-stop-B20-ALECU-exceed}{27.2}
\resdef{n15-s030-stop-B20-ALECU-qmax}{90}
\resdef{n15-s030-stop-B20-ALECU-qtot}{47.5}
\resdef{n15-s030-stop-B20-ALECU-errearly}{19.6}
\resdef{n15-s030-stop-B20-ALECU-errci}{[16.9, 22.7]}
\resdef{n15-s030-stop-B20-ALECU-regret}{1.97}
\resdef{n15-s030-stop-B20-ALECU-save2n}{-8}
\resdef{n15-s030-stop-B10-ALECU-s06-q}{59.4}
\resdef{n15-s030-stop-B10-ALECU-s06-hit}{76.2}
\resdef{n15-s030-stop-B10-ALECU-s06-agree}{99.9}
\resdef{n15-s030-stop-B10-ALECU-s06-early}{42.3}
\resdef{n15-s030-stop-B10-ALECU-s06-exceed}{23.8}
\resdef{n15-s030-stop-B10-ALECU-s06-qmax}{90}
\resdef{n15-s030-stop-B10-ALECU-s06-qtot}{74.4}
\resdef{n15-s030-stop-B10-ALECU-s06-errearly}{5.9}
\resdef{n15-s030-stop-B10-ALECU-s06-errci}{[4.0, 8.6]}
\resdef{n15-s030-stop-B10-ALECU-s06-regret}{1.54}
\resdef{n15-s030-stop-B10-ALECU-s06-save2n}{-98}
\resdef{n15-s030-stop-B10-ALECU-cap-q}{20.5}
\resdef{n15-s030-stop-B10-ALECU-cap-hit}{74.2}
\resdef{n15-s030-stop-B10-ALECU-cap-agree}{87.6}
\resdef{n15-s030-stop-B10-ALECU-cap-early}{49.1}
\resdef{n15-s030-stop-B10-ALECU-cap-exceed}{25.8}
\resdef{n15-s030-stop-B10-ALECU-cap-qmax}{30}
\resdef{n15-s030-stop-B10-ALECU-cap-qtot}{35.5}
\resdef{n15-s030-stop-B10-ALECU-cap-errearly}{10.4}
\resdef{n15-s030-stop-B10-ALECU-cap-errci}{[8.0, 13.4]}
\resdef{n15-s030-stop-B10-ALECU-cap-regret}{1.77}
\resdef{n15-s030-stop-B10-ALECU-cap-qround}{20}
\resdef{n15-s030-stop-B10-ALECU-cap-fixedhit}{73.4}
\resdef{n15-s030-stop-B10-ALECU-cap-fixedreg}{1.93}
\resdef{n15-s030-stop-B10-ALECU-cap-save2n}{32}
\resdef{n15-s030-stop-B10-ALECU-oracle-q}{47.0}
\resdef{n15-s030-stop-B10-ALECU-oracle-hit}{75.1}
\resdef{n15-s030-stop-B10-ALECU-oracle-agree}{96.5}
\resdef{n15-s030-stop-B10-ALECU-oracle-early}{53.4}
\resdef{n15-s030-stop-B10-ALECU-oracle-exceed}{24.9}
\resdef{n15-s030-stop-B10-ALECU-oracle-qmax}{90}
\resdef{n15-s030-stop-B10-ALECU-oracle-qtot}{62.0}
\resdef{n15-s030-stop-B10-ALECU-oracle-errearly}{13.7}
\resdef{n15-s030-stop-B10-ALECU-oracle-errci}{[11.0, 16.8]}
\resdef{n15-s030-stop-B10-ALECU-oracle-regret}{1.68}
\resdef{n15-s030-stop-B10-ALECU-oracle-save2n}{-57}
\resdef{n15-s030-stop-E02-ALEC-q}{50.6}
\resdef{n15-s030-stop-E02-ALEC-hit}{75.6}
\resdef{n15-s030-stop-E02-ALEC-agree}{98.6}
\resdef{n15-s030-stop-E02-ALEC-early}{53.0}
\resdef{n15-s030-stop-E02-ALEC-exceed}{19.6}
\resdef{n15-s030-stop-E02-ALEC-qmax}{90}
\resdef{n15-s030-stop-E02-ALEC-qtot}{65.6}
\resdef{n15-s030-stop-E02-ALEC-errearly}{8.1}
\resdef{n15-s030-stop-E02-ALEC-errci}{[6.1, 10.8]}
\resdef{n15-s030-stop-E02-ALEC-regret}{1.64}
\resdef{n15-s030-stop-E02-ALEC-save2n}{-69}
\resdef{n15-s030-stop-E05-ALEC-q}{41.4}
\resdef{n15-s030-stop-E05-ALEC-hit}{75.2}
\resdef{n15-s030-stop-E05-ALEC-agree}{97.1}
\resdef{n15-s030-stop-E05-ALEC-early}{65.0}
\resdef{n15-s030-stop-E05-ALEC-exceed}{13.7}
\resdef{n15-s030-stop-E05-ALEC-qmax}{90}
\resdef{n15-s030-stop-E05-ALEC-qtot}{56.4}
\resdef{n15-s030-stop-E05-ALEC-errearly}{7.2}
\resdef{n15-s030-stop-E05-ALEC-errci}{[5.5, 9.5]}
\resdef{n15-s030-stop-E05-ALEC-regret}{1.71}
\resdef{n15-s030-stop-E05-ALEC-save2n}{-38}
\resdef{n15-s030-stop-B10-ALEC-cap-q}{22.9}
\resdef{n15-s030-stop-B10-ALEC-cap-hit}{74.5}
\resdef{n15-s030-stop-B10-ALEC-cap-agree}{95.6}
\resdef{n15-s030-stop-B10-ALEC-cap-early}{42.8}
\resdef{n15-s030-stop-B10-ALEC-cap-exceed}{25.5}
\resdef{n15-s030-stop-B10-ALEC-cap-qmax}{30}
\resdef{n15-s030-stop-B10-ALEC-cap-qtot}{37.9}
\resdef{n15-s030-stop-B10-ALEC-cap-errearly}{9.8}
\resdef{n15-s030-stop-B10-ALEC-cap-errci}{[7.3, 13.0]}
\resdef{n15-s030-stop-B10-ALEC-cap-regret}{1.76}
\resdef{n15-s030-stop-B10-ALEC-cap-qround}{23}
\resdef{n15-s030-stop-B10-ALEC-cap-fixedhit}{75.2}
\resdef{n15-s030-stop-B10-ALEC-cap-fixedreg}{1.65}
\resdef{n15-s030-stop-B10-ALEC-cap-save2n}{24}
\resdef{n15-s030-stop-E02-ALEC-cap-q}{21.9}
\resdef{n15-s030-stop-E02-ALEC-cap-hit}{74.6}
\resdef{n15-s030-stop-E02-ALEC-cap-agree}{95.3}
\resdef{n15-s030-stop-E02-ALEC-cap-early}{49.0}
\resdef{n15-s030-stop-E02-ALEC-cap-exceed}{20.6}
\resdef{n15-s030-stop-E02-ALEC-cap-qmax}{30}
\resdef{n15-s030-stop-E02-ALEC-cap-qtot}{36.9}
\resdef{n15-s030-stop-E02-ALEC-cap-errearly}{7.8}
\resdef{n15-s030-stop-E02-ALEC-cap-errci}{[5.7, 10.5]}
\resdef{n15-s030-stop-E02-ALEC-cap-regret}{1.77}
\resdef{n15-s030-stop-E02-ALEC-cap-qround}{22}
\resdef{n15-s030-stop-E02-ALEC-cap-fixedhit}{74.4}
\resdef{n15-s030-stop-E02-ALEC-cap-fixedreg}{1.74}
\resdef{n15-s030-stop-E02-ALEC-cap-save2n}{27}
\resdef{n15-s030-stop-E05-ALEC-cap-q}{19.8}
\resdef{n15-s030-stop-E05-ALEC-cap-hit}{74.2}
\resdef{n15-s030-stop-E05-ALEC-cap-agree}{93.8}
\resdef{n15-s030-stop-E05-ALEC-cap-early}{59.4}
\resdef{n15-s030-stop-E05-ALEC-cap-exceed}{14.7}
\resdef{n15-s030-stop-E05-ALEC-cap-qmax}{30}
\resdef{n15-s030-stop-E05-ALEC-cap-qtot}{34.8}
\resdef{n15-s030-stop-E05-ALEC-cap-errearly}{7.2}
\resdef{n15-s030-stop-E05-ALEC-cap-errci}{[5.4, 9.6]}
\resdef{n15-s030-stop-E05-ALEC-cap-regret}{1.84}
\resdef{n15-s030-stop-E05-ALEC-cap-qround}{20}
\resdef{n15-s030-stop-E05-ALEC-cap-fixedhit}{73.4}
\resdef{n15-s030-stop-E05-ALEC-cap-fixedreg}{1.93}
\resdef{n15-s030-stop-E05-ALEC-cap-save2n}{34}
\resdef{n15-s030-stop-E05-ALEC-cap-s06-q}{23.7}
\resdef{n15-s030-stop-E05-ALEC-cap-s06-hit}{74.9}
\resdef{n15-s030-stop-E05-ALEC-cap-s06-agree}{96.1}
\resdef{n15-s030-stop-E05-ALEC-cap-s06-early}{48.8}
\resdef{n15-s030-stop-E05-ALEC-cap-s06-exceed}{14.2}
\resdef{n15-s030-stop-E05-ALEC-cap-s06-qmax}{30}
\resdef{n15-s030-stop-E05-ALEC-cap-s06-qtot}{38.7}
\resdef{n15-s030-stop-E05-ALEC-cap-s06-errearly}{3.5}
\resdef{n15-s030-stop-E05-ALEC-cap-s06-errci}{[2.2, 5.5]}
\resdef{n15-s030-stop-E05-ALEC-cap-s06-regret}{1.71}
\resdef{n15-s030-stop-E05-ALEC-cap-s06-qround}{24}
\resdef{n15-s030-stop-E05-ALEC-cap-s06-fixedhit}{75.1}
\resdef{n15-s030-stop-E05-ALEC-cap-s06-fixedreg}{1.65}
\resdef{n15-s030-stop-E05-ALEC-cap-s06-save2n}{21}
\resdef{n15-s050-full-hit}{67.5}
\resdef{n15-s050-tot-h}{23}
\resdef{n15-s050-frac-h}{22}
\resdef{n15-s050-RANDOM-hit-h}{40.6}
\resdef{n15-s050-RANDOM-hitllsm-h}{40.6}
\resdef{n15-s050-RANDOM-tau-h}{0.688}
\resdef{n15-s050-RANDOM-reg-h}{11.70}
\resdef{n15-s050-RANDOM-agree-h}{46.5}
\resdef{n15-s050-p-RANDOM-h}{$p<0.001$}
\resdef{n15-s050-p-RANDOM-h-llsm}{$p<0.001$}
\resdef{n15-s050-HARKER-PG-hit-h}{40.0}
\resdef{n15-s050-HARKER-PG-hitllsm-h}{43.6}
\resdef{n15-s050-HARKER-PG-tau-h}{0.647}
\resdef{n15-s050-HARKER-PG-reg-h}{12.63}
\resdef{n15-s050-HARKER-PG-agree-h}{45.3}
\resdef{n15-s050-p-HARKER-PG-h}{$p<0.001$}
\resdef{n15-s050-p-HARKER-PG-h-llsm}{$p<0.001$}
\resdef{n15-s050-HARKER-hit-h}{45.1}
\resdef{n15-s050-HARKER-hitllsm-h}{45.2}
\resdef{n15-s050-HARKER-tau-h}{0.649}
\resdef{n15-s050-HARKER-reg-h}{9.96}
\resdef{n15-s050-HARKER-agree-h}{51.1}
\resdef{n15-s050-p-HARKER-h}{$p<0.001$}
\resdef{n15-s050-p-HARKER-h-llsm}{$p<0.001$}
\resdef{n15-s050-HARKER-LOG-hit-h}{50.9}
\resdef{n15-s050-HARKER-LOG-hitllsm-h}{52.2}
\resdef{n15-s050-HARKER-LOG-tau-h}{0.587}
\resdef{n15-s050-HARKER-LOG-reg-h}{7.57}
\resdef{n15-s050-HARKER-LOG-agree-h}{62.1}
\resdef{n15-s050-p-HARKER-LOG-h}{$p=0.566$}
\resdef{n15-s050-p-HARKER-LOG-h-llsm}{$p=0.894$}
\resdef{n15-s050-AOPT-hit-h}{43.0}
\resdef{n15-s050-AOPT-hitllsm-h}{43.0}
\resdef{n15-s050-AOPT-tau-h}{0.711}
\resdef{n15-s050-AOPT-reg-h}{10.82}
\resdef{n15-s050-AOPT-agree-h}{48.9}
\resdef{n15-s050-p-AOPT-h}{$p<0.001$}
\resdef{n15-s050-p-AOPT-h-llsm}{$p<0.001$}
\resdef{n15-s050-FROZEN-hit-h}{47.6}
\resdef{n15-s050-FROZEN-hitllsm-h}{47.6}
\resdef{n15-s050-FROZEN-tau-h}{0.666}
\resdef{n15-s050-FROZEN-reg-h}{9.21}
\resdef{n15-s050-FROZEN-agree-h}{54.3}
\resdef{n15-s050-p-FROZEN-h}{$p=0.005$}
\resdef{n15-s050-p-FROZEN-h-llsm}{$p=0.005$}
\resdef{n15-s050-STAR-hit-h}{47.6}
\resdef{n15-s050-STAR-hitllsm-h}{47.6}
\resdef{n15-s050-STAR-tau-h}{0.666}
\resdef{n15-s050-STAR-reg-h}{9.21}
\resdef{n15-s050-STAR-agree-h}{54.3}
\resdef{n15-s050-p-STAR-h}{$p=0.005$}
\resdef{n15-s050-p-STAR-h-llsm}{$p=0.005$}
\resdef{n15-s050-ALEC-hit-h}{51.9}
\resdef{n15-s050-ALEC-hitllsm-h}{51.9}
\resdef{n15-s050-ALEC-tau-h}{0.658}
\resdef{n15-s050-ALEC-reg-h}{7.01}
\resdef{n15-s050-ALEC-agree-h}{61.7}
\resdef{n15-s050-ALEC-U-hit-h}{51.6}
\resdef{n15-s050-ALEC-U-hitllsm-h}{51.6}
\resdef{n15-s050-ALEC-U-tau-h}{0.662}
\resdef{n15-s050-ALEC-U-reg-h}{6.82}
\resdef{n15-s050-ALEC-U-agree-h}{63.4}
\resdef{n15-s050-p-ALEC-U-h}{$p=0.890$}
\resdef{n15-s050-p-ALEC-U-h-llsm}{$p=0.890$}
\resdef{n15-s050-tot-n}{30}
\resdef{n15-s050-frac-n}{29}
\resdef{n15-s050-RANDOM-hit-n}{43.8}
\resdef{n15-s050-RANDOM-hitllsm-n}{43.8}
\resdef{n15-s050-RANDOM-tau-n}{0.739}
\resdef{n15-s050-RANDOM-reg-n}{10.04}
\resdef{n15-s050-RANDOM-agree-n}{53.4}
\resdef{n15-s050-p-RANDOM-n}{$p<0.001$}
\resdef{n15-s050-p-RANDOM-n-llsm}{$p<0.001$}
\resdef{n15-s050-HARKER-PG-hit-n}{49.6}
\resdef{n15-s050-HARKER-PG-hitllsm-n}{50.4}
\resdef{n15-s050-HARKER-PG-tau-n}{0.706}
\resdef{n15-s050-HARKER-PG-reg-n}{8.93}
\resdef{n15-s050-HARKER-PG-agree-n}{56.6}
\resdef{n15-s050-p-HARKER-PG-n}{$p<0.001$}
\resdef{n15-s050-p-HARKER-PG-n-llsm}{$p=0.004$}
\resdef{n15-s050-HARKER-hit-n}{51.4}
\resdef{n15-s050-HARKER-hitllsm-n}{51.1}
\resdef{n15-s050-HARKER-tau-n}{0.701}
\resdef{n15-s050-HARKER-reg-n}{7.86}
\resdef{n15-s050-HARKER-agree-n}{58.4}
\resdef{n15-s050-p-HARKER-n}{$p=0.019$}
\resdef{n15-s050-p-HARKER-n-llsm}{$p=0.013$}
\resdef{n15-s050-HARKER-LOG-hit-n}{53.4}
\resdef{n15-s050-HARKER-LOG-hitllsm-n}{55.8}
\resdef{n15-s050-HARKER-LOG-tau-n}{0.628}
\resdef{n15-s050-HARKER-LOG-reg-n}{6.16}
\resdef{n15-s050-HARKER-LOG-agree-n}{68.9}
\resdef{n15-s050-p-HARKER-LOG-n}{$p=0.109$}
\resdef{n15-s050-p-HARKER-LOG-n-llsm}{$p=0.948$}
\resdef{n15-s050-AOPT-hit-n}{47.9}
\resdef{n15-s050-AOPT-hitllsm-n}{47.9}
\resdef{n15-s050-AOPT-tau-n}{0.757}
\resdef{n15-s050-AOPT-reg-n}{8.42}
\resdef{n15-s050-AOPT-agree-n}{57.6}
\resdef{n15-s050-p-AOPT-n}{$p<0.001$}
\resdef{n15-s050-p-AOPT-n-llsm}{$p<0.001$}
\resdef{n15-s050-FROZEN-hit-n}{49.9}
\resdef{n15-s050-FROZEN-hitllsm-n}{49.9}
\resdef{n15-s050-FROZEN-tau-n}{0.726}
\resdef{n15-s050-FROZEN-reg-n}{7.83}
\resdef{n15-s050-FROZEN-agree-n}{60.9}
\resdef{n15-s050-p-FROZEN-n}{$p<0.001$}
\resdef{n15-s050-p-FROZEN-n-llsm}{$p<0.001$}
\resdef{n15-s050-STAR-hit-n}{55.0}
\resdef{n15-s050-STAR-hitllsm-n}{55.0}
\resdef{n15-s050-STAR-tau-n}{0.724}
\resdef{n15-s050-STAR-reg-n}{6.28}
\resdef{n15-s050-STAR-agree-n}{66.5}
\resdef{n15-s050-p-STAR-n}{$p=0.505$}
\resdef{n15-s050-p-STAR-n-llsm}{$p=0.505$}
\resdef{n15-s050-ALEC-hit-n}{56.0}
\resdef{n15-s050-ALEC-hitllsm-n}{56.0}
\resdef{n15-s050-ALEC-tau-n}{0.720}
\resdef{n15-s050-ALEC-reg-n}{5.71}
\resdef{n15-s050-ALEC-agree-n}{69.5}
\resdef{n15-s050-ALEC-U-hit-n}{57.0}
\resdef{n15-s050-ALEC-U-hitllsm-n}{57.0}
\resdef{n15-s050-ALEC-U-tau-n}{0.712}
\resdef{n15-s050-ALEC-U-reg-n}{5.21}
\resdef{n15-s050-ALEC-U-agree-n}{72.4}
\resdef{n15-s050-p-ALEC-U-n}{$p=0.548$}
\resdef{n15-s050-p-ALEC-U-n-llsm}{$p=0.548$}
\resdef{n15-s050-tot-2n}{45}
\resdef{n15-s050-frac-2n}{43}
\resdef{n15-s050-RANDOM-hit-2n}{52.5}
\resdef{n15-s050-RANDOM-hitllsm-2n}{52.5}
\resdef{n15-s050-RANDOM-tau-2n}{0.799}
\resdef{n15-s050-RANDOM-reg-2n}{7.26}
\resdef{n15-s050-RANDOM-agree-2n}{64.5}
\resdef{n15-s050-p-RANDOM-2n}{$p<0.001$}
\resdef{n15-s050-p-RANDOM-2n-llsm}{$p<0.001$}
\resdef{n15-s050-HARKER-PG-hit-2n}{--}
\resdef{n15-s050-HARKER-PG-hitllsm-2n}{--}
\resdef{n15-s050-HARKER-PG-tau-2n}{--}
\resdef{n15-s050-HARKER-PG-reg-2n}{--}
\resdef{n15-s050-HARKER-PG-agree-2n}{--}
\resdef{n15-s050-HARKER-hit-2n}{56.3}
\resdef{n15-s050-HARKER-hitllsm-2n}{57.1}
\resdef{n15-s050-HARKER-tau-2n}{0.760}
\resdef{n15-s050-HARKER-reg-2n}{5.55}
\resdef{n15-s050-HARKER-agree-2n}{67.3}
\resdef{n15-s050-p-HARKER-2n}{$p<0.001$}
\resdef{n15-s050-p-HARKER-2n-llsm}{$p<0.001$}
\resdef{n15-s050-HARKER-LOG-hit-2n}{59.5}
\resdef{n15-s050-HARKER-LOG-hitllsm-2n}{59.9}
\resdef{n15-s050-HARKER-LOG-tau-2n}{0.713}
\resdef{n15-s050-HARKER-LOG-reg-2n}{4.56}
\resdef{n15-s050-HARKER-LOG-agree-2n}{77.6}
\resdef{n15-s050-p-HARKER-LOG-2n}{$p<0.001$}
\resdef{n15-s050-p-HARKER-LOG-2n-llsm}{$p<0.001$}
\resdef{n15-s050-AOPT-hit-2n}{54.3}
\resdef{n15-s050-AOPT-hitllsm-2n}{54.3}
\resdef{n15-s050-AOPT-tau-2n}{0.805}
\resdef{n15-s050-AOPT-reg-2n}{6.09}
\resdef{n15-s050-AOPT-agree-2n}{66.1}
\resdef{n15-s050-p-AOPT-2n}{$p<0.001$}
\resdef{n15-s050-p-AOPT-2n-llsm}{$p<0.001$}
\resdef{n15-s050-FROZEN-hit-2n}{60.1}
\resdef{n15-s050-FROZEN-hitllsm-2n}{60.1}
\resdef{n15-s050-FROZEN-tau-2n}{0.785}
\resdef{n15-s050-FROZEN-reg-2n}{4.85}
\resdef{n15-s050-FROZEN-agree-2n}{76.0}
\resdef{n15-s050-p-FROZEN-2n}{$p<0.001$}
\resdef{n15-s050-p-FROZEN-2n-llsm}{$p<0.001$}
\resdef{n15-s050-STAR-hit-2n}{63.7}
\resdef{n15-s050-STAR-hitllsm-2n}{63.7}
\resdef{n15-s050-STAR-tau-2n}{0.785}
\resdef{n15-s050-STAR-reg-2n}{3.85}
\resdef{n15-s050-STAR-agree-2n}{87.5}
\resdef{n15-s050-p-STAR-2n}{$p=0.119$}
\resdef{n15-s050-p-STAR-2n-llsm}{$p=0.119$}
\resdef{n15-s050-ALEC-hit-2n}{65.2}
\resdef{n15-s050-ALEC-hitllsm-2n}{65.2}
\resdef{n15-s050-ALEC-tau-2n}{0.783}
\resdef{n15-s050-ALEC-reg-2n}{3.51}
\resdef{n15-s050-ALEC-agree-2n}{93.0}
\resdef{n15-s050-ALEC-U-hit-2n}{63.0}
\resdef{n15-s050-ALEC-U-hitllsm-2n}{63.0}
\resdef{n15-s050-ALEC-U-tau-2n}{0.780}
\resdef{n15-s050-ALEC-U-reg-2n}{3.84}
\resdef{n15-s050-ALEC-U-agree-2n}{82.5}
\resdef{n15-s050-p-ALEC-U-2n}{$p=0.073$}
\resdef{n15-s050-p-ALEC-U-2n-llsm}{$p=0.073$}
\resdef{n15-s050-stop-LC-q}{15.8}
\resdef{n15-s050-stop-LC-hit}{59.2}
\resdef{n15-s050-stop-LC-agree}{78.3}
\resdef{n15-s050-stop-LC-early}{100.0}
\resdef{n15-s050-stop-LC-exceed}{40.8}
\resdef{n15-s050-stop-LC-qmax}{36}
\resdef{n15-s050-stop-LC-qtot}{30.8}
\resdef{n15-s050-stop-LC-errearly}{40.8}
\resdef{n15-s050-stop-LC-errci}{[37.8, 43.9]}
\resdef{n15-s050-stop-LC-regret}{4.62}
\resdef{n15-s050-stop-LC-qround}{16}
\resdef{n15-s050-stop-LC-fixedhit}{56.7}
\resdef{n15-s050-stop-LC-fixedreg}{5.36}
\resdef{n15-s050-stop-LC-save2n}{47}
\resdef{n15-s050-stop-B10-ALEC-q}{63.9}
\resdef{n15-s050-stop-B10-ALEC-hit}{65.6}
\resdef{n15-s050-stop-B10-ALEC-agree}{96.7}
\resdef{n15-s050-stop-B10-ALEC-early}{36.0}
\resdef{n15-s050-stop-B10-ALEC-exceed}{34.4}
\resdef{n15-s050-stop-B10-ALEC-qmax}{90}
\resdef{n15-s050-stop-B10-ALEC-qtot}{78.9}
\resdef{n15-s050-stop-B10-ALEC-errearly}{22.8}
\resdef{n15-s050-stop-B10-ALEC-errci}{[18.7, 27.4]}
\resdef{n15-s050-stop-B10-ALEC-regret}{3.39}
\resdef{n15-s050-stop-B10-ALEC-save2n}{-113}
\resdef{n15-s050-stop-B05-ALECU-q}{70.2}
\resdef{n15-s050-stop-B05-ALECU-hit}{66.9}
\resdef{n15-s050-stop-B05-ALECU-agree}{98.3}
\resdef{n15-s050-stop-B05-ALECU-early}{26.7}
\resdef{n15-s050-stop-B05-ALECU-exceed}{33.1}
\resdef{n15-s050-stop-B05-ALECU-qmax}{90}
\resdef{n15-s050-stop-B05-ALECU-qtot}{85.2}
\resdef{n15-s050-stop-B05-ALECU-errearly}{15.0}
\resdef{n15-s050-stop-B05-ALECU-errci}{[11.2, 19.8]}
\resdef{n15-s050-stop-B05-ALECU-regret}{3.21}
\resdef{n15-s050-stop-B05-ALECU-save2n}{-134}
\resdef{n15-s050-stop-B10-ALECU-q}{58.2}
\resdef{n15-s050-stop-B10-ALECU-hit}{66.3}
\resdef{n15-s050-stop-B10-ALECU-agree}{96.5}
\resdef{n15-s050-stop-B10-ALECU-early}{41.5}
\resdef{n15-s050-stop-B10-ALECU-exceed}{33.7}
\resdef{n15-s050-stop-B10-ALECU-qmax}{90}
\resdef{n15-s050-stop-B10-ALECU-qtot}{73.2}
\resdef{n15-s050-stop-B10-ALECU-errearly}{21.4}
\resdef{n15-s050-stop-B10-ALECU-errci}{[17.8, 25.6]}
\resdef{n15-s050-stop-B10-ALECU-regret}{3.35}
\resdef{n15-s050-stop-B10-ALECU-save2n}{-94}
\resdef{n15-s050-stop-B20-ALECU-q}{39.8}
\resdef{n15-s050-stop-B20-ALECU-hit}{63.4}
\resdef{n15-s050-stop-B20-ALECU-agree}{89.2}
\resdef{n15-s050-stop-B20-ALECU-early}{63.6}
\resdef{n15-s050-stop-B20-ALECU-exceed}{36.6}
\resdef{n15-s050-stop-B20-ALECU-qmax}{90}
\resdef{n15-s050-stop-B20-ALECU-qtot}{54.8}
\resdef{n15-s050-stop-B20-ALECU-errearly}{30.7}
\resdef{n15-s050-stop-B20-ALECU-errci}{[27.2, 34.4]}
\resdef{n15-s050-stop-B20-ALECU-regret}{4.08}
\resdef{n15-s050-stop-B20-ALECU-save2n}{-33}
\resdef{n15-s050-stop-B10-ALECU-s06-q}{68.1}
\resdef{n15-s050-stop-B10-ALECU-s06-hit}{67.2}
\resdef{n15-s050-stop-B10-ALECU-s06-agree}{99.5}
\resdef{n15-s050-stop-B10-ALECU-s06-early}{30.9}
\resdef{n15-s050-stop-B10-ALECU-s06-exceed}{32.8}
\resdef{n15-s050-stop-B10-ALECU-s06-qmax}{90}
\resdef{n15-s050-stop-B10-ALECU-s06-qtot}{83.1}
\resdef{n15-s050-stop-B10-ALECU-s06-errearly}{12.3}
\resdef{n15-s050-stop-B10-ALECU-s06-errci}{[9.1, 16.4]}
\resdef{n15-s050-stop-B10-ALECU-s06-regret}{3.10}
\resdef{n15-s050-stop-B10-ALECU-s06-save2n}{-127}
\resdef{n15-s050-stop-B10-ALECU-cap-q}{22.7}
\resdef{n15-s050-stop-B10-ALECU-cap-hit}{62.4}
\resdef{n15-s050-stop-B10-ALECU-cap-agree}{81.4}
\resdef{n15-s050-stop-B10-ALECU-cap-early}{37.4}
\resdef{n15-s050-stop-B10-ALECU-cap-exceed}{37.6}
\resdef{n15-s050-stop-B10-ALECU-cap-qmax}{30}
\resdef{n15-s050-stop-B10-ALECU-cap-qtot}{37.7}
\resdef{n15-s050-stop-B10-ALECU-cap-errearly}{22.5}
\resdef{n15-s050-stop-B10-ALECU-cap-errci}{[18.5, 27.0]}
\resdef{n15-s050-stop-B10-ALECU-cap-regret}{4.02}
\resdef{n15-s050-stop-B10-ALECU-cap-qround}{23}
\resdef{n15-s050-stop-B10-ALECU-cap-fixedhit}{61.8}
\resdef{n15-s050-stop-B10-ALECU-cap-fixedreg}{4.16}
\resdef{n15-s050-stop-B10-ALECU-cap-save2n}{24}
\resdef{n15-s050-stop-B10-ALECU-oracle-q}{65.0}
\resdef{n15-s050-stop-B10-ALECU-oracle-hit}{67.2}
\resdef{n15-s050-stop-B10-ALECU-oracle-agree}{98.6}
\resdef{n15-s050-stop-B10-ALECU-oracle-early}{33.0}
\resdef{n15-s050-stop-B10-ALECU-oracle-exceed}{32.8}
\resdef{n15-s050-stop-B10-ALECU-oracle-qmax}{90}
\resdef{n15-s050-stop-B10-ALECU-oracle-qtot}{80.0}
\resdef{n15-s050-stop-B10-ALECU-oracle-errearly}{17.9}
\resdef{n15-s050-stop-B10-ALECU-oracle-errci}{[14.1, 22.4]}
\resdef{n15-s050-stop-B10-ALECU-oracle-regret}{3.21}
\resdef{n15-s050-stop-B10-ALECU-oracle-save2n}{-117}
\resdef{n15-s050-stop-E02-ALEC-q}{60.7}
\resdef{n15-s050-stop-E02-ALEC-hit}{65.7}
\resdef{n15-s050-stop-E02-ALEC-agree}{96.1}
\resdef{n15-s050-stop-E02-ALEC-early}{40.6}
\resdef{n15-s050-stop-E02-ALEC-exceed}{28.6}
\resdef{n15-s050-stop-E02-ALEC-qmax}{90}
\resdef{n15-s050-stop-E02-ALEC-qtot}{75.7}
\resdef{n15-s050-stop-E02-ALEC-errearly}{20.9}
\resdef{n15-s050-stop-E02-ALEC-errci}{[17.3, 25.2]}
\resdef{n15-s050-stop-E02-ALEC-regret}{3.46}
\resdef{n15-s050-stop-E02-ALEC-save2n}{-102}
\resdef{n15-s050-stop-E05-ALEC-q}{55.2}
\resdef{n15-s050-stop-E05-ALEC-hit}{65.3}
\resdef{n15-s050-stop-E05-ALEC-agree}{94.9}
\resdef{n15-s050-stop-E05-ALEC-early}{47.5}
\resdef{n15-s050-stop-E05-ALEC-exceed}{23.5}
\resdef{n15-s050-stop-E05-ALEC-qmax}{90}
\resdef{n15-s050-stop-E05-ALEC-qtot}{70.2}
\resdef{n15-s050-stop-E05-ALEC-errearly}{16.8}
\resdef{n15-s050-stop-E05-ALEC-errci}{[13.7, 20.5]}
\resdef{n15-s050-stop-E05-ALEC-regret}{3.55}
\resdef{n15-s050-stop-E05-ALEC-save2n}{-84}
\resdef{n15-s050-stop-B10-ALEC-cap-q}{24.7}
\resdef{n15-s050-stop-B10-ALEC-cap-hit}{63.1}
\resdef{n15-s050-stop-B10-ALEC-cap-agree}{90.3}
\resdef{n15-s050-stop-B10-ALEC-cap-early}{29.6}
\resdef{n15-s050-stop-B10-ALEC-cap-exceed}{36.9}
\resdef{n15-s050-stop-B10-ALEC-cap-qmax}{30}
\resdef{n15-s050-stop-B10-ALEC-cap-qtot}{39.7}
\resdef{n15-s050-stop-B10-ALEC-cap-errearly}{24.3}
\resdef{n15-s050-stop-B10-ALEC-cap-errci}{[19.8, 29.5]}
\resdef{n15-s050-stop-B10-ALEC-cap-regret}{3.86}
\resdef{n15-s050-stop-B10-ALEC-cap-qround}{25}
\resdef{n15-s050-stop-B10-ALEC-cap-fixedhit}{62.9}
\resdef{n15-s050-stop-B10-ALEC-cap-fixedreg}{3.91}
\resdef{n15-s050-stop-B10-ALEC-cap-save2n}{18}
\resdef{n15-s050-stop-E02-ALEC-cap-q}{24.0}
\resdef{n15-s050-stop-E02-ALEC-cap-hit}{63.2}
\resdef{n15-s050-stop-E02-ALEC-cap-agree}{89.7}
\resdef{n15-s050-stop-E02-ALEC-cap-early}{33.5}
\resdef{n15-s050-stop-E02-ALEC-cap-exceed}{31.0}
\resdef{n15-s050-stop-E02-ALEC-cap-qmax}{30}
\resdef{n15-s050-stop-E02-ALEC-cap-qtot}{39.0}
\resdef{n15-s050-stop-E02-ALEC-cap-errearly}{20.9}
\resdef{n15-s050-stop-E02-ALEC-cap-errci}{[16.9, 25.6]}
\resdef{n15-s050-stop-E02-ALEC-cap-regret}{3.93}
\resdef{n15-s050-stop-E02-ALEC-cap-qround}{24}
\resdef{n15-s050-stop-E02-ALEC-cap-fixedhit}{63.3}
\resdef{n15-s050-stop-E02-ALEC-cap-fixedreg}{3.83}
\resdef{n15-s050-stop-E02-ALEC-cap-save2n}{20}
\resdef{n15-s050-stop-E05-ALEC-cap-q}{22.8}
\resdef{n15-s050-stop-E05-ALEC-cap-hit}{63.0}
\resdef{n15-s050-stop-E05-ALEC-cap-agree}{88.8}
\resdef{n15-s050-stop-E05-ALEC-cap-early}{39.6}
\resdef{n15-s050-stop-E05-ALEC-cap-exceed}{25.4}
\resdef{n15-s050-stop-E05-ALEC-cap-qmax}{30}
\resdef{n15-s050-stop-E05-ALEC-cap-qtot}{37.8}
\resdef{n15-s050-stop-E05-ALEC-cap-errearly}{17.7}
\resdef{n15-s050-stop-E05-ALEC-cap-errci}{[14.2, 21.7]}
\resdef{n15-s050-stop-E05-ALEC-cap-regret}{3.99}
\resdef{n15-s050-stop-E05-ALEC-cap-qround}{23}
\resdef{n15-s050-stop-E05-ALEC-cap-fixedhit}{61.8}
\resdef{n15-s050-stop-E05-ALEC-cap-fixedreg}{4.16}
\resdef{n15-s050-stop-E05-ALEC-cap-save2n}{24}
\resdef{n15-s050-stop-E05-ALEC-cap-s06-q}{26.6}
\resdef{n15-s050-stop-E05-ALEC-cap-s06-hit}{64.9}
\resdef{n15-s050-stop-E05-ALEC-cap-s06-agree}{92.2}
\resdef{n15-s050-stop-E05-ALEC-cap-s06-early}{28.3}
\resdef{n15-s050-stop-E05-ALEC-cap-s06-exceed}{23.3}
\resdef{n15-s050-stop-E05-ALEC-cap-s06-qmax}{30}
\resdef{n15-s050-stop-E05-ALEC-cap-s06-qtot}{41.6}
\resdef{n15-s050-stop-E05-ALEC-cap-s06-errearly}{9.5}
\resdef{n15-s050-stop-E05-ALEC-cap-s06-errci}{[6.6, 13.5]}
\resdef{n15-s050-stop-E05-ALEC-cap-s06-regret}{3.57}
\resdef{n15-s050-stop-E05-ALEC-cap-s06-qround}{27}
\resdef{n15-s050-stop-E05-ALEC-cap-s06-fixedhit}{64.4}
\resdef{n15-s050-stop-E05-ALEC-cap-s06-fixedreg}{3.72}
\resdef{n15-s050-stop-E05-ALEC-cap-s06-save2n}{11}
\resdef{n11-s030-saaty-full-hit}{75.1}
\resdef{n11-s030-saaty-tot-h}{17}
\resdef{n11-s030-saaty-frac-h}{31}
\resdef{n11-s030-saaty-RANDOM-hit-h}{60.8}
\resdef{n11-s030-saaty-RANDOM-hitllsm-h}{60.8}
\resdef{n11-s030-saaty-RANDOM-tau-h}{0.786}
\resdef{n11-s030-saaty-RANDOM-reg-h}{5.78}
\resdef{n11-s030-saaty-RANDOM-agree-h}{64.9}
\resdef{n11-s030-saaty-p-RANDOM-h}{$p<0.001$}
\resdef{n11-s030-saaty-p-RANDOM-h-llsm}{$p<0.001$}
\resdef{n11-s030-saaty-HARKER-PG-hit-h}{60.4}
\resdef{n11-s030-saaty-HARKER-PG-hitllsm-h}{61.0}
\resdef{n11-s030-saaty-HARKER-PG-tau-h}{0.767}
\resdef{n11-s030-saaty-HARKER-PG-reg-h}{5.93}
\resdef{n11-s030-saaty-HARKER-PG-agree-h}{63.5}
\resdef{n11-s030-saaty-p-HARKER-PG-h}{$p<0.001$}
\resdef{n11-s030-saaty-p-HARKER-PG-h-llsm}{$p<0.001$}
\resdef{n11-s030-saaty-HARKER-hit-h}{62.3}
\resdef{n11-s030-saaty-HARKER-hitllsm-h}{62.0}
\resdef{n11-s030-saaty-HARKER-tau-h}{0.769}
\resdef{n11-s030-saaty-HARKER-reg-h}{5.10}
\resdef{n11-s030-saaty-HARKER-agree-h}{65.5}
\resdef{n11-s030-saaty-p-HARKER-h}{$p=0.004$}
\resdef{n11-s030-saaty-p-HARKER-h-llsm}{$p=0.002$}
\resdef{n11-s030-saaty-HARKER-LOG-hit-h}{64.3}
\resdef{n11-s030-saaty-HARKER-LOG-hitllsm-h}{64.7}
\resdef{n11-s030-saaty-HARKER-LOG-tau-h}{0.744}
\resdef{n11-s030-saaty-HARKER-LOG-reg-h}{4.17}
\resdef{n11-s030-saaty-HARKER-LOG-agree-h}{74.4}
\resdef{n11-s030-saaty-p-HARKER-LOG-h}{$p=0.050$}
\resdef{n11-s030-saaty-p-HARKER-LOG-h-llsm}{$p=0.084$}
\resdef{n11-s030-saaty-AOPT-hit-h}{59.6}
\resdef{n11-s030-saaty-AOPT-hitllsm-h}{59.6}
\resdef{n11-s030-saaty-AOPT-tau-h}{0.801}
\resdef{n11-s030-saaty-AOPT-reg-h}{5.88}
\resdef{n11-s030-saaty-AOPT-agree-h}{63.2}
\resdef{n11-s030-saaty-p-AOPT-h}{$p<0.001$}
\resdef{n11-s030-saaty-p-AOPT-h-llsm}{$p<0.001$}
\resdef{n11-s030-saaty-FROZEN-hit-h}{64.6}
\resdef{n11-s030-saaty-FROZEN-hitllsm-h}{64.6}
\resdef{n11-s030-saaty-FROZEN-tau-h}{0.779}
\resdef{n11-s030-saaty-FROZEN-reg-h}{4.75}
\resdef{n11-s030-saaty-FROZEN-agree-h}{70.2}
\resdef{n11-s030-saaty-p-FROZEN-h}{$p=0.022$}
\resdef{n11-s030-saaty-p-FROZEN-h-llsm}{$p=0.022$}
\resdef{n11-s030-saaty-STAR-hit-h}{64.6}
\resdef{n11-s030-saaty-STAR-hitllsm-h}{64.6}
\resdef{n11-s030-saaty-STAR-tau-h}{0.779}
\resdef{n11-s030-saaty-STAR-reg-h}{4.74}
\resdef{n11-s030-saaty-STAR-agree-h}{70.3}
\resdef{n11-s030-saaty-p-STAR-h}{$p=0.022$}
\resdef{n11-s030-saaty-p-STAR-h-llsm}{$p=0.022$}
\resdef{n11-s030-saaty-ALEC-hit-h}{67.2}
\resdef{n11-s030-saaty-ALEC-hitllsm-h}{67.2}
\resdef{n11-s030-saaty-ALEC-tau-h}{0.779}
\resdef{n11-s030-saaty-ALEC-reg-h}{3.72}
\resdef{n11-s030-saaty-ALEC-agree-h}{74.8}
\resdef{n11-s030-saaty-ALEC-U-hit-h}{66.3}
\resdef{n11-s030-saaty-ALEC-U-hitllsm-h}{66.3}
\resdef{n11-s030-saaty-ALEC-U-tau-h}{0.777}
\resdef{n11-s030-saaty-ALEC-U-reg-h}{3.95}
\resdef{n11-s030-saaty-ALEC-U-agree-h}{74.3}
\resdef{n11-s030-saaty-p-ALEC-U-h}{$p=0.523$}
\resdef{n11-s030-saaty-p-ALEC-U-h-llsm}{$p=0.523$}
\resdef{n11-s030-saaty-tot-n}{22}
\resdef{n11-s030-saaty-frac-n}{40}
\resdef{n11-s030-saaty-RANDOM-hit-n}{62.7}
\resdef{n11-s030-saaty-RANDOM-hitllsm-n}{62.7}
\resdef{n11-s030-saaty-RANDOM-tau-n}{0.818}
\resdef{n11-s030-saaty-RANDOM-reg-n}{4.72}
\resdef{n11-s030-saaty-RANDOM-agree-n}{71.8}
\resdef{n11-s030-saaty-p-RANDOM-n}{$p<0.001$}
\resdef{n11-s030-saaty-p-RANDOM-n-llsm}{$p<0.001$}
\resdef{n11-s030-saaty-HARKER-PG-hit-n}{65.8}
\resdef{n11-s030-saaty-HARKER-PG-hitllsm-n}{66.9}
\resdef{n11-s030-saaty-HARKER-PG-tau-n}{0.808}
\resdef{n11-s030-saaty-HARKER-PG-reg-n}{4.39}
\resdef{n11-s030-saaty-HARKER-PG-agree-n}{72.0}
\resdef{n11-s030-saaty-p-HARKER-PG-n}{$p<0.001$}
\resdef{n11-s030-saaty-p-HARKER-PG-n-llsm}{$p=0.006$}
\resdef{n11-s030-saaty-HARKER-hit-n}{68.0}
\resdef{n11-s030-saaty-HARKER-hitllsm-n}{67.9}
\resdef{n11-s030-saaty-HARKER-tau-n}{0.808}
\resdef{n11-s030-saaty-HARKER-reg-n}{3.88}
\resdef{n11-s030-saaty-HARKER-agree-n}{72.7}
\resdef{n11-s030-saaty-p-HARKER-n}{$p=0.047$}
\resdef{n11-s030-saaty-p-HARKER-n-llsm}{$p=0.037$}
\resdef{n11-s030-saaty-HARKER-LOG-hit-n}{69.2}
\resdef{n11-s030-saaty-HARKER-LOG-hitllsm-n}{71.3}
\resdef{n11-s030-saaty-HARKER-LOG-tau-n}{0.787}
\resdef{n11-s030-saaty-HARKER-LOG-reg-n}{3.07}
\resdef{n11-s030-saaty-HARKER-LOG-agree-n}{81.2}
\resdef{n11-s030-saaty-p-HARKER-LOG-n}{$p=0.078$}
\resdef{n11-s030-saaty-p-HARKER-LOG-n-llsm}{$p=1.000$}
\resdef{n11-s030-saaty-AOPT-hit-n}{65.2}
\resdef{n11-s030-saaty-AOPT-hitllsm-n}{65.2}
\resdef{n11-s030-saaty-AOPT-tau-n}{0.832}
\resdef{n11-s030-saaty-AOPT-reg-n}{4.38}
\resdef{n11-s030-saaty-AOPT-agree-n}{70.2}
\resdef{n11-s030-saaty-p-AOPT-n}{$p<0.001$}
\resdef{n11-s030-saaty-p-AOPT-n-llsm}{$p<0.001$}
\resdef{n11-s030-saaty-FROZEN-hit-n}{67.4}
\resdef{n11-s030-saaty-FROZEN-hitllsm-n}{67.4}
\resdef{n11-s030-saaty-FROZEN-tau-n}{0.815}
\resdef{n11-s030-saaty-FROZEN-reg-n}{3.85}
\resdef{n11-s030-saaty-FROZEN-agree-n}{78.8}
\resdef{n11-s030-saaty-p-FROZEN-n}{$p<0.001$}
\resdef{n11-s030-saaty-p-FROZEN-n-llsm}{$p<0.001$}
\resdef{n11-s030-saaty-STAR-hit-n}{70.7}
\resdef{n11-s030-saaty-STAR-hitllsm-n}{70.7}
\resdef{n11-s030-saaty-STAR-tau-n}{0.813}
\resdef{n11-s030-saaty-STAR-reg-n}{3.14}
\resdef{n11-s030-saaty-STAR-agree-n}{82.6}
\resdef{n11-s030-saaty-p-STAR-n}{$p=0.539$}
\resdef{n11-s030-saaty-p-STAR-n-llsm}{$p=0.539$}
\resdef{n11-s030-saaty-ALEC-hit-n}{71.3}
\resdef{n11-s030-saaty-ALEC-hitllsm-n}{71.3}
\resdef{n11-s030-saaty-ALEC-tau-n}{0.815}
\resdef{n11-s030-saaty-ALEC-reg-n}{2.83}
\resdef{n11-s030-saaty-ALEC-agree-n}{85.4}
\resdef{n11-s030-saaty-ALEC-U-hit-n}{70.1}
\resdef{n11-s030-saaty-ALEC-U-hitllsm-n}{70.1}
\resdef{n11-s030-saaty-ALEC-U-tau-n}{0.811}
\resdef{n11-s030-saaty-ALEC-U-reg-n}{2.83}
\resdef{n11-s030-saaty-ALEC-U-agree-n}{83.8}
\resdef{n11-s030-saaty-p-ALEC-U-n}{$p=0.331$}
\resdef{n11-s030-saaty-p-ALEC-U-n-llsm}{$p=0.331$}
\resdef{n11-s030-saaty-tot-2n}{33}
\resdef{n11-s030-saaty-frac-2n}{60}
\resdef{n11-s030-saaty-RANDOM-hit-2n}{68.6}
\resdef{n11-s030-saaty-RANDOM-hitllsm-2n}{68.6}
\resdef{n11-s030-saaty-RANDOM-tau-2n}{0.860}
\resdef{n11-s030-saaty-RANDOM-reg-2n}{3.10}
\resdef{n11-s030-saaty-RANDOM-agree-2n}{80.4}
\resdef{n11-s030-saaty-p-RANDOM-2n}{$p<0.001$}
\resdef{n11-s030-saaty-p-RANDOM-2n-llsm}{$p<0.001$}
\resdef{n11-s030-saaty-HARKER-PG-hit-2n}{--}
\resdef{n11-s030-saaty-HARKER-PG-hitllsm-2n}{--}
\resdef{n11-s030-saaty-HARKER-PG-tau-2n}{--}
\resdef{n11-s030-saaty-HARKER-PG-reg-2n}{--}
\resdef{n11-s030-saaty-HARKER-PG-agree-2n}{--}
\resdef{n11-s030-saaty-HARKER-hit-2n}{71.7}
\resdef{n11-s030-saaty-HARKER-hitllsm-2n}{72.0}
\resdef{n11-s030-saaty-HARKER-tau-2n}{0.859}
\resdef{n11-s030-saaty-HARKER-reg-2n}{2.64}
\resdef{n11-s030-saaty-HARKER-agree-2n}{83.6}
\resdef{n11-s030-saaty-p-HARKER-2n}{$p=0.004$}
\resdef{n11-s030-saaty-p-HARKER-2n-llsm}{$p=0.007$}
\resdef{n11-s030-saaty-HARKER-LOG-hit-2n}{73.9}
\resdef{n11-s030-saaty-HARKER-LOG-hitllsm-2n}{74.6}
\resdef{n11-s030-saaty-HARKER-LOG-tau-2n}{0.850}
\resdef{n11-s030-saaty-HARKER-LOG-reg-2n}{2.25}
\resdef{n11-s030-saaty-HARKER-LOG-agree-2n}{91.1}
\resdef{n11-s030-saaty-p-HARKER-LOG-2n}{$p=0.177$}
\resdef{n11-s030-saaty-p-HARKER-LOG-2n-llsm}{$p=0.532$}
\resdef{n11-s030-saaty-AOPT-hit-2n}{71.2}
\resdef{n11-s030-saaty-AOPT-hitllsm-2n}{71.2}
\resdef{n11-s030-saaty-AOPT-tau-2n}{0.865}
\resdef{n11-s030-saaty-AOPT-reg-2n}{2.82}
\resdef{n11-s030-saaty-AOPT-agree-2n}{82.0}
\resdef{n11-s030-saaty-p-AOPT-2n}{$p=0.002$}
\resdef{n11-s030-saaty-p-AOPT-2n-llsm}{$p=0.002$}
\resdef{n11-s030-saaty-FROZEN-hit-2n}{74.1}
\resdef{n11-s030-saaty-FROZEN-hitllsm-2n}{74.1}
\resdef{n11-s030-saaty-FROZEN-tau-2n}{0.859}
\resdef{n11-s030-saaty-FROZEN-reg-2n}{2.28}
\resdef{n11-s030-saaty-FROZEN-agree-2n}{93.1}
\resdef{n11-s030-saaty-p-FROZEN-2n}{$p=0.185$}
\resdef{n11-s030-saaty-p-FROZEN-2n-llsm}{$p=0.185$}
\resdef{n11-s030-saaty-STAR-hit-2n}{74.7}
\resdef{n11-s030-saaty-STAR-hitllsm-2n}{74.7}
\resdef{n11-s030-saaty-STAR-tau-2n}{0.860}
\resdef{n11-s030-saaty-STAR-reg-2n}{2.15}
\resdef{n11-s030-saaty-STAR-agree-2n}{98.3}
\resdef{n11-s030-saaty-p-STAR-2n}{$p=0.180$}
\resdef{n11-s030-saaty-p-STAR-2n-llsm}{$p=0.180$}
\resdef{n11-s030-saaty-ALEC-hit-2n}{75.2}
\resdef{n11-s030-saaty-ALEC-hitllsm-2n}{75.2}
\resdef{n11-s030-saaty-ALEC-tau-2n}{0.860}
\resdef{n11-s030-saaty-ALEC-reg-2n}{2.06}
\resdef{n11-s030-saaty-ALEC-agree-2n}{99.3}
\resdef{n11-s030-saaty-ALEC-U-hit-2n}{73.9}
\resdef{n11-s030-saaty-ALEC-U-hitllsm-2n}{73.9}
\resdef{n11-s030-saaty-ALEC-U-tau-2n}{0.861}
\resdef{n11-s030-saaty-ALEC-U-reg-2n}{2.24}
\resdef{n11-s030-saaty-ALEC-U-agree-2n}{94.3}
\resdef{n11-s030-saaty-p-ALEC-U-2n}{$p=0.092$}
\resdef{n11-s030-saaty-p-ALEC-U-2n-llsm}{$p=0.092$}
\resdef{n11-s030-saaty-stop-LC-q}{9.2}
\resdef{n11-s030-saaty-stop-LC-hit}{70.4}
\resdef{n11-s030-saaty-stop-LC-agree}{82.4}
\resdef{n11-s030-saaty-stop-LC-early}{100.0}
\resdef{n11-s030-saaty-stop-LC-exceed}{29.6}
\resdef{n11-s030-saaty-stop-LC-qmax}{20}
\resdef{n11-s030-saaty-stop-LC-qtot}{20.2}
\resdef{n11-s030-saaty-stop-LC-errearly}{29.6}
\resdef{n11-s030-saaty-stop-LC-errci}{[26.9, 32.5]}
\resdef{n11-s030-saaty-stop-LC-regret}{3.02}
\resdef{n11-s030-saaty-stop-LC-qround}{9}
\resdef{n11-s030-saaty-stop-LC-fixedhit}{68.4}
\resdef{n11-s030-saaty-stop-LC-fixedreg}{3.49}
\resdef{n11-s030-saaty-stop-LC-save2n}{58}
\resdef{n11-s030-saaty-stop-B10-ALEC-q}{27.8}
\resdef{n11-s030-saaty-stop-B10-ALEC-hit}{74.9}
\resdef{n11-s030-saaty-stop-B10-ALEC-agree}{99.1}
\resdef{n11-s030-saaty-stop-B10-ALEC-early}{45.9}
\resdef{n11-s030-saaty-stop-B10-ALEC-exceed}{25.1}
\resdef{n11-s030-saaty-stop-B10-ALEC-qmax}{44}
\resdef{n11-s030-saaty-stop-B10-ALEC-qtot}{38.8}
\resdef{n11-s030-saaty-stop-B10-ALEC-errearly}{8.3}
\resdef{n11-s030-saaty-stop-B10-ALEC-errci}{[6.1, 11.2]}
\resdef{n11-s030-saaty-stop-B10-ALEC-regret}{2.16}
\resdef{n11-s030-saaty-stop-B10-ALEC-save2n}{-26}
\resdef{n11-s030-saaty-stop-B05-ALECU-q}{30.6}
\resdef{n11-s030-saaty-stop-B05-ALECU-hit}{75.0}
\resdef{n11-s030-saaty-stop-B05-ALECU-agree}{99.9}
\resdef{n11-s030-saaty-stop-B05-ALECU-early}{37.6}
\resdef{n11-s030-saaty-stop-B05-ALECU-exceed}{25.0}
\resdef{n11-s030-saaty-stop-B05-ALECU-qmax}{44}
\resdef{n11-s030-saaty-stop-B05-ALECU-qtot}{41.6}
\resdef{n11-s030-saaty-stop-B05-ALECU-errearly}{5.9}
\resdef{n11-s030-saaty-stop-B05-ALECU-errci}{[3.9, 8.7]}
\resdef{n11-s030-saaty-stop-B05-ALECU-regret}{2.07}
\resdef{n11-s030-saaty-stop-B05-ALECU-save2n}{-39}
\resdef{n11-s030-saaty-stop-B10-ALECU-q}{26.4}
\resdef{n11-s030-saaty-stop-B10-ALECU-hit}{74.9}
\resdef{n11-s030-saaty-stop-B10-ALECU-agree}{99.3}
\resdef{n11-s030-saaty-stop-B10-ALECU-early}{48.0}
\resdef{n11-s030-saaty-stop-B10-ALECU-exceed}{25.1}
\resdef{n11-s030-saaty-stop-B10-ALECU-qmax}{44}
\resdef{n11-s030-saaty-stop-B10-ALECU-qtot}{37.4}
\resdef{n11-s030-saaty-stop-B10-ALECU-errearly}{8.5}
\resdef{n11-s030-saaty-stop-B10-ALECU-errci}{[6.4, 11.4]}
\resdef{n11-s030-saaty-stop-B10-ALECU-regret}{2.10}
\resdef{n11-s030-saaty-stop-B10-ALECU-save2n}{-20}
\resdef{n11-s030-saaty-stop-B20-ALECU-q}{19.3}
\resdef{n11-s030-saaty-stop-B20-ALECU-hit}{73.4}
\resdef{n11-s030-saaty-stop-B20-ALECU-agree}{95.7}
\resdef{n11-s030-saaty-stop-B20-ALECU-early}{66.2}
\resdef{n11-s030-saaty-stop-B20-ALECU-exceed}{26.6}
\resdef{n11-s030-saaty-stop-B20-ALECU-qmax}{44}
\resdef{n11-s030-saaty-stop-B20-ALECU-qtot}{30.3}
\resdef{n11-s030-saaty-stop-B20-ALECU-errearly}{16.5}
\resdef{n11-s030-saaty-stop-B20-ALECU-errci}{[13.8, 19.5]}
\resdef{n11-s030-saaty-stop-B20-ALECU-regret}{2.40}
\resdef{n11-s030-saaty-stop-B20-ALECU-qround}{19}
\resdef{n11-s030-saaty-stop-B20-ALECU-fixedhit}{75.0}
\resdef{n11-s030-saaty-stop-B20-ALECU-fixedreg}{2.13}
\resdef{n11-s030-saaty-stop-B20-ALECU-save2n}{12}
\resdef{n11-s030-saaty-stop-B10-ALECU-s06-q}{30.8}
\resdef{n11-s030-saaty-stop-B10-ALECU-s06-hit}{75.1}
\resdef{n11-s030-saaty-stop-B10-ALECU-s06-agree}{100.0}
\resdef{n11-s030-saaty-stop-B10-ALECU-s06-early}{39.3}
\resdef{n11-s030-saaty-stop-B10-ALECU-s06-exceed}{24.9}
\resdef{n11-s030-saaty-stop-B10-ALECU-s06-qmax}{44}
\resdef{n11-s030-saaty-stop-B10-ALECU-s06-qtot}{41.8}
\resdef{n11-s030-saaty-stop-B10-ALECU-s06-errearly}{5.1}
\resdef{n11-s030-saaty-stop-B10-ALECU-s06-errci}{[3.3, 7.7]}
\resdef{n11-s030-saaty-stop-B10-ALECU-s06-regret}{2.06}
\resdef{n11-s030-saaty-stop-B10-ALECU-s06-save2n}{-40}
\resdef{n11-s030-saaty-stop-B10-ALECU-cap-q}{14.9}
\resdef{n11-s030-saaty-stop-B10-ALECU-cap-hit}{73.8}
\resdef{n11-s030-saaty-stop-B10-ALECU-cap-agree}{93.7}
\resdef{n11-s030-saaty-stop-B10-ALECU-cap-early}{46.3}
\resdef{n11-s030-saaty-stop-B10-ALECU-cap-exceed}{26.2}
\resdef{n11-s030-saaty-stop-B10-ALECU-cap-qmax}{22}
\resdef{n11-s030-saaty-stop-B10-ALECU-cap-qtot}{25.9}
\resdef{n11-s030-saaty-stop-B10-ALECU-cap-errearly}{8.9}
\resdef{n11-s030-saaty-stop-B10-ALECU-cap-errci}{[6.6, 11.8]}
\resdef{n11-s030-saaty-stop-B10-ALECU-cap-regret}{2.24}
\resdef{n11-s030-saaty-stop-B10-ALECU-cap-qround}{15}
\resdef{n11-s030-saaty-stop-B10-ALECU-cap-fixedhit}{73.2}
\resdef{n11-s030-saaty-stop-B10-ALECU-cap-fixedreg}{2.33}
\resdef{n11-s030-saaty-stop-B10-ALECU-cap-save2n}{32}
\resdef{n11-s030-saaty-stop-B10-ALECU-oracle-q}{21.5}
\resdef{n11-s030-saaty-stop-B10-ALECU-oracle-hit}{73.0}
\resdef{n11-s030-saaty-stop-B10-ALECU-oracle-agree}{96.0}
\resdef{n11-s030-saaty-stop-B10-ALECU-oracle-early}{57.3}
\resdef{n11-s030-saaty-stop-B10-ALECU-oracle-exceed}{27.0}
\resdef{n11-s030-saaty-stop-B10-ALECU-oracle-qmax}{44}
\resdef{n11-s030-saaty-stop-B10-ALECU-oracle-qtot}{32.5}
\resdef{n11-s030-saaty-stop-B10-ALECU-oracle-errearly}{15.4}
\resdef{n11-s030-saaty-stop-B10-ALECU-oracle-errci}{[12.6, 18.5]}
\resdef{n11-s030-saaty-stop-B10-ALECU-oracle-regret}{2.55}
\resdef{n11-s030-saaty-stop-B10-ALECU-oracle-qround}{22}
\resdef{n11-s030-saaty-stop-B10-ALECU-oracle-fixedhit}{75.2}
\resdef{n11-s030-saaty-stop-B10-ALECU-oracle-fixedreg}{2.06}
\resdef{n11-s030-saaty-stop-B10-ALECU-oracle-save2n}{2}
\resdef{n11-s030-saaty-stop-E02-ALEC-q}{25.9}
\resdef{n11-s030-saaty-stop-E02-ALEC-hit}{74.7}
\resdef{n11-s030-saaty-stop-E02-ALEC-agree}{98.9}
\resdef{n11-s030-saaty-stop-E02-ALEC-early}{51.5}
\resdef{n11-s030-saaty-stop-E02-ALEC-exceed}{20.3}
\resdef{n11-s030-saaty-stop-E02-ALEC-qmax}{44}
\resdef{n11-s030-saaty-stop-E02-ALEC-qtot}{36.9}
\resdef{n11-s030-saaty-stop-E02-ALEC-errearly}{8.2}
\resdef{n11-s030-saaty-stop-E02-ALEC-errci}{[6.1, 10.8]}
\resdef{n11-s030-saaty-stop-E02-ALEC-regret}{2.18}
\resdef{n11-s030-saaty-stop-E02-ALEC-save2n}{-18}
\resdef{n11-s030-saaty-stop-E05-ALEC-q}{23.0}
\resdef{n11-s030-saaty-stop-E05-ALEC-hit}{74.7}
\resdef{n11-s030-saaty-stop-E05-ALEC-agree}{98.1}
\resdef{n11-s030-saaty-stop-E05-ALEC-early}{60.1}
\resdef{n11-s030-saaty-stop-E05-ALEC-exceed}{14.3}
\resdef{n11-s030-saaty-stop-E05-ALEC-qmax}{44}
\resdef{n11-s030-saaty-stop-E05-ALEC-qtot}{34.0}
\resdef{n11-s030-saaty-stop-E05-ALEC-errearly}{6.7}
\resdef{n11-s030-saaty-stop-E05-ALEC-errci}{[4.9, 8.9]}
\resdef{n11-s030-saaty-stop-E05-ALEC-regret}{2.25}
\resdef{n11-s030-saaty-stop-E05-ALEC-save2n}{-5}
\resdef{n11-s030-saaty-stop-B10-ALEC-cap-q}{15.7}
\resdef{n11-s030-saaty-stop-B10-ALEC-cap-hit}{75.0}
\resdef{n11-s030-saaty-stop-B10-ALEC-cap-agree}{98.4}
\resdef{n11-s030-saaty-stop-B10-ALEC-cap-early}{43.8}
\resdef{n11-s030-saaty-stop-B10-ALEC-cap-exceed}{25.0}
\resdef{n11-s030-saaty-stop-B10-ALEC-cap-qmax}{22}
\resdef{n11-s030-saaty-stop-B10-ALEC-cap-qtot}{26.7}
\resdef{n11-s030-saaty-stop-B10-ALEC-cap-errearly}{8.2}
\resdef{n11-s030-saaty-stop-B10-ALEC-cap-errci}{[6.0, 11.2]}
\resdef{n11-s030-saaty-stop-B10-ALEC-cap-regret}{2.15}
\resdef{n11-s030-saaty-stop-B10-ALEC-cap-qround}{16}
\resdef{n11-s030-saaty-stop-B10-ALEC-cap-fixedhit}{73.6}
\resdef{n11-s030-saaty-stop-B10-ALEC-cap-fixedreg}{2.28}
\resdef{n11-s030-saaty-stop-B10-ALEC-cap-save2n}{29}
\resdef{n11-s030-saaty-stop-E02-ALEC-cap-q}{15.1}
\resdef{n11-s030-saaty-stop-E02-ALEC-cap-hit}{74.8}
\resdef{n11-s030-saaty-stop-E02-ALEC-cap-agree}{98.2}
\resdef{n11-s030-saaty-stop-E02-ALEC-cap-early}{48.7}
\resdef{n11-s030-saaty-stop-E02-ALEC-cap-exceed}{20.2}
\resdef{n11-s030-saaty-stop-E02-ALEC-cap-qmax}{22}
\resdef{n11-s030-saaty-stop-E02-ALEC-cap-qtot}{26.1}
\resdef{n11-s030-saaty-stop-E02-ALEC-cap-errearly}{8.0}
\resdef{n11-s030-saaty-stop-E02-ALEC-cap-errci}{[5.9, 10.8]}
\resdef{n11-s030-saaty-stop-E02-ALEC-cap-regret}{2.17}
\resdef{n11-s030-saaty-stop-E02-ALEC-cap-qround}{15}
\resdef{n11-s030-saaty-stop-E02-ALEC-cap-fixedhit}{73.2}
\resdef{n11-s030-saaty-stop-E02-ALEC-cap-fixedreg}{2.33}
\resdef{n11-s030-saaty-stop-E02-ALEC-cap-save2n}{32}
\resdef{n11-s030-saaty-stop-E05-ALEC-cap-q}{14.1}
\resdef{n11-s030-saaty-stop-E05-ALEC-cap-hit}{74.8}
\resdef{n11-s030-saaty-stop-E05-ALEC-cap-agree}{97.4}
\resdef{n11-s030-saaty-stop-E05-ALEC-cap-early}{56.6}
\resdef{n11-s030-saaty-stop-E05-ALEC-cap-exceed}{14.2}
\resdef{n11-s030-saaty-stop-E05-ALEC-cap-qmax}{22}
\resdef{n11-s030-saaty-stop-E05-ALEC-cap-qtot}{25.1}
\resdef{n11-s030-saaty-stop-E05-ALEC-cap-errearly}{6.5}
\resdef{n11-s030-saaty-stop-E05-ALEC-cap-errci}{[4.8, 8.9]}
\resdef{n11-s030-saaty-stop-E05-ALEC-cap-regret}{2.24}
\resdef{n11-s030-saaty-stop-E05-ALEC-cap-qround}{14}
\resdef{n11-s030-saaty-stop-E05-ALEC-cap-fixedhit}{73.6}
\resdef{n11-s030-saaty-stop-E05-ALEC-cap-fixedreg}{2.25}
\resdef{n11-s030-saaty-stop-E05-ALEC-cap-save2n}{36}
\resdef{n11-s030-saaty-stop-E05-ALEC-cap-s06-q}{16.8}
\resdef{n11-s030-saaty-stop-E05-ALEC-cap-s06-hit}{75.2}
\resdef{n11-s030-saaty-stop-E05-ALEC-cap-s06-agree}{99.3}
\resdef{n11-s030-saaty-stop-E05-ALEC-cap-s06-early}{45.1}
\resdef{n11-s030-saaty-stop-E05-ALEC-cap-s06-exceed}{13.6}
\resdef{n11-s030-saaty-stop-E05-ALEC-cap-s06-qmax}{22}
\resdef{n11-s030-saaty-stop-E05-ALEC-cap-s06-qtot}{27.8}
\resdef{n11-s030-saaty-stop-E05-ALEC-cap-s06-errearly}{2.4}
\resdef{n11-s030-saaty-stop-E05-ALEC-cap-s06-errci}{[1.4, 4.3]}
\resdef{n11-s030-saaty-stop-E05-ALEC-cap-s06-regret}{2.06}
\resdef{n11-s030-saaty-stop-E05-ALEC-cap-s06-qround}{17}
\resdef{n11-s030-saaty-stop-E05-ALEC-cap-s06-fixedhit}{74.0}
\resdef{n11-s030-saaty-stop-E05-ALEC-cap-s06-fixedreg}{2.24}
\resdef{n11-s030-saaty-stop-E05-ALEC-cap-s06-save2n}{24}
\resdef{n11-s030-dir-full-hit}{88.6}
\resdef{n11-s030-dir-tot-h}{17}
\resdef{n11-s030-dir-frac-h}{31}
\resdef{n11-s030-dir-RANDOM-hit-h}{77.6}
\resdef{n11-s030-dir-RANDOM-hitllsm-h}{77.6}
\resdef{n11-s030-dir-RANDOM-tau-h}{0.883}
\resdef{n11-s030-dir-RANDOM-reg-h}{3.80}
\resdef{n11-s030-dir-RANDOM-agree-h}{79.4}
\resdef{n11-s030-dir-p-RANDOM-h}{$p<0.001$}
\resdef{n11-s030-dir-p-RANDOM-h-llsm}{$p<0.001$}
\resdef{n11-s030-dir-HARKER-PG-hit-h}{79.4}
\resdef{n11-s030-dir-HARKER-PG-hitllsm-h}{80.4}
\resdef{n11-s030-dir-HARKER-PG-tau-h}{0.875}
\resdef{n11-s030-dir-HARKER-PG-reg-h}{3.30}
\resdef{n11-s030-dir-HARKER-PG-agree-h}{82.1}
\resdef{n11-s030-dir-p-HARKER-PG-h}{$p<0.001$}
\resdef{n11-s030-dir-p-HARKER-PG-h-llsm}{$p=0.006$}
\resdef{n11-s030-dir-HARKER-hit-h}{78.6}
\resdef{n11-s030-dir-HARKER-hitllsm-h}{78.8}
\resdef{n11-s030-dir-HARKER-tau-h}{0.873}
\resdef{n11-s030-dir-HARKER-reg-h}{3.22}
\resdef{n11-s030-dir-HARKER-agree-h}{82.0}
\resdef{n11-s030-dir-p-HARKER-h}{$p<0.001$}
\resdef{n11-s030-dir-p-HARKER-h-llsm}{$p<0.001$}
\resdef{n11-s030-dir-HARKER-LOG-hit-h}{83.7}
\resdef{n11-s030-dir-HARKER-LOG-hitllsm-h}{83.5}
\resdef{n11-s030-dir-HARKER-LOG-tau-h}{0.869}
\resdef{n11-s030-dir-HARKER-LOG-reg-h}{1.91}
\resdef{n11-s030-dir-HARKER-LOG-agree-h}{88.6}
\resdef{n11-s030-dir-p-HARKER-LOG-h}{$p=0.764$}
\resdef{n11-s030-dir-p-HARKER-LOG-h-llsm}{$p=0.610$}
\resdef{n11-s030-dir-AOPT-hit-h}{80.3}
\resdef{n11-s030-dir-AOPT-hitllsm-h}{80.3}
\resdef{n11-s030-dir-AOPT-tau-h}{0.893}
\resdef{n11-s030-dir-AOPT-reg-h}{2.74}
\resdef{n11-s030-dir-AOPT-agree-h}{83.2}
\resdef{n11-s030-dir-p-AOPT-h}{$p=0.003$}
\resdef{n11-s030-dir-p-AOPT-h-llsm}{$p=0.003$}
\resdef{n11-s030-dir-FROZEN-hit-h}{83.2}
\resdef{n11-s030-dir-FROZEN-hitllsm-h}{83.2}
\resdef{n11-s030-dir-FROZEN-tau-h}{0.882}
\resdef{n11-s030-dir-FROZEN-reg-h}{2.13}
\resdef{n11-s030-dir-FROZEN-agree-h}{85.1}
\resdef{n11-s030-dir-p-FROZEN-h}{$p=0.328$}
\resdef{n11-s030-dir-p-FROZEN-h-llsm}{$p=0.328$}
\resdef{n11-s030-dir-STAR-hit-h}{83.3}
\resdef{n11-s030-dir-STAR-hitllsm-h}{83.3}
\resdef{n11-s030-dir-STAR-tau-h}{0.882}
\resdef{n11-s030-dir-STAR-reg-h}{2.12}
\resdef{n11-s030-dir-STAR-agree-h}{85.1}
\resdef{n11-s030-dir-p-STAR-h}{$p=0.389$}
\resdef{n11-s030-dir-p-STAR-h-llsm}{$p=0.389$}
\resdef{n11-s030-dir-ALEC-hit-h}{84.1}
\resdef{n11-s030-dir-ALEC-hitllsm-h}{84.1}
\resdef{n11-s030-dir-ALEC-tau-h}{0.880}
\resdef{n11-s030-dir-ALEC-reg-h}{1.82}
\resdef{n11-s030-dir-ALEC-agree-h}{87.4}
\resdef{n11-s030-dir-ALEC-U-hit-h}{84.9}
\resdef{n11-s030-dir-ALEC-U-hitllsm-h}{84.9}
\resdef{n11-s030-dir-ALEC-U-tau-h}{0.878}
\resdef{n11-s030-dir-ALEC-U-reg-h}{1.68}
\resdef{n11-s030-dir-ALEC-U-agree-h}{89.3}
\resdef{n11-s030-dir-p-ALEC-U-h}{$p=0.403$}
\resdef{n11-s030-dir-p-ALEC-U-h-llsm}{$p=0.403$}
\resdef{n11-s030-dir-tot-n}{22}
\resdef{n11-s030-dir-frac-n}{40}
\resdef{n11-s030-dir-RANDOM-hit-n}{81.3}
\resdef{n11-s030-dir-RANDOM-hitllsm-n}{81.3}
\resdef{n11-s030-dir-RANDOM-tau-n}{0.905}
\resdef{n11-s030-dir-RANDOM-reg-n}{2.54}
\resdef{n11-s030-dir-RANDOM-agree-n}{83.9}
\resdef{n11-s030-dir-p-RANDOM-n}{$p<0.001$}
\resdef{n11-s030-dir-p-RANDOM-n-llsm}{$p<0.001$}
\resdef{n11-s030-dir-HARKER-PG-hit-n}{82.8}
\resdef{n11-s030-dir-HARKER-PG-hitllsm-n}{82.6}
\resdef{n11-s030-dir-HARKER-PG-tau-n}{0.898}
\resdef{n11-s030-dir-HARKER-PG-reg-n}{2.20}
\resdef{n11-s030-dir-HARKER-PG-agree-n}{85.8}
\resdef{n11-s030-dir-p-HARKER-PG-n}{$p=0.012$}
\resdef{n11-s030-dir-p-HARKER-PG-n-llsm}{$p=0.008$}
\resdef{n11-s030-dir-HARKER-hit-n}{83.2}
\resdef{n11-s030-dir-HARKER-hitllsm-n}{82.7}
\resdef{n11-s030-dir-HARKER-tau-n}{0.898}
\resdef{n11-s030-dir-HARKER-reg-n}{2.05}
\resdef{n11-s030-dir-HARKER-agree-n}{87.1}
\resdef{n11-s030-dir-p-HARKER-n}{$p=0.030$}
\resdef{n11-s030-dir-p-HARKER-n-llsm}{$p=0.010$}
\resdef{n11-s030-dir-HARKER-LOG-hit-n}{86.6}
\resdef{n11-s030-dir-HARKER-LOG-hitllsm-n}{87.3}
\resdef{n11-s030-dir-HARKER-LOG-tau-n}{0.895}
\resdef{n11-s030-dir-HARKER-LOG-reg-n}{1.35}
\resdef{n11-s030-dir-HARKER-LOG-agree-n}{92.5}
\resdef{n11-s030-dir-p-HARKER-LOG-n}{$p=0.477$}
\resdef{n11-s030-dir-p-HARKER-LOG-n-llsm}{$p=0.109$}
\resdef{n11-s030-dir-AOPT-hit-n}{84.3}
\resdef{n11-s030-dir-AOPT-hitllsm-n}{84.3}
\resdef{n11-s030-dir-AOPT-tau-n}{0.911}
\resdef{n11-s030-dir-AOPT-reg-n}{1.92}
\resdef{n11-s030-dir-AOPT-agree-n}{85.6}
\resdef{n11-s030-dir-p-AOPT-n}{$p=0.185$}
\resdef{n11-s030-dir-p-AOPT-n-llsm}{$p=0.185$}
\resdef{n11-s030-dir-FROZEN-hit-n}{86.6}
\resdef{n11-s030-dir-FROZEN-hitllsm-n}{86.6}
\resdef{n11-s030-dir-FROZEN-tau-n}{0.901}
\resdef{n11-s030-dir-FROZEN-reg-n}{1.36}
\resdef{n11-s030-dir-FROZEN-agree-n}{91.4}
\resdef{n11-s030-dir-p-FROZEN-n}{$p=0.464$}
\resdef{n11-s030-dir-p-FROZEN-n-llsm}{$p=0.464$}
\resdef{n11-s030-dir-STAR-hit-n}{86.3}
\resdef{n11-s030-dir-STAR-hitllsm-n}{86.3}
\resdef{n11-s030-dir-STAR-tau-n}{0.902}
\resdef{n11-s030-dir-STAR-reg-n}{1.35}
\resdef{n11-s030-dir-STAR-agree-n}{91.8}
\resdef{n11-s030-dir-p-STAR-n}{$p=0.618$}
\resdef{n11-s030-dir-p-STAR-n-llsm}{$p=0.618$}
\resdef{n11-s030-dir-ALEC-hit-n}{85.9}
\resdef{n11-s030-dir-ALEC-hitllsm-n}{85.9}
\resdef{n11-s030-dir-ALEC-tau-n}{0.901}
\resdef{n11-s030-dir-ALEC-reg-n}{1.41}
\resdef{n11-s030-dir-ALEC-agree-n}{92.7}
\resdef{n11-s030-dir-ALEC-U-hit-n}{87.0}
\resdef{n11-s030-dir-ALEC-U-hitllsm-n}{87.0}
\resdef{n11-s030-dir-ALEC-U-tau-n}{0.899}
\resdef{n11-s030-dir-ALEC-U-reg-n}{1.31}
\resdef{n11-s030-dir-ALEC-U-agree-n}{93.1}
\resdef{n11-s030-dir-p-ALEC-U-n}{$p=0.228$}
\resdef{n11-s030-dir-p-ALEC-U-n-llsm}{$p=0.228$}
\resdef{n11-s030-dir-tot-2n}{33}
\resdef{n11-s030-dir-frac-2n}{60}
\resdef{n11-s030-dir-RANDOM-hit-2n}{86.6}
\resdef{n11-s030-dir-RANDOM-hitllsm-2n}{86.6}
\resdef{n11-s030-dir-RANDOM-tau-2n}{0.929}
\resdef{n11-s030-dir-RANDOM-reg-2n}{1.45}
\resdef{n11-s030-dir-RANDOM-agree-2n}{89.3}
\resdef{n11-s030-dir-p-RANDOM-2n}{$p=0.054$}
\resdef{n11-s030-dir-p-RANDOM-2n-llsm}{$p=0.054$}
\resdef{n11-s030-dir-HARKER-PG-hit-2n}{--}
\resdef{n11-s030-dir-HARKER-PG-hitllsm-2n}{--}
\resdef{n11-s030-dir-HARKER-PG-tau-2n}{--}
\resdef{n11-s030-dir-HARKER-PG-reg-2n}{--}
\resdef{n11-s030-dir-HARKER-PG-agree-2n}{--}
\resdef{n11-s030-dir-HARKER-hit-2n}{87.3}
\resdef{n11-s030-dir-HARKER-hitllsm-2n}{86.8}
\resdef{n11-s030-dir-HARKER-tau-2n}{0.924}
\resdef{n11-s030-dir-HARKER-reg-2n}{1.28}
\resdef{n11-s030-dir-HARKER-agree-2n}{92.9}
\resdef{n11-s030-dir-p-HARKER-2n}{$p=0.136$}
\resdef{n11-s030-dir-p-HARKER-2n-llsm}{$p=0.030$}
\resdef{n11-s030-dir-HARKER-LOG-hit-2n}{87.5}
\resdef{n11-s030-dir-HARKER-LOG-hitllsm-2n}{88.2}
\resdef{n11-s030-dir-HARKER-LOG-tau-2n}{0.921}
\resdef{n11-s030-dir-HARKER-LOG-reg-2n}{1.13}
\resdef{n11-s030-dir-HARKER-LOG-agree-2n}{95.9}
\resdef{n11-s030-dir-p-HARKER-LOG-2n}{$p=0.099$}
\resdef{n11-s030-dir-p-HARKER-LOG-2n-llsm}{$p=0.585$}
\resdef{n11-s030-dir-AOPT-hit-2n}{85.0}
\resdef{n11-s030-dir-AOPT-hitllsm-2n}{85.0}
\resdef{n11-s030-dir-AOPT-tau-2n}{0.927}
\resdef{n11-s030-dir-AOPT-reg-2n}{1.47}
\resdef{n11-s030-dir-AOPT-agree-2n}{89.7}
\resdef{n11-s030-dir-p-AOPT-2n}{$p<0.001$}
\resdef{n11-s030-dir-p-AOPT-2n-llsm}{$p<0.001$}
\resdef{n11-s030-dir-FROZEN-hit-2n}{88.2}
\resdef{n11-s030-dir-FROZEN-hitllsm-2n}{88.2}
\resdef{n11-s030-dir-FROZEN-tau-2n}{0.923}
\resdef{n11-s030-dir-FROZEN-reg-2n}{0.98}
\resdef{n11-s030-dir-FROZEN-agree-2n}{98.1}
\resdef{n11-s030-dir-p-FROZEN-2n}{$p=0.388$}
\resdef{n11-s030-dir-p-FROZEN-2n-llsm}{$p=0.388$}
\resdef{n11-s030-dir-STAR-hit-2n}{88.5}
\resdef{n11-s030-dir-STAR-hitllsm-2n}{88.5}
\resdef{n11-s030-dir-STAR-tau-2n}{0.923}
\resdef{n11-s030-dir-STAR-reg-2n}{0.97}
\resdef{n11-s030-dir-STAR-agree-2n}{99.6}
\resdef{n11-s030-dir-p-STAR-2n}{$p=1.000$}
\resdef{n11-s030-dir-p-STAR-2n-llsm}{$p=1.000$}
\resdef{n11-s030-dir-ALEC-hit-2n}{88.6}
\resdef{n11-s030-dir-ALEC-hitllsm-2n}{88.6}
\resdef{n11-s030-dir-ALEC-tau-2n}{0.923}
\resdef{n11-s030-dir-ALEC-reg-2n}{0.96}
\resdef{n11-s030-dir-ALEC-agree-2n}{99.7}
\resdef{n11-s030-dir-ALEC-U-hit-2n}{88.7}
\resdef{n11-s030-dir-ALEC-U-hitllsm-2n}{88.7}
\resdef{n11-s030-dir-ALEC-U-tau-2n}{0.925}
\resdef{n11-s030-dir-ALEC-U-reg-2n}{0.95}
\resdef{n11-s030-dir-ALEC-U-agree-2n}{98.3}
\resdef{n11-s030-dir-p-ALEC-U-2n}{$p=1.000$}
\resdef{n11-s030-dir-p-ALEC-U-2n-llsm}{$p=1.000$}
\resdef{n11-s030-dir-stop-LC-q}{8.5}
\resdef{n11-s030-dir-stop-LC-hit}{85.6}
\resdef{n11-s030-dir-stop-LC-agree}{90.8}
\resdef{n11-s030-dir-stop-LC-early}{100.0}
\resdef{n11-s030-dir-stop-LC-exceed}{14.4}
\resdef{n11-s030-dir-stop-LC-qmax}{18}
\resdef{n11-s030-dir-stop-LC-qtot}{19.5}
\resdef{n11-s030-dir-stop-LC-errearly}{14.4}
\resdef{n11-s030-dir-stop-LC-errci}{[12.4, 16.7]}
\resdef{n11-s030-dir-stop-LC-regret}{1.47}
\resdef{n11-s030-dir-stop-LC-qround}{9}
\resdef{n11-s030-dir-stop-LC-fixedhit}{85.4}
\resdef{n11-s030-dir-stop-LC-fixedreg}{1.59}
\resdef{n11-s030-dir-stop-LC-save2n}{61}
\resdef{n11-s030-dir-stop-B10-ALEC-q}{15.3}
\resdef{n11-s030-dir-stop-B10-ALEC-hit}{88.3}
\resdef{n11-s030-dir-stop-B10-ALEC-agree}{99.1}
\resdef{n11-s030-dir-stop-B10-ALEC-early}{74.3}
\resdef{n11-s030-dir-stop-B10-ALEC-exceed}{11.7}
\resdef{n11-s030-dir-stop-B10-ALEC-qmax}{44}
\resdef{n11-s030-dir-stop-B10-ALEC-qtot}{26.3}
\resdef{n11-s030-dir-stop-B10-ALEC-errearly}{3.1}
\resdef{n11-s030-dir-stop-B10-ALEC-errci}{[2.1, 4.6]}
\resdef{n11-s030-dir-stop-B10-ALEC-regret}{1.01}
\resdef{n11-s030-dir-stop-B10-ALEC-qround}{15}
\resdef{n11-s030-dir-stop-B10-ALEC-fixedhit}{89.0}
\resdef{n11-s030-dir-stop-B10-ALEC-fixedreg}{0.94}
\resdef{n11-s030-dir-stop-B10-ALEC-save2n}{30}
\resdef{n11-s030-dir-stop-B05-ALECU-q}{17.8}
\resdef{n11-s030-dir-stop-B05-ALECU-hit}{88.6}
\resdef{n11-s030-dir-stop-B05-ALECU-agree}{99.6}
\resdef{n11-s030-dir-stop-B05-ALECU-early}{68.6}
\resdef{n11-s030-dir-stop-B05-ALECU-exceed}{11.4}
\resdef{n11-s030-dir-stop-B05-ALECU-qmax}{44}
\resdef{n11-s030-dir-stop-B05-ALECU-qtot}{28.8}
\resdef{n11-s030-dir-stop-B05-ALECU-errearly}{1.9}
\resdef{n11-s030-dir-stop-B05-ALECU-errci}{[1.1, 3.2]}
\resdef{n11-s030-dir-stop-B05-ALECU-regret}{0.92}
\resdef{n11-s030-dir-stop-B05-ALECU-qround}{18}
\resdef{n11-s030-dir-stop-B05-ALECU-fixedhit}{89.0}
\resdef{n11-s030-dir-stop-B05-ALECU-fixedreg}{0.91}
\resdef{n11-s030-dir-stop-B05-ALECU-save2n}{19}
\resdef{n11-s030-dir-stop-B10-ALECU-q}{13.8}
\resdef{n11-s030-dir-stop-B10-ALECU-hit}{88.2}
\resdef{n11-s030-dir-stop-B10-ALECU-agree}{99.0}
\resdef{n11-s030-dir-stop-B10-ALECU-early}{76.2}
\resdef{n11-s030-dir-stop-B10-ALECU-exceed}{11.8}
\resdef{n11-s030-dir-stop-B10-ALECU-qmax}{44}
\resdef{n11-s030-dir-stop-B10-ALECU-qtot}{24.8}
\resdef{n11-s030-dir-stop-B10-ALECU-errearly}{3.9}
\resdef{n11-s030-dir-stop-B10-ALECU-errci}{[2.8, 5.6]}
\resdef{n11-s030-dir-stop-B10-ALECU-regret}{0.98}
\resdef{n11-s030-dir-stop-B10-ALECU-qround}{14}
\resdef{n11-s030-dir-stop-B10-ALECU-fixedhit}{88.1}
\resdef{n11-s030-dir-stop-B10-ALECU-fixedreg}{1.01}
\resdef{n11-s030-dir-stop-B10-ALECU-save2n}{37}
\resdef{n11-s030-dir-stop-B20-ALECU-q}{8.3}
\resdef{n11-s030-dir-stop-B20-ALECU-hit}{86.5}
\resdef{n11-s030-dir-stop-B20-ALECU-agree}{95.8}
\resdef{n11-s030-dir-stop-B20-ALECU-early}{88.1}
\resdef{n11-s030-dir-stop-B20-ALECU-exceed}{13.5}
\resdef{n11-s030-dir-stop-B20-ALECU-qmax}{44}
\resdef{n11-s030-dir-stop-B20-ALECU-qtot}{19.3}
\resdef{n11-s030-dir-stop-B20-ALECU-errearly}{8.1}
\resdef{n11-s030-dir-stop-B20-ALECU-errci}{[6.4, 10.0]}
\resdef{n11-s030-dir-stop-B20-ALECU-regret}{1.22}
\resdef{n11-s030-dir-stop-B20-ALECU-qround}{8}
\resdef{n11-s030-dir-stop-B20-ALECU-fixedhit}{84.8}
\resdef{n11-s030-dir-stop-B20-ALECU-fixedreg}{1.71}
\resdef{n11-s030-dir-stop-B20-ALECU-save2n}{62}
\resdef{n11-s030-dir-stop-B10-ALECU-s06-q}{18.3}
\resdef{n11-s030-dir-stop-B10-ALECU-s06-hit}{88.5}
\resdef{n11-s030-dir-stop-B10-ALECU-s06-agree}{99.9}
\resdef{n11-s030-dir-stop-B10-ALECU-s06-early}{70.2}
\resdef{n11-s030-dir-stop-B10-ALECU-s06-exceed}{11.5}
\resdef{n11-s030-dir-stop-B10-ALECU-s06-qmax}{44}
\resdef{n11-s030-dir-stop-B10-ALECU-s06-qtot}{29.3}
\resdef{n11-s030-dir-stop-B10-ALECU-s06-errearly}{1.7}
\resdef{n11-s030-dir-stop-B10-ALECU-s06-errci}{[1.0, 3.0]}
\resdef{n11-s030-dir-stop-B10-ALECU-s06-regret}{0.95}
\resdef{n11-s030-dir-stop-B10-ALECU-s06-qround}{18}
\resdef{n11-s030-dir-stop-B10-ALECU-s06-fixedhit}{89.0}
\resdef{n11-s030-dir-stop-B10-ALECU-s06-fixedreg}{0.91}
\resdef{n11-s030-dir-stop-B10-ALECU-s06-save2n}{17}
\resdef{n11-s030-dir-stop-B10-ALECU-cap-q}{8.6}
\resdef{n11-s030-dir-stop-B10-ALECU-cap-hit}{88.3}
\resdef{n11-s030-dir-stop-B10-ALECU-cap-agree}{97.5}
\resdef{n11-s030-dir-stop-B10-ALECU-cap-early}{75.7}
\resdef{n11-s030-dir-stop-B10-ALECU-cap-exceed}{11.7}
\resdef{n11-s030-dir-stop-B10-ALECU-cap-qmax}{22}
\resdef{n11-s030-dir-stop-B10-ALECU-cap-qtot}{19.6}
\resdef{n11-s030-dir-stop-B10-ALECU-cap-errearly}{4.0}
\resdef{n11-s030-dir-stop-B10-ALECU-cap-errci}{[2.8, 5.6]}
\resdef{n11-s030-dir-stop-B10-ALECU-cap-regret}{1.04}
\resdef{n11-s030-dir-stop-B10-ALECU-cap-qround}{9}
\resdef{n11-s030-dir-stop-B10-ALECU-cap-fixedhit}{85.4}
\resdef{n11-s030-dir-stop-B10-ALECU-cap-fixedreg}{1.59}
\resdef{n11-s030-dir-stop-B10-ALECU-cap-save2n}{61}
\resdef{n11-s030-dir-stop-B10-ALECU-oracle-q}{12.0}
\resdef{n11-s030-dir-stop-B10-ALECU-oracle-hit}{87.3}
\resdef{n11-s030-dir-stop-B10-ALECU-oracle-agree}{97.9}
\resdef{n11-s030-dir-stop-B10-ALECU-oracle-early}{77.9}
\resdef{n11-s030-dir-stop-B10-ALECU-oracle-exceed}{12.7}
\resdef{n11-s030-dir-stop-B10-ALECU-oracle-qmax}{44}
\resdef{n11-s030-dir-stop-B10-ALECU-oracle-qtot}{23.0}
\resdef{n11-s030-dir-stop-B10-ALECU-oracle-errearly}{5.9}
\resdef{n11-s030-dir-stop-B10-ALECU-oracle-errci}{[4.5, 7.8]}
\resdef{n11-s030-dir-stop-B10-ALECU-oracle-regret}{1.12}
\resdef{n11-s030-dir-stop-B10-ALECU-oracle-qround}{12}
\resdef{n11-s030-dir-stop-B10-ALECU-oracle-fixedhit}{86.8}
\resdef{n11-s030-dir-stop-B10-ALECU-oracle-fixedreg}{1.25}
\resdef{n11-s030-dir-stop-B10-ALECU-oracle-save2n}{45}
\resdef{n11-s030-dir-stop-E02-ALEC-q}{13.7}
\resdef{n11-s030-dir-stop-E02-ALEC-hit}{88.1}
\resdef{n11-s030-dir-stop-E02-ALEC-agree}{98.7}
\resdef{n11-s030-dir-stop-E02-ALEC-early}{78.4}
\resdef{n11-s030-dir-stop-E02-ALEC-exceed}{9.8}
\resdef{n11-s030-dir-stop-E02-ALEC-qmax}{44}
\resdef{n11-s030-dir-stop-E02-ALEC-qtot}{24.7}
\resdef{n11-s030-dir-stop-E02-ALEC-errearly}{3.3}
\resdef{n11-s030-dir-stop-E02-ALEC-errci}{[2.3, 4.8]}
\resdef{n11-s030-dir-stop-E02-ALEC-regret}{1.04}
\resdef{n11-s030-dir-stop-E02-ALEC-qround}{14}
\resdef{n11-s030-dir-stop-E02-ALEC-fixedhit}{88.1}
\resdef{n11-s030-dir-stop-E02-ALEC-fixedreg}{1.01}
\resdef{n11-s030-dir-stop-E02-ALEC-save2n}{38}
\resdef{n11-s030-dir-stop-E05-ALEC-q}{11.4}
\resdef{n11-s030-dir-stop-E05-ALEC-hit}{88.2}
\resdef{n11-s030-dir-stop-E05-ALEC-agree}{98.1}
\resdef{n11-s030-dir-stop-E05-ALEC-early}{84.6}
\resdef{n11-s030-dir-stop-E05-ALEC-exceed}{7.7}
\resdef{n11-s030-dir-stop-E05-ALEC-qmax}{44}
\resdef{n11-s030-dir-stop-E05-ALEC-qtot}{22.4}
\resdef{n11-s030-dir-stop-E05-ALEC-errearly}{3.5}
\resdef{n11-s030-dir-stop-E05-ALEC-errci}{[2.5, 5.0]}
\resdef{n11-s030-dir-stop-E05-ALEC-regret}{1.03}
\resdef{n11-s030-dir-stop-E05-ALEC-qround}{11}
\resdef{n11-s030-dir-stop-E05-ALEC-fixedhit}{85.9}
\resdef{n11-s030-dir-stop-E05-ALEC-fixedreg}{1.41}
\resdef{n11-s030-dir-stop-E05-ALEC-save2n}{48}
\resdef{n11-s030-dir-stop-B10-ALEC-cap-q}{9.6}
\resdef{n11-s030-dir-stop-B10-ALEC-cap-hit}{88.3}
\resdef{n11-s030-dir-stop-B10-ALEC-cap-agree}{98.8}
\resdef{n11-s030-dir-stop-B10-ALEC-cap-early}{73.4}
\resdef{n11-s030-dir-stop-B10-ALEC-cap-exceed}{11.7}
\resdef{n11-s030-dir-stop-B10-ALEC-cap-qmax}{22}
\resdef{n11-s030-dir-stop-B10-ALEC-cap-qtot}{20.6}
\resdef{n11-s030-dir-stop-B10-ALEC-cap-errearly}{3.1}
\resdef{n11-s030-dir-stop-B10-ALEC-cap-errci}{[2.1, 4.7]}
\resdef{n11-s030-dir-stop-B10-ALEC-cap-regret}{1.04}
\resdef{n11-s030-dir-stop-B10-ALEC-cap-qround}{10}
\resdef{n11-s030-dir-stop-B10-ALEC-cap-fixedhit}{85.4}
\resdef{n11-s030-dir-stop-B10-ALEC-cap-fixedreg}{1.56}
\resdef{n11-s030-dir-stop-B10-ALEC-cap-save2n}{56}
\resdef{n11-s030-dir-stop-E02-ALEC-cap-q}{8.9}
\resdef{n11-s030-dir-stop-E02-ALEC-cap-hit}{88.1}
\resdef{n11-s030-dir-stop-E02-ALEC-cap-agree}{98.5}
\resdef{n11-s030-dir-stop-E02-ALEC-cap-early}{77.1}
\resdef{n11-s030-dir-stop-E02-ALEC-cap-exceed}{9.8}
\resdef{n11-s030-dir-stop-E02-ALEC-cap-qmax}{22}
\resdef{n11-s030-dir-stop-E02-ALEC-cap-qtot}{19.9}
\resdef{n11-s030-dir-stop-E02-ALEC-cap-errearly}{3.4}
\resdef{n11-s030-dir-stop-E02-ALEC-cap-errci}{[2.3, 4.9]}
\resdef{n11-s030-dir-stop-E02-ALEC-cap-regret}{1.06}
\resdef{n11-s030-dir-stop-E02-ALEC-cap-qround}{9}
\resdef{n11-s030-dir-stop-E02-ALEC-cap-fixedhit}{85.4}
\resdef{n11-s030-dir-stop-E02-ALEC-cap-fixedreg}{1.59}
\resdef{n11-s030-dir-stop-E02-ALEC-cap-save2n}{60}
\resdef{n11-s030-dir-stop-E05-ALEC-cap-q}{7.9}
\resdef{n11-s030-dir-stop-E05-ALEC-cap-hit}{88.2}
\resdef{n11-s030-dir-stop-E05-ALEC-cap-agree}{97.9}
\resdef{n11-s030-dir-stop-E05-ALEC-cap-early}{82.9}
\resdef{n11-s030-dir-stop-E05-ALEC-cap-exceed}{7.7}
\resdef{n11-s030-dir-stop-E05-ALEC-cap-qmax}{22}
\resdef{n11-s030-dir-stop-E05-ALEC-cap-qtot}{18.9}
\resdef{n11-s030-dir-stop-E05-ALEC-cap-errearly}{3.4}
\resdef{n11-s030-dir-stop-E05-ALEC-cap-errci}{[2.3, 4.8]}
\resdef{n11-s030-dir-stop-E05-ALEC-cap-regret}{1.05}
\resdef{n11-s030-dir-stop-E05-ALEC-cap-qround}{8}
\resdef{n11-s030-dir-stop-E05-ALEC-cap-fixedhit}{84.8}
\resdef{n11-s030-dir-stop-E05-ALEC-cap-fixedreg}{1.71}
\resdef{n11-s030-dir-stop-E05-ALEC-cap-save2n}{64}
\resdef{n11-s030-dir-stop-E05-ALEC-cap-s06-q}{11.2}
\resdef{n11-s030-dir-stop-E05-ALEC-cap-s06-hit}{88.6}
\resdef{n11-s030-dir-stop-E05-ALEC-cap-s06-agree}{99.5}
\resdef{n11-s030-dir-stop-E05-ALEC-cap-s06-early}{74.2}
\resdef{n11-s030-dir-stop-E05-ALEC-cap-s06-exceed}{7.1}
\resdef{n11-s030-dir-stop-E05-ALEC-cap-s06-qmax}{22}
\resdef{n11-s030-dir-stop-E05-ALEC-cap-s06-qtot}{22.2}
\resdef{n11-s030-dir-stop-E05-ALEC-cap-s06-errearly}{1.6}
\resdef{n11-s030-dir-stop-E05-ALEC-cap-s06-errci}{[0.9, 2.8]}
\resdef{n11-s030-dir-stop-E05-ALEC-cap-s06-regret}{0.97}
\resdef{n11-s030-dir-stop-E05-ALEC-cap-s06-qround}{11}
\resdef{n11-s030-dir-stop-E05-ALEC-cap-s06-fixedhit}{85.9}
\resdef{n11-s030-dir-stop-E05-ALEC-cap-s06-fixedreg}{1.41}
\resdef{n11-s030-dir-stop-E05-ALEC-cap-s06-save2n}{49}
\resdef{n11-s030-path-full-hit}{76.6}
\resdef{n11-s030-path-tot-h}{17}
\resdef{n11-s030-path-frac-h}{31}
\resdef{n11-s030-path-RANDOM-hit-h}{61.6}
\resdef{n11-s030-path-RANDOM-hitllsm-h}{61.6}
\resdef{n11-s030-path-RANDOM-tau-h}{0.807}
\resdef{n11-s030-path-RANDOM-reg-h}{4.83}
\resdef{n11-s030-path-RANDOM-agree-h}{69.5}
\resdef{n11-s030-path-p-RANDOM-h}{$p<0.001$}
\resdef{n11-s030-path-p-RANDOM-h-llsm}{$p<0.001$}
\resdef{n11-s030-path-HARKER-PG-hit-h}{63.0}
\resdef{n11-s030-path-HARKER-PG-hitllsm-h}{65.7}
\resdef{n11-s030-path-HARKER-PG-tau-h}{0.796}
\resdef{n11-s030-path-HARKER-PG-reg-h}{4.74}
\resdef{n11-s030-path-HARKER-PG-agree-h}{69.4}
\resdef{n11-s030-path-p-HARKER-PG-h}{$p<0.001$}
\resdef{n11-s030-path-p-HARKER-PG-h-llsm}{$p<0.001$}
\resdef{n11-s030-path-HARKER-hit-h}{63.3}
\resdef{n11-s030-path-HARKER-hitllsm-h}{63.8}
\resdef{n11-s030-path-HARKER-tau-h}{0.799}
\resdef{n11-s030-path-HARKER-reg-h}{4.51}
\resdef{n11-s030-path-HARKER-agree-h}{70.4}
\resdef{n11-s030-path-p-HARKER-h}{$p<0.001$}
\resdef{n11-s030-path-p-HARKER-h-llsm}{$p<0.001$}
\resdef{n11-s030-path-HARKER-LOG-hit-h}{70.4}
\resdef{n11-s030-path-HARKER-LOG-hitllsm-h}{71.0}
\resdef{n11-s030-path-HARKER-LOG-tau-h}{0.781}
\resdef{n11-s030-path-HARKER-LOG-reg-h}{3.16}
\resdef{n11-s030-path-HARKER-LOG-agree-h}{76.3}
\resdef{n11-s030-path-p-HARKER-LOG-h}{$p=0.204$}
\resdef{n11-s030-path-p-HARKER-LOG-h-llsm}{$p=0.425$}
\resdef{n11-s030-path-AOPT-hit-h}{65.3}
\resdef{n11-s030-path-AOPT-hitllsm-h}{65.3}
\resdef{n11-s030-path-AOPT-tau-h}{0.828}
\resdef{n11-s030-path-AOPT-reg-h}{4.13}
\resdef{n11-s030-path-AOPT-agree-h}{72.1}
\resdef{n11-s030-path-p-AOPT-h}{$p<0.001$}
\resdef{n11-s030-path-p-AOPT-h-llsm}{$p<0.001$}
\resdef{n11-s030-path-FROZEN-hit-h}{69.0}
\resdef{n11-s030-path-FROZEN-hitllsm-h}{69.0}
\resdef{n11-s030-path-FROZEN-tau-h}{0.804}
\resdef{n11-s030-path-FROZEN-reg-h}{3.34}
\resdef{n11-s030-path-FROZEN-agree-h}{76.4}
\resdef{n11-s030-path-p-FROZEN-h}{$p=0.015$}
\resdef{n11-s030-path-p-FROZEN-h-llsm}{$p=0.015$}
\resdef{n11-s030-path-STAR-hit-h}{69.0}
\resdef{n11-s030-path-STAR-hitllsm-h}{69.0}
\resdef{n11-s030-path-STAR-tau-h}{0.804}
\resdef{n11-s030-path-STAR-reg-h}{3.34}
\resdef{n11-s030-path-STAR-agree-h}{76.4}
\resdef{n11-s030-path-p-STAR-h}{$p=0.015$}
\resdef{n11-s030-path-p-STAR-h-llsm}{$p=0.015$}
\resdef{n11-s030-path-ALEC-hit-h}{72.1}
\resdef{n11-s030-path-ALEC-hitllsm-h}{72.1}
\resdef{n11-s030-path-ALEC-tau-h}{0.799}
\resdef{n11-s030-path-ALEC-reg-h}{2.83}
\resdef{n11-s030-path-ALEC-agree-h}{78.9}
\resdef{n11-s030-path-ALEC-U-hit-h}{71.4}
\resdef{n11-s030-path-ALEC-U-hitllsm-h}{71.4}
\resdef{n11-s030-path-ALEC-U-tau-h}{0.805}
\resdef{n11-s030-path-ALEC-U-reg-h}{2.78}
\resdef{n11-s030-path-ALEC-U-agree-h}{80.2}
\resdef{n11-s030-path-p-ALEC-U-h}{$p=0.600$}
\resdef{n11-s030-path-p-ALEC-U-h-llsm}{$p=0.600$}
\resdef{n11-s030-path-tot-n}{22}
\resdef{n11-s030-path-frac-n}{40}
\resdef{n11-s030-path-RANDOM-hit-n}{65.2}
\resdef{n11-s030-path-RANDOM-hitllsm-n}{65.2}
\resdef{n11-s030-path-RANDOM-tau-n}{0.843}
\resdef{n11-s030-path-RANDOM-reg-n}{3.92}
\resdef{n11-s030-path-RANDOM-agree-n}{72.8}
\resdef{n11-s030-path-p-RANDOM-n}{$p<0.001$}
\resdef{n11-s030-path-p-RANDOM-n-llsm}{$p<0.001$}
\resdef{n11-s030-path-HARKER-PG-hit-n}{70.9}
\resdef{n11-s030-path-HARKER-PG-hitllsm-n}{70.9}
\resdef{n11-s030-path-HARKER-PG-tau-n}{0.834}
\resdef{n11-s030-path-HARKER-PG-reg-n}{2.99}
\resdef{n11-s030-path-HARKER-PG-agree-n}{77.8}
\resdef{n11-s030-path-p-HARKER-PG-n}{$p=0.009$}
\resdef{n11-s030-path-p-HARKER-PG-n-llsm}{$p=0.008$}
\resdef{n11-s030-path-HARKER-hit-n}{70.5}
\resdef{n11-s030-path-HARKER-hitllsm-n}{69.9}
\resdef{n11-s030-path-HARKER-tau-n}{0.833}
\resdef{n11-s030-path-HARKER-reg-n}{2.77}
\resdef{n11-s030-path-HARKER-agree-n}{77.8}
\resdef{n11-s030-path-p-HARKER-n}{$p=0.004$}
\resdef{n11-s030-path-p-HARKER-n-llsm}{$p<0.001$}
\resdef{n11-s030-path-HARKER-LOG-hit-n}{73.7}
\resdef{n11-s030-path-HARKER-LOG-hitllsm-n}{73.8}
\resdef{n11-s030-path-HARKER-LOG-tau-n}{0.824}
\resdef{n11-s030-path-HARKER-LOG-reg-n}{2.19}
\resdef{n11-s030-path-HARKER-LOG-agree-n}{83.3}
\resdef{n11-s030-path-p-HARKER-LOG-n}{$p=0.359$}
\resdef{n11-s030-path-p-HARKER-LOG-n-llsm}{$p=0.391$}
\resdef{n11-s030-path-AOPT-hit-n}{69.5}
\resdef{n11-s030-path-AOPT-hitllsm-n}{69.5}
\resdef{n11-s030-path-AOPT-tau-n}{0.855}
\resdef{n11-s030-path-AOPT-reg-n}{3.01}
\resdef{n11-s030-path-AOPT-agree-n}{77.3}
\resdef{n11-s030-path-p-AOPT-n}{$p<0.001$}
\resdef{n11-s030-path-p-AOPT-n-llsm}{$p<0.001$}
\resdef{n11-s030-path-FROZEN-hit-n}{73.1}
\resdef{n11-s030-path-FROZEN-hitllsm-n}{73.1}
\resdef{n11-s030-path-FROZEN-tau-n}{0.836}
\resdef{n11-s030-path-FROZEN-reg-n}{2.42}
\resdef{n11-s030-path-FROZEN-agree-n}{81.5}
\resdef{n11-s030-path-p-FROZEN-n}{$p=0.155$}
\resdef{n11-s030-path-p-FROZEN-n-llsm}{$p=0.155$}
\resdef{n11-s030-path-STAR-hit-n}{74.9}
\resdef{n11-s030-path-STAR-hitllsm-n}{74.9}
\resdef{n11-s030-path-STAR-tau-n}{0.835}
\resdef{n11-s030-path-STAR-reg-n}{2.04}
\resdef{n11-s030-path-STAR-agree-n}{84.9}
\resdef{n11-s030-path-p-STAR-n}{$p=1.000$}
\resdef{n11-s030-path-p-STAR-n-llsm}{$p=1.000$}
\resdef{n11-s030-path-ALEC-hit-n}{74.8}
\resdef{n11-s030-path-ALEC-hitllsm-n}{74.8}
\resdef{n11-s030-path-ALEC-tau-n}{0.833}
\resdef{n11-s030-path-ALEC-reg-n}{1.86}
\resdef{n11-s030-path-ALEC-agree-n}{86.7}
\resdef{n11-s030-path-ALEC-U-hit-n}{74.7}
\resdef{n11-s030-path-ALEC-U-hitllsm-n}{74.7}
\resdef{n11-s030-path-ALEC-U-tau-n}{0.835}
\resdef{n11-s030-path-ALEC-U-reg-n}{2.12}
\resdef{n11-s030-path-ALEC-U-agree-n}{85.1}
\resdef{n11-s030-path-p-ALEC-U-n}{$p=1.000$}
\resdef{n11-s030-path-p-ALEC-U-n-llsm}{$p=1.000$}
\resdef{n11-s030-path-tot-2n}{33}
\resdef{n11-s030-path-frac-2n}{60}
\resdef{n11-s030-path-RANDOM-hit-2n}{71.3}
\resdef{n11-s030-path-RANDOM-hitllsm-2n}{71.3}
\resdef{n11-s030-path-RANDOM-tau-2n}{0.881}
\resdef{n11-s030-path-RANDOM-reg-2n}{2.54}
\resdef{n11-s030-path-RANDOM-agree-2n}{81.2}
\resdef{n11-s030-path-p-RANDOM-2n}{$p<0.001$}
\resdef{n11-s030-path-p-RANDOM-2n-llsm}{$p<0.001$}
\resdef{n11-s030-path-HARKER-PG-hit-2n}{--}
\resdef{n11-s030-path-HARKER-PG-hitllsm-2n}{--}
\resdef{n11-s030-path-HARKER-PG-tau-2n}{--}
\resdef{n11-s030-path-HARKER-PG-reg-2n}{--}
\resdef{n11-s030-path-HARKER-PG-agree-2n}{--}
\resdef{n11-s030-path-HARKER-hit-2n}{73.2}
\resdef{n11-s030-path-HARKER-hitllsm-2n}{73.2}
\resdef{n11-s030-path-HARKER-tau-2n}{0.875}
\resdef{n11-s030-path-HARKER-reg-2n}{2.09}
\resdef{n11-s030-path-HARKER-agree-2n}{86.0}
\resdef{n11-s030-path-p-HARKER-2n}{$p=0.003$}
\resdef{n11-s030-path-p-HARKER-2n-llsm}{$p=0.004$}
\resdef{n11-s030-path-HARKER-LOG-hit-2n}{75.2}
\resdef{n11-s030-path-HARKER-LOG-hitllsm-2n}{75.4}
\resdef{n11-s030-path-HARKER-LOG-tau-2n}{0.870}
\resdef{n11-s030-path-HARKER-LOG-reg-2n}{1.75}
\resdef{n11-s030-path-HARKER-LOG-agree-2n}{90.8}
\resdef{n11-s030-path-p-HARKER-LOG-2n}{$p=0.165$}
\resdef{n11-s030-path-p-HARKER-LOG-2n-llsm}{$p=0.222$}
\resdef{n11-s030-path-AOPT-hit-2n}{73.5}
\resdef{n11-s030-path-AOPT-hitllsm-2n}{73.5}
\resdef{n11-s030-path-AOPT-tau-2n}{0.883}
\resdef{n11-s030-path-AOPT-reg-2n}{2.22}
\resdef{n11-s030-path-AOPT-agree-2n}{83.7}
\resdef{n11-s030-path-p-AOPT-2n}{$p=0.013$}
\resdef{n11-s030-path-p-AOPT-2n-llsm}{$p=0.013$}
\resdef{n11-s030-path-FROZEN-hit-2n}{76.2}
\resdef{n11-s030-path-FROZEN-hitllsm-2n}{76.2}
\resdef{n11-s030-path-FROZEN-tau-2n}{0.872}
\resdef{n11-s030-path-FROZEN-reg-2n}{1.68}
\resdef{n11-s030-path-FROZEN-agree-2n}{95.2}
\resdef{n11-s030-path-p-FROZEN-2n}{$p=0.749$}
\resdef{n11-s030-path-p-FROZEN-2n-llsm}{$p=0.749$}
\resdef{n11-s030-path-STAR-hit-2n}{76.5}
\resdef{n11-s030-path-STAR-hitllsm-2n}{76.5}
\resdef{n11-s030-path-STAR-tau-2n}{0.873}
\resdef{n11-s030-path-STAR-reg-2n}{1.58}
\resdef{n11-s030-path-STAR-agree-2n}{99.1}
\resdef{n11-s030-path-p-STAR-2n}{$p=1.000$}
\resdef{n11-s030-path-p-STAR-2n-llsm}{$p=1.000$}
\resdef{n11-s030-path-ALEC-hit-2n}{76.5}
\resdef{n11-s030-path-ALEC-hitllsm-2n}{76.5}
\resdef{n11-s030-path-ALEC-tau-2n}{0.872}
\resdef{n11-s030-path-ALEC-reg-2n}{1.59}
\resdef{n11-s030-path-ALEC-agree-2n}{99.4}
\resdef{n11-s030-path-ALEC-U-hit-2n}{75.8}
\resdef{n11-s030-path-ALEC-U-hitllsm-2n}{75.8}
\resdef{n11-s030-path-ALEC-U-tau-2n}{0.874}
\resdef{n11-s030-path-ALEC-U-reg-2n}{1.64}
\resdef{n11-s030-path-ALEC-U-agree-2n}{94.8}
\resdef{n11-s030-path-p-ALEC-U-2n}{$p=0.360$}
\resdef{n11-s030-path-p-ALEC-U-2n-llsm}{$p=0.360$}
\resdef{n11-s030-path-stop-LC-q}{9.2}
\resdef{n11-s030-path-stop-LC-hit}{75.0}
\resdef{n11-s030-path-stop-LC-agree}{85.7}
\resdef{n11-s030-path-stop-LC-early}{100.0}
\resdef{n11-s030-path-stop-LC-exceed}{25.0}
\resdef{n11-s030-path-stop-LC-qmax}{21}
\resdef{n11-s030-path-stop-LC-qtot}{20.2}
\resdef{n11-s030-path-stop-LC-errearly}{25.0}
\resdef{n11-s030-path-stop-LC-errci}{[22.4, 27.8]}
\resdef{n11-s030-path-stop-LC-regret}{2.05}
\resdef{n11-s030-path-stop-LC-qround}{9}
\resdef{n11-s030-path-stop-LC-fixedhit}{72.5}
\resdef{n11-s030-path-stop-LC-fixedreg}{2.67}
\resdef{n11-s030-path-stop-LC-save2n}{58}
\resdef{n11-s030-path-stop-B10-ALEC-q}{25.5}
\resdef{n11-s030-path-stop-B10-ALEC-hit}{77.1}
\resdef{n11-s030-path-stop-B10-ALEC-agree}{98.9}
\resdef{n11-s030-path-stop-B10-ALEC-early}{52.0}
\resdef{n11-s030-path-stop-B10-ALEC-exceed}{22.9}
\resdef{n11-s030-path-stop-B10-ALEC-qmax}{44}
\resdef{n11-s030-path-stop-B10-ALEC-qtot}{36.5}
\resdef{n11-s030-path-stop-B10-ALEC-errearly}{7.1}
\resdef{n11-s030-path-stop-B10-ALEC-errci}{[5.2, 9.7]}
\resdef{n11-s030-path-stop-B10-ALEC-regret}{1.55}
\resdef{n11-s030-path-stop-B10-ALEC-save2n}{-16}
\resdef{n11-s030-path-stop-B05-ALECU-q}{29.1}
\resdef{n11-s030-path-stop-B05-ALECU-hit}{76.8}
\resdef{n11-s030-path-stop-B05-ALECU-agree}{99.6}
\resdef{n11-s030-path-stop-B05-ALECU-early}{42.0}
\resdef{n11-s030-path-stop-B05-ALECU-exceed}{23.2}
\resdef{n11-s030-path-stop-B05-ALECU-qmax}{44}
\resdef{n11-s030-path-stop-B05-ALECU-qtot}{40.1}
\resdef{n11-s030-path-stop-B05-ALECU-errearly}{5.2}
\resdef{n11-s030-path-stop-B05-ALECU-errci}{[3.5, 7.8]}
\resdef{n11-s030-path-stop-B05-ALECU-regret}{1.54}
\resdef{n11-s030-path-stop-B05-ALECU-save2n}{-32}
\resdef{n11-s030-path-stop-B10-ALECU-q}{23.3}
\resdef{n11-s030-path-stop-B10-ALECU-hit}{77.1}
\resdef{n11-s030-path-stop-B10-ALECU-agree}{98.8}
\resdef{n11-s030-path-stop-B10-ALECU-early}{56.6}
\resdef{n11-s030-path-stop-B10-ALECU-exceed}{22.9}
\resdef{n11-s030-path-stop-B10-ALECU-qmax}{44}
\resdef{n11-s030-path-stop-B10-ALECU-qtot}{34.3}
\resdef{n11-s030-path-stop-B10-ALECU-errearly}{8.1}
\resdef{n11-s030-path-stop-B10-ALECU-errci}{[6.1, 10.7]}
\resdef{n11-s030-path-stop-B10-ALECU-regret}{1.52}
\resdef{n11-s030-path-stop-B10-ALECU-save2n}{-6}
\resdef{n11-s030-path-stop-B20-ALECU-q}{15.8}
\resdef{n11-s030-path-stop-B20-ALECU-hit}{76.3}
\resdef{n11-s030-path-stop-B20-ALECU-agree}{95.3}
\resdef{n11-s030-path-stop-B20-ALECU-early}{73.9}
\resdef{n11-s030-path-stop-B20-ALECU-exceed}{23.7}
\resdef{n11-s030-path-stop-B20-ALECU-qmax}{44}
\resdef{n11-s030-path-stop-B20-ALECU-qtot}{26.8}
\resdef{n11-s030-path-stop-B20-ALECU-errearly}{15.0}
\resdef{n11-s030-path-stop-B20-ALECU-errci}{[12.6, 17.8]}
\resdef{n11-s030-path-stop-B20-ALECU-regret}{1.69}
\resdef{n11-s030-path-stop-B20-ALECU-qround}{16}
\resdef{n11-s030-path-stop-B20-ALECU-fixedhit}{76.5}
\resdef{n11-s030-path-stop-B20-ALECU-fixedreg}{1.64}
\resdef{n11-s030-path-stop-B20-ALECU-save2n}{28}
\resdef{n11-s030-path-stop-B10-ALECU-s06-q}{29.6}
\resdef{n11-s030-path-stop-B10-ALECU-s06-hit}{76.6}
\resdef{n11-s030-path-stop-B10-ALECU-s06-agree}{100.0}
\resdef{n11-s030-path-stop-B10-ALECU-s06-early}{43.5}
\resdef{n11-s030-path-stop-B10-ALECU-s06-exceed}{23.4}
\resdef{n11-s030-path-stop-B10-ALECU-s06-qmax}{44}
\resdef{n11-s030-path-stop-B10-ALECU-s06-qtot}{40.6}
\resdef{n11-s030-path-stop-B10-ALECU-s06-errearly}{3.4}
\resdef{n11-s030-path-stop-B10-ALECU-s06-errci}{[2.1, 5.6]}
\resdef{n11-s030-path-stop-B10-ALECU-s06-regret}{1.55}
\resdef{n11-s030-path-stop-B10-ALECU-s06-save2n}{-35}
\resdef{n11-s030-path-stop-B10-ALECU-cap-q}{13.7}
\resdef{n11-s030-path-stop-B10-ALECU-cap-hit}{76.3}
\resdef{n11-s030-path-stop-B10-ALECU-cap-agree}{93.8}
\resdef{n11-s030-path-stop-B10-ALECU-cap-early}{55.0}
\resdef{n11-s030-path-stop-B10-ALECU-cap-exceed}{23.7}
\resdef{n11-s030-path-stop-B10-ALECU-cap-qmax}{22}
\resdef{n11-s030-path-stop-B10-ALECU-cap-qtot}{24.7}
\resdef{n11-s030-path-stop-B10-ALECU-cap-errearly}{8.4}
\resdef{n11-s030-path-stop-B10-ALECU-cap-errci}{[6.3, 11.0]}
\resdef{n11-s030-path-stop-B10-ALECU-cap-regret}{1.61}
\resdef{n11-s030-path-stop-B10-ALECU-cap-qround}{14}
\resdef{n11-s030-path-stop-B10-ALECU-cap-fixedhit}{74.9}
\resdef{n11-s030-path-stop-B10-ALECU-cap-fixedreg}{1.83}
\resdef{n11-s030-path-stop-B10-ALECU-cap-save2n}{38}
\resdef{n11-s030-path-stop-B10-ALECU-oracle-q}{21.8}
\resdef{n11-s030-path-stop-B10-ALECU-oracle-hit}{77.2}
\resdef{n11-s030-path-stop-B10-ALECU-oracle-agree}{97.9}
\resdef{n11-s030-path-stop-B10-ALECU-oracle-early}{57.3}
\resdef{n11-s030-path-stop-B10-ALECU-oracle-exceed}{22.8}
\resdef{n11-s030-path-stop-B10-ALECU-oracle-qmax}{44}
\resdef{n11-s030-path-stop-B10-ALECU-oracle-qtot}{32.8}
\resdef{n11-s030-path-stop-B10-ALECU-oracle-errearly}{9.4}
\resdef{n11-s030-path-stop-B10-ALECU-oracle-errci}{[7.3, 12.1]}
\resdef{n11-s030-path-stop-B10-ALECU-oracle-regret}{1.55}
\resdef{n11-s030-path-stop-B10-ALECU-oracle-qround}{22}
\resdef{n11-s030-path-stop-B10-ALECU-oracle-fixedhit}{76.5}
\resdef{n11-s030-path-stop-B10-ALECU-oracle-fixedreg}{1.59}
\resdef{n11-s030-path-stop-B10-ALECU-oracle-save2n}{1}
\resdef{n11-s030-path-stop-E02-ALEC-q}{23.3}
\resdef{n11-s030-path-stop-E02-ALEC-hit}{76.7}
\resdef{n11-s030-path-stop-E02-ALEC-agree}{98.4}
\resdef{n11-s030-path-stop-E02-ALEC-early}{58.1}
\resdef{n11-s030-path-stop-E02-ALEC-exceed}{18.6}
\resdef{n11-s030-path-stop-E02-ALEC-qmax}{44}
\resdef{n11-s030-path-stop-E02-ALEC-qtot}{34.3}
\resdef{n11-s030-path-stop-E02-ALEC-errearly}{7.1}
\resdef{n11-s030-path-stop-E02-ALEC-errci}{[5.2, 9.4]}
\resdef{n11-s030-path-stop-E02-ALEC-regret}{1.61}
\resdef{n11-s030-path-stop-E02-ALEC-save2n}{-6}
\resdef{n11-s030-path-stop-E05-ALEC-q}{19.9}
\resdef{n11-s030-path-stop-E05-ALEC-hit}{76.8}
\resdef{n11-s030-path-stop-E05-ALEC-agree}{97.7}
\resdef{n11-s030-path-stop-E05-ALEC-early}{67.9}
\resdef{n11-s030-path-stop-E05-ALEC-exceed}{11.6}
\resdef{n11-s030-path-stop-E05-ALEC-qmax}{44}
\resdef{n11-s030-path-stop-E05-ALEC-qtot}{30.9}
\resdef{n11-s030-path-stop-E05-ALEC-errearly}{4.6}
\resdef{n11-s030-path-stop-E05-ALEC-errci}{[3.2, 6.4]}
\resdef{n11-s030-path-stop-E05-ALEC-regret}{1.60}
\resdef{n11-s030-path-stop-E05-ALEC-qround}{20}
\resdef{n11-s030-path-stop-E05-ALEC-fixedhit}{76.5}
\resdef{n11-s030-path-stop-E05-ALEC-fixedreg}{1.61}
\resdef{n11-s030-path-stop-E05-ALEC-save2n}{9}
\resdef{n11-s030-path-stop-B10-ALEC-cap-q}{14.8}
\resdef{n11-s030-path-stop-B10-ALEC-cap-hit}{77.0}
\resdef{n11-s030-path-stop-B10-ALEC-cap-agree}{98.3}
\resdef{n11-s030-path-stop-B10-ALEC-cap-early}{50.1}
\resdef{n11-s030-path-stop-B10-ALEC-cap-exceed}{23.0}
\resdef{n11-s030-path-stop-B10-ALEC-cap-qmax}{22}
\resdef{n11-s030-path-stop-B10-ALEC-cap-qtot}{25.8}
\resdef{n11-s030-path-stop-B10-ALEC-cap-errearly}{7.2}
\resdef{n11-s030-path-stop-B10-ALEC-cap-errci}{[5.2, 9.8]}
\resdef{n11-s030-path-stop-B10-ALEC-cap-regret}{1.59}
\resdef{n11-s030-path-stop-B10-ALEC-cap-qround}{15}
\resdef{n11-s030-path-stop-B10-ALEC-cap-fixedhit}{75.9}
\resdef{n11-s030-path-stop-B10-ALEC-cap-fixedreg}{1.75}
\resdef{n11-s030-path-stop-B10-ALEC-cap-save2n}{33}
\resdef{n11-s030-path-stop-E02-ALEC-cap-q}{14.0}
\resdef{n11-s030-path-stop-E02-ALEC-cap-hit}{76.6}
\resdef{n11-s030-path-stop-E02-ALEC-cap-agree}{97.8}
\resdef{n11-s030-path-stop-E02-ALEC-cap-early}{56.2}
\resdef{n11-s030-path-stop-E02-ALEC-cap-exceed}{18.8}
\resdef{n11-s030-path-stop-E02-ALEC-cap-qmax}{22}
\resdef{n11-s030-path-stop-E02-ALEC-cap-qtot}{25.0}
\resdef{n11-s030-path-stop-E02-ALEC-cap-errearly}{7.1}
\resdef{n11-s030-path-stop-E02-ALEC-cap-errci}{[5.3, 9.5]}
\resdef{n11-s030-path-stop-E02-ALEC-cap-regret}{1.65}
\resdef{n11-s030-path-stop-E02-ALEC-cap-qround}{14}
\resdef{n11-s030-path-stop-E02-ALEC-cap-fixedhit}{74.9}
\resdef{n11-s030-path-stop-E02-ALEC-cap-fixedreg}{1.83}
\resdef{n11-s030-path-stop-E02-ALEC-cap-save2n}{36}
\resdef{n11-s030-path-stop-E05-ALEC-cap-q}{12.7}
\resdef{n11-s030-path-stop-E05-ALEC-cap-hit}{76.7}
\resdef{n11-s030-path-stop-E05-ALEC-cap-agree}{97.1}
\resdef{n11-s030-path-stop-E05-ALEC-cap-early}{65.3}
\resdef{n11-s030-path-stop-E05-ALEC-cap-exceed}{11.9}
\resdef{n11-s030-path-stop-E05-ALEC-cap-qmax}{22}
\resdef{n11-s030-path-stop-E05-ALEC-cap-qtot}{23.7}
\resdef{n11-s030-path-stop-E05-ALEC-cap-errearly}{4.4}
\resdef{n11-s030-path-stop-E05-ALEC-cap-errci}{[3.1, 6.3]}
\resdef{n11-s030-path-stop-E05-ALEC-cap-regret}{1.64}
\resdef{n11-s030-path-stop-E05-ALEC-cap-qround}{13}
\resdef{n11-s030-path-stop-E05-ALEC-cap-fixedhit}{75.9}
\resdef{n11-s030-path-stop-E05-ALEC-cap-fixedreg}{1.87}
\resdef{n11-s030-path-stop-E05-ALEC-cap-save2n}{42}
\resdef{n11-s030-path-stop-E05-ALEC-cap-s06-q}{16.2}
\resdef{n11-s030-path-stop-E05-ALEC-cap-s06-hit}{76.7}
\resdef{n11-s030-path-stop-E05-ALEC-cap-s06-agree}{99.2}
\resdef{n11-s030-path-stop-E05-ALEC-cap-s06-early}{51.6}
\resdef{n11-s030-path-stop-E05-ALEC-cap-s06-exceed}{11.8}
\resdef{n11-s030-path-stop-E05-ALEC-cap-s06-qmax}{22}
\resdef{n11-s030-path-stop-E05-ALEC-cap-s06-qtot}{27.2}
\resdef{n11-s030-path-stop-E05-ALEC-cap-s06-errearly}{1.9}
\resdef{n11-s030-path-stop-E05-ALEC-cap-s06-errci}{[1.1, 3.5]}
\resdef{n11-s030-path-stop-E05-ALEC-cap-s06-regret}{1.59}
\resdef{n11-s030-path-stop-E05-ALEC-cap-s06-qround}{16}
\resdef{n11-s030-path-stop-E05-ALEC-cap-s06-fixedhit}{76.5}
\resdef{n11-s030-path-stop-E05-ALEC-cap-s06-fixedreg}{1.64}
\resdef{n11-s030-path-stop-E05-ALEC-cap-s06-save2n}{26}
\resdef{n11-s030-star-full-hit}{76.6}
\resdef{n11-s030-star-tot-h}{17}
\resdef{n11-s030-star-frac-h}{31}
\resdef{n11-s030-star-RANDOM-hit-h}{59.3}
\resdef{n11-s030-star-RANDOM-hitllsm-h}{59.3}
\resdef{n11-s030-star-RANDOM-tau-h}{0.803}
\resdef{n11-s030-star-RANDOM-reg-h}{5.56}
\resdef{n11-s030-star-RANDOM-agree-h}{65.1}
\resdef{n11-s030-star-p-RANDOM-h}{$p<0.001$}
\resdef{n11-s030-star-p-RANDOM-h-llsm}{$p<0.001$}
\resdef{n11-s030-star-HARKER-PG-hit-h}{59.6}
\resdef{n11-s030-star-HARKER-PG-hitllsm-h}{60.8}
\resdef{n11-s030-star-HARKER-PG-tau-h}{0.777}
\resdef{n11-s030-star-HARKER-PG-reg-h}{5.68}
\resdef{n11-s030-star-HARKER-PG-agree-h}{64.2}
\resdef{n11-s030-star-p-HARKER-PG-h}{$p<0.001$}
\resdef{n11-s030-star-p-HARKER-PG-h-llsm}{$p<0.001$}
\resdef{n11-s030-star-HARKER-hit-h}{59.5}
\resdef{n11-s030-star-HARKER-hitllsm-h}{59.9}
\resdef{n11-s030-star-HARKER-tau-h}{0.791}
\resdef{n11-s030-star-HARKER-reg-h}{5.41}
\resdef{n11-s030-star-HARKER-agree-h}{67.1}
\resdef{n11-s030-star-p-HARKER-h}{$p<0.001$}
\resdef{n11-s030-star-p-HARKER-h-llsm}{$p<0.001$}
\resdef{n11-s030-star-HARKER-LOG-hit-h}{70.3}
\resdef{n11-s030-star-HARKER-LOG-hitllsm-h}{70.1}
\resdef{n11-s030-star-HARKER-LOG-tau-h}{0.790}
\resdef{n11-s030-star-HARKER-LOG-reg-h}{3.21}
\resdef{n11-s030-star-HARKER-LOG-agree-h}{76.9}
\resdef{n11-s030-star-p-HARKER-LOG-h}{$p=0.526$}
\resdef{n11-s030-star-p-HARKER-LOG-h-llsm}{$p=0.616$}
\resdef{n11-s030-star-AOPT-hit-h}{61.3}
\resdef{n11-s030-star-AOPT-hitllsm-h}{61.3}
\resdef{n11-s030-star-AOPT-tau-h}{0.804}
\resdef{n11-s030-star-AOPT-reg-h}{5.29}
\resdef{n11-s030-star-AOPT-agree-h}{66.9}
\resdef{n11-s030-star-p-AOPT-h}{$p<0.001$}
\resdef{n11-s030-star-p-AOPT-h-llsm}{$p<0.001$}
\resdef{n11-s030-star-FROZEN-hit-h}{67.3}
\resdef{n11-s030-star-FROZEN-hitllsm-h}{67.3}
\resdef{n11-s030-star-FROZEN-tau-h}{0.799}
\resdef{n11-s030-star-FROZEN-reg-h}{3.96}
\resdef{n11-s030-star-FROZEN-agree-h}{73.2}
\resdef{n11-s030-star-p-FROZEN-h}{$p=0.080$}
\resdef{n11-s030-star-p-FROZEN-h-llsm}{$p=0.080$}
\resdef{n11-s030-star-STAR-hit-h}{67.6}
\resdef{n11-s030-star-STAR-hitllsm-h}{67.6}
\resdef{n11-s030-star-STAR-tau-h}{0.799}
\resdef{n11-s030-star-STAR-reg-h}{3.88}
\resdef{n11-s030-star-STAR-agree-h}{73.5}
\resdef{n11-s030-star-p-STAR-h}{$p=0.123$}
\resdef{n11-s030-star-p-STAR-h-llsm}{$p=0.123$}
\resdef{n11-s030-star-ALEC-hit-h}{69.4}
\resdef{n11-s030-star-ALEC-hitllsm-h}{69.4}
\resdef{n11-s030-star-ALEC-tau-h}{0.799}
\resdef{n11-s030-star-ALEC-reg-h}{3.11}
\resdef{n11-s030-star-ALEC-agree-h}{78.2}
\resdef{n11-s030-star-ALEC-U-hit-h}{70.7}
\resdef{n11-s030-star-ALEC-U-hitllsm-h}{70.7}
\resdef{n11-s030-star-ALEC-U-tau-h}{0.797}
\resdef{n11-s030-star-ALEC-U-reg-h}{2.86}
\resdef{n11-s030-star-ALEC-U-agree-h}{78.0}
\resdef{n11-s030-star-p-ALEC-U-h}{$p=0.294$}
\resdef{n11-s030-star-p-ALEC-U-h-llsm}{$p=0.294$}
\resdef{n11-s030-star-tot-n}{22}
\resdef{n11-s030-star-frac-n}{40}
\resdef{n11-s030-star-RANDOM-hit-n}{64.5}
\resdef{n11-s030-star-RANDOM-hitllsm-n}{64.5}
\resdef{n11-s030-star-RANDOM-tau-n}{0.835}
\resdef{n11-s030-star-RANDOM-reg-n}{4.11}
\resdef{n11-s030-star-RANDOM-agree-n}{71.6}
\resdef{n11-s030-star-p-RANDOM-n}{$p<0.001$}
\resdef{n11-s030-star-p-RANDOM-n-llsm}{$p<0.001$}
\resdef{n11-s030-star-HARKER-PG-hit-n}{65.5}
\resdef{n11-s030-star-HARKER-PG-hitllsm-n}{64.9}
\resdef{n11-s030-star-HARKER-PG-tau-n}{0.817}
\resdef{n11-s030-star-HARKER-PG-reg-n}{4.09}
\resdef{n11-s030-star-HARKER-PG-agree-n}{72.3}
\resdef{n11-s030-star-p-HARKER-PG-n}{$p<0.001$}
\resdef{n11-s030-star-p-HARKER-PG-n-llsm}{$p<0.001$}
\resdef{n11-s030-star-HARKER-hit-n}{67.6}
\resdef{n11-s030-star-HARKER-hitllsm-n}{66.6}
\resdef{n11-s030-star-HARKER-tau-n}{0.824}
\resdef{n11-s030-star-HARKER-reg-n}{3.50}
\resdef{n11-s030-star-HARKER-agree-n}{75.0}
\resdef{n11-s030-star-p-HARKER-n}{$p<0.001$}
\resdef{n11-s030-star-p-HARKER-n-llsm}{$p<0.001$}
\resdef{n11-s030-star-HARKER-LOG-hit-n}{71.7}
\resdef{n11-s030-star-HARKER-LOG-hitllsm-n}{71.3}
\resdef{n11-s030-star-HARKER-LOG-tau-n}{0.817}
\resdef{n11-s030-star-HARKER-LOG-reg-n}{2.46}
\resdef{n11-s030-star-HARKER-LOG-agree-n}{82.3}
\resdef{n11-s030-star-p-HARKER-LOG-n}{$p=0.145$}
\resdef{n11-s030-star-p-HARKER-LOG-n-llsm}{$p=0.064$}
\resdef{n11-s030-star-AOPT-hit-n}{65.5}
\resdef{n11-s030-star-AOPT-hitllsm-n}{65.5}
\resdef{n11-s030-star-AOPT-tau-n}{0.843}
\resdef{n11-s030-star-AOPT-reg-n}{3.85}
\resdef{n11-s030-star-AOPT-agree-n}{73.7}
\resdef{n11-s030-star-p-AOPT-n}{$p<0.001$}
\resdef{n11-s030-star-p-AOPT-n-llsm}{$p<0.001$}
\resdef{n11-s030-star-FROZEN-hit-n}{71.7}
\resdef{n11-s030-star-FROZEN-hitllsm-n}{71.7}
\resdef{n11-s030-star-FROZEN-tau-n}{0.831}
\resdef{n11-s030-star-FROZEN-reg-n}{2.81}
\resdef{n11-s030-star-FROZEN-agree-n}{81.5}
\resdef{n11-s030-star-p-FROZEN-n}{$p=0.113$}
\resdef{n11-s030-star-p-FROZEN-n-llsm}{$p=0.113$}
\resdef{n11-s030-star-STAR-hit-n}{72.7}
\resdef{n11-s030-star-STAR-hitllsm-n}{72.7}
\resdef{n11-s030-star-STAR-tau-n}{0.829}
\resdef{n11-s030-star-STAR-reg-n}{2.48}
\resdef{n11-s030-star-STAR-agree-n}{82.5}
\resdef{n11-s030-star-p-STAR-n}{$p=0.380$}
\resdef{n11-s030-star-p-STAR-n-llsm}{$p=0.380$}
\resdef{n11-s030-star-ALEC-hit-n}{73.6}
\resdef{n11-s030-star-ALEC-hitllsm-n}{73.6}
\resdef{n11-s030-star-ALEC-tau-n}{0.829}
\resdef{n11-s030-star-ALEC-reg-n}{2.29}
\resdef{n11-s030-star-ALEC-agree-n}{84.8}
\resdef{n11-s030-star-ALEC-U-hit-n}{72.5}
\resdef{n11-s030-star-ALEC-U-hitllsm-n}{72.5}
\resdef{n11-s030-star-ALEC-U-tau-n}{0.825}
\resdef{n11-s030-star-ALEC-U-reg-n}{2.37}
\resdef{n11-s030-star-ALEC-U-agree-n}{82.8}
\resdef{n11-s030-star-p-ALEC-U-n}{$p=0.371$}
\resdef{n11-s030-star-p-ALEC-U-n-llsm}{$p=0.371$}
\resdef{n11-s030-star-tot-2n}{33}
\resdef{n11-s030-star-frac-2n}{60}
\resdef{n11-s030-star-RANDOM-hit-2n}{71.0}
\resdef{n11-s030-star-RANDOM-hitllsm-2n}{71.0}
\resdef{n11-s030-star-RANDOM-tau-2n}{0.878}
\resdef{n11-s030-star-RANDOM-reg-2n}{2.59}
\resdef{n11-s030-star-RANDOM-agree-2n}{80.1}
\resdef{n11-s030-star-p-RANDOM-2n}{$p<0.001$}
\resdef{n11-s030-star-p-RANDOM-2n-llsm}{$p<0.001$}
\resdef{n11-s030-star-HARKER-PG-hit-2n}{--}
\resdef{n11-s030-star-HARKER-PG-hitllsm-2n}{--}
\resdef{n11-s030-star-HARKER-PG-tau-2n}{--}
\resdef{n11-s030-star-HARKER-PG-reg-2n}{--}
\resdef{n11-s030-star-HARKER-PG-agree-2n}{--}
\resdef{n11-s030-star-HARKER-hit-2n}{74.0}
\resdef{n11-s030-star-HARKER-hitllsm-2n}{74.4}
\resdef{n11-s030-star-HARKER-tau-2n}{0.872}
\resdef{n11-s030-star-HARKER-reg-2n}{2.04}
\resdef{n11-s030-star-HARKER-agree-2n}{85.0}
\resdef{n11-s030-star-p-HARKER-2n}{$p=0.100$}
\resdef{n11-s030-star-p-HARKER-2n-llsm}{$p=0.201$}
\resdef{n11-s030-star-HARKER-LOG-hit-2n}{76.1}
\resdef{n11-s030-star-HARKER-LOG-hitllsm-2n}{76.7}
\resdef{n11-s030-star-HARKER-LOG-tau-2n}{0.866}
\resdef{n11-s030-star-HARKER-LOG-reg-2n}{1.69}
\resdef{n11-s030-star-HARKER-LOG-agree-2n}{90.6}
\resdef{n11-s030-star-p-HARKER-LOG-2n}{$p=1.000$}
\resdef{n11-s030-star-p-HARKER-LOG-2n-llsm}{$p=0.470$}
\resdef{n11-s030-star-AOPT-hit-2n}{73.7}
\resdef{n11-s030-star-AOPT-hitllsm-2n}{73.7}
\resdef{n11-s030-star-AOPT-tau-2n}{0.885}
\resdef{n11-s030-star-AOPT-reg-2n}{2.21}
\resdef{n11-s030-star-AOPT-agree-2n}{82.6}
\resdef{n11-s030-star-p-AOPT-2n}{$p=0.077$}
\resdef{n11-s030-star-p-AOPT-2n-llsm}{$p=0.077$}
\resdef{n11-s030-star-FROZEN-hit-2n}{75.6}
\resdef{n11-s030-star-FROZEN-hitllsm-2n}{75.6}
\resdef{n11-s030-star-FROZEN-tau-2n}{0.872}
\resdef{n11-s030-star-FROZEN-reg-2n}{1.73}
\resdef{n11-s030-star-FROZEN-agree-2n}{94.5}
\resdef{n11-s030-star-p-FROZEN-2n}{$p=0.672$}
\resdef{n11-s030-star-p-FROZEN-2n-llsm}{$p=0.672$}
\resdef{n11-s030-star-STAR-hit-2n}{75.6}
\resdef{n11-s030-star-STAR-hitllsm-2n}{75.6}
\resdef{n11-s030-star-STAR-tau-2n}{0.871}
\resdef{n11-s030-star-STAR-reg-2n}{1.63}
\resdef{n11-s030-star-STAR-agree-2n}{98.3}
\resdef{n11-s030-star-p-STAR-2n}{$p=0.289$}
\resdef{n11-s030-star-p-STAR-2n-llsm}{$p=0.289$}
\resdef{n11-s030-star-ALEC-hit-2n}{76.0}
\resdef{n11-s030-star-ALEC-hitllsm-2n}{76.0}
\resdef{n11-s030-star-ALEC-tau-2n}{0.871}
\resdef{n11-s030-star-ALEC-reg-2n}{1.62}
\resdef{n11-s030-star-ALEC-agree-2n}{98.6}
\resdef{n11-s030-star-ALEC-U-hit-2n}{75.5}
\resdef{n11-s030-star-ALEC-U-hitllsm-2n}{75.5}
\resdef{n11-s030-star-ALEC-U-tau-2n}{0.872}
\resdef{n11-s030-star-ALEC-U-reg-2n}{1.67}
\resdef{n11-s030-star-ALEC-U-agree-2n}{93.7}
\resdef{n11-s030-star-p-ALEC-U-2n}{$p=0.603$}
\resdef{n11-s030-star-p-ALEC-U-2n-llsm}{$p=0.603$}
\resdef{n11-s030-star-stop-LC-q}{8.8}
\resdef{n11-s030-star-stop-LC-hit}{72.2}
\resdef{n11-s030-star-stop-LC-agree}{83.0}
\resdef{n11-s030-star-stop-LC-early}{100.0}
\resdef{n11-s030-star-stop-LC-exceed}{27.8}
\resdef{n11-s030-star-stop-LC-qmax}{24}
\resdef{n11-s030-star-stop-LC-qtot}{19.8}
\resdef{n11-s030-star-stop-LC-errearly}{27.8}
\resdef{n11-s030-star-stop-LC-errci}{[25.1, 30.7]}
\resdef{n11-s030-star-stop-LC-regret}{2.41}
\resdef{n11-s030-star-stop-LC-qround}{9}
\resdef{n11-s030-star-stop-LC-fixedhit}{71.5}
\resdef{n11-s030-star-stop-LC-fixedreg}{2.62}
\resdef{n11-s030-star-stop-LC-save2n}{60}
\resdef{n11-s030-star-stop-B10-ALEC-q}{25.4}
\resdef{n11-s030-star-stop-B10-ALEC-hit}{76.7}
\resdef{n11-s030-star-stop-B10-ALEC-agree}{98.7}
\resdef{n11-s030-star-stop-B10-ALEC-early}{52.5}
\resdef{n11-s030-star-stop-B10-ALEC-exceed}{23.3}
\resdef{n11-s030-star-stop-B10-ALEC-qmax}{44}
\resdef{n11-s030-star-stop-B10-ALEC-qtot}{36.4}
\resdef{n11-s030-star-stop-B10-ALEC-errearly}{8.2}
\resdef{n11-s030-star-stop-B10-ALEC-errci}{[6.1, 10.9]}
\resdef{n11-s030-star-stop-B10-ALEC-regret}{1.59}
\resdef{n11-s030-star-stop-B10-ALEC-save2n}{-15}
\resdef{n11-s030-star-stop-B05-ALECU-q}{29.1}
\resdef{n11-s030-star-stop-B05-ALECU-hit}{76.5}
\resdef{n11-s030-star-stop-B05-ALECU-agree}{99.7}
\resdef{n11-s030-star-stop-B05-ALECU-early}{42.1}
\resdef{n11-s030-star-stop-B05-ALECU-exceed}{23.5}
\resdef{n11-s030-star-stop-B05-ALECU-qmax}{44}
\resdef{n11-s030-star-stop-B05-ALECU-qtot}{40.1}
\resdef{n11-s030-star-stop-B05-ALECU-errearly}{5.0}
\resdef{n11-s030-star-stop-B05-ALECU-errci}{[3.3, 7.5]}
\resdef{n11-s030-star-stop-B05-ALECU-regret}{1.57}
\resdef{n11-s030-star-stop-B05-ALECU-save2n}{-32}
\resdef{n11-s030-star-stop-B10-ALECU-q}{24.2}
\resdef{n11-s030-star-stop-B10-ALECU-hit}{76.8}
\resdef{n11-s030-star-stop-B10-ALECU-agree}{98.6}
\resdef{n11-s030-star-stop-B10-ALECU-early}{54.2}
\resdef{n11-s030-star-stop-B10-ALECU-exceed}{23.2}
\resdef{n11-s030-star-stop-B10-ALECU-qmax}{44}
\resdef{n11-s030-star-stop-B10-ALECU-qtot}{35.2}
\resdef{n11-s030-star-stop-B10-ALECU-errearly}{8.5}
\resdef{n11-s030-star-stop-B10-ALECU-errci}{[6.4, 11.1]}
\resdef{n11-s030-star-stop-B10-ALECU-regret}{1.58}
\resdef{n11-s030-star-stop-B10-ALECU-save2n}{-10}
\resdef{n11-s030-star-stop-B20-ALECU-q}{16.0}
\resdef{n11-s030-star-stop-B20-ALECU-hit}{75.3}
\resdef{n11-s030-star-stop-B20-ALECU-agree}{94.4}
\resdef{n11-s030-star-stop-B20-ALECU-early}{73.9}
\resdef{n11-s030-star-stop-B20-ALECU-exceed}{24.7}
\resdef{n11-s030-star-stop-B20-ALECU-qmax}{44}
\resdef{n11-s030-star-stop-B20-ALECU-qtot}{27.0}
\resdef{n11-s030-star-stop-B20-ALECU-errearly}{16.2}
\resdef{n11-s030-star-stop-B20-ALECU-errci}{[13.8, 19.1]}
\resdef{n11-s030-star-stop-B20-ALECU-regret}{1.93}
\resdef{n11-s030-star-stop-B20-ALECU-qround}{16}
\resdef{n11-s030-star-stop-B20-ALECU-fixedhit}{75.7}
\resdef{n11-s030-star-stop-B20-ALECU-fixedreg}{1.61}
\resdef{n11-s030-star-stop-B20-ALECU-save2n}{27}
\resdef{n11-s030-star-stop-B10-ALECU-s06-q}{29.7}
\resdef{n11-s030-star-stop-B10-ALECU-s06-hit}{76.6}
\resdef{n11-s030-star-stop-B10-ALECU-s06-agree}{100.0}
\resdef{n11-s030-star-stop-B10-ALECU-s06-early}{43.4}
\resdef{n11-s030-star-stop-B10-ALECU-s06-exceed}{23.4}
\resdef{n11-s030-star-stop-B10-ALECU-s06-qmax}{44}
\resdef{n11-s030-star-stop-B10-ALECU-s06-qtot}{40.7}
\resdef{n11-s030-star-stop-B10-ALECU-s06-errearly}{3.2}
\resdef{n11-s030-star-stop-B10-ALECU-s06-errci}{[1.9, 5.3]}
\resdef{n11-s030-star-stop-B10-ALECU-s06-regret}{1.55}
\resdef{n11-s030-star-stop-B10-ALECU-s06-save2n}{-35}
\resdef{n11-s030-star-stop-B10-ALECU-cap-q}{14.0}
\resdef{n11-s030-star-stop-B10-ALECU-cap-hit}{75.6}
\resdef{n11-s030-star-stop-B10-ALECU-cap-agree}{92.8}
\resdef{n11-s030-star-stop-B10-ALECU-cap-early}{52.5}
\resdef{n11-s030-star-stop-B10-ALECU-cap-exceed}{24.4}
\resdef{n11-s030-star-stop-B10-ALECU-cap-qmax}{22}
\resdef{n11-s030-star-stop-B10-ALECU-cap-qtot}{25.0}
\resdef{n11-s030-star-stop-B10-ALECU-cap-errearly}{8.6}
\resdef{n11-s030-star-stop-B10-ALECU-cap-errci}{[6.5, 11.3]}
\resdef{n11-s030-star-stop-B10-ALECU-cap-regret}{1.70}
\resdef{n11-s030-star-stop-B10-ALECU-cap-qround}{14}
\resdef{n11-s030-star-stop-B10-ALECU-cap-fixedhit}{74.2}
\resdef{n11-s030-star-stop-B10-ALECU-cap-fixedreg}{1.93}
\resdef{n11-s030-star-stop-B10-ALECU-cap-save2n}{37}
\resdef{n11-s030-star-stop-B10-ALECU-oracle-q}{22.9}
\resdef{n11-s030-star-stop-B10-ALECU-oracle-hit}{76.5}
\resdef{n11-s030-star-stop-B10-ALECU-oracle-agree}{97.5}
\resdef{n11-s030-star-stop-B10-ALECU-oracle-early}{54.2}
\resdef{n11-s030-star-stop-B10-ALECU-oracle-exceed}{23.5}
\resdef{n11-s030-star-stop-B10-ALECU-oracle-qmax}{44}
\resdef{n11-s030-star-stop-B10-ALECU-oracle-qtot}{33.9}
\resdef{n11-s030-star-stop-B10-ALECU-oracle-errearly}{9.0}
\resdef{n11-s030-star-stop-B10-ALECU-oracle-errci}{[6.9, 11.8]}
\resdef{n11-s030-star-stop-B10-ALECU-oracle-regret}{1.71}
\resdef{n11-s030-star-stop-B10-ALECU-oracle-save2n}{-4}
\resdef{n11-s030-star-stop-E02-ALEC-q}{23.1}
\resdef{n11-s030-star-stop-E02-ALEC-hit}{76.6}
\resdef{n11-s030-star-stop-E02-ALEC-agree}{98.0}
\resdef{n11-s030-star-stop-E02-ALEC-early}{59.0}
\resdef{n11-s030-star-stop-E02-ALEC-exceed}{18.5}
\resdef{n11-s030-star-stop-E02-ALEC-qmax}{44}
\resdef{n11-s030-star-stop-E02-ALEC-qtot}{34.1}
\resdef{n11-s030-star-stop-E02-ALEC-errearly}{7.6}
\resdef{n11-s030-star-stop-E02-ALEC-errci}{[5.7, 10.1]}
\resdef{n11-s030-star-stop-E02-ALEC-regret}{1.61}
\resdef{n11-s030-star-stop-E02-ALEC-save2n}{-5}
\resdef{n11-s030-star-stop-E05-ALEC-q}{19.7}
\resdef{n11-s030-star-stop-E05-ALEC-hit}{76.5}
\resdef{n11-s030-star-stop-E05-ALEC-agree}{97.1}
\resdef{n11-s030-star-stop-E05-ALEC-early}{68.0}
\resdef{n11-s030-star-stop-E05-ALEC-exceed}{12.1}
\resdef{n11-s030-star-stop-E05-ALEC-qmax}{44}
\resdef{n11-s030-star-stop-E05-ALEC-qtot}{30.7}
\resdef{n11-s030-star-stop-E05-ALEC-errearly}{5.7}
\resdef{n11-s030-star-stop-E05-ALEC-errci}{[4.2, 7.7]}
\resdef{n11-s030-star-stop-E05-ALEC-regret}{1.68}
\resdef{n11-s030-star-stop-E05-ALEC-qround}{20}
\resdef{n11-s030-star-stop-E05-ALEC-fixedhit}{76.0}
\resdef{n11-s030-star-stop-E05-ALEC-fixedreg}{1.61}
\resdef{n11-s030-star-stop-E05-ALEC-save2n}{11}
\resdef{n11-s030-star-stop-B10-ALEC-cap-q}{14.8}
\resdef{n11-s030-star-stop-B10-ALEC-cap-hit}{76.1}
\resdef{n11-s030-star-stop-B10-ALEC-cap-agree}{97.3}
\resdef{n11-s030-star-stop-B10-ALEC-cap-early}{50.7}
\resdef{n11-s030-star-stop-B10-ALEC-cap-exceed}{23.9}
\resdef{n11-s030-star-stop-B10-ALEC-cap-qmax}{22}
\resdef{n11-s030-star-stop-B10-ALEC-cap-qtot}{25.8}
\resdef{n11-s030-star-stop-B10-ALEC-cap-errearly}{8.5}
\resdef{n11-s030-star-stop-B10-ALEC-cap-errci}{[6.4, 11.2]}
\resdef{n11-s030-star-stop-B10-ALEC-cap-regret}{1.66}
\resdef{n11-s030-star-stop-B10-ALEC-cap-qround}{15}
\resdef{n11-s030-star-stop-B10-ALEC-cap-fixedhit}{75.3}
\resdef{n11-s030-star-stop-B10-ALEC-cap-fixedreg}{1.70}
\resdef{n11-s030-star-stop-B10-ALEC-cap-save2n}{33}
\resdef{n11-s030-star-stop-E02-ALEC-cap-q}{13.9}
\resdef{n11-s030-star-stop-E02-ALEC-cap-hit}{76.0}
\resdef{n11-s030-star-stop-E02-ALEC-cap-agree}{96.6}
\resdef{n11-s030-star-stop-E02-ALEC-cap-early}{56.1}
\resdef{n11-s030-star-stop-E02-ALEC-cap-exceed}{19.2}
\resdef{n11-s030-star-stop-E02-ALEC-cap-qmax}{22}
\resdef{n11-s030-star-stop-E02-ALEC-cap-qtot}{24.9}
\resdef{n11-s030-star-stop-E02-ALEC-cap-errearly}{8.0}
\resdef{n11-s030-star-stop-E02-ALEC-cap-errci}{[6.0, 10.6]}
\resdef{n11-s030-star-stop-E02-ALEC-cap-regret}{1.67}
\resdef{n11-s030-star-stop-E02-ALEC-cap-qround}{14}
\resdef{n11-s030-star-stop-E02-ALEC-cap-fixedhit}{74.2}
\resdef{n11-s030-star-stop-E02-ALEC-cap-fixedreg}{1.93}
\resdef{n11-s030-star-stop-E02-ALEC-cap-save2n}{37}
\resdef{n11-s030-star-stop-E05-ALEC-cap-q}{12.5}
\resdef{n11-s030-star-stop-E05-ALEC-cap-hit}{75.9}
\resdef{n11-s030-star-stop-E05-ALEC-cap-agree}{95.7}
\resdef{n11-s030-star-stop-E05-ALEC-cap-early}{65.3}
\resdef{n11-s030-star-stop-E05-ALEC-cap-exceed}{12.9}
\resdef{n11-s030-star-stop-E05-ALEC-cap-qmax}{22}
\resdef{n11-s030-star-stop-E05-ALEC-cap-qtot}{23.5}
\resdef{n11-s030-star-stop-E05-ALEC-cap-errearly}{5.8}
\resdef{n11-s030-star-stop-E05-ALEC-cap-errci}{[4.3, 7.9]}
\resdef{n11-s030-star-stop-E05-ALEC-cap-regret}{1.75}
\resdef{n11-s030-star-stop-E05-ALEC-cap-qround}{12}
\resdef{n11-s030-star-stop-E05-ALEC-cap-fixedhit}{73.2}
\resdef{n11-s030-star-stop-E05-ALEC-cap-fixedreg}{2.13}
\resdef{n11-s030-star-stop-E05-ALEC-cap-save2n}{43}
\resdef{n11-s030-star-stop-E05-ALEC-cap-s06-q}{16.4}
\resdef{n11-s030-star-stop-E05-ALEC-cap-s06-hit}{76.1}
\resdef{n11-s030-star-stop-E05-ALEC-cap-s06-agree}{98.5}
\resdef{n11-s030-star-stop-E05-ALEC-cap-s06-early}{50.3}
\resdef{n11-s030-star-stop-E05-ALEC-cap-s06-exceed}{12.2}
\resdef{n11-s030-star-stop-E05-ALEC-cap-s06-qmax}{22}
\resdef{n11-s030-star-stop-E05-ALEC-cap-s06-qtot}{27.4}
\resdef{n11-s030-star-stop-E05-ALEC-cap-s06-errearly}{2.2}
\resdef{n11-s030-star-stop-E05-ALEC-cap-s06-errci}{[1.2, 3.9]}
\resdef{n11-s030-star-stop-E05-ALEC-cap-s06-regret}{1.61}
\resdef{n11-s030-star-stop-E05-ALEC-cap-s06-qround}{16}
\resdef{n11-s030-star-stop-E05-ALEC-cap-s06-fixedhit}{75.7}
\resdef{n11-s030-star-stop-E05-ALEC-cap-s06-fixedreg}{1.61}
\resdef{n11-s030-star-stop-E05-ALEC-cap-s06-save2n}{25}
\resdef{count-RANDOM-h-win}{9}
\resdef{count-RANDOM-h-lose}{0}
\resdef{count-RANDOM-h-tie}{0}
\resdef{count-RANDOM-h-holm}{9}
\resdef{count-RANDOM-h-holmlose}{0}
\resdef{count-RANDOM-h-of}{9}
\resdef{ci-RANDOM-h-lo-min}{0.5}
\resdef{ci-RANDOM-h-hi-max}{14.8}
\resdef{gap-RANDOM-h-min}{$3.3$}
\resdef{gap-RANDOM-h-max}{$11.3$}
\resdef{count-RANDOM-h-llsm-win}{9}
\resdef{count-RANDOM-h-llsm-lose}{0}
\resdef{count-RANDOM-h-llsm-tie}{0}
\resdef{count-RANDOM-h-llsm-holm}{9}
\resdef{count-RANDOM-h-llsm-holmlose}{0}
\resdef{count-RANDOM-h-llsm-of}{9}
\resdef{ci-RANDOM-h-llsm-lo-min}{0.5}
\resdef{ci-RANDOM-h-llsm-hi-max}{14.8}
\resdef{gap-RANDOM-h-llsm-min}{$3.3$}
\resdef{gap-RANDOM-h-llsm-max}{$11.3$}
\resdef{count-RANDOM-n-win}{9}
\resdef{count-RANDOM-n-lose}{0}
\resdef{count-RANDOM-n-tie}{0}
\resdef{count-RANDOM-n-holm}{9}
\resdef{count-RANDOM-n-holmlose}{0}
\resdef{count-RANDOM-n-of}{9}
\resdef{ci-RANDOM-n-lo-min}{1.3}
\resdef{ci-RANDOM-n-hi-max}{15.7}
\resdef{gap-RANDOM-n-min}{$4.0$}
\resdef{gap-RANDOM-n-max}{$12.2$}
\resdef{count-RANDOM-n-llsm-win}{9}
\resdef{count-RANDOM-n-llsm-lose}{0}
\resdef{count-RANDOM-n-llsm-tie}{0}
\resdef{count-RANDOM-n-llsm-holm}{9}
\resdef{count-RANDOM-n-llsm-holmlose}{0}
\resdef{count-RANDOM-n-llsm-of}{9}
\resdef{ci-RANDOM-n-llsm-lo-min}{1.3}
\resdef{ci-RANDOM-n-llsm-hi-max}{15.7}
\resdef{gap-RANDOM-n-llsm-min}{$4.0$}
\resdef{gap-RANDOM-n-llsm-max}{$12.2$}
\resdef{count-RANDOM-2n-win}{6}
\resdef{count-RANDOM-2n-lose}{0}
\resdef{count-RANDOM-2n-tie}{3}
\resdef{count-RANDOM-2n-holm}{6}
\resdef{count-RANDOM-2n-holmlose}{0}
\resdef{count-RANDOM-2n-of}{9}
\resdef{ci-RANDOM-2n-lo-min}{-0.4}
\resdef{ci-RANDOM-2n-hi-max}{15.9}
\resdef{gap-RANDOM-2n-min}{$0.0$}
\resdef{gap-RANDOM-2n-max}{$12.7$}
\resdef{count-RANDOM-2n-llsm-win}{6}
\resdef{count-RANDOM-2n-llsm-lose}{0}
\resdef{count-RANDOM-2n-llsm-tie}{3}
\resdef{count-RANDOM-2n-llsm-holm}{6}
\resdef{count-RANDOM-2n-llsm-holmlose}{0}
\resdef{count-RANDOM-2n-llsm-of}{9}
\resdef{ci-RANDOM-2n-llsm-lo-min}{-0.4}
\resdef{ci-RANDOM-2n-llsm-hi-max}{15.9}
\resdef{gap-RANDOM-2n-llsm-min}{$0.0$}
\resdef{gap-RANDOM-2n-llsm-max}{$12.7$}
\resdef{count-HARKER-PG-h-win}{6}
\resdef{count-HARKER-PG-h-lose}{0}
\resdef{count-HARKER-PG-h-tie}{3}
\resdef{count-HARKER-PG-h-holm}{6}
\resdef{count-HARKER-PG-h-holmlose}{0}
\resdef{count-HARKER-PG-h-of}{9}
\resdef{ci-HARKER-PG-h-lo-min}{-1.2}
\resdef{ci-HARKER-PG-h-hi-max}{15.5}
\resdef{gap-HARKER-PG-h-min}{$0.7$}
\resdef{gap-HARKER-PG-h-max}{$11.9$}
\resdef{count-HARKER-PG-h-llsm-win}{6}
\resdef{count-HARKER-PG-h-llsm-lose}{0}
\resdef{count-HARKER-PG-h-llsm-tie}{3}
\resdef{count-HARKER-PG-h-llsm-holm}{6}
\resdef{count-HARKER-PG-h-llsm-holmlose}{0}
\resdef{count-HARKER-PG-h-llsm-of}{9}
\resdef{ci-HARKER-PG-h-llsm-lo-min}{-1.1}
\resdef{ci-HARKER-PG-h-llsm-hi-max}{11.9}
\resdef{gap-HARKER-PG-h-llsm-min}{$1.3$}
\resdef{gap-HARKER-PG-h-llsm-max}{$8.6$}
\resdef{count-HARKER-PG-n-win}{6}
\resdef{count-HARKER-PG-n-lose}{0}
\resdef{count-HARKER-PG-n-tie}{3}
\resdef{count-HARKER-PG-n-holm}{5}
\resdef{count-HARKER-PG-n-holmlose}{0}
\resdef{count-HARKER-PG-n-of}{9}
\resdef{ci-HARKER-PG-n-lo-min}{-0.2}
\resdef{ci-HARKER-PG-n-hi-max}{10.1}
\resdef{gap-HARKER-PG-n-min}{$1.2$}
\resdef{gap-HARKER-PG-n-max}{$6.4$}
\resdef{count-HARKER-PG-n-llsm-win}{5}
\resdef{count-HARKER-PG-n-llsm-lose}{0}
\resdef{count-HARKER-PG-n-llsm-tie}{4}
\resdef{count-HARKER-PG-n-llsm-holm}{3}
\resdef{count-HARKER-PG-n-llsm-holmlose}{0}
\resdef{count-HARKER-PG-n-llsm-of}{9}
\resdef{ci-HARKER-PG-n-llsm-lo-min}{-0.4}
\resdef{ci-HARKER-PG-n-llsm-hi-max}{9.3}
\resdef{gap-HARKER-PG-n-llsm-min}{$1.1$}
\resdef{gap-HARKER-PG-n-llsm-max}{$5.6$}
\resdef{count-HARKER-PG-2n-win}{0}
\resdef{count-HARKER-PG-2n-lose}{0}
\resdef{count-HARKER-PG-2n-tie}{3}
\resdef{count-HARKER-PG-2n-holm}{0}
\resdef{count-HARKER-PG-2n-holmlose}{0}
\resdef{count-HARKER-PG-2n-of}{3}
\resdef{ci-HARKER-PG-2n-lo-min}{-1.9}
\resdef{ci-HARKER-PG-2n-hi-max}{1.3}
\resdef{gap-HARKER-PG-2n-min}{$-0.9$}
\resdef{gap-HARKER-PG-2n-max}{$0.0$}
\resdef{count-HARKER-PG-2n-llsm-win}{0}
\resdef{count-HARKER-PG-2n-llsm-lose}{0}
\resdef{count-HARKER-PG-2n-llsm-tie}{3}
\resdef{count-HARKER-PG-2n-llsm-holm}{0}
\resdef{count-HARKER-PG-2n-llsm-holmlose}{0}
\resdef{count-HARKER-PG-2n-llsm-of}{3}
\resdef{ci-HARKER-PG-2n-llsm-lo-min}{-0.4}
\resdef{ci-HARKER-PG-2n-llsm-hi-max}{0.4}
\resdef{gap-HARKER-PG-2n-llsm-min}{$0.0$}
\resdef{gap-HARKER-PG-2n-llsm-max}{$0.0$}
\resdef{count-HARKER-h-win}{7}
\resdef{count-HARKER-h-lose}{0}
\resdef{count-HARKER-h-tie}{2}
\resdef{count-HARKER-h-holm}{6}
\resdef{count-HARKER-h-holmlose}{0}
\resdef{count-HARKER-h-of}{9}
\resdef{ci-HARKER-h-lo-min}{-2.6}
\resdef{ci-HARKER-h-hi-max}{11.2}
\resdef{gap-HARKER-h-min}{$-0.2$}
\resdef{gap-HARKER-h-max}{$7.9$}
\resdef{count-HARKER-h-llsm-win}{6}
\resdef{count-HARKER-h-llsm-lose}{0}
\resdef{count-HARKER-h-llsm-tie}{3}
\resdef{count-HARKER-h-llsm-holm}{6}
\resdef{count-HARKER-h-llsm-holmlose}{0}
\resdef{count-HARKER-h-llsm-of}{9}
\resdef{ci-HARKER-h-llsm-lo-min}{-2.9}
\resdef{ci-HARKER-h-llsm-hi-max}{12.3}
\resdef{gap-HARKER-h-llsm-min}{$-0.5$}
\resdef{gap-HARKER-h-llsm-max}{$9.0$}
\resdef{count-HARKER-n-win}{8}
\resdef{count-HARKER-n-lose}{0}
\resdef{count-HARKER-n-tie}{1}
\resdef{count-HARKER-n-holm}{2}
\resdef{count-HARKER-n-holmlose}{0}
\resdef{count-HARKER-n-of}{9}
\resdef{ci-HARKER-n-lo-min}{-0.0}
\resdef{ci-HARKER-n-hi-max}{9.3}
\resdef{gap-HARKER-n-min}{$2.0$}
\resdef{gap-HARKER-n-max}{$6.1$}
\resdef{count-HARKER-n-llsm-win}{8}
\resdef{count-HARKER-n-llsm-lose}{0}
\resdef{count-HARKER-n-llsm-tie}{1}
\resdef{count-HARKER-n-llsm-holm}{2}
\resdef{count-HARKER-n-llsm-holmlose}{0}
\resdef{count-HARKER-n-llsm-of}{9}
\resdef{ci-HARKER-n-llsm-lo-min}{-0.9}
\resdef{ci-HARKER-n-llsm-hi-max}{8.7}
\resdef{gap-HARKER-n-llsm-min}{$1.5$}
\resdef{gap-HARKER-n-llsm-max}{$4.9$}
\resdef{count-HARKER-2n-win}{5}
\resdef{count-HARKER-2n-lose}{0}
\resdef{count-HARKER-2n-tie}{4}
\resdef{count-HARKER-2n-holm}{5}
\resdef{count-HARKER-2n-holmlose}{0}
\resdef{count-HARKER-2n-of}{9}
\resdef{ci-HARKER-2n-lo-min}{-1.9}
\resdef{ci-HARKER-2n-hi-max}{12.0}
\resdef{gap-HARKER-2n-min}{$-0.9$}
\resdef{gap-HARKER-2n-max}{$8.9$}
\resdef{count-HARKER-2n-llsm-win}{5}
\resdef{count-HARKER-2n-llsm-lose}{0}
\resdef{count-HARKER-2n-llsm-tie}{4}
\resdef{count-HARKER-2n-llsm-holm}{5}
\resdef{count-HARKER-2n-llsm-holmlose}{0}
\resdef{count-HARKER-2n-llsm-of}{9}
\resdef{ci-HARKER-2n-llsm-lo-min}{-0.5}
\resdef{ci-HARKER-2n-llsm-hi-max}{11.1}
\resdef{gap-HARKER-2n-llsm-min}{$0.0$}
\resdef{gap-HARKER-2n-llsm-max}{$8.1$}
\resdef{count-HARKER-LOG-h-win}{1}
\resdef{count-HARKER-LOG-h-lose}{0}
\resdef{count-HARKER-LOG-h-tie}{8}
\resdef{count-HARKER-LOG-h-holm}{0}
\resdef{count-HARKER-LOG-h-holmlose}{0}
\resdef{count-HARKER-LOG-h-of}{9}
\resdef{ci-HARKER-LOG-h-lo-min}{-2.8}
\resdef{ci-HARKER-LOG-h-hi-max}{5.4}
\resdef{gap-HARKER-LOG-h-min}{$-1.1$}
\resdef{gap-HARKER-LOG-h-max}{$2.8$}
\resdef{count-HARKER-LOG-h-llsm-win}{0}
\resdef{count-HARKER-LOG-h-llsm-lose}{0}
\resdef{count-HARKER-LOG-h-llsm-tie}{9}
\resdef{count-HARKER-LOG-h-llsm-holm}{0}
\resdef{count-HARKER-LOG-h-llsm-holmlose}{0}
\resdef{count-HARKER-LOG-h-llsm-of}{9}
\resdef{ci-HARKER-LOG-h-llsm-lo-min}{-3.4}
\resdef{ci-HARKER-LOG-h-llsm-hi-max}{4.8}
\resdef{gap-HARKER-LOG-h-llsm-min}{$-1.5$}
\resdef{gap-HARKER-LOG-h-llsm-max}{$2.3$}
\resdef{count-HARKER-LOG-n-win}{0}
\resdef{count-HARKER-LOG-n-lose}{0}
\resdef{count-HARKER-LOG-n-tie}{9}
\resdef{count-HARKER-LOG-n-holm}{0}
\resdef{count-HARKER-LOG-n-holmlose}{0}
\resdef{count-HARKER-LOG-n-of}{9}
\resdef{ci-HARKER-LOG-n-lo-min}{-5.0}
\resdef{ci-HARKER-LOG-n-hi-max}{5.7}
\resdef{gap-HARKER-LOG-n-min}{$-2.5$}
\resdef{gap-HARKER-LOG-n-max}{$2.6$}
\resdef{count-HARKER-LOG-n-llsm-win}{0}
\resdef{count-HARKER-LOG-n-llsm-lose}{2}
\resdef{count-HARKER-LOG-n-llsm-tie}{7}
\resdef{count-HARKER-LOG-n-llsm-holm}{0}
\resdef{count-HARKER-LOG-n-llsm-holmlose}{0}
\resdef{count-HARKER-LOG-n-llsm-of}{9}
\resdef{ci-HARKER-LOG-n-llsm-lo-min}{-6.1}
\resdef{ci-HARKER-LOG-n-llsm-hi-max}{3.2}
\resdef{gap-HARKER-LOG-n-llsm-min}{$-3.5$}
\resdef{gap-HARKER-LOG-n-llsm-max}{$1.2$}
\resdef{count-HARKER-LOG-2n-win}{2}
\resdef{count-HARKER-LOG-2n-lose}{0}
\resdef{count-HARKER-LOG-2n-tie}{7}
\resdef{count-HARKER-LOG-2n-holm}{2}
\resdef{count-HARKER-LOG-2n-holmlose}{0}
\resdef{count-HARKER-LOG-2n-of}{9}
\resdef{ci-HARKER-LOG-2n-lo-min}{-1.9}
\resdef{ci-HARKER-LOG-2n-hi-max}{8.4}
\resdef{gap-HARKER-LOG-2n-min}{$-0.9$}
\resdef{gap-HARKER-LOG-2n-max}{$5.7$}
\resdef{count-HARKER-LOG-2n-llsm-win}{2}
\resdef{count-HARKER-LOG-2n-llsm-lose}{0}
\resdef{count-HARKER-LOG-2n-llsm-tie}{7}
\resdef{count-HARKER-LOG-2n-llsm-holm}{2}
\resdef{count-HARKER-LOG-2n-llsm-holmlose}{0}
\resdef{count-HARKER-LOG-2n-llsm-of}{9}
\resdef{ci-HARKER-LOG-2n-llsm-lo-min}{-1.3}
\resdef{ci-HARKER-LOG-2n-llsm-hi-max}{7.8}
\resdef{gap-HARKER-LOG-2n-llsm-min}{$0.0$}
\resdef{gap-HARKER-LOG-2n-llsm-max}{$5.3$}
\resdef{count-AOPT-h-win}{6}
\resdef{count-AOPT-h-lose}{0}
\resdef{count-AOPT-h-tie}{3}
\resdef{count-AOPT-h-holm}{6}
\resdef{count-AOPT-h-holmlose}{0}
\resdef{count-AOPT-h-of}{9}
\resdef{ci-AOPT-h-lo-min}{-1.0}
\resdef{ci-AOPT-h-hi-max}{12.4}
\resdef{gap-AOPT-h-min}{$1.1$}
\resdef{gap-AOPT-h-max}{$9.1$}
\resdef{count-AOPT-h-llsm-win}{6}
\resdef{count-AOPT-h-llsm-lose}{0}
\resdef{count-AOPT-h-llsm-tie}{3}
\resdef{count-AOPT-h-llsm-holm}{6}
\resdef{count-AOPT-h-llsm-holmlose}{0}
\resdef{count-AOPT-h-llsm-of}{9}
\resdef{ci-AOPT-h-llsm-lo-min}{-1.0}
\resdef{ci-AOPT-h-llsm-hi-max}{12.4}
\resdef{gap-AOPT-h-llsm-min}{$1.1$}
\resdef{gap-AOPT-h-llsm-max}{$9.1$}
\resdef{count-AOPT-n-win}{9}
\resdef{count-AOPT-n-lose}{0}
\resdef{count-AOPT-n-tie}{0}
\resdef{count-AOPT-n-holm}{6}
\resdef{count-AOPT-n-holmlose}{0}
\resdef{count-AOPT-n-of}{9}
\resdef{ci-AOPT-n-lo-min}{0.4}
\resdef{ci-AOPT-n-hi-max}{11.6}
\resdef{gap-AOPT-n-min}{$2.0$}
\resdef{gap-AOPT-n-max}{$8.1$}
\resdef{count-AOPT-n-llsm-win}{9}
\resdef{count-AOPT-n-llsm-lose}{0}
\resdef{count-AOPT-n-llsm-tie}{0}
\resdef{count-AOPT-n-llsm-holm}{6}
\resdef{count-AOPT-n-llsm-holmlose}{0}
\resdef{count-AOPT-n-llsm-of}{9}
\resdef{ci-AOPT-n-llsm-lo-min}{0.4}
\resdef{ci-AOPT-n-llsm-hi-max}{11.6}
\resdef{gap-AOPT-n-llsm-min}{$2.0$}
\resdef{gap-AOPT-n-llsm-max}{$8.1$}
\resdef{count-AOPT-2n-win}{6}
\resdef{count-AOPT-2n-lose}{0}
\resdef{count-AOPT-2n-tie}{3}
\resdef{count-AOPT-2n-holm}{5}
\resdef{count-AOPT-2n-holmlose}{0}
\resdef{count-AOPT-2n-of}{9}
\resdef{ci-AOPT-2n-lo-min}{-0.4}
\resdef{ci-AOPT-2n-hi-max}{14.0}
\resdef{gap-AOPT-2n-min}{$0.0$}
\resdef{gap-AOPT-2n-max}{$10.9$}
\resdef{count-AOPT-2n-llsm-win}{6}
\resdef{count-AOPT-2n-llsm-lose}{0}
\resdef{count-AOPT-2n-llsm-tie}{3}
\resdef{count-AOPT-2n-llsm-holm}{5}
\resdef{count-AOPT-2n-llsm-holmlose}{0}
\resdef{count-AOPT-2n-llsm-of}{9}
\resdef{ci-AOPT-2n-llsm-lo-min}{-0.4}
\resdef{ci-AOPT-2n-llsm-hi-max}{14.0}
\resdef{gap-AOPT-2n-llsm-min}{$0.0$}
\resdef{gap-AOPT-2n-llsm-max}{$10.9$}
\resdef{count-FROZEN-h-win}{5}
\resdef{count-FROZEN-h-lose}{0}
\resdef{count-FROZEN-h-tie}{4}
\resdef{count-FROZEN-h-holm}{5}
\resdef{count-FROZEN-h-holmlose}{0}
\resdef{count-FROZEN-h-of}{9}
\resdef{ci-FROZEN-h-lo-min}{-1.9}
\resdef{ci-FROZEN-h-hi-max}{9.0}
\resdef{gap-FROZEN-h-min}{$-0.4$}
\resdef{gap-FROZEN-h-max}{$6.6$}
\resdef{count-FROZEN-h-llsm-win}{5}
\resdef{count-FROZEN-h-llsm-lose}{0}
\resdef{count-FROZEN-h-llsm-tie}{4}
\resdef{count-FROZEN-h-llsm-holm}{5}
\resdef{count-FROZEN-h-llsm-holmlose}{0}
\resdef{count-FROZEN-h-llsm-of}{9}
\resdef{ci-FROZEN-h-llsm-lo-min}{-1.9}
\resdef{ci-FROZEN-h-llsm-hi-max}{9.0}
\resdef{gap-FROZEN-h-llsm-min}{$-0.4$}
\resdef{gap-FROZEN-h-llsm-max}{$6.6$}
\resdef{count-FROZEN-n-win}{1}
\resdef{count-FROZEN-n-lose}{0}
\resdef{count-FROZEN-n-tie}{8}
\resdef{count-FROZEN-n-holm}{1}
\resdef{count-FROZEN-n-holmlose}{0}
\resdef{count-FROZEN-n-of}{9}
\resdef{ci-FROZEN-n-lo-min}{-1.5}
\resdef{ci-FROZEN-n-hi-max}{9.3}
\resdef{gap-FROZEN-n-min}{$-0.3$}
\resdef{gap-FROZEN-n-max}{$6.1$}
\resdef{count-FROZEN-n-llsm-win}{1}
\resdef{count-FROZEN-n-llsm-lose}{0}
\resdef{count-FROZEN-n-llsm-tie}{8}
\resdef{count-FROZEN-n-llsm-holm}{1}
\resdef{count-FROZEN-n-llsm-holmlose}{0}
\resdef{count-FROZEN-n-llsm-of}{9}
\resdef{ci-FROZEN-n-llsm-lo-min}{-1.5}
\resdef{ci-FROZEN-n-llsm-hi-max}{9.3}
\resdef{gap-FROZEN-n-llsm-min}{$-0.3$}
\resdef{gap-FROZEN-n-llsm-max}{$6.1$}
\resdef{count-FROZEN-2n-win}{3}
\resdef{count-FROZEN-2n-lose}{0}
\resdef{count-FROZEN-2n-tie}{6}
\resdef{count-FROZEN-2n-holm}{3}
\resdef{count-FROZEN-2n-holmlose}{0}
\resdef{count-FROZEN-2n-of}{9}
\resdef{ci-FROZEN-2n-lo-min}{-1.4}
\resdef{ci-FROZEN-2n-hi-max}{7.6}
\resdef{gap-FROZEN-2n-min}{$-0.1$}
\resdef{gap-FROZEN-2n-max}{$5.1$}
\resdef{count-FROZEN-2n-llsm-win}{3}
\resdef{count-FROZEN-2n-llsm-lose}{0}
\resdef{count-FROZEN-2n-llsm-tie}{6}
\resdef{count-FROZEN-2n-llsm-holm}{3}
\resdef{count-FROZEN-2n-llsm-holmlose}{0}
\resdef{count-FROZEN-2n-llsm-of}{9}
\resdef{ci-FROZEN-2n-llsm-lo-min}{-1.4}
\resdef{ci-FROZEN-2n-llsm-hi-max}{7.6}
\resdef{gap-FROZEN-2n-llsm-min}{$-0.1$}
\resdef{gap-FROZEN-2n-llsm-max}{$5.1$}
\resdef{count-STAR-h-win}{5}
\resdef{count-STAR-h-lose}{0}
\resdef{count-STAR-h-tie}{4}
\resdef{count-STAR-h-holm}{5}
\resdef{count-STAR-h-holmlose}{0}
\resdef{count-STAR-h-of}{9}
\resdef{ci-STAR-h-lo-min}{-1.8}
\resdef{ci-STAR-h-hi-max}{8.9}
\resdef{gap-STAR-h-min}{$-0.2$}
\resdef{gap-STAR-h-max}{$6.5$}
\resdef{count-STAR-h-llsm-win}{5}
\resdef{count-STAR-h-llsm-lose}{0}
\resdef{count-STAR-h-llsm-tie}{4}
\resdef{count-STAR-h-llsm-holm}{5}
\resdef{count-STAR-h-llsm-holmlose}{0}
\resdef{count-STAR-h-llsm-of}{9}
\resdef{ci-STAR-h-llsm-lo-min}{-1.8}
\resdef{ci-STAR-h-llsm-hi-max}{8.9}
\resdef{gap-STAR-h-llsm-min}{$-0.2$}
\resdef{gap-STAR-h-llsm-max}{$6.5$}
\resdef{count-STAR-n-win}{0}
\resdef{count-STAR-n-lose}{0}
\resdef{count-STAR-n-tie}{9}
\resdef{count-STAR-n-holm}{0}
\resdef{count-STAR-n-holmlose}{0}
\resdef{count-STAR-n-of}{9}
\resdef{ci-STAR-n-lo-min}{-3.6}
\resdef{ci-STAR-n-hi-max}{3.6}
\resdef{gap-STAR-n-min}{$-1.4$}
\resdef{gap-STAR-n-max}{$1.0$}
\resdef{count-STAR-n-llsm-win}{0}
\resdef{count-STAR-n-llsm-lose}{0}
\resdef{count-STAR-n-llsm-tie}{9}
\resdef{count-STAR-n-llsm-holm}{0}
\resdef{count-STAR-n-llsm-holmlose}{0}
\resdef{count-STAR-n-llsm-of}{9}
\resdef{ci-STAR-n-llsm-lo-min}{-3.6}
\resdef{ci-STAR-n-llsm-hi-max}{3.6}
\resdef{gap-STAR-n-llsm-min}{$-1.4$}
\resdef{gap-STAR-n-llsm-max}{$1.0$}
\resdef{count-STAR-2n-win}{0}
\resdef{count-STAR-2n-lose}{0}
\resdef{count-STAR-2n-tie}{9}
\resdef{count-STAR-2n-holm}{0}
\resdef{count-STAR-2n-holmlose}{0}
\resdef{count-STAR-2n-of}{9}
\resdef{ci-STAR-2n-lo-min}{-1.0}
\resdef{ci-STAR-2n-hi-max}{3.3}
\resdef{gap-STAR-2n-min}{$0.0$}
\resdef{gap-STAR-2n-max}{$1.5$}
\resdef{count-STAR-2n-llsm-win}{0}
\resdef{count-STAR-2n-llsm-lose}{0}
\resdef{count-STAR-2n-llsm-tie}{9}
\resdef{count-STAR-2n-llsm-holm}{0}
\resdef{count-STAR-2n-llsm-holmlose}{0}
\resdef{count-STAR-2n-llsm-of}{9}
\resdef{ci-STAR-2n-llsm-lo-min}{-1.0}
\resdef{ci-STAR-2n-llsm-hi-max}{3.3}
\resdef{gap-STAR-2n-llsm-min}{$0.0$}
\resdef{gap-STAR-2n-llsm-max}{$1.5$}
\resdef{count-ALEC-U-h-win}{0}
\resdef{count-ALEC-U-h-lose}{0}
\resdef{count-ALEC-U-h-tie}{9}
\resdef{count-ALEC-U-h-holm}{0}
\resdef{count-ALEC-U-h-holmlose}{0}
\resdef{count-ALEC-U-h-of}{9}
\resdef{ci-ALEC-U-h-lo-min}{-4.5}
\resdef{ci-ALEC-U-h-hi-max}{4.8}
\resdef{gap-ALEC-U-h-min}{$-2.0$}
\resdef{gap-ALEC-U-h-max}{$2.4$}
\resdef{count-ALEC-U-h-llsm-win}{0}
\resdef{count-ALEC-U-h-llsm-lose}{0}
\resdef{count-ALEC-U-h-llsm-tie}{9}
\resdef{count-ALEC-U-h-llsm-holm}{0}
\resdef{count-ALEC-U-h-llsm-holmlose}{0}
\resdef{count-ALEC-U-h-llsm-of}{9}
\resdef{ci-ALEC-U-h-llsm-lo-min}{-4.5}
\resdef{ci-ALEC-U-h-llsm-hi-max}{4.8}
\resdef{gap-ALEC-U-h-llsm-min}{$-2.0$}
\resdef{gap-ALEC-U-h-llsm-max}{$2.4$}
\resdef{count-ALEC-U-n-win}{0}
\resdef{count-ALEC-U-n-lose}{0}
\resdef{count-ALEC-U-n-tie}{9}
\resdef{count-ALEC-U-n-holm}{0}
\resdef{count-ALEC-U-n-holmlose}{0}
\resdef{count-ALEC-U-n-of}{9}
\resdef{ci-ALEC-U-n-lo-min}{-3.9}
\resdef{ci-ALEC-U-n-hi-max}{2.9}
\resdef{gap-ALEC-U-n-min}{$-1.5$}
\resdef{gap-ALEC-U-n-max}{$1.2$}
\resdef{count-ALEC-U-n-llsm-win}{0}
\resdef{count-ALEC-U-n-llsm-lose}{0}
\resdef{count-ALEC-U-n-llsm-tie}{9}
\resdef{count-ALEC-U-n-llsm-holm}{0}
\resdef{count-ALEC-U-n-llsm-holmlose}{0}
\resdef{count-ALEC-U-n-llsm-of}{9}
\resdef{ci-ALEC-U-n-llsm-lo-min}{-3.9}
\resdef{ci-ALEC-U-n-llsm-hi-max}{2.9}
\resdef{gap-ALEC-U-n-llsm-min}{$-1.5$}
\resdef{gap-ALEC-U-n-llsm-max}{$1.2$}
\resdef{count-ALEC-U-2n-win}{0}
\resdef{count-ALEC-U-2n-lose}{0}
\resdef{count-ALEC-U-2n-tie}{9}
\resdef{count-ALEC-U-2n-holm}{0}
\resdef{count-ALEC-U-2n-holmlose}{0}
\resdef{count-ALEC-U-2n-of}{9}
\resdef{ci-ALEC-U-2n-lo-min}{-1.4}
\resdef{ci-ALEC-U-2n-hi-max}{4.5}
\resdef{gap-ALEC-U-2n-min}{$0.0$}
\resdef{gap-ALEC-U-2n-max}{$2.2$}
\resdef{count-ALEC-U-2n-llsm-win}{0}
\resdef{count-ALEC-U-2n-llsm-lose}{0}
\resdef{count-ALEC-U-2n-llsm-tie}{9}
\resdef{count-ALEC-U-2n-llsm-holm}{0}
\resdef{count-ALEC-U-2n-llsm-holmlose}{0}
\resdef{count-ALEC-U-2n-llsm-of}{9}
\resdef{ci-ALEC-U-2n-llsm-lo-min}{-1.4}
\resdef{ci-ALEC-U-2n-llsm-hi-max}{4.5}
\resdef{gap-ALEC-U-2n-llsm-min}{$0.0$}
\resdef{gap-ALEC-U-2n-llsm-max}{$2.2$}
\resdef{R-syn}{1000}
\resdef{n15-s015-budget-dfrozen-2n}{1.1}
\resdef{n15-s015-budget-dfrozen-max}{4.2}
\resdef{n15-s015-budget-hlog-lead-max}{2.3}
\resdef{n15-s030-budget-dfrozen-2n}{3.2}
\resdef{n15-s030-budget-dfrozen-max}{7.2}
\resdef{n15-s030-budget-hlog-lead-max}{2.5}
\resdef{n15-s050-budget-dfrozen-2n}{5.1}
\resdef{n15-s050-budget-dfrozen-max}{7.5}
\resdef{n15-s050-budget-hlog-lead-max}{0.0}
\resdef{budget-hlog-lead-max}{2.5}
\makeatletter
\def\xr@Mkeys{}\def\xr@Skeys{}
\newcommand{\xrdef}[3]{% #1 = M (article) or S (supplement), #2 = label, #3 = number
  \expandafter\def\csname xr@#1@#2\endcsname{#3}%
  \expandafter\g@addto@macro\csname xr@#1keys\endcsname{\xr@do{#1}{#2}}}
\newcommand{\xr@get}[2]{\ifcsname xr@#1@#2\endcsname\csname xr@#1@#2\endcsname\else\textbf{\textcolor{red}{??#2??}}\fi}
\newcommand{\Mref}[1]{\xr@get{M}{#1}}
\newcommand{\Sref}[1]{\xr@get{S}{#1}}
\newcommand{\xrverify}[1]{% compare the numbers of document #1 (M or S) with its own labels
  \AtEndDocument{%
    \def\xr@do##1##2{%
      \@ifundefined{r@##2}{\PackageWarningNoLine{crossrefs}{Label ##2 not found (yet)}}{%
        \protected@edef\xr@have{\expandafter\expandafter\expandafter\@car\csname r@##2\endcsname\@nil}%
        \protected@edef\xr@want{\csname xr@##1@##2\endcsname}%
        \ifx\xr@have\xr@want\else
          \PackageWarningNoLine{crossrefs}{##2 is numbered \xr@have\space but crossrefs.tex says \xr@want}%
        \fi}}%
    \csname xr@#1keys\endcsname}}
\makeatother

\xrdef{M}{sec:posterior}{4.2}
\xrdef{M}{sec:stopping}{4.4}
\xrdef{M}{sec:benchmarks}{4.6}
\xrdef{M}{sec:example}{5}
\xrdef{M}{sec:design}{6.1}
\xrdef{M}{sec:ref-results}{6.2}
\xrdef{M}{sec:fixed}{6.3}
\xrdef{M}{sec:stop-results}{6.4}
\xrdef{M}{sec:robust}{6.5}
\xrdef{M}{prop:resistance}{1}
\xrdef{M}{prop:direct}{2}
\xrdef{M}{cor:eps}{1}
\xrdef{M}{eq:rankone}{3}
\xrdef{M}{eq:betaeps}{5}
\xrdef{M}{prop:posterior}{3}
\xrdef{M}{prop:winner}{4}
\xrdef{M}{eq:beta}{4}
\xrdef{M}{tab:example}{1}
\xrdef{M}{tab:hit-n}{2}
\xrdef{M}{tab:protocol}{3}
\xrdef{M}{tab:samestop}{4}

\xrdef{S}{sec:S-example}{S1}
\xrdef{S}{sec:S-ref}{S2}
\xrdef{S}{sec:S-fixed}{S3}
\xrdef{S}{sec:S-stop}{S4}
\xrdef{S}{sec:S-calib}{S5}
\xrdef{S}{sec:S-robust}{S6}
\xrdef{S}{fig:example}{S1}
\xrdef{S}{fig:curves}{S2}
\xrdef{S}{fig:stopping}{S3}
\xrdef{S}{tab:esmaeili}{S1}
\xrdef{S}{tab:beverage}{S2}
\xrdef{S}{tab:hit-h}{S3}
\xrdef{S}{tab:hit-2n}{S4}
\xrdef{S}{tab:hit-llsm}{S5}
\xrdef{S}{tab:hit-llsm-2n}{S6}
\xrdef{S}{tab:regret}{S7}
\xrdef{S}{tab:tau}{S8}
\xrdef{S}{tab:ci-hlog}{S9}
\xrdef{S}{tab:timing}{S10}
\xrdef{S}{tab:gap}{S11}
\xrdef{S}{tab:savings}{S12}
\xrdef{S}{tab:stop}{S13}
\xrdef{S}{tab:stop015}{S14}
\xrdef{S}{tab:stop050}{S15}
\xrdef{S}{tab:winner}{S16}
\xrdef{S}{tab:eps}{S17}
\xrdef{S}{tab:protocol-full}{S18}
\xrdef{S}{tab:fixed-full}{S19}
\xrdef{S}{tab:samestop-full}{S20}
\xrdef{S}{tab:bayes}{S21}
\xrdef{S}{tab:calib-alpha}{S22}
\xrdef{S}{tab:calib-sigma}{S23}
\xrdef{S}{tab:robust}{S24}
\xrdef{S}{tab:saaty-stop}{S25}
\xrdef{S}{sec:S-technical}{S7}   % numbers of the items of the supplementary material (\Sref)
\xrverify{M}

\newtheorem{proposition}{Proposition}

\newtheorem{corollary}{Corollary}
\theoremstyle{definition}

\theoremstyle{remark}
\newtheorem{remark}{Remark}

\newcommand{\R}{\mathbb{R}}
\newcommand{\E}{\mathcal{E}}
\newcommand{\G}{\mathcal{G}}
\newcommand{\one}{\mathbf{1}}
\newcommand{\xhat}{\hat{\bm{x}}}
\newcommand{\Lp}{L^{+}}
\newcommand{\Reff}{R}
\newcommand{\coderepo}{\href{https://github.com/LuisAngelCalvoPascual/ALEC-code-and-replication-repository}{ALEC code and replication repository}}
\DeclareMathOperator*{\argmax}{arg\,max}

\journal{Decision Support Systems}

\begin{document}
%\linenumbers   % uncomment for submission if the journal asks for line numbers

\begin{frontmatter}

\title{Adaptive pairwise preference elicitation for decision support: identifying the best alternative with fewer judgments}

%% ------------------------------------------------------------------ authors
%% Affiliations taken from public profiles; authors to confirm departments.
%% L. A. Calvo Pascual is the first and corresponding author (only the corresponding author's e-mail is listed).
\author[a]{Luis \'Angel Calvo Pascual\corref{cor1}}
\ead{luisangelcalvopascual@gmail.com}
\cortext[cor1]{Corresponding author.}
\author[b]{Ji\v{r}\'{i} Mazurek}
\affiliation[a]{organization={Universidad Pontificia Comillas}, city={Madrid}, country={Spain}}
\affiliation[b]{organization={School of Business Administration in Karviná, Silesian University in Opava}, city={Karviná}, country={Czech Republic}}

\begin{abstract}
Interactive decision support systems often ask users for repeated preference judgments. When only the best alternative is needed, eliciting all $n(n-1)/2$ pairwise comparisons imposes cognitive and interaction costs that the decision may not justify. We propose ALEC (Adaptive Lexicographic Elicitation of Comparisons), an elicitation mechanism that re-ranks the alternatives after each judgment and asks the first unobserved pair in the lexicographic order of the current ranking, with the Analytic Hierarchy Process as the concrete setting. Under a log-normal error model, the direct comparison of a contrast is the single question that most reduces its effective resistance, and posterior bounds on the probability of a wrong choice or of exceeding a relative-loss tolerance yield stopping rules that remain valid under adaptive questioning; a consistency-and-margin result gives a deterministic certificate of agreement with the full-matrix winner. In simulations, re-ranking is most useful in larger problems with restricted budgets, and ALEC is significantly more accurate than Harker's original rule in most scenarios with $n\ge11$ while remaining competitive with a stronger variant. With a 5\% loss tolerance and a cap of $2n$ additional comparisons, the stopping rule saved \res{prot-s06-saveadd-min}--\res{prot-s06-saveadd-max}\% of them; among runs stopped by the bound, the observed rate of losses above the tolerance was at most \res{s06-errearly-max}\%, against a nominal 10\%. ALEC can therefore serve as an adaptive elicitation module for decision support systems, selecting the next query and stopping when further user effort is not justified by the required decision precision.
\end{abstract}


\end{frontmatter}

%% ============================================================================
\section{Introduction}\label{sec:intro}
%% ============================================================================

Many decision support systems (DSSs) obtain their preference information interactively: the system asks, the decision maker answers, and the recommendation is built from the answers. Each judgment costs time and attention, and the effort a decision aid demands shapes how decision makers use it and how good the resulting decisions are \citep{Payne1993,ToddBenbasat1999}. Traditional preference elicitation is known to be time-consuming and tedious for users \citep{ChenPu2004}, so how much to ask, and what, is a design question for the DSS itself: the system should request the information that matters for the decision and stop when further questions are no longer worth the user's effort. This is especially pressing when the decision requires only the best alternative rather than a complete preference model.

Pairwise comparisons are a natural interface for such judgments, because each question involves only two alternatives, and they are the elementary building block of the Analytic Hierarchy Process (AHP) \citep{Saaty1977,Saaty1980}. Their drawback is volume: a complete pairwise-comparison matrix of order $n$ requires $n(n-1)/2$ judgments, 105 for 15 alternatives. We use the AHP as the concrete setting and ask how a DSS should manage this elicitation adaptively.

Existing approaches take two broad forms. Some prescribe in advance which comparisons to collect, as in the Best--Worst Method \citep{Rezaei2015} and related reduced-comparison methods \citep{Leal2020,Haseli2020,Abastante2019}, spanning-tree designs \citep{Siraj2012,Lundy2017,MazurekKulakowski2022} and optimal filling patterns \citep{Szadoczki2022,SzadoczkiBozoki2025}. Others estimate priorities from whichever comparisons are available \citep{HarkerIncomplete1987,HarkerQuestioning1987,Bozoki2010,Kulakowski2020,AgostonCsato2024}; see \citet{CsatoBozoki2026} for a recent overview. The question connecting these two strands is adaptive: \emph{given the answers observed so far, which comparison should be asked next?} Harker proposed a sensitivity-based rule, later extended to hierarchies \citep{HarkerIncomplete1987,MilletHarker1990}, and other sequential or adaptive strategies have been studied both within and beyond multi-criteria decision analysis \citep{FedrizziGiove2013,Ciomek2017,Meissner2010,Branke2017,Vasconcelos2019}.

We study this adaptive question for the \emph{choice} problematic: identifying the best alternative rather than estimating the entire priority vector equally well. In the terms of the DSS literature, the answer is a form of adaptive decisional guidance \citep{Silver1991,Fazlollahi1997}: the system uses the information elicited so far to decide which question is most useful to pose next and whether to pose any more. The contribution is organised around four questions.
\begin{enumerate}[leftmargin=*,label=(\arabic*)]
\item \textbf{Which comparison should be asked next?} ALEC (Adaptive Lexicographic Elicitation of Comparisons) re-ranks after every answer and asks the first missing pair in the lexicographic order induced by the current ranking. Under the log-normal model behind LLSM \citep{CrawfordWilliams1985}, estimated contrasts have variances governed by effective resistance, and Proposition~\ref{prop:direct} shows that the direct comparison of a specified contrast is locally the most informative single question. This motivates ALEC without claiming global optimality.
\item \textbf{When does adaptation help?} Comparing ALEC with a frozen initial lexicographic order isolates the value of re-ranking. We also compare ALEC with Harker's original rule, an orientation-invariant Harker variant, random questioning and a global A-optimal design.
\item \textbf{When can elicitation stop, and what can be certified?} Proposition~\ref{prop:posterior} and Corollary~\ref{cor:eps} give posterior bounds for the probability of a wrong choice or of exceeding a relative-loss tolerance under adaptive questioning. Proposition~\ref{prop:winner} complements these probabilistic criteria with a sharp consistency-and-margin certificate for agreement with the complete-matrix LLSM winner.
\item \textbf{How can adaptive elicitation be embedded in a decision support system?} Section~\ref{sec:dss} proposes an implementable architecture: a connected seed design followed by an elicitation loop with separate components for the preference store, state estimation, stopping and query selection (Figure~\ref{fig:loop}). The user sets only the loss tolerance, and the computation between two answers is small at the problem sizes studied.
\end{enumerate}

The simulations support three main messages: adaptive re-ranking is most useful in larger problems at small budgets; tolerance-based stopping can reduce elicitation substantially when its working model is reasonably calibrated; and, under a common stopping rule, ALEC and the strongest Harker variant exhibit a clear accuracy--cost trade-off rather than universal dominance. Detailed sensitivity, calibration and robustness results are reported in the supplementary material.

Section~\ref{sec:lit} reviews related work, Sections~\ref{sec:model}--\ref{sec:methods} develop the model, the elicitation rules and their implementation in a DSS, Section~\ref{sec:example} gives an illustration, Section~\ref{sec:experiments} reports the computational study, and Sections~\ref{sec:discussion}--\ref{sec:conclusion} discuss implications and conclude. All proofs, benchmark definitions, further derivations and secondary experiments are in the supplementary material; the proofs are collected in Supplementary Section~\Sref{sec:S-technical}.

%% ============================================================================
\section{Related literature}\label{sec:lit}
%% ============================================================================

\paragraph{Adaptive decision support and interactive elicitation}
Adaptive decision support systems adjust the support they provide to the decision maker, the problem and the context while leaving the user in control \citep{Fazlollahi1997}. A related stream studies decisional guidance, the way a DSS steers its users, deliberately or not, as they make judgments and choose their next steps \citep{Silver1991}; \citet{Morana2017} review the design features through which guidance is delivered and its effects on users. Behavioural research on decision aids shows that users trade decision quality against cognitive effort, and that a DSS changes this trade-off by changing the effort that each decision strategy requires \citep{Payne1993,ToddBenbasat1999}. In interactive preference elicitation the cost of questioning the user is a recognised bottleneck \citep{ChenPu2004}, and for utility-based decision support, elicitation strategies have been designed that choose queries to reduce minimax regret and stop once a (near-)optimal decision can be identified \citep{Boutilier2006}. ALEC belongs to this family of adaptive guidance mechanisms, with the pairwise-comparison matrix of the AHP as the preference model: it uses the comparisons elicited so far to decide which question is most useful next and whether further questions are warranted, and it reports how reliable the recommendation is when it stops.

\paragraph{Priorities from incomplete pairwise comparisons}
An incomplete PC matrix is conveniently represented by an undirected graph whose edges are the observed comparisons; priorities are identified (up to scale) if and only if the graph is connected. \citet{HarkerIncomplete1987} completed each missing entry by the geometric mean of the intensities of all paths that connect the two alternatives and applied the eigenvector method to the completed matrix. \citet{HarkerQuestioning1987} instead replaced each missing entry by the ratio of the (unknown) weights, which leads to the Perron vector of an auxiliary matrix with zeros in the missing positions and a modified diagonal. \citet{Bozoki2010} showed that the incomplete LLSM has a unique solution if and only if the graph is connected and that it solves a linear system whose matrix is the graph Laplacian. \citet{Lundy2017} proved that, for complete matrices, the row geometric mean equals the geometric mean of the weight vectors of all spanning trees, and \citet{BozokiTsyganok2019} extended this result to incomplete matrices. We use the incomplete LLSM for ALEC because of its transparent statistical interpretation.

\paragraph{Which comparisons to ask}
\citet{Wedley1993} studied starting rules for incomplete comparisons, and \citet{Carmone1997} evaluated by simulation how much information is lost when entries are deleted. \citet{Szadoczki2022} showed that regular or quasi-regular comparison graphs with minimal diameter are the best filling patterns for a given number of comparisons, and \citet{SzadoczkiBozoki2025} determined optimal patterns and optimal \emph{sequences} of patterns for $n\le 6$; an experiment with human participants supported these results \citep{Szadoczki2025exp}. These designs are non-adaptive, and the adaptive rules for AHP elicitation \citep{HarkerIncomplete1987,MilletHarker1990,FedrizziGiove2013} aim at the whole priority vector. Choosing the best alternative calls for a different allocation of questions, because comparisons among clearly inferior alternatives are of little value.

\paragraph{Adaptive comparisons outside multi-criteria decision analysis}
Identifying the best item from noisy pairwise comparisons is central to active ranking and dueling bandits \citep{Yue2012,Heckel2019}, where comparisons are usually binary and may be repeated. Our problem fits the framework of \emph{transductive} linear bandits \citep{Fiez2019}, in which the vectors that can be measured, here $\bm e_i-\bm e_j$, differ from the alternatives among which one must choose, here the coordinates of $\bm u$; the AHP setting adds ratio judgments on a bounded scale and the restriction that each pair is asked at most once. The connection between log-least-squares error and the effective resistance of the comparison graph is known \citep{Hendrickx2019}, and graph designs with small effective resistance have been studied by \citet{Ghosh2008} and, for statistical rankings, by \citet{Osting2014}. Finally, the problem is related to Bayesian ranking and selection, e.g.\ optimal computing budget allocation \citep{ChenOCBA2000} and the knowledge gradient for correlated beliefs \citep{Frazier2009}. Our contribution relative to these literatures is an interpretable rule tailored to AHP elicitation, stopping criteria with a direct decision meaning, and systematic benchmarking against AHP-specific competitors.

%% ============================================================================
\section{Model}\label{sec:model}
%% ============================================================================

\subsection{Incomplete comparisons and the LLSM}
Let $\{1,\dots,n\}$ index the alternatives. A judgment $a_{ij}>0$ encodes the relative preference of alternative $i$ over $j$, with $a_{ji}=1/a_{ij}$. The set of observed unordered pairs is $\E$, and $\G=(\{1,\dots,n\},\E)$ is the \emph{comparison graph}. We write $y_{ij}=\log a_{ij}$ for $\{i,j\}\in\E$, oriented with $i<j$. The incomplete LLSM solves
\begin{equation}\label{eq:llsm}
\min_{\bm{x}\in\R^n}\ \sum_{\{i,j\}\in\E}\bigl(y_{ij}-(x_i-x_j)\bigr)^2,
\end{equation}
and the normalised priorities are $w_i=\exp(x_i)/\sum_k \exp(x_k)$. Let $B\in\R^{|\E|\times n}$ be the oriented incidence matrix of $\G$ and $L=B^\top B$ the graph Laplacian. The normal equations are $L\bm{x}=B^\top\bm{y}$. If $\G$ is connected, the solution normalised by $\one^\top\bm{x}=0$ is $\xhat=\Lp B^\top\bm{y}$, where $\Lp$ is the Moore--Penrose inverse of $L$ \citep{Bozoki2010}. On a spanning tree, $\xhat$ reproduces every observed ratio exactly.

\subsection{Statistical model and effective resistance}
Following \citet{CrawfordWilliams1985} we assume multiplicative log-normal errors around true log-priorities $\bm{u}$, normalised by $\one^\top\bm{u}=0$:
\begin{equation}\label{eq:model}
y_{ij}=u_i-u_j+\varepsilon_{ij},\qquad \varepsilon_{ij}\ \text{i.i.d.}\ \mathcal{N}(0,\sigma^2),\qquad \{i,j\}\in\E.
\end{equation}
The best alternative is $b^\ast=\argmax_i u_i$. Throughout, each pair is asked at most once, as is customary in AHP elicitation. This is a restriction of the present study rather than of the model: a repeated question would add a parallel edge to $\G$, and all formulas below would remain valid with $L$ the corresponding multigraph Laplacian. The \emph{effective resistance} between $i$ and $j$ \citep{KleinRandic1993} is $\Reff_{ij}=(\bm e_i-\bm e_j)^\top \Lp(\bm e_i-\bm e_j)$, the electrical resistance between $i$ and $j$ when every comparison is a unit resistor \citep{DoyleSnell1984}. On a tree, $\Reff_{ij}$ is the length of the path between $i$ and $j$; for an observed pair, $\Reff_{ij}\le 1$; and adding a comparison never increases any effective resistance.

\begin{proposition}[Precision of LLSM contrasts]\label{prop:resistance}
Under \eqref{eq:model}, for a connected comparison graph fixed independently of the answers, $\xhat\sim\mathcal{N}(\bm u,\sigma^2\Lp)$ on the subspace $\one^\top\bm x=0$. In particular, $\hat x_i-\hat x_j\sim\mathcal{N}\bigl(u_i-u_j,\ \sigma^2\Reff_{ij}\bigr)$, and $\hat x_i-\hat x_j$ is the best linear unbiased estimator of $u_i-u_j$.
\end{proposition}

Proposition~\ref{prop:resistance} is a direct consequence of the Gauss--Markov theorem for the linear model \eqref{eq:model}. We state it because the rest of the paper builds on it; we do not claim it as a new result \citep[cf.][]{Hendrickx2019}. When the graph is built adaptively, which pairs are observed depends on earlier answers and the sampling distribution of $\xhat$ changes. This is why our stopping rules rest on a posterior argument (Section~\ref{sec:posterior}), which adaptivity does not affect.

\begin{remark}[Geometric consistency index]\label{rem:gci}
For a fixed design, $s^2=\mathrm{RSS}/(|\E|-n+1)$, with $\mathrm{RSS}$ the residual sum of squares of \eqref{eq:llsm}, is an unbiased estimator of $\sigma^2$. For a complete matrix, LLSM coincides with the row geometric mean and $|\E|-n+1=(n-1)(n-2)/2$, so $s^2$ equals the geometric consistency index (GCI) of \citet{CrawfordWilliams1985} and \citet{AguaronMoreno2003}. Thus $\mathrm{RSS}/(|\E|-n+1)$ is the natural incomplete analogue of the GCI; a regularised version of it enters the stopping rules.
\end{remark}

\subsection{The value of one more comparison}
If the pair $\{p,q\}\notin\E$ is added, the Laplacian changes by a rank-one term, and the Sherman--Morrison formula gives (proof in Supplementary Section~\Sref{sec:S-technical})
\begin{equation}\label{eq:rankone}
\Reff'_{ij}=\Reff_{ij}-\frac{c_{ij,pq}^2}{1+\Reff_{pq}},\qquad c_{ij,pq}=(\bm e_i-\bm e_j)^\top\Lp(\bm e_p-\bm e_q).
\end{equation}

\begin{proposition}[Direct comparisons are locally the most informative]\label{prop:direct}
Let $\G$ be connected and $\{i,j\}\notin\E$. Among all pairs $\{p,q\}\notin\E$, the reduction of $\Reff_{ij}$ in \eqref{eq:rankone} is uniquely maximised by $\{p,q\}=\{i,j\}$, in which case $\Reff'_{ij}=\Reff_{ij}/(1+\Reff_{ij})$.
\end{proposition}

Proposition~\ref{prop:direct} answers a \emph{local} question: once we know which contrast we want to sharpen, the best single question is its direct comparison. It does not say which contrast should be targeted, so ALEC, which targets the contrasts of the current leader, is a deliberately simple heuristic. A concrete counterexample and the uncertainty-based ALEC-U variant are given in Supplementary Section~\Sref{sec:S-technical}.

%% ============================================================================
\section{Elicitation procedures}\label{sec:methods}
%% ============================================================================

\subsection{ALEC}\label{sec:alec}
ALEC starts from an initial connected set of comparisons, for instance a spanning tree chosen by the decision maker plus one further pair, as in our experiments. After each answer it recomputes the LLSM priorities, ranks the alternatives (ties broken by a fixed index order), and forms the lexicographic list of pairs induced by the ranking $\pi^{(t)}$:
\[
\{\pi_1,\pi_2\},\{\pi_1,\pi_3\},\dots,\{\pi_1,\pi_n\},\{\pi_2,\pi_3\},\dots,\{\pi_{n-1},\pi_n\}.
\]
The next question is the first pair in this list that has not been asked yet (Algorithm~\ref{alg:alec}). As long as the leader $\ell=\pi_1$ has not been compared with every other alternative, ALEC therefore asks $\ell$ against the highest-ranked alternative not yet compared with it; once $\ell$ is \emph{leader-complete}, i.e.\ has been compared directly with every other alternative, it proceeds to the second-ranked alternative, and so on. Because the list is rebuilt after every answer, both the leader and the order of its challengers can change during elicitation.

\begin{algorithm}[t]
\caption{ALEC: Adaptive Lexicographic Elicitation of Comparisons}\label{alg:alec}
\begin{algorithmic}[1]
\Require initial connected comparison set $\E_0$ with judgments; stopping rule (budget $m$, leader-complete, or $\beta_\varepsilon\le\alpha$)
\State $\E\gets\E_0$;\quad $\xhat\gets\textsc{LLSM}(\E)$;\quad $\pi\gets\textsc{Rank}(\xhat)$
\While{the stopping criterion is not met, the budget is not exhausted and $\E$ is not complete}
  \State $\{i,j\}\gets$ first pair of $\textsc{LexPairs}(\pi)$ with $\{i,j\}\notin\E$
  \State elicit $a_{ij}$;\quad $\E\gets\E\cup\{\{i,j\}\}$
  \State $\xhat\gets\textsc{LLSM}(\E)$;\quad $\pi\gets\textsc{Rank}(\xhat)$
\EndWhile
\State \Return leader $\pi_1$, ranking $\pi$, priorities $\bm w\propto\exp(\xhat)$, final bound $\beta_\varepsilon$, and the reason for stopping (criterion met, budget exhausted, or matrix complete)
\end{algorithmic}
\end{algorithm}

\begin{remark}[Relation to the Best--Worst Method]\label{rem:bwm}
When ALEC is run until the leader is complete, $\G$ contains the star centred at the leader, i.e.\ a full vector of ``best-to-others'' judgments as in the Best--Worst Method \citep{Rezaei2015}. ALEC can thus be read as an adaptive way of choosing the reference alternative that BWM fixes in advance.
\end{remark}

\subsection{Posterior bounds for the best alternative}\label{sec:posterior}
We place a scaled-inverse-$\chi^2$ prior with $\nu_0$ degrees of freedom and scale $\sigma_0^2$ on $\sigma^2$, and a flat prior on $\bm u$ in the subspace $\one^\top\bm u=0$.

\begin{proposition}[Posterior error bound under adaptive questioning]\label{prop:posterior}
Let the initial comparisons be chosen independently of $(\bm u,\sigma)$ and the subsequent ones sequentially, each choice depending only on the previous questions and answers (and possibly on randomisation independent of $(\bm u,\sigma)$), and let elicitation stop at a time determined by the same information. Let ``data'' denote the questions and answers obtained and $\G$ the resulting connected comparison graph. Let $\nu=\nu_0+|\E|-n+1>0$ and $s^2=(\nu_0\sigma_0^2+\mathrm{RSS})/\nu$. Then, for every pair $i,j$,
\[
\frac{(u_i-u_j)-(\hat x_i-\hat x_j)}{s\sqrt{\Reff_{ij}}}\ \Big|\ \text{data}\ \sim\ t_\nu ,
\]
and, for the current leader $\ell$,
\begin{equation}\label{eq:beta}
\Pr(\ell\neq b^\ast\mid\text{data})\ \le\ \beta:=\sum_{j\neq\ell}F_{\nu}\!\left(-\frac{\hat x_\ell-\hat x_j}{s\sqrt{\Reff_{\ell j}}}\right),
\end{equation}
where $F_\nu$ is the Student-$t$ distribution function with $\nu$ degrees of freedom.
\end{proposition}

When two alternatives are almost tied, certifying the exact winner may require many questions although choosing either would cost little. Let $\rho=(w_{b^\ast}-w_\ell)/w_{b^\ast}=1-\exp(u_\ell-\max_j u_j)$ be the relative loss of choosing the leader, and for a tolerance $0\le\varepsilon<1$ let $\delta_\varepsilon=-\log(1-\varepsilon)$.

\begin{corollary}[Tolerance bound]\label{cor:eps}
Under the conditions of Proposition~\ref{prop:posterior},
\begin{equation}\label{eq:betaeps}
\Pr(\rho>\varepsilon\mid\text{data})\ \le\ \beta_\varepsilon:=\sum_{j\neq\ell}F_{\nu}\!\left(-\frac{\hat x_\ell-\hat x_j+\delta_\varepsilon}{s\sqrt{\Reff_{\ell j}}}\right),
\end{equation}
and $\beta_0=\beta$.
\end{corollary}

Stopping when $\beta_\varepsilon\le\alpha$ thus states that, under the posterior, the relative loss of the chosen alternative does not exceed $\varepsilon$ with probability at least $1-\alpha$. In general, $\max\{0,1-\beta_\varepsilon\}$ is a lower bound on the posterior probability $\Pr(\rho\le\varepsilon\mid\text{data})$, not that probability itself: $\beta_\varepsilon$ adds marginal probabilities and can exceed one. This distinction affects how the result should be reported, but not the stopping condition, since $\alpha<1$. The adaptivity of the questions does not invalidate these statements, because the questions depend on $\bm u$ only through the observed answers and the likelihood is therefore the same as for a fixed design \citep{BergerWolpert1988}. The union bound makes $\beta_\varepsilon$ conservative when several challengers are close. Because \eqref{eq:beta} and \eqref{eq:betaeps} are Bayesian statements, their frequentist calibration depends on how well the priors match the process generating the weights and the noise; Section~\ref{sec:stop-results} examines matched and mismatched cases. With $\nu_0=0$ (the improper prior $p(\sigma^2)\propto\sigma^{-2}$), $s^2$ is the incomplete GCI of Remark~\ref{rem:gci}; with $\nu_0>0$ the rules can be used from the first step, when a spanning tree provides no residual degrees of freedom. We use $\nu_0=2$ with two prior scales: $\sigma_0=0.6$, the setting we recommend in practice, motivated by the GCI acceptability threshold of about $0.37$ for $n>4$ \citep{AguaronMoreno2003}, since $\sqrt{0.37}\approx0.61$; and the less prudent $\sigma_0=0.3$, used in the illustrative example and in the sensitivity analyses. The GCI of an observed matrix, the noise parameter $\sigma$ and the prior scale $\sigma_0$ are different objects: the threshold makes $\sigma_0=0.6$ an interpretable working scale, not an upper bound on the true noise.

\subsection{ALEC-U: an uncertainty-reduction variant}\label{sec:alecu}
For sensitivity analysis we also consider ALEC-U, which chooses the missing pair predicted to reduce the posterior error bound most when the current estimates are held fixed. It combines estimated gaps with effective resistances, but it is not a full Bayesian look-ahead because the next answer may change the estimates, leader and noise estimate. Its exact criterion and a counterexample separating it from ALEC are given in Supplementary Section~\Sref{sec:S-technical}; its results are reported with those of the other rules (Table~\ref{tab:hit-n} and Supplementary Sections~\Sref{sec:S-ref}--\Sref{sec:S-calib}).

\subsection{Stopping rules}\label{sec:stopping}
We consider (a) a fixed budget of $m$ additional comparisons; (b) the \emph{leader-complete} rule, which stops when the current leader has been compared directly with all other alternatives; (c) the \emph{posterior} rule $\beta\le\alpha$; (d) the posterior rule capped at $2n$ additional comparisons; and (e) the \emph{tolerance} rule $\beta_\varepsilon\le\alpha$ with $\varepsilon>0$, with or without the same cap. Combining ALEC's questioning rule with the capped tolerance rule gives the protocol that we recommend in Section~\ref{sec:discussion}; Section~\ref{sec:stop-results} evaluates exactly this combination. Rules (c)--(e) adapt the number of questions to the instance. Reaching the cap in (d) or exhausting a budget does not mean that the posterior threshold has been met, which is why Algorithm~\ref{alg:alec} reports the reason for stopping.

The leader-complete rule raises a question that does not involve the posterior: if the complete matrix were elicited, would a leader that has won all its direct comparisons remain its LLSM winner?

\begin{proposition}[Direct winners and the complete matrix]\label{prop:winner}
Let $A$ be a complete PC matrix of order $n\ge3$ with log-judgments $y_{ij}$, LLSM log-priorities $\bm x$, residuals $e_{ij}=y_{ij}-(x_i-x_j)$ and $\mathrm{GCI}(A)=\frac{2}{(n-1)(n-2)}\sum_{i<j}e_{ij}^2$, and let $\gamma_n(A)=(n-2)\sqrt{(n-1)\,\mathrm{GCI}(A)/(2n)}$.
\begin{enumerate}[label=(\alph*)]
\item For every $k\neq i$, $x_k-x_i=y_{ki}-e_{ki}$; hence $k$ is the unique LLSM winner of $A$ if and only if $e_{ki}<y_{ki}$ for every $i\neq k$.
\item $|e_{ij}|\le\gamma_n(A)$ for every pair $\{i,j\}$, and equality is possible.
\item If $k$ wins every direct comparison with a margin $\delta_k=\min_{i\neq k}y_{ki}>\gamma_n(A)$, then $k$ is the unique LLSM winner of $A$. The threshold cannot be lowered: for every $g>0$ and every $\delta\in\bigl(0,(n-2)\sqrt{(n-1)g/(2n)}\,\bigr]$, some complete matrix with $\mathrm{GCI}=g$ has an alternative that wins every direct comparison with margin $\delta$ and is not the unique LLSM winner.
\end{enumerate}
\end{proposition}

With a tie ($y_{ki}=0$), part (a) shows that $k$ stays ahead of $i$ only if $e_{ki}<0$; the GCI bounds the size of the residuals but not their sign, so no condition on the GCI alone guarantees the outcome. The condition in (c) involves the GCI of the complete matrix, which is not observed during elicitation. It can therefore serve only as a worst-case certificate based on an assumed upper bound for that GCI, and the resulting requirement is demanding. For $n=11$, even a leader that wins every comparison 9 to 1 is certified only if $\mathrm{GCI}(A)<\res{dw-gci-n11-a9}$, about a third of the acceptability threshold of $0.37$ \citep{AguaronMoreno2003}; a margin of 3 to 1 requires $\mathrm{GCI}(A)<\res{dw-gci-n11-a3}$. In our simulations the certificate applies to at most \res{dw-cert-max}\% of the direct winners and almost never for $n\ge11$, although a direct winner is usually the LLSM winner of the complete matrix (Section~\ref{sec:stop-results}). This is one reason to base stopping on posterior bounds rather than on worst-case guarantees.

\subsection{Decision-support implementation}\label{sec:dss}
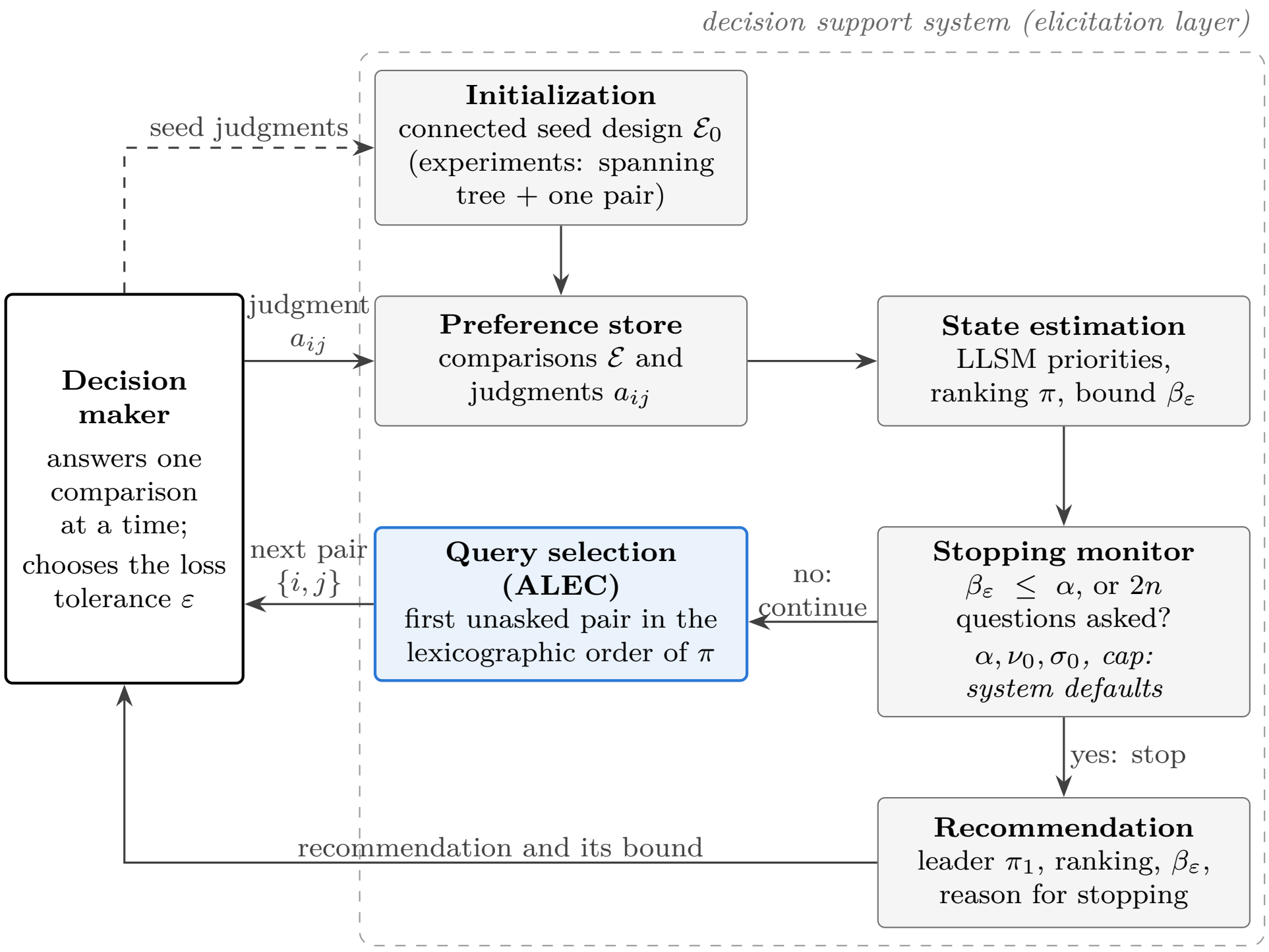
\begin{figure}[t]
\centering
\begin{tikzpicture}[font=\scriptsize, node distance=10mm and 13mm, execute at begin node={\hyphenpenalty=10000\exhyphenpenalty=10000},
  comp/.style={draw=black!55, rounded corners=2pt, fill=black!4, align=center, text width=3.45cm, inner sep=4pt, minimum height=13mm},
  alec/.style={comp, draw=alecblue, fill=alecblue!10, line width=0.8pt},
  user/.style={draw=black, rounded corners=2pt, align=center, text width=2.1cm, inner sep=4pt, line width=0.8pt},
  flow/.style={-{Stealth[length=2.2mm]}, black!75, line width=0.6pt},
  lab/.style={font=\scriptsize, align=center, inner sep=1.5pt}]
\node[comp] (store) {\textbf{Preference store}\\ comparisons $\E$ and judgments $a_{ij}$};
\node[comp, above=7mm of store, minimum height=9mm] (init) {\textbf{Initialization}\\ connected seed design $\E_0$\\ (experiments: spanning tree + one pair)};
\node[comp, right=of store] (est) {\textbf{State estimation}\\ LLSM priorities,\\ ranking $\pi$, bound $\beta_\varepsilon$};
\node[comp, below=of est] (stop) {\textbf{Stopping monitor}\\ $\beta_\varepsilon\le\alpha$, or $2n$\\ questions asked?\\[1pt] {\itshape $\alpha,\nu_0,\sigma_0$, cap: system defaults}};
\node[alec, below=of store] (query) {\textbf{Query selection (ALEC)}\\ first unasked pair in the\\ lexicographic order of $\pi$};
\node[comp, below=8mm of stop] (rec) {\textbf{Recommendation}\\ leader $\pi_1$, ranking, $\beta_\varepsilon$,\\ reason for stopping};
\path (store.north west) -- (query.south west) coordinate[pos=0.5] (mid);
\node[user, left=of mid, minimum height=39mm] (dm) {\textbf{Decision maker}\\[3pt] answers one\\ comparison at a time;\\[2pt] chooses the loss\\ tolerance $\varepsilon$};
\draw[flow] (dm.east |- store) -- node[lab, above] {judgment\\ $a_{ij}$} (store.west);
\draw[flow] (init) -- (store);
\draw[flow, dashed] (dm.north) |- node[lab, above, pos=0.75] {seed judgments} (init.west);
\draw[flow] (store) -- (est);
\draw[flow] (est) -- (stop);
\draw[flow] (stop.west) -- node[lab, above] {no:\\ continue} (stop.west -| query.east);
\draw[flow] (query.west) -- node[lab, above] {next pair\\ $\{i,j\}$} (dm.east |- query);
\draw[flow] (stop) -- node[lab, right] {yes: stop} (rec);
\draw[flow] (rec.west) -| node[lab, above, pos=0.25] {recommendation and its bound} (dm.south);
\begin{scope}[on background layer]
\node[draw=black!40, dashed, rounded corners=4pt, inner sep=5pt, inner xsep=4pt, fit=(init) (store) (est) (stop) (query) (rec),
      label={[font=\scriptsize\itshape, text=black!60, anchor=south east]north east:decision support system (elicitation layer)}] {};
\end{scope}
\end{tikzpicture}
\caption{ALEC as an adaptive decision-support elicitation loop. The loop starts from a connected seed design $\E_0$ answered by the decision maker, who then answers one pairwise comparison at a time and chooses the loss tolerance $\varepsilon$; the risk level $\alpha$, the prior parameters $\nu_0,\sigma_0$ and the budget cap are set by the system designer, with defaults. After every answer the system updates its estimate, checks whether further questions are warranted and, if so, selects the next pair. Only the query-selection component is specific to ALEC; replacing it by another questioning rule, such as Harker-log, leaves the rest of the loop unchanged.}\label{fig:loop}
\end{figure}

Figure~\ref{fig:loop} shows ALEC as the elicitation layer of a DSS. Before the adaptive loop starts, the system elicits a connected seed design $\E_0$, so that all priorities are identified; in our experiments this consists of a spanning tree plus one additional comparison (Section~\ref{sec:design}), and in an application the tree can also be chosen by the decision maker or the analyst (Section~\ref{sec:alec}). A preference store then keeps the comparisons answered so far. After each answer, a state-estimation component recomputes the LLSM priorities, the ranking and the bound $\beta_\varepsilon$; a stopping monitor compares the bound with the risk level $\alpha$ and counts the questions asked against the cap; and, if elicitation continues, the query-selection component returns the next question, the first unasked pair in ALEC's lexicographic order. When the monitor stops, the system reports the leader, the ranking, the final bound and the reason for stopping, as in Algorithm~\ref{alg:alec}. The components interact only through the store and the current estimate, so the questioning rule or the stopping criterion can be replaced, for instance by Harker-log or by a fixed budget, without changing the rest of the system; this is how ALEC and Harker-log are compared under the same stopping rule in Section~\ref{sec:stop-results}. The loss tolerance $\varepsilon$ is the user-facing preference parameter (``an alternative whose priority is at most 5\% below that of the best is acceptable''), whereas $\alpha$, $\nu_0$, $\sigma_0$ and the maximum budget can be configured by the DSS designer or provided as defaults, such as the values recommended in Section~\ref{sec:discussion}.

The computational cost per question is small. In our Python implementation, selecting the next question and updating the estimates takes \res{time-ALEC-min}--\res{time-ALEC-max}\,ms per question for $n=7$--$15$ (Supplementary Table~\Sref{tab:timing}), and the bound adds only the pseudoinverse of an $n\times n$ Laplacian and $n-1$ evaluations of a Student-$t$ distribution function. These timings indicate that computational latency is unlikely to be a bottleneck at the problem sizes studied, making ALEC suitable for real-time interactive implementation.

\subsection{Benchmark questioning rules}\label{sec:benchmarks}
We use six benchmark families. \textbf{Harker-PG} is Harker's original path-completion and sensitivity procedure \citep{HarkerIncomplete1987}. \textbf{Harker-H} applies the same sensitivity idea to Harker's auxiliary-matrix estimator \citep{HarkerQuestioning1987}; this combination is our adaptation. \textbf{Harker-log} ranks missing pairs by sensitivity with respect to $\log a_{kl}$, removing the arbitrary dependence on the orientation of a reciprocal pair. \textbf{Frozen} follows the initial lexicographic order throughout and therefore isolates the value of re-ranking. \textbf{Random} chooses a missing pair uniformly, whereas \textbf{A-opt} selects the comparison giving the largest reduction in $\operatorname{tr}(L^+)$ and targets global estimation accuracy. We also report Star-ALEC, ALEC-H and ALEC-U as sensitivity variants.

The exact Harker derivatives, implementation details and the full ALEC-U criterion are provided in Supplementary Section~\Sref{sec:S-technical}. Because exact path enumeration becomes expensive as the comparison graph gets denser, Harker-PG is run up to $m=2n$ additional comparisons for $n=7$ and up to $m=n$ for $n\ge11$. To separate question selection from the final estimator, selected experiments re-evaluate every rule's acquired comparisons with a common final LLSM estimator.

%% ============================================================================
%% ============================================================================
\section{Illustrative example}\label{sec:example}
%% ============================================================================

We use Saaty's beverage-consumption matrix \citep{Saaty2008} with seven alternatives (Coffee, Wine, Tea, Beer, Soda, Milk, Water) and consistency ratio $\mathrm{CR}=\res{ex-cr}$. The eigenvector method (EVM) applied to the complete matrix ranks them as \res{ex-rank-evm}. The decision maker starts with the spanning tree
\begin{multline*}
\E_0=\bigl\{\{\text{Water},\text{Coffee}\},\{\text{Water},\text{Soda}\},\{\text{Coffee},\text{Tea}\},\\
\{\text{Tea},\text{Beer}\},\{\text{Beer},\text{Milk}\},\{\text{Milk},\text{Wine}\}\bigr\}
\end{multline*}
with judgments $2,2,2,1,1,3$. On a tree, LLSM reproduces the judgments exactly: $w^{(0)}_{\text{Water}}=\res{ex-w0-Water}$, $w^{(0)}_{\text{Coffee}}=w^{(0)}_{\text{Soda}}=\res{ex-w0-Coffee}$, $w^{(0)}_{\text{Tea}}=w^{(0)}_{\text{Beer}}=w^{(0)}_{\text{Milk}}=\res{ex-w0-Tea}$ and $w^{(0)}_{\text{Wine}}=\res{ex-w0-Wine}$, with ties broken by the index order of the matrix.

Table~\ref{tab:example} shows the sequence, and Figure~\Sref{fig:example} in the supplementary material plots the weights. Water vs Coffee and Water vs Soda are already known, so ALEC asks Water against Tea, Beer, Milk and Wine in turn, re-ranking after each answer. After three additional comparisons the LLSM ranking coincides with the ranking of the complete matrix. The leader-complete rule stops after $m=\res{ex-stop-lc}$ additional comparisons, i.e.\ 10 of the 21 judgments. The bound $\beta$ falls from \res{ex-beta-0} on the tree to \res{ex-beta-2} after two answers, so the posterior rule with $\alpha=0.10$ stops at $m=\res{ex-stop-b10}$, with Water identified from 8 judgments.

\begin{table}[t]
\centering\footnotesize\setlength{\tabcolsep}{4pt}
\caption{ALEC on Saaty's beverage matrix from the spanning tree $\E_0$. Abbreviations: Wa Water, C Coffee, S Soda, T Tea, B Beer, M Milk, Wi Wine. $\beta$ is the bound \eqref{eq:beta} ($\nu_0=2$, $\sigma_0=0.3$); LC: the leader has been compared with all others; ``Next pair'' and $a$ are the question asked after step $m$ and its answer, oriented from the higher-ranked alternative.}\label{tab:example}
\begin{tabular}{@{}r l c c c l c@{}}
\toprule
$m$ & Ranking & $w_{\text{Water}}$ & $\beta$ & LC & Next pair & $a$ \\
\midrule
0 & Wa$\succ$C$\succ$S$\succ$T$\succ$B$\succ$M$\succ$Wi & 0.353 & 0.353 & no & Water--Tea & 3 \\
1 & Wa$\succ$C$\succ$S$\succ$T$\succ$B$\succ$M$\succ$Wi & 0.327 & 0.143 & no & Water--Beer & 3 \\
2 & Wa$\succ$C$\succ$S$\succ$B$\succ$M$\succ$T$\succ$Wi & 0.321 & 0.049 & no & Water--Milk & 2 \\
3 & Wa$\succ$C$\succ$S$\succ$M$\succ$B$\succ$T$\succ$Wi & 0.302 & 0.037 & no & Water--Wine & 9 \\
4 & Wa$\succ$C$\succ$S$\succ$M$\succ$B$\succ$T$\succ$Wi & 0.307 & 0.021 & yes & Coffee--Soda & 1 \\
5 & Wa$\succ$C$\succ$S$\succ$M$\succ$B$\succ$T$\succ$Wi & 0.307 & 0.006 & yes & Coffee--Milk & 1 \\
6 & Wa$\succ$C$\succ$S$\succ$M$\succ$B$\succ$T$\succ$Wi & 0.309 & 0.002 & yes & -- & -- \\
\bottomrule
\end{tabular}
\end{table}

%% ============================================================================
\section{Computational study}\label{sec:experiments}
%% ============================================================================

\subsection{Design}\label{sec:design}
All rules start from the same $n$ comparisons: a uniformly random labelled spanning tree plus one random additional pair. Every rule sees the same instances and initial comparisons (common random numbers), so all comparisons between rules are paired. We study $n\in\{7,11,15\}$ and budgets of $m=0,\dots,2n$ additional comparisons; for $n=7$, $m=2n$ completes the matrix. Two published complete matrices serve as reference examples: the beverage matrix of Section~\ref{sec:example} and the $11\times11$ matrix of safety indicators of \citet{Esmaeili2023} ($\mathrm{CR}=0.092$); rules are scored against the leader of the complete matrix (\res{R-ref} runs, seed 20261013). In the main simulation, true log-weights are drawn as $u_i\sim\mathcal U(0,\log 9)$ and judgments follow \eqref{eq:model} with $\sigma\in\{0.15,0.3,0.5\}$ (\res{R-syn} replications per scenario). Judgments are continuous and are not clipped; complete matrices generated in this way have mean consistency ratios of about \res{cr-s015-mean}, \res{cr-s030-mean} and \res{cr-s050-mean} for the three noise levels. Robustness scenarios for $n=11$ and $\sigma=0.3$ vary one feature at a time: judgments are rounded to Saaty's scale, weights are drawn from a flat Dirichlet distribution, or the initial graph is a random Hamiltonian path or a star. For the rounding experiment, log-judgments are clipped to $[-\log 9,\log 9]$ and rounded to the nearest scale value; the stopping-rule check also includes $\sigma=0.5$. The rounded scenario and the two alternative initial designs reuse the true weights and continuous judgments of the baseline, so their comparisons with it are paired.

The primary outcome is the hit rate, the percentage of runs in which the estimated leader is truly best; the complete matrix is itself noisy, so its hit rate is a reference rather than an upper bound. We also study relative regret $\rho$, Kendall's $\tau_b$ and agreement with the complete-matrix leader. Paired hit rates are compared with exact McNemar tests \citep{McNemar1947}; families of main comparisons use a Holm correction, and paired confidence intervals use Newcombe's score method \citep{Newcombe1998}. Common-estimator analyses and robustness tables are in the supplementary material.

\subsection{Reference matrices}\label{sec:ref-results}
On the two published matrices, ALEC reaches high agreement with the complete-matrix leader at small budgets and improves clearly on Harker's original procedure. Harker-log is a substantially stronger benchmark. On the matrix of \citet{Esmaeili2023} at the largest reported budget, $m=10$, ALEC finds the leader in \res{esm-ALEC-hit-10}\% of the runs, compared with \res{esm-ALEC-U-hit-10}\% for ALEC-U, \res{esm-HARKER-LOG-hit-10}\% for Harker-log and \res{esm-STAR-hit-10}\% for Star-ALEC. The complete trajectories and weight-reconstruction metrics are reported in Tables~\Sref{tab:esmaeili}--\Sref{tab:beverage} of the supplementary material.

\subsection{Ground truth: fixed budgets}\label{sec:fixed}
\begin{table}[t]
\centering\scriptsize\setlength{\tabcolsep}{2.6pt}
\caption{Hit rate (\%) after $m=n$ additional comparisons, i.e.\ $2n$ of the $n(n-1)/2$ comparisons, \res{R-syn} replications per row. ``Full'': hit rate of the complete matrix (LLSM). Harker rules use their own estimators. Best rule per row in bold. Difference from ALEC (exact McNemar test): $^\dagger$ significant after a Holm correction over the nine rows of the column; $^\ast$ $p<0.05$ without correction only. Abbreviations: Rand.\ Random, H-PG Harker-PG, H-H Harker-H, H-log Harker-log, Star-A Star-ALEC.}\label{tab:hit-n}
\begin{tabular}{@{}r r r rrrrrrrrr@{}}
\toprule
$n$ & $\sigma$ & Full & Rand. & H-PG & H-H & H-log & A-opt & Frozen & Star-A & ALEC & ALEC-U \\
\midrule
7 & 0.15 & 91.4 & 86.8$^{\dagger}$ & 90.2 & 89.3$^{*}$ & 90.4 & 89.4$^{*}$ & \textbf{91.7} & 91.6 & 91.4 & 90.2 \\
7 & 0.30 & 80.8 & 74.9$^{\dagger}$ & 77.6 & 77.5$^{*}$ & 78.5 & 77.1$^{*}$ & 79.8 & \textbf{80.2} & 79.5 & 79.2 \\
7 & 0.50 & 72.5 & 67.2$^{\dagger}$ & 68.5$^{*}$ & 68.7 & 69.6 & 67.7$^{\dagger}$ & 70.0 & 70.4 & \textbf{71.2} & 70.2 \\
\addlinespace
11 & 0.15 & 88.9 & 80.2$^{\dagger}$ & 83.0 & 82.3$^{*}$ & 85.1 & 82.7$^{*}$ & 85.3 & \textbf{85.9} & 85.3 & 85.8 \\
11 & 0.30 & 76.6 & 66.3$^{\dagger}$ & 68.7$^{\dagger}$ & 71.0$^{*}$ & \textbf{74.5} & 69.0$^{\dagger}$ & 73.0 & 73.7 & 74.3 & 74.2 \\
11 & 0.50 & 67.7 & 52.3$^{\dagger}$ & 54.7$^{\dagger}$ & 56.2$^{\dagger}$ & \textbf{63.0} & 56.3$^{\dagger}$ & 59.1 & 60.7 & 61.1 & 62.2 \\
\addlinespace
15 & 0.15 & 85.5 & 74.2$^{\dagger}$ & 77.7$^{\dagger}$ & 78.7$^{*}$ & 82.4 & 77.5$^{\dagger}$ & 79.9 & 82.1 & 82.1 & \textbf{83.3} \\
15 & 0.30 & 76.1 & 61.4$^{\dagger}$ & 61.9$^{\dagger}$ & 62.1$^{\dagger}$ & \textbf{70.7} & 60.6$^{\dagger}$ & 66.5 & 69.6 & 68.2 & 69.7 \\
15 & 0.50 & 67.5 & 43.8$^{\dagger}$ & 49.6$^{\dagger}$ & 51.4$^{*}$ & 53.4 & 47.9$^{\dagger}$ & 49.9$^{\dagger}$ & 55.0 & 56.0 & \textbf{57.0} \\
\bottomrule
\end{tabular}
\end{table}

\begin{figure}[t]
\centering
\includegraphics[width=\textwidth]{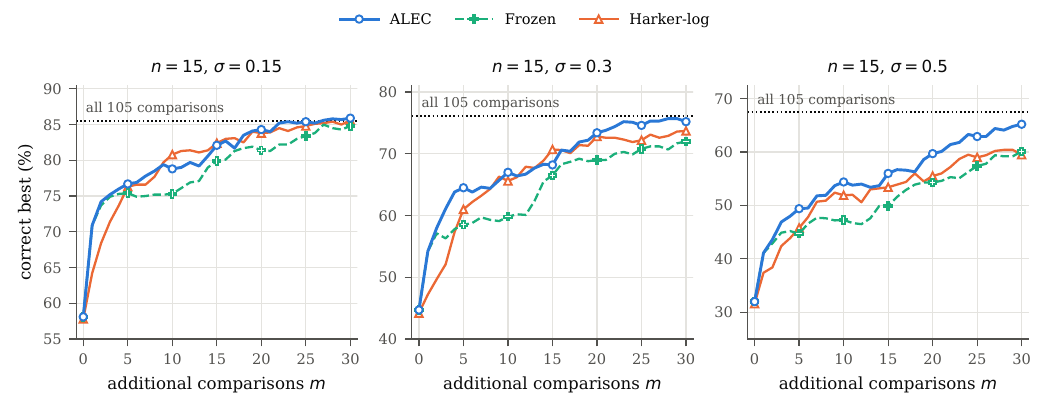}
\caption{Hit rate as a function of the number of additional comparisons $m$ for $n=15$, under the three noise levels of the main simulation. ALEC is compared with the frozen lexicographic order and with Harker-log, the strongest Harker benchmark; the dotted horizontal line is the hit rate of the complete matrix. At $m=2n=30$, the $n$ initial and $m$ additional comparisons are 45 of the 105 possible ones. \res{R-syn} paired replications per scenario; Figure~\Sref{fig:curves} in the supplementary material shows all sizes and rules.}\label{fig:budget}
\end{figure}

Table~\ref{tab:hit-n} gives the central fixed-budget comparison. After a Holm correction, ALEC is significantly more accurate than Harker-PG in \res{count-HARKER-PG-n-holm} of the nine scenarios, all with $n\ge11$. This advantage is not explained solely by the final estimation method: when every rule's comparisons are evaluated with a common LLSM estimator, ALEC's hit rate remains above Harker-PG's in all six scenarios with $n\ge11$, significantly in \res{count-HARKER-PG-n-llsm-holm} of them (Table~\Sref{tab:hit-llsm}). Harker-log is the strongest Harker benchmark: at the smaller budgets neither rule dominates, and the paired intervals are wide enough that we describe them as competitive rather than equivalent (Table~\Sref{tab:ci-hlog}).

Re-ranking is most useful in larger problems at small budgets. Compared with Frozen, ALEC shows its clearest gains for $n=11$ and $15$, whereas the difference is negligible for $n=7$. Pooling the three noise levels, re-ranking raises the hit rate for $n=15$ by \res{pool-n15-h-diff} percentage points at $m=\lceil n/2\rceil$ (paired 95\% score interval \res{pool-n15-h-ci} over the $3\times\res{R-syn}$ runs; hit rates by noise level in Table~\Sref{tab:hit-h}). Stratification by the true gap between the top two alternatives suggests that adaptation helps most in intermediate cases, although the between-stratum permutation test is not conclusive; full results are in Table~\Sref{tab:gap}. A-opt performs well on global ranking metrics but less well on identifying the winner, illustrating that a good design for the whole priority vector need not be a good design for the choice problem.

Figure~\ref{fig:budget} follows the hit rate along the whole budget for $n=15$. ALEC and Frozen ask the same first question, so their curves coincide at the start; from the second additional question on, ALEC is above Frozen at every budget and noise level. With $\sigma=0.15$ the gap peaks at \res{n15-s015-budget-dfrozen-max} percentage points and narrows to \res{n15-s015-budget-dfrozen-2n} at $m=2n$, whereas with $\sigma=0.3$ and $0.5$ it is still \res{n15-s030-budget-dfrozen-2n} and \res{n15-s050-budget-dfrozen-2n} points at $m=2n$, both significant after a Holm correction (Table~\Sref{tab:hit-2n}). Harker-log trails ALEC over the first few questions, is slightly ahead (by at most \res{budget-hlog-lead-max} points) at some budgets around $m=n$ when $\sigma\le0.3$, and stays below ALEC at every budget when $\sigma=0.5$. These curves are point estimates; paired tests at $m=\lceil n/2\rceil$, $n$ and $2n$ are reported in Table~\Sref{tab:hit-h}, Table~\ref{tab:hit-n} and Table~\Sref{tab:hit-2n}.

At $m=2n$, ALEC has used $15+30=45$ of the 105 possible comparisons for $n=15$ (\res{n15-s030-frac-2n}\%) and finds the best alternative in \res{n15-s015-ALEC-hit-2n}\%, \res{n15-s030-ALEC-hit-2n}\% and \res{n15-s050-ALEC-hit-2n}\% of the runs, against \res{n15-s015-full-hit}\%, \res{n15-s030-full-hit}\% and \res{n15-s050-full-hit}\% with all 105. Among the scenarios with $n\ge11$, the shortfall of \res{alec2n-diff-min-abs} points for $n=15$ and $\sigma=0.5$ is the only difference from the complete matrix whose paired interval excludes zero (Table~\Sref{tab:savings}, which also reports the regret).

\subsection{Stopping rules}\label{sec:stop-results}
The leader-complete rule is inexpensive but does not generally certify the full-matrix winner. In the nine main scenarios, among runs where the leader has won every direct comparison, it is the unique LLSM winner of the complete matrix in \res{dw-agree-min}--\res{dw-agree-max}\% of the runs, yet the certificate of Proposition~\ref{prop:winner}, which is sharp in the worst case, applies to at most \res{dw-cert-max}\% of these direct winners, and to at most \res{dw-cert-max-n11}\% for $n\ge11$ (Table~\Sref{tab:winner}). This motivates the posterior tolerance rule used in the recommended protocol.

\begin{table}[tb]
\centering\footnotesize\setlength{\tabcolsep}{2.2pt}
\caption{Recommended protocol: ALEC questions from the $n$ initial comparisons, stop when $\beta_\varepsilon\le0.10$ or after $2n$ additional comparisons, with $\varepsilon=0.05$ and $\sigma_0=0.6$. $\bar q$: mean number of additional comparisons; Regret: mean relative loss $\rho$; Stop/Cap: shares stopped by the bound or by the budget; last column: share with $\rho>\varepsilon$ among runs stopped by the bound, with a Wilson 95\% interval. Full sensitivity results are in Table~\Sref{tab:protocol-full}.}\label{tab:protocol}
\begin{tabular}{@{}r r r r r r r l@{}}
\toprule
$n$ & $\sigma$ & $\bar q$ & Regret & $\rho>\varepsilon$ & Stop & Cap & $\rho>\varepsilon$ if stopped \\
 & & & (\%) & (\% runs) & (\%) & (\%) & (\% runs) [95\% CI] \\
\midrule
11 & 0.15 & 14.0 & 0.43 & 3.2 & 68 & 32 & 0.4 [0.1, 1.3] \\
11 & 0.30 & 16.4 & 1.59 & 11.6 & 51 & 49 & 1.8 [0.9, 3.3] \\
11 & 0.50 & 18.2 & 3.47 & 22.9 & 36 & 64 & 7.6 [5.3, 10.8] \\
15 & 0.15 & 20.2 & 0.42 & 3.0 & 70 & 30 & 0.1 [0.0, 0.8] \\
15 & 0.30 & 23.7 & 1.71 & 14.2 & 49 & 51 & 3.5 [2.2, 5.5] \\
15 & 0.50 & 26.6 & 3.57 & 23.3 & 28 & 72 & 9.5 [6.6, 13.5] \\
\bottomrule
\end{tabular}
\end{table}

\begin{figure}[t]
\centering
\includegraphics[width=\textwidth]{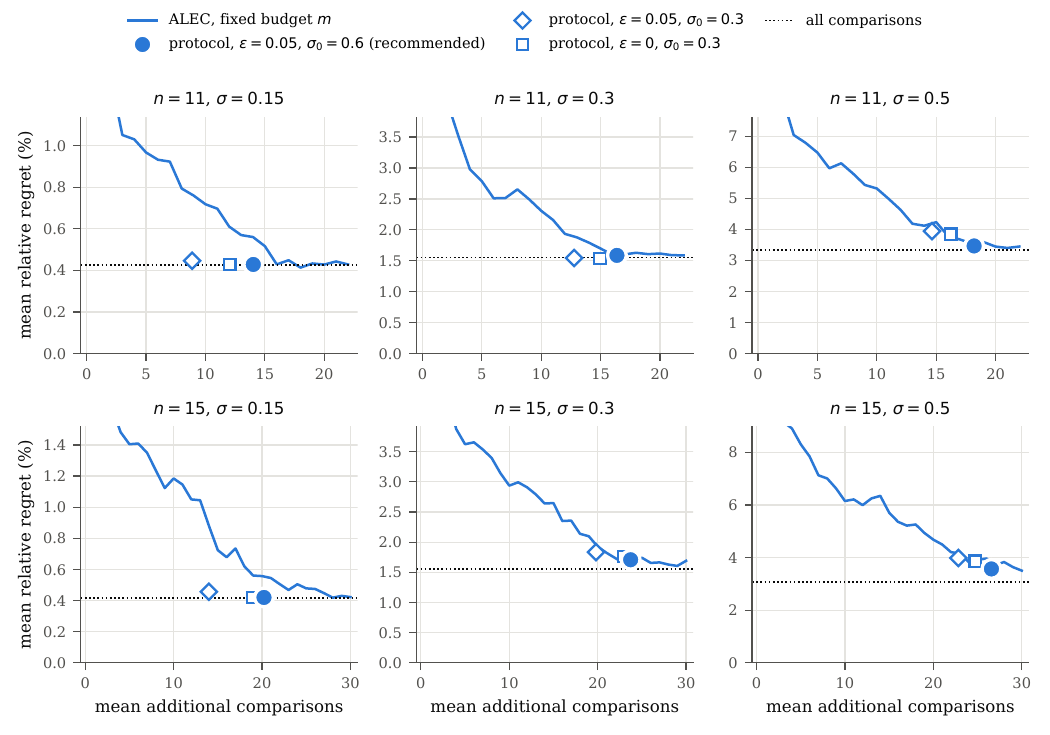}
\caption{Mean relative regret against the mean number of additional comparisons. Line: ALEC with a fixed budget $m=0,\dots,2n$; markers: the protocol (ALEC, stop when $\beta_\varepsilon\le0.10$ or after $2n$ additional comparisons) with the recommended setting (filled circle) and two settings with $\sigma_0=0.3$; dotted line: the complete matrix. \res{R-syn} replications per scenario.}\label{fig:protocol}
\end{figure}

Table~\ref{tab:protocol} evaluates the recommended protocol. It saves on average \res{prot-s06-saveadd-min}--\res{prot-s06-saveadd-max}\% of the $2n$ additional comparisons (\res{prot-s06-savetot-min}--\res{prot-s06-savetot-max}\% of the $3n$ comparisons in total), and \res{prot-s06-stop-min}--\res{prot-s06-stop-max}\% of the runs are stopped by the bound; the others reach the cap and have not met the posterior criterion. Among the runs stopped by the bound, the point estimates of the rate of losses above $\varepsilon$ are at or below the nominal 10\% in every scenario: at most \res{errearly-E05-ALEC-cap-s06-s030-max}\% for $\sigma\le0.3$ and \res{errearly-E05-ALEC-cap-s06-s050-min}--\res{errearly-E05-ALEC-cap-s06-s050-max}\% for $\sigma=0.5$, where the Wilson intervals extend above 10\%. The less prudent prior scale $\sigma_0=0.3$ stops earlier, but for $\sigma=0.5$ this rate rises to \res{errearly-E05-ALEC-cap-s050-min}--\res{errearly-E05-ALEC-cap-s050-max}\%. Figure~\ref{fig:protocol} shows the cost--regret trade-off, and Table~\Sref{tab:protocol-full} gives all tolerance and prior-scale settings.

Against a fixed ALEC budget with at least the same mean cost, hit-rate differences are generally small: after a Holm correction, \res{fix-hit-holm} of the \res{fix-rows} paired contrasts (six scenarios, two prior scales; \res{fix-s06-hit-holm} of the six with the recommended $\sigma_0=0.6$) favour the protocol, all with $\sigma=0.15$, and none favours the fixed budget (Supplementary Table~\Sref{tab:fixed-full}). One plausible interpretation, which the experiment does not isolate, is that under low noise the variable-length protocol stops earlier on easier instances and continues longer on harder ones.

\begin{table}[tb]
\centering\footnotesize\setlength{\tabcolsep}{3.0pt}
\caption{ALEC minus Harker-log under the same stopping rule ($\varepsilon=0.05$, $\alpha=0.10$, $\nu_0=2$, $\sigma_0=0.6$, cap $2n$): paired differences with 95\% intervals. $^\dagger$: significant after a Holm correction across the six scenarios. Full method-specific results are in Supplementary Table~\Sref{tab:samestop-full}.}\label{tab:samestop}
\begin{tabular}{@{}r r l r r r r@{}}
\toprule
$n$ & $\sigma$ & & $\Delta\bar q$ & $\Delta$Hit & $\Delta$Regret & $\Delta(\rho>\varepsilon)$ \\
 & & & (questions) & (points) & (points) & (points) \\
\midrule
11 & 0.15 & difference & $+0.10$ & $+1.9$$^{\dagger}$ & $-0.08$ & $-0.4$ \\
 & & 95\% CI & [$-0.04$, $+0.25$] & [$+0.8$, $+3.1$] & [$-0.15$, $-0.01$] & [$-1.2$, $+0.4$] \\
\addlinespace
11 & 0.30 & difference & $+0.53$$^{\dagger}$ & $+0.1$ & $-0.03$ & $-0.1$ \\
 & & 95\% CI & [$+0.30$, $+0.76$] & [$-1.4$, $+1.6$] & [$-0.21$, $+0.14$] & [$-1.5$, $+1.3$] \\
\addlinespace
11 & 0.50 & difference & $+0.71$$^{\dagger}$ & $+1.4$ & $-0.31$ & $-1.7$ \\
 & & 95\% CI & [$+0.40$, $+1.02$] & [$-0.5$, $+3.3$] & [$-0.64$, $+0.02$] & [$-3.5$, $+0.1$] \\
\addlinespace
15 & 0.15 & difference & $-0.02$ & $+0.7$ & $-0.04$ & $-0.3$ \\
 & & 95\% CI & [$-0.28$, $+0.24$] & [$-1.0$, $+2.4$] & [$-0.14$, $+0.06$] & [$-1.4$, $+0.8$] \\
\addlinespace
15 & 0.30 & difference & $+0.95$$^{\dagger}$ & $+0.8$ & $-0.17$ & $-1.3$ \\
 & & 95\% CI & [$+0.59$, $+1.32$] & [$-1.3$, $+2.9$] & [$-0.39$, $+0.05$] & [$-3.2$, $+0.5$] \\
\addlinespace
15 & 0.50 & difference & $+1.13$$^{\dagger}$ & $+5.0$$^{\dagger}$ & $-0.71$$^{\dagger}$ & $-4.1$$^{\dagger}$ \\
 & & 95\% CI & [$+0.73$, $+1.53$] & [$+2.5$, $+7.5$] & [$-1.14$, $-0.27$] & [$-6.3$, $-1.9$] \\
\bottomrule
\end{tabular}
\end{table}

Under the same stopping rule, ALEC asks slightly more questions (significantly more in \res{ss-holm-q} of the 6 scenarios after a Holm correction) but has higher hit-rate point estimates in all six, by \res{ss-dhit-min} to \res{ss-dhit-max} percentage points, with the clearest gain in the largest high-noise setting (Table~\ref{tab:samestop}). The full table also reports regret, stopping frequencies and conditional loss rates (Supplementary Table~\Sref{tab:samestop-full}). Thus the comparison is an accuracy--cost trade-off rather than a dominance result.

Calibration experiments in Supplementary Section~\Sref{sec:S-calib} show that the posterior bound behaves conservatively under a matched model with known noise, while model mismatch and an optimistic prior scale can raise the frequentist error rate. This distinction is important: the bound is a posterior statement under the working model, not a model-free frequency guarantee.

\subsection{Robustness}\label{sec:robust}
With alternative initial graphs, Dirichlet-distributed weights or judgments rounded to Saaty's scale, ALEC remains significantly more accurate than Harker-PG and random questioning in the unadjusted tests, while no significant difference from Harker-log is detected (Supplementary Table~\Sref{tab:robust}). Rounding reduces the accuracy of all rules, as expected; with Dirichlet weights, which more often produce a dominant alternative, the gain from re-ranking over Frozen disappears. Applying the recommended stopping rule to rounded judgments gives similar point estimates of the conditional loss rate in the two scenarios examined (Supplementary Table~\Sref{tab:saaty-stop}); larger problems and other discretisations remain open.

%% ============================================================================
\section{Discussion}\label{sec:discussion}
%% ============================================================================

\paragraph{Design implications for decision support systems}
From a DSS design perspective, ALEC separates the elicitation layer from the underlying decision model. The decision maker remains responsible for the judgments and for choosing the acceptable loss tolerance $\varepsilon$, the one parameter with a direct decision meaning; the risk level $\alpha$, the prior parameters and the budget cap are design choices of the DSS, with defaults such as those recommended below. The system determines which information to request next and whether further elicitation is warranted. This architecture makes the trade-off between user effort and decision reliability explicit rather than fixing the number of questions ex ante: with the same tolerance, the protocol stopped before the cap in \res{prot-s06-stop-min}--\res{prot-s06-stop-max}\% of the runs, and its mean number of questions grew with the noise in the judgments (Table~\ref{tab:protocol}). The final bound and the reason for stopping are part of the guidance the system gives, because they quantify the posterior support for the recommendation under the working model.

\paragraph{Practical guidance}
For a choice problem with $n$ alternatives, we recommend starting from a connected set of $n$ comparisons, updating the LLSM ranking after every answer, and asking the first unobserved pair in ALEC's current lexicographic order. Elicitation stops when $\beta_\varepsilon\le\alpha$ for a decision-maker-specified loss tolerance $\varepsilon$ or after $2n$ additional comparisons; the reported output should include the leader, the final bound and the reason for stopping. The configuration studied most extensively, and the one we recommend, uses $\varepsilon=0.05$, $\alpha=0.10$, $\nu_0=2$, $\sigma_0=0.6$ and this cap of $2n$; in a DSS, $\alpha$, $\nu_0$, $\sigma_0$ and the cap can be fixed at these values, and $\varepsilon=0.05$ offered as the default tolerance.

\paragraph{Interpreting the stopping rule}
The prior scale controls a genuine cost--reliability trade-off. We regard $\sigma_0=0.6$ as a prudent working value when consistency is unknown, not as a guarantee. Moreover, $\max\{0,1-\beta_\varepsilon\}$ is a lower bound on a posterior probability under the working model, not a model-free confidence level. A run that reaches the budget cap has not met the posterior criterion. Proposition~\ref{prop:winner} provides a complementary deterministic guarantee, but its worst-case consistency requirement becomes demanding as $n$ grows.

\paragraph{What the benchmark comparisons show}
Harker-log removes the orientation dependence of the raw sensitivity parametrisation and is substantially stronger than the original Harker procedure in our experiments. ALEC remains attractive because its next question is transparent and requires no derivative calculation; with a common stopping rule, it trades a modest increase in mean questioning cost for higher hit-rate point estimates in the scenarios studied. These findings should be interpreted as comparative evidence, not as proof of universal dominance.

\paragraph{Limitations and further research}
The analysis assumes independent homoscedastic log-normal errors, one judgment per pair and a single criterion. Real judgments are bounded, discretised and may exhibit heteroscedasticity, fatigue, anchoring or order effects. The rounding and alternative-design experiments address only part of this mismatch. Natural extensions include a hierarchy-level version of ALEC, repeated comparisons, richer priors or a full Bayesian look-ahead, and comparisons with further adaptive benchmarks such as \citet{FedrizziGiove2013}. Above all, the simulations measure the number of questions but not the user's experience: an experiment in which participants use an interactive implementation, comparing ALEC with a fixed questioning sequence in terms of the number of questions, time, accuracy and perceived effort, would test the user-side benefits that motivate the method. Data with objective magnitudes would permit external correctness checks, whereas subjective preference data would instead assess agreement with complete elicitation or direct ratings.

%% ============================================================================
\section{Conclusions}\label{sec:conclusion}
%% ============================================================================
ALEC is a simple adaptive elicitation mechanism for decision support with pairwise comparisons, for problems in which the goal is to identify the best alternative rather than estimate the entire priority vector uniformly well. Effective resistance gives a local statistical rationale for direct comparisons, while the simulation study shows that re-ranking is most valuable in larger problems with restricted budgets. ALEC was significantly more accurate than Harker's original procedure in most scenarios with $n\ge11$ and remained competitive with the stronger orientation-invariant Harker-log benchmark.

The posterior tolerance rule turns ALEC into a complete adaptive elicitation protocol that reduced the number of additional questions by \res{prot-s06-saveadd-min}--\res{prot-s06-saveadd-max}\% in our scenarios while making the decision risk explicit under a working model. The deterministic full-matrix certificate of Proposition~\ref{prop:winner} clarifies why a direct winner is not automatically a global one and why worst-case guarantees can be too conservative. Together, these results support an interpretable accuracy--cost framework for incomplete AHP elicitation. ALEC can therefore serve as an adaptive elicitation module for decision support systems, selecting the next query and stopping when further user effort is not justified by the required decision precision; validation with human decision makers using such a system is the main next step.

%% ============================================================================
\section*{Data and code availability}
The code, simulation results and scripts used to generate the tables, figures and numerical results in the article and supplementary material are available in the public \coderepo\ (release v1.3.2). Running \texttt{python code/run\_all.py} reproduces all results except the machine-dependent timings; the version used for this manuscript, its random seeds and the software versions are recorded in \texttt{code/VERSION}.

\section*{Declaration of generative AI and AI-assisted technologies in the writing process}
During the preparation of this work, the authors used ChatGPT (OpenAI) for exploratory calculations and an initial derivation of the conditions in Proposition~\ref{prop:winner}, and Claude (Anthropic) to assist with writing and testing the simulation code, checking and refining derivations, and editing the manuscript. All results reported in the paper were regenerated with the accompanying code, which includes unit tests of the computations and of the propositions. After using these tools, the authors reviewed and edited the content as needed and take full responsibility for the content of the published article.

\section*{CRediT authorship contribution statement}
%% Proposed by the drafting assistant from the history of the project; authors to confirm.
{\sloppy
\noindent\textbf{Luis \'Angel Calvo Pascual:} Methodology, Software, Formal analysis, Validation, Visualization, Writing -- review \& editing.

\noindent\textbf{Ji\v{r}\'{i} Mazurek:} Conceptualization, Methodology, Formal analysis, Writing -- original draft, Writing -- review \& editing.\par}

\section*{Declaration of competing interest}
The authors declare that they have no known competing financial interests or personal relationships that could have appeared to influence the work reported in this paper.

\section*{Acknowledgments}
This work was partially supported by the Ministry of Education, Youth and Sports Czech Republic within the Institutional Support for Long-term Development of a Research Organization in 2026.

\FloatBarrier

\end{document}